\documentclass[english,12pt]{paper}
\usepackage{geometry} 
\usepackage{imakeidx}

\makeindex[name=t, title=Index]
\makeindex[name=n, title=Index of Notations]

\usepackage{pgfplots}
\usepackage{psfrag}
\usepackage{amsmath}
\usepackage{calc}
\usepackage{systeme}
\usepackage{amsfonts}
\usepackage{amsthm}
\usepackage{enumerate}
\usepackage {amssymb}
\usepackage{algorithm2e}
\usepackage[latin1]{inputenc}

\usepackage[T1]{fontenc}
\usepackage{lmodern}

\usepackage{bbold}
\usepackage{graphicx}
\usepackage{caption}
\usepackage{subcaption}
\usepackage{color,overpic,cite,pdfsync,url}
\usepackage{amsmath,amsfonts,amssymb,amsthm,bbm,bm}
\usepackage{epstopdf,float}

\usepackage{tikz}
\usetikzlibrary{arrows.meta, positioning, calc, matrix}

\usepackage[numbers]{natbib}

\newtheorem{Theorem}{Theorem}[section]
\newtheorem{Lemma}[Theorem]{Lemma}
\newtheorem{Rem}[Theorem]{Remark}
\newtheorem{Def}[Theorem]{Definition}
\newtheorem{Prop}[Theorem]{Proposition}
\newtheorem{Cor}[Theorem]{Corollary}
\newtheorem{Ex}[Theorem]{Example}
\newtheorem{hyp}[Theorem]{Hypothesis}

\newcommand{\be}{\begin{equation}}
\newcommand{\ee}{\end{equation}}
\newcommand{\ba}{\begin{aligned}}
\newcommand{\ea}{\end{aligned}}

\newcommand{\R}{\mathbb{R}}
\newcommand{\T}{\mathbb{T}}

\newcommand{\M}{\mathcal{M}}

\newcommand{\mptd}{\mathcal{P}(\mathbb{T}^d)}
\newcommand{\mpo}{\mathcal{P}(\mathcal O)}
\newcommand{\mppo}{\mathcal{P}_p(\mathcal O)}
\newcommand{\mpprd}{\mathcal{P}_p(\R^d)}

\newcommand{\mpt}{\mathcal{P}_2(\R^d)}

\newcommand{\mprd}{\mathcal{P}(\mathbb{R}^d)}

\newcommand{\mo}{\mathcal{O}}
\newcommand{\eps}{\epsilon}

\newcommand{\mmo}{\mathcal{M}(\mathcal O)}

\newcommand{\mmpo}{\mathcal{M}_+(\mathcal O)}

\newcommand{\motd}{\mathcal{M}_1(\mathbb{T}^d)}
\newcommand{\mmuo}{\mathcal{M}_1(\mathcal O)}

\usepackage{appendix}
\usepackage{chngcntr}
\usepackage{etoolbox}
\usepackage{xcolor} 
\usepackage{hyperref}

\hypersetup{
    colorlinks=true, 
    linkcolor=blue,  
    citecolor=blue,  
    urlcolor=blue    
}
\AtBeginEnvironment{subappendices}{%
\section*{Appendix}
\addcontentsline{toc}{section}{Appendices}
\counterwithin{figure}{section}
\counterwithin{table}{section}
}

\numberwithin{equation}{section}

\title{Analysis on spaces of measures: Mean field game master equations and first order mean field Hamilton-Jacobi equations}
\author{Charles Bertucci}

\begin{document}
\maketitle
\setcounter{tocdepth}{2}

\tableofcontents
\newpage
\section*{Introduction}
I will not spare the reader this now classical introduction to a new research book on mathematics: an almost apologetic justification of its existence. An easy, and slightly dishonest explanation could be to say that this (short) book simply originates from the Peccot course I gave at Collège de France in 2023, where I presented some results I obtained on mean field game master equations and Hamilton-Jacobi equations, both set on spaces of measures. However, this does not justify at all the first half of what is about to follow...  In fact, while discussing the possibility of publishing the lecture notes of this course with some colleagues, in particular with my great friend Sylvain Sorin, I got convinced that a new presentation of basic aspects of analysis on spaces of measures might be of interest, in particular to young or non-specialist (in this topic) mathematicians, especially around the notions of derivatives. I insist upon the plural for the term notion, as several of them are used in the literature: Lions, flat, Wasserstein or linear derivatives can be encountered... It is often the largest source of confusion among inexperienced mathematicians in this field, and, if it is true that some of these terms refer to the same objects, there are also co-existing notions which have very little to do with one another. One of my main objectives when writing this book is to reach a more understandable presentation of these concepts, or at least to give my point of view on these questions.\\

The literature on analysis on spaces of measures, in particular on Wasserstein spaces, is growing at an incredible pace. It started with topics of dynamics: gradient flows, questions about what should be the equivalent of ordinary differential equations in this context,... and now concerns also partial differential equations (PDEs) on spaces of measures. Most of those mathematical objects arise from quite applied problems: the dynamics of optimal transport have proved to be of interest in image processing or economics for instance, mean field games are natural models in economics, mean field control problems naturally appear in engineering when dealing with decentralized protocols of large fleet of devices, contract problems in financial mathematics sometimes yield similar objects, and the fluctuations of models in statistical physics (or more generally mean field models) leads to the study of Hamilton-Jacobi-Bellman equations on spaces of probability measures.\\

Several excellent references already define derivatives of a function of a (probability) measure, such as Part II of the textbook of Ambrosio, Gigli and Savaré \citep{ags}, the presentation of Otto's calculus in Villani's textbook on optimal transport \citep{villani}, Lions' lectures at Collège de France in 2007-2008, Chapter V in Carmona and Delarue's textbook \citep{carmona2017probabilistic} or the short book of Cardaliaguet, Delarue, Lasry and Lions on mean field game master equations \citep{cdll}. The first two references do not include the more recent developments that were initiated by Lions to study mean field games, and the last two references have a presentation which is quite oriented toward the study of master equations, and may discourage researchers not interested in mean field games. In the present work, I try to give a self-motivated, and thus quite abstract, study of functions of a measure argument. I will focus on basic notions such as convexity, sub-differentiability, optimality conditions, and keep a safe distance from second order derivatives.\\

This presentation is not an end supposed to close a discussion, as several questions remain open, nor a beginning, because of the already existing literature, but a step aside, trying to offer a new perspective on what may appear as confused developments. It aims at presenting in a concise way a point of view that turned out to be quite helpful in my research, and which I hope will help others as well: there are at least two ways to differentiate a function of a measure, the two approaches are extremely detached from one another, and both are helpful in different situations. The first way, which I shall call vertical, consists in adding or taking away mass of the current measure. The second one, which I shall call horizontal, consists in moving around the mass of the current measure. I will be quite caricatural in this distinction, insisting that the notion of vertical derivative needs almost nothing on the underlying space on which measures are built, and I hope that this exaggeration will accentuate the inherent distinction between the two notions, rather than discourage applied mathematicians from working with unnecessary abstract notions.\\

Once these general developments are given in Part I, Part II of the book presents a study of mean field game master equations, which are non-linear transport equations written on spaces of measures. I will define weak notions of solutions, and show how to obtain existence, uniqueness and stability results for them, in what I believe is a rather concise presentation. Part II can serve as a comprehensive introduction to master equations, for people having already a small background on mean field games, namely on the forward-backward systems characterizing Nash equilibria. I will omit developments on second order master equations here (but still treat cases with common noise!).\\

Part III is concerned with mean field control problems and the typical Hamilton-Jacobi-Bellman equations that arise from them. I will derive such equations from control problems, obtain results of uniqueness of viscosity solutions through comparison principles, and prove results of existence by using, among other techniques, the representation of the solution with the value function of the problem. I then present a stochastic case, which motivates the study of such equations. While I give a rather complete study of deterministic problems, I will not talk about equations with singular terms, such as the ones arising from so called idiosyncratic noises. In my opinion Part III can serve as a general introduction to researchers interested in such questions.\\

The three parts are rather mathematically independent of one another, and only a handful of results of Part I will be used in Parts II and III (results which I will highlight). The main reason for this independence is that the typical solutions of the partial differential equations at hand in Parts II and III are not, in general, sufficiently regular to be differentiable. However, the distinction between the tools used in Part II and III serves as an illustration of the co-existence between the two notions of derivatives introduced in Part I, and where the ideas might come from. I insist upon this because, at a formal level, all partial differential equations of interest can be written with either one of those notions, but only the vertical point of view shall be used in Part II while only tools from the horizontal one will help in Part III.\\

Concerning the writing of this book, my first thanks go to Sylvain Sorin for his all-around help in this redaction, and mainly for pushing me to write what constitutes the first part of this work. As usual, I have made significant progress when trying to put on paper what I thought I understood. For this, I cannot thank him enough. I feel extremely grateful to Pierre-Louis Lions, on too many points to mention them all here, but especially for agreeing to write a preface of this book, for his valuable input on an early version of the manuscript and for the amazing lectures and seminars he gave during his tenure at Coll\`ege de France. I learned so much from them that I feel obliged to give the following advice: (mathematical) French is not that hard, it is worth it. Several discussions with Jean-Michel Lasry, especially on the Lagrangian point of view influenced me a lot, and for this I thank him here. Leaving typos or small gaps is one of my specialty and to remove them, I receive great help from Valentin Pesce, Charles Meynard and Giacomo Ceccherini Silberstein. There is no doubt that some remain, but they helped me to remove many of them. A great deal of this book was written in the kind atmosphere of the CEREMADE which I thank for welcoming me. A significant portion of this work has also been realized in the calm Lagrange Mathematics and Computing Research Center to which I am also grateful for its hospitality. The ERC is funding me this last year under the project PArtial DIfferential Equations on SEts of Measures. Last but not least, I want to thank \'Emilie and Georges for their love and patience.

 \begin{flushright} Charles Feral-Bertucci,
 
 Paris, July 2026.\end{flushright}

\newpage

\section*{Notation}
Here are some general notations. Several more are to be introduced in the core of the book.\\

The $d$ dimensional Euclidean space is $\R^d$ while the flat torus is $\T^d$, which is defined as the quotient space $\R^d/ \mathbb Z^d$, where $\mathbb Z$ is the set of integers. The flat torus $\T^d$ is equipped with the distance $d(x,y) = \inf_{z \in \mathbb Z^d}|x-y + z|$, for $|\cdot|$ the Euclidean norm on $\R^d$. The scalar product between $x,y \in \R^d$ is denoted $x \cdot y$. The set of real $d\times d$ matrices is denoted by $\mathcal M_d(\R)$.\\

Functional spaces are denoted as follows:
\begin{itemize}
\item Given a Banach space $E$, a probability space $(\Omega, \mathcal A,µ)$ and $p \in [1,\infty]$, the $p$ Lebesgue space of equivalence classes of functions $f : \Omega \mapsto E$ for the relation "being equal $µ$ almost everywhere" is denoted by $L^p((\Omega,µ),E)$. When $\Omega \subset \R^d$ and $µ$ is the Lebesgue measure, I simply write $L^p(E)$, and the $E$ can also be dropped when the context is clear.
\item Given a Banach space $E$, a measurable space $(\Omega, \mathcal A)$, the space of bounded measurable functions $f : \Omega\mapsto E$ is denoted by $\mathcal B(\Omega,E)$. In particular elements of $\mathcal B$ are defined everywhere while elements of $L^\infty((\Omega,µ),E)$ are only defined $µ$ almost everywhere.
\item For $X,Y$ two metric spaces, the set of continuous functions $f:X\mapsto Y$ is $\mathcal C(X,Y)$, the set of $\alpha$ H\"older continuous functions $f:X\mapsto Y$, for $\alpha \in (0,1)$ is $\mathcal C^\alpha(X,Y)$. The Lipschitz constant of a function $f: X \to Y$ is $Lip(f) := \sup_{x\ne y}\frac{d_Y(f(x),f(y))}{d_X(x,y)}$.
\item For $k \geq 1$ integer and $p \in [1,\infty]$, the Sobolev spaces of $\R^d$ valued functions over a domain $\mo \subset \R^d$ with derivatives up to order $k$ in $L^p$ are denoted as usual $W^{k,p}(\mo,\R^d)$, and $H^k(\mo,\R^d)$ when $p=2$.
\item For $\mo$ a domain of $\R^d$, the space of continuously differentiable functions taking values in $\R^k$ is $\mathcal C^1(\mo,\R^k)$. The set of such functions which are infinitely many times continuously differentiable is $\mathcal C^\infty(\mo,\R^k)$. When such functions are restricted to have compact support, I write $\mathcal C_c^1(\mo,\R^k)$ and $\mathcal C_c^\infty(\mo,\R^k)$, and the same goes for continuous functions with compact support: $\mathcal C_c(\mo,\R^k)$.
\item The set of real bounded uniformly continuous functions on a metric space $(E,d)$ is denoted by $BUC(E)$.
\item When the space into which the function is valued is $\R$, it is not written. For instance $\mathcal C(X)$ is the space of continuous real functions on $X$.
\end{itemize}

A modulus of continuity is a continuous concave non-decreasing map $\omega: \R_+\mapsto \R_+$, smooth on $(0,\infty)$ such that $\omega(0) = 0$. In particular, it is sub-additive.\\

For a function $u: \R^d \mapsto \R\cup\{-\infty; + \infty\}$, the sub-differential $\partial^-u(x)$ of $u$ at $x \in \R^d$ such that $u(x) \in \R$ is the set of elements $p$ of $\R^d$ such that there exists a modulus of continuity $\omega$ such that for all $y \in \R^d$
$$
u(y) \geq u(x) + p\cdot(y-x) -|y-x|\omega(|y-x|).
$$

For $E$ a topological space, the upper and lower semi continuous envelope of a locally bounded function $u: E \mapsto \R$ are denoted respectively by $u^*$ and $u_*$. They are defined as $u^*(x) = \limsup_{y \to x} u(y)$ and $u_*(x) = \liminf_{y \to x} u(y)$. Moreover, upper semi continuity and lower semi continuity are abbreviated respectively usc and lsc.\\

A correspondence $A$ between two sets $X$ and $Y$ is a map defined on $X$ taking values in the subsets of $Y$, it is noted $A: X\rightrightarrows Y$.\\

Given a tuple $x=(x_1,x_2,\dots,x_n)$ and $1 \leq i \leq n$, the map $\pi_i$ is the projection on the $i$-th coordinate, that is $\pi_i(x) = x_i$.\\

For $t \in \R$, the map $e_t$ is the evaluation mapping at $t$. Hence, if $I \subset \R$, $\theta : I \mapsto X$ and $t \in I$, then $e_t(\theta) = \theta(t)$.\\

For $h \in \R^d$, $\tau_h : \R^d \mapsto \R^d, x \mapsto x +h$.\\

The topological dual of a topological vector space $X$ is denoted by $X'$. The duality product will be noted by $\langle \cdot,\cdot\rangle_{E',E}$, or simply $\langle\cdot,\cdot\rangle $ when the context is clear. In particular, when dealing with function spaces, it will always refer to the extension of the $L^2$ scalar product. Given $Y \subset X$, a mapping $f: Y \mapsto E'$ is monotone if for all $x,y \in Y$,
$$
\langle f(x) -f(y),x-y\rangle_{E',E} \geq 0.
$$
The mapping $f$ is strictly monotone if the previous inequality is strict as soon as $x \ne y$.\\

In a metric space $(E,d)$, the closed ball of radius $R \geq 0$ centered at $x_0 \in E$ is denoted by $B(x_0,R)$.\\

The notation $o(\cdot)$ stands for the usual Landau notation.\\

For $p \in (1,\infty)$, $p' = \frac{p}{p-1}$ is its conjugate exponent.

\newpage

\part{Analysis on spaces of measures}

In this first part, I present several properties of functions $U: X \mapsto \R$ for $X$ a certain subset of the space of signed measures $\mmo$ on a measurable space $\mo$. I will concentrate on cases in which $X$ is either a subset of the space of signed measures equipped with total variation, or in cases where $X$ is some Wasserstein space. I recall some standard definitions, as well as metric and topological properties of spaces of measures and their implications on the study of continuous functions.

The core of this first part is concerned with notions of derivatives of a function $U: X \mapsto \R$. In practice, such a notion is most of the time used to compute
\be\label{quotient}
\lim_{t \to 0^+}\frac{U(m_t)-U(m_0)}{t},
\ee
for some variation, or path, $(m_t)_{t \in [0,1]}$ valued in $X$. As typical paths of interest can be irregular and of different natures, it is mathematically helpful to consider a more abstract notion of derivative than just giving a definition for the previous limits. That is why I will focus on defining notions of derivatives by finding equivalents of first order Taylor expansions of a map $U$ near a measure $µ$, quite in the spirit of the Fréchet derivative.\\

In order to present the different problems at hand, I will list some classical variations $(m_t)_{t \in [0,1]}$ for which it might be interesting to be able to compute \eqref{quotient}. I will then introduce two notions of derivatives, in two different settings: a vertical notion which requires nothing on the base space $\mo$ but to be a measurable space, and a horizontal notion which shall rely on the possibility to consider differential calculus on the base space $\mo$ itself. I then present different notions of convexity, of super-differentiability, $\mathcal C^{1,\alpha}$ regularity, as well as some optimality conditions and results of perturbed optimization.\\

I feel the need to recall here that this part is quite independent of the rest of the book. Indeed, as the two next parts are concerned with the study of non-differentiable solutions to PDEs, differentiability shall play no role in what follows. Hence, for people only interested in precise statements on MFG master equations or on mean field optimal control, they can skip this Part a bit. However, I do not particularly encourage it, as it is always a good idea to learn to walk before trying to run, which here takes the form of understanding usual notions of derivatives before going to weaker notions. In any case, only Section \ref{sec:optim} will be of use in both the next two parts. Section \ref{sec:flat} allows one to understand the terms in MFG master equations in Part II, while Section \ref{sec:horiz} helps to understand the HJB equations of Part III. The notion of super-differentials of Section \ref{sec:superdiff} is also fundamental in Part III.

\section{Distances, topologies and continuous functions on spaces of measures}
In this part, I define the main objects of interest for the rest of the book. I also recall some well-known results about various topologies and distances on sets of measures.

\subsection{Measures and total variation distance}
In what follows, a \textit{measure} $µ$ on a measurable space $(\mo,\mathcal A)$ is a $\sigma$-additive real function on the $\sigma$-algebra $\mathcal A$ of $\mo$, which gives value $0$ to the empty set. When $\mo$ is a topological space, $\mathcal A$ is always taken as the Borel $\sigma$-algebra, that is the smallest $\sigma$-algebra containing all open sets of $\mo$. The positive and negative parts $µ^+$ and $µ^-$ of a measure $µ$ are given by Hahn-Jordan decomposition theorem and given by $µ^+(A) = \sup_{B \in \mathcal A, B \subset A} µ(B)$ and by $µ^-(A) = -\inf_{B \in \mathcal A, B \subset A}µ(B)$. The total variation $µ_{TV}$ of $µ$ is equal to $µ^+ + µ^-$. The \textit{norm of a measure} $µ$ is the quantity 
$$
|µ| := µ_{TV}(\mo).
$$
The set of \textit{finite measures}, that is measures $µ$ such that $|µ| < \infty$, is denoted by $\mathcal M(\mo)$. In particular, $(\mathcal M(\mo),|\cdot|)$ is a normed vector space, which is complete since $\R$ is. A measure $µ$ is non-negative if $µ^- = 0$. The set of non-negative finite measures is denoted by $\mathcal M_+(\mo)$, the set of non-negative measures with mass at most one is denoted by $\mathcal M_1(\mo)$ and the set of probability measures, that is non-negative measures such that $µ(\mo) = 1$, by $\mathcal P(\mo)$. Observe that $\mathcal M_+(\mo), \mathcal M_1(\mo)$ and $\mathcal P(\mo)$ are closed convex subsets of $\mmo$. 

For $µ,\nu \in \mathcal M(\mo)$, the \textit{total variation distance} between $µ$ and $\nu$ is the norm of the difference $|µ-\nu|$. The Dirac mass at $x \in \mo$ is denoted by $\delta_x$.\\

For $\phi$ a measurable real valued function on $\mo$, the Lebesgue integral, when it is well defined, of $\phi$ against $µ$ is denoted by both
$$
\langle \phi,µ\rangle = \int_\mo\phi(x)µ(dx).
$$
Note in particular that if $\phi \in \mathcal B(\mo)$, the set of bounded measurable functions over $\mo$, then for any $µ \in \mmo$, the integral is well defined and $\langle \phi,µ\rangle \leq \|\phi\|_\infty|µ|$. In particular, this implies that we can identify $\mathcal B(\mo)$ with a subspace of $(\mmo)'$, the topological dual of $\mmo$. The space $(\mmo)'$ is much larger than $\mathcal B(\mo)$, for instance consider the case $\mo = [0,1]$ (equipped with the Borel $\sigma$-algebra) and consider the map $\mathcal T: µ \mapsto µ_s(\mo)$, where $µ_s$ is the singular part of $µ$ given by the Lebesgue decomposition Theorem. Clearly, $\mathcal T$ is linear, and $|T(µ)| \leq |µ|$, hence it is also continuous. However, if it were given by an element $\phi \in \mathcal B(\mo)$, we should have $\phi = 1$, which is not possible.\\

If $\phi$ is valued in a normed vector space, the notation $\langle \phi,µ\rangle$ will by convention refer to the Bochner integral of $\phi$ against $µ$, see II.2 of \citep{diestel}.\\

For a measurable mapping $f: \mo \mapsto Y$, for $Y$ another measurable space and $µ \in \mathcal M(\mo)$, the \textit{image measure} $f_\#µ \in \mathcal M(Y)$ is defined, for all $B$ measurable set in $Y$, by
$$
f_\#µ(B) = µ(f^{-1}(B)).
$$

\subsection{Weak convergence of measures}
The previous development does not require any topology on $\mo$, and thus sees all points of $\mo$ as infinitely far away from each other. In this section, $\mo$ is a topological space and I want to consider topologies on $\mmo$ which formalize the idea that two measures are close if they put similar mass on points which are close to one another.

The most frequently studied topology which formalizes this idea is the so-called \textit{weak topology} of measures. It is induced by the following notion of convergence: the sequence $(µ_n)_{n \geq 0}$ valued in $\mmo$ weakly converges toward $µ \in \mmo$ if for every function $\phi \in \mathcal C_b(\mo)$, the space of bounded continuous functions on $\mo$, 
$$
\lim_{n \to \infty} \langle \phi,µ_n\rangle = \langle \phi,µ\rangle.
$$
For measures with finite norms, convergence in total variation implies weak convergence, as for any $\phi \in \mathcal C_b(\mo)$,
$$
|\langle \phi,µ_n - µ\rangle | \leq \|\phi\|_\infty |µ_n -µ|.
$$
The converse is false since $(\delta_{x_n})_{n \geq 0}$ converges in total variation toward $\delta_x$ if and only if $(x_n)_{n \geq 0}$ is equal to $x$ after a certain rank, while it converges weakly as soon as $(x_n)_{n \geq 0}$ converges toward $x$ in $\mo$.

Complete separable metric spaces are of particular interest in measure theory, because all probability measures on them are in fact Radon measures (Chapter 7 in \citep{bogachev}), and thus quite regular. Such spaces are called \textit{Polish spaces}. When $\mo$ is a Polish space, there exist metrics which give the same notion of convergence as the one of weak convergence on $\mpo$, this is for instance the case of the Lévy-Prokhorov metric (Theorem 8.2.2 in \citep{bogachev}). Furthermore, Prokhorov's Theorem gives an easy way to verify pre-compactness of subsets of $\mpo$ for the weak topology.
\begin{Theorem}\label{thm:proh}
Let $\mo$ be a Polish space and $A \subset \mpo$. The closure of $A$ in $\mpo$ (for the weak convergence) is sequentially compact if and only if, for any $\eps > 0$, there exists a compact set $K_\eps \subset  \mo$ such that $\inf \{ µ(K_\eps) | µ \in A\} \geq 1- \eps$.
\end{Theorem}
In particular, if $\mo$ is a Polish compact set, then $\mpo$ is also compact.

\subsection{Couplings, gluing and disintegration}
Let $(\mo,d)$ be a Polish space. For $Y$ another Polish space, $µ \in \mathcal P(\mo), \nu \in \mathcal P(Y)$, the set of \textit{couplings} $\Pi(µ,\nu)$ between $µ$ and $\nu$ is defined by
$$
\Pi(µ,\nu):= \{ \gamma \in \mathcal P(\mo\times Y)| (\pi_1)_\#\gamma = µ, (\pi_2)_\#\gamma = \nu\}.
$$
It is always non-empty, as we can always consider the independent coupling

\noindent$µ(dx)\nu(dy) \in \Pi(µ,\nu)$. A coupling $\gamma$ is called \textit{deterministic} if there exists a measurable map $T: \mo \mapsto Y$ such that $\gamma = (Id,T)_\#µ$.

The following result is called \textit{disintegration} theorem, and it yields the existence of so-called conditional probability measures. It relies on a notion of measurability for a mapping $Y \mapsto \mathcal P(\mo)$, which is a priori not canonical since I did not fix a $\sigma$-algebra, or even a topology, on $\mpo$. By convention, for the rest of the book, it is chosen as follows. A mapping $\psi: Y \mapsto \mpo$ is called measurable if for any measurable set $B \subset \mo$, the map $y \mapsto \psi_y(B)$ is measurable.
\begin{Theorem}\label{thm:disint}
Let $Y$ be a Polish space, $µ \in \mathcal P(\mo)$ and $T: \mo \mapsto Y$ be measurable. Let $\nu = T_\#µ \in \mathcal P(Y)$. Then, there exists a measurable mapping $Y \mapsto \mpo$, $y \mapsto µ_y$, such that 
\begin{itemize}
\item $\nu$ almost surely, $µ_y(T^{-1}(\{y\})) = 1$, 
\item For any measurable map $f: \mo \mapsto [0,\infty)$,
$$
\int_\mo f(x)µ(dx) = \int_Y \int_{T^{-1}(\{y\})}f(x)µ_y(dx)\nu(dy).
$$
\end{itemize}
Furthermore, two such mappings are equal $\nu$ almost everywhere.
\end{Theorem}

In what follows this result shall be used in two main arguments. In the first one, when given a coupling $\gamma \in \Pi(µ,\nu)$, for $µ \in \mathcal P(\mo)$ and $\nu \in \mathcal P(Y)$, I shall use it to write $\gamma(dx,dy) = µ(dx)\psi_x(dy)$, where $\psi: \mo \mapsto \mathcal P(Y)$ is a measurable map given by the disintegration Theorem along the map $\pi_1$. The second one will be in the proof of the Gluing Lemma which follows.

\begin{Lemma}\label{lemma:gluing}
Let $\mo,Y,Z$ be Polish spaces, $µ,\nu,\eta$ in respectively $\mathcal P(\mo),\mathcal P(Y),\mathcal P(Z)$. Let $\gamma \in \Pi(µ,\nu), \tilde \gamma \in \Pi(µ,\eta)$. Then, there exists a $\Gamma \in \mathcal P(\mo\times Y\times Z)$ such that $(\pi_1,\pi_2)_\#\Gamma = \gamma$, $(\pi_1,\pi_3)_\#\Gamma = \tilde \gamma$. Furthermore, if either $\gamma$ or $\tilde\gamma$ is deterministic, then there exists a unique such $\Gamma$.
\end{Lemma}
\begin{proof}
By the disintegration Theorem, we can write $\gamma(dx,dy) = µ(dx)\psi_x(dy)$. Writing $\Gamma(dx,dy,dz) = \tilde \gamma(dx,dz)\psi_x(dy)$ produces a required coupling.
The uniqueness part follows from the following: if $\Gamma$ is a gluing of $\gamma$ and $\tilde \gamma$, and, denoting by $\Gamma_x,\gamma_x$ and $\tilde \gamma_x$ the disintegration of the three previous probability measures along $\pi_1$, then, $µ$ almost everywhere, $\Gamma_x \in \Pi(\gamma_x,\tilde \gamma_x)$.

 First, observe that this last property clearly implies the uniqueness of the coupling. Indeed, if $\gamma$ (or $\tilde \gamma$) is deterministic, it implies that $µ$ almost everywhere, $\gamma_x$ (or $\tilde \gamma_x$) is a Dirac mass. Since the set of couplings between Dirac masses and anything else is reduced to a singleton, it follows that $\Gamma_x$ is uniquely defined $µ$ almost everywhere which yields the result.
 
 Now this last property is implied by the fact that $(\pi_1,\pi_2)_\#\Gamma$ can be computed in two ways. On one hand it is given by $\gamma$, and on the other hand, it is also given by $((\pi_1)_\#\Gamma_x(dy,dz))µ(dx)$. Hence, the uniqueness of the disintegration implies that necessarily, $µ$ almost everywhere, $(\pi_1)_\#\Gamma_x = \gamma_x$. The same obviously holds for the other component.

\end{proof}

\subsection{Wasserstein spaces and first look at optimal transport}
To metricize weak convergence, the attention is often put on subsets $\mathcal P(\mo)$, equipped with so-called Wasserstein metrics. Here, $(\mo,d)$ is still a Polish space. Let $ p \in [1,\infty)$. The set $\mathcal P_{ p}(\mo)$ is defined as 
$$
\mathcal P_{ p}(\mo) := \left\{ µ \in \mathcal P(\mo), \exists x_0 \in \mo, \int_\mo d(x,x_0)^pµ(dx) < \infty\right\}.
$$
This definition does not depend on the choice of $x_0$, as if it is satisfied for one, it is also satisfied for all elements of $\mo$. The previous definition has an analogue when $p = \infty$, given by
$$
\mathcal P_{ \infty}(\mo) := \left\{ µ \in \mathcal P(\mo), \exists R > 0, x_0 \in \mo, µ(B(x_0,R)) =1 \right\}.
$$
When $\mo$ is bounded, the restriction of $\mpo$ to $\mppo$ is useless as the sets are then equal. For $p \in [1,\infty)$, $µ,\nu \in \mathcal P_{ p}(\mo)$ and $\gamma \in \Pi(µ,\nu)$, we have $\gamma \in \mathcal P_p(\mo\times \mo)$ and I define
$$
C_p(\gamma) := \left(\int_{\mo\times \mo}d(x,y)^p\gamma(dx,dy)  \right)^\frac{1}{p}.
$$
For $p \in [1,\infty)$, the \textit{Wasserstein distance} on $\mathcal P_{ p}(\mo)$ is defined by
$$
\mathcal W_p(µ,\nu) := \inf_{\gamma \in \Pi(µ,\nu)}C_p(\gamma).
$$
I leave as an exercise to verify, using the Gluing Lemma, that the previous quantity indeed defines a distance. The spaces $(\mppo,\mathcal W_p)$ are called \textit{Wasserstein spaces}. Note that $\mathcal P_p(\mo) \subset \mathcal P_{q}(\mo)$ and $\mathcal W_{q} \leq \mathcal W_p$, for $q\leq p$. One can also define a distance $\mathcal W_\infty$ on $\mathcal P_\infty(\mo)\times \mathcal P_\infty(\mo)$ as follows
$$
\mathcal W_\infty(µ,\nu) := \inf_{ \gamma \in \Pi(µ,\nu)}\|d(x,y)\|_{L^\infty((\mo^2,\gamma))}.
$$

The following result is probably the most important one on Wasserstein spaces.
\begin{Theorem}\label{thm:Wasserstein}
For $p\in[1,\infty)$, the space $(\mathcal P_p(\mo),\mathcal W_p)$ is a complete separable metric space. For any $µ,\nu \in \mathcal P_p(\mo)$, the infimum in $\mathcal W_p(µ,\nu)$ is reached. 
\end{Theorem}

I will also use several times the following result.
\begin{Theorem}\label{thm:glivenko}
Let $p \in [1,\infty)$ and $µ \in \mppo$, there exists a sequence $(x^n)_{n \geq 1}$ such that $x^n \in \mo^n$ for all $n \geq 1$ and 
$$\lim_{n \to \infty}\mathcal W_p\left(µ, \frac1n \sum_{i=1}^n \delta_{x^n_i}\right) = 0.$$
\end{Theorem}

The space $\mathcal P_1(\mo)$ is rather special, namely because of the following property that satisfies its metric, which in fact comes from a norm as the next result shows. This result is called Kantorovich's duality, and it is a classical example of duality result in constrained optimization.
\begin{Theorem}\label{thm:kantorovich}
Let $µ,\nu \in \mathcal P_1(\mo)$, then 
$$
\ba
\mathcal W_1(µ,\nu) &= \sup_{Lip(\phi) \leq 1} \langle \phi,µ-\nu\rangle,\\
&=: \|µ-\nu\|_{W^{-1,1}}.
\ea
$$
\end{Theorem}
The last equality can be taken as the definition of the norm of $W^{-1,1}$, which is by definition the topological dual of the space of Lipschitz functions equipped with an appropriate norm. Finally, as mentioned above, Wasserstein metrics are naturally associated with weak convergence of measures, as the next result shows.
\begin{Theorem}\label{thm:topowass}
Let $p \in [1,\infty)$, $µ \in \mppo$ and $(µ_n)_{n \geq 0}$ be a sequence valued in $\mppo$. Then the two following statements are equivalent.
\begin{itemize}
\item $\lim_{n \to \infty} \mathcal{W}_p(µ_n,µ) = 0$,
\item $(µ_n)_{n \geq 0}$ converges weakly toward $µ$ and, for some $x_0 \in \mo$,
$$
\lim_{n \to \infty} \int_{\mo}d(x,x_0)^pµ_n(dx) = \int_{\mo}d(x,x_0)^pµ(dx).
$$
\end{itemize}
\end{Theorem}
As an immediate corollary of the previous and Theorem \ref{thm:proh}, we have the following.
\begin{Prop}
Let $p,q \in [1,\infty)$ with $p < q$. Then, bounded sets of $\mathcal P_q(\R^d)$ are relatively compact in $\mpprd$.
\end{Prop}

The Wasserstein metrics are the distances of optimal transport, the problem of transporting optimally a measure $µ$ toward a measure $\nu$ for a given cost. Entire textbooks are written on optimal transport and I will simply hint at why such a link exists here. Associated with a coupling $\gamma \in \Pi(µ,\nu)$ is the idea of pairing elements of mass $dx$ and $dy$ by giving a mass $\gamma(dx,dy)$ to the pair. If the cost of transporting an element of mass from $x$ to $y$ is given by $d(x,y)^p$ then, the average cost associated to all the pairings induced by $\gamma$ is clearly $C_p^p(\gamma)$. Thus $\mathcal W^p_p(µ,\nu)$ is the infimum of all cost of the transports given by couplings.\\

When $\mo = \R^d$, there are other functionals, closely related to Wasserstein distances, that will be helpful in what follows. For $p \in (1,\infty)$, the functional $\mathcal I_p$ is defined on $\mathcal P_p(\R^d)\times \mathcal P_{p'}(\R^d)$ by
$$
\mathcal I_p(µ,\nu) := \inf_{\gamma \in \Pi(µ,\nu)} -\int_{\R^d\times \R^d}x\cdot y\quad \gamma(dx,dy).
$$
Just as it was the case in the previous result, the infimum is always reached here. For $p \in (1,\infty), µ \in \mpprd$, we denote by $M_p(µ) := \int_{\R^d}|x|^pµ(dx)$. We then have the identity $\mathcal{W}_2^2(µ,\nu) = M_2(µ) + M_2(\nu) + 2\mathcal I_2(µ,\nu)$.\\

I shall also use the notation $\mathcal P_{pp'}(X\times Y):= \{\gamma \in \Pi(µ,\nu) | µ \in \mathcal P_p(X), \nu \in \mathcal P_{p'}(Y)\}$.

\subsection{Probabilistic formulation of problems in $\mathcal P_p(\mo)$}
Sometimes, it will be useful in what follows to pass from a problem posed on the set of probability measures, to a problem posed on the set of random variables. This procedure is called a \textit{lift}, since the space of random variables is significantly larger than the set of their laws. Here, I make precise some concepts in this direction. In this section, $(\mo,d)$ is a Polish space.\\

To consider a lift, we first need a probability space $(\Omega, \mathcal A,\mathbb P)$, whose expectation is denoted $\mathbb E$. This choice is mainly free. The reader should choose the one she or he wants. The only requirements are: that $\Omega$ is Polish, and that it is atomless (i.e. without $\omega \in \Omega$ such that $\mathbb P(\{\omega\})> 0$) to apply Theorem \ref{thm:existrv} below. I will call such probability spaces \textit{standard}. Measurable functions from this probability space to $\mo$ will be called random variables. For a random variable $X$, its law, denoted by $\mathcal L(X)$ is the probability measure $X_\#\mathbb P$. I write $X \sim µ$ when $X$ is of law $µ$.\\

For $p \in [1,\infty)$ and $X \in L^p(\Omega,\mo)$, its law $\mathcal L(X)$ belongs to $\mppo$. The $p$ Wasserstein distance can also be formulated through 
$$
\mathcal W_p^p(µ,\nu)= \inf_{X \sim µ,Y\sim \nu} \mathbb E[d(X,Y)^p].
$$
Note that when $\mo \subset \R^d$, $d(X,Y)^p$ denotes the Euclidean norm of $X-Y$ to the power $p$, and not some $p$-norm in $\R^d$.\\

When using lifting methods, two classical probability results shall be useful.
\begin{Theorem}\label{thm:existrv}
Let $(\Omega,\mathcal A,\mathbb P)$ be a standard probability space and $µ \in \mpo$. Then, there exists a random variable $X: \Omega\mapsto \mo$ such that $\mathcal L(X) = µ$.
\end{Theorem}
\begin{proof}
I sketch a proof, in the case in which $\mo$ is uncountable. The first step consists in constructing one-to-one mappings $f: [0,1] \mapsto \mo$ and $g: [0,1]\mapsto \Omega$ such that $f,f^{-1},g$ and $g^{-1}$ are all measurable ($[0,1]$ is equipped with the Lebesgue measure), see Theorem 13.1.1 in \citep{dudley}. We have $(f^{-1})_\#µ \in \mathcal P([0,1])$. We can consider $\xi:([0,1],g^{-1}_\# \mathbb P)\mapsto [0,1]$ whose law is $(f^{-1})_\#µ$ by using the cumulative distribution functions, because $g^{-1}_\# \mathbb P$ is atomless. Observe now that $f(\xi(g^{-1}))$ is a suitable random variable.
\end{proof}

\begin{Theorem}\label{thm:proba}
Let $(\Omega,\mathcal A,\mathbb P)$ be a standard probability space, and consider two random variables $X,X'$ which have the same law. Then, for any $\eps > 0$, there exists a one-to-one measurable mapping $\tau: \Omega\mapsto \Omega$ such that $\tau^{-1}$ is measurable, $\tau_\#\mathbb P = (\tau^{-1})_\#\mathbb P = \mathbb P$ and 
$$
\|d(X\circ \tau,X')\|_\infty \leq \eps.
$$
\end{Theorem}

For $\mo \subset \R^d$, $p \in [1,\infty)$, $X \in L^p(\Omega,\mo)$ of law $µ$, we will consider functions of $X$, that is to consider $f(X)$, for $f \in L^p_µ:= L^p(\mo,µ)$, the usual Lebesgue space of order $p$. These types of $L^p_µ$ space have to be manipulated with care as they are simply equivalence classes for the equivalence relation: "being equal $µ$ almost everywhere". In particular, they can then contain elements which are not Borel measurable. As a convention, in what follows, I will always restrict my attention to elements of $L^p_µ$ which are Borel-measurable.\\

As a convention, $L^p(\Omega,\R^d)$ shall be denoted by $\mathbb L^p$, and $\mathbb L^2 = \mathbb H$. Recall that it is a Banach space, when equipped with the norm $\|X\|_p = (\mathbb E[|X|^p])^\frac1p$, and that $\mathbb H$ is a Hilbert space. Furthermore, I define the \textit{lift} of a function $U : \mpprd \mapsto \R$ as the function $\mathcal U: \mathbb L^p \mapsto \R$ defined by
$$
\forall X \in \mathbb L^p,\, \,\mathcal U(X) = U(\mathcal L(X)).
$$

\subsection{Spaces of continuous functions over spaces of measures}
The Wasserstein and total variation metrics allow us to define naturally several regularity spaces for functions over $\mppo$ or $\mmo$, such as the regularity spaces $\mathcal C$, $\mathcal C^{\alpha}$ for $\alpha \in (0,1)$, or spaces of Lipschitz continuous functions for instance. In addition to the standard properties of such spaces, I present some approximation of continuous functions of measures and the link with the continuity of the lift.\\

\textbf{Polynomials as approximations of continuous functions of measures}\\
I assume in this section that $\mo$ is a compact set of $\R^d$ and I consider $\mmuo$ equipped with the weak topology. I say that a function $U: \mmuo \mapsto \R$ is a \textit{monomial} when either it is constant, or either there exist $n \geq 1$ and $n$ continuous functions $(\phi_i)_{i =1,\dots,n}$ such that 
$$
U(µ) = \prod_{i=1}^n \int_\mo \phi_i(x)µ(dx).
$$
A \textit{polynomial} is a real linear combination of monomials. It is clear that all polynomials are continuous with respect to the topology of weak convergence on $\mmo$, that the set of polynomials forms an algebra which separates points, and for all $µ$, there exists a polynomial which does not vanish on $µ$. Hence, by the Stone-Weierstrass Theorem, we have the following.
\begin{Theorem}\label{thm:poly}
For any $\eps > 0$ and $U: \mmuo \mapsto \R$ continuous for the weak convergence of measures, there exists a polynomial $P: \mmuo \mapsto \R$ such that
$$
\sup_{µ \in \mmuo}|U(µ)-P(µ)| \leq \eps.
$$
\end{Theorem}
Remark that, if we equip $\mmuo$ with the total variation distance, it is no longer compact, and thus unsuitable for the use of Stone-Weierstrass's Theorem.\\

\textbf{Approximation by Dirac masses}\\
Consider a compact subset $\mo$ of $\R^d$. There is a link between functions $U: \mpo \mapsto \R$ and \textit{symmetric functions} $U^N : \mo^N \mapsto \R$ that I now illustrate.

Given a permutation $\sigma: \{1,\dots,N\}\mapsto\{1,\dots,N\}$, it acts on $\mo^N$ through $\sigma (x) = (x_{\sigma(1)},\dots,x_{\sigma(N)})$. A function $U^N: \mo^N \mapsto \R$ is called symmetric if, for any permutation $\sigma$, $U^N(x) = U^N(\sigma(x))$. Because of this symmetry, $U^N(x)$ is in fact a function not of $x$, but say of $\{x_1;\dots;x_N\}$. This set can be represented by the \textit{empirical measure} of $x$ defined by $µ^N_{x}:=N^{-1}\sum_{i = 1}^N\delta_{x_i}$. Hence, we can view symmetric functions on $\mo^N$ as functions of the empirical measure. Thus, as $N$ grows large, $U^N$ should resemble a function $U: \mpo \mapsto \R$. Of course a sort of uniform continuity assumption has to be satisfied by the $U^N$ for this sort of limit to make sense in $\mathcal C (\mpo)$. The following result makes this idea precise. I refer to the notation section for a precise definition of a \textit{modulus of continuity}.

\begin{Theorem}\label{thm:ascoli}
Let $p \in [1,\infty)$, $(U^N)_{N \geq 1}$ be a sequence of symmetric functions each defined on $\mo^N$, $C > 0$ and $\omega$ a modulus of continuity such that, for all $N \geq 1$, $U^N: \mo^N \mapsto \R$ satisfies 
 $$
\forall x,y \in \mo^N,\quad |U^N(x) - U^N(y)| \leq \omega(\mathcal{W}_p(µ^N_x,µ^N_y)),
$$
 $$ \|U^N\|_\infty \leq C.$$
Then, extracting a subsequence if necessary, there exists a continuous function $U: \mpo\mapsto \R$, with modulus of continuity $\omega$, such that
$$
\lim_{N \to \infty}\sup_{x \in \mo^N}|U(µ^N_x)-U^N(x)| = 0.
$$
\end{Theorem}
\begin{proof}
The main idea is classical and consists in considering an extension of $U^N$ defined on the whole $\mpo$, which has the same modulus of continuity. Let $V^N : \mpo \mapsto \R$ be defined by
$$
V^N(µ) = \inf_{x \in \mo^N}\{U^N(µ^N_x) + \omega(\mathcal{W}_p(µ^N_x,µ))\}.
$$
The sequence $(V^N)_{N \geq 1}$ is equi-bounded on $\mpo$. It is also equi-continuous of modulus $\omega$. Indeed, let $µ,µ' \in \mpo$, and consider $x \in \mo^N$ which is $\eps$-optimal for $V^N(µ')$. Then 
$$\ba
V^N(µ) &\leq U^N(µ^N_x) + \omega(\mathcal{W}_p(µ,µ^N_x))\\
&\leq U^N(µ^N_x) + \omega(\mathcal{W}_p(µ',µ^N_x)) + \omega(\mathcal{W}_p(µ,µ^N_x)) - \omega(\mathcal{W}_p(µ^N_x,µ'))\\
& \leq V^N(µ') + \eps + \omega(\mathcal{W}_p(µ,µ') + \mathcal{W}_p(µ',µ^N_x)) - \omega(\mathcal{W}_p(µ^N_x,µ'))\\
& \leq V^N(µ') + \eps + \omega(\mathcal{W}_p(µ,µ'))
\ea
$$
where I used the sub-additivity of $\omega$. Hence, using Arzela-Ascoli Theorem, it follows that, up to a subsequence, $(V^N)$ converges uniformly to some function $U: \mpo\mapsto\R$ with modulus of continuity $\omega$. To conclude, it is sufficient to remark that $U^N$ and $V^N$ agree on $\{ µ^N_x | x \in \mo^N\}$.
\end{proof}

\textbf{Continuity of a function and continuity of its lift}\\
There is a natural link between topological properties of a function on $\mpprd$ and its lift.
\begin{Prop}\label{prop:samecont}
 Let $p \in [1,\infty)$, $U: \mpprd \mapsto \R$ and its lift $\mathcal U: \mathbb L^p \mapsto \R$. The function $U$ is usc (resp. lsc, resp. continuous) on $\mpprd$ if and only if $\mathcal U$ is usc (resp. lsc, resp. continuous) on $\mathbb L^p$.
 \end{Prop}
 \begin{proof}
 Let $p \in [1,\infty)$ and $U: \mpprd \mapsto \R$ be usc. Let $(X_n)_{n \geq 0}$ be converging toward some $X \in \mathbb L^p$. This convergence implies that 
 $$
 \mathcal W_p(\mathcal L(X_n),\mathcal L(X)) \leq \|X_n-X\|_p \underset{n \to \infty}{\longrightarrow} 0.
 $$
 Hence, by upper semi continuity of $U$, $\limsup_{n \to \infty} U(\mathcal L(X_n)) \leq  U(\mathcal L(X))$ and thus $\mathcal U$ is also usc.
 
 Conversely, let $\mathcal U$ be usc and consider $(µ_n)_{n \geq 0}$ converging toward $µ \in \mpprd$. Consider $X$ of law $µ$. By Theorem \ref{thm:proba}, for any $n \geq 0$, there exists $X_n$ of law $µ_n$ such that $\|X_n-X\|_p \leq \mathcal W_p(µ_n,µ) + \frac{1}{n+1}$. In particular, $(X_n)_{n \geq 0}$ converges toward $X$ in $\mathbb L^p$ and thus $\limsup_{n \to \infty}\mathcal U(X_n) \leq \mathcal U(X)$. Hence $U$ is usc as well. The lsc and the continuity properties can be proven in the same way.
 \end{proof}

\subsection*{Bibliographical comments}
There are several excellent textbooks on measure theory, such as Bogachev \citep{bogachev} for instance. More modern approaches, especially concerning Wasserstein spaces, can be found in the books of Ambrosio, Gigli and Savaré \citep{ags}, Villani \citep{villani} or Santambrogio \citep{santambrogio}. A statement and proof of Prokhorov's Theorem can be found in \citep{bogachev}, Theorem 8.6.2. The disintegration Theorem is proven by Dellacherie and Meyer in \cite{dellacherie-meyer}, III-70. This version of the Gluing Lemma is presented in Lemma 5.3.2 in \citep{ags}. Theorems \ref{thm:Wasserstein}, \ref{thm:glivenko}, \ref{thm:kantorovich} and \ref{thm:topowass} can be found in Theorems 6.18, 5.10 and 6.9 in \citep{villani}. Theorem \ref{thm:proba} is Lemma 5.23 in the book of Carmona and Delarue \citep{carmona2017probabilistic}. Theorems \ref{thm:poly} and \ref{thm:ascoli}, as well as the generalization of lifting techniques, are due to Lions and were presented in his lectures at Collège de France \citep{lions20072008}. Proposition \ref{prop:samecont} can be found in the work of Gangbo and Tudorascu \citep{gangbotudorascu}.

\newpage

\section{Variations on spaces of measures}\label{sec:variations}

Before giving precise notions of derivatives in spaces of measures, I want to introduce several typical variations of measures of interest, both in terms of applications and for theoretical purposes. It is typically along such variations $(m_t)_{t \in [0,1]}$ in $\mmo$ that one is interested in limits of the form 
\be\label{quotient2}
\lim_{t \to 0}\frac{U(m_t)-U(m_0)}{t},
\ee
for $U : \mmo \mapsto \R$.

As a definition, a \textit{variation} (or a curve, or a path) on a topological space $X$ is a path $(x_t)_{t \in [0,1]}\in\mathcal C([0,1],X)$. Usually, we would like to take limits in \eqref{quotient2} using some sort of chain rule, which should hint at restricting our attention to more regular paths than simply continuous ones. Because in several cases of interest such extra regularity will not hold, I will forget it for the moment.  As we shall see, there are several natural types of variations on spaces of measures, which can be quite different from each other. 

The aim of the present section is essentially to fix some vocabulary and to give some intuition to the reader unused to spaces of measures. The more informed reader is invited to skip or pass quickly over this section.

\subsection{Linear-like combinations}\label{sec:linearvar}
The first variation that I want to insist upon may be the most natural one, as it does not rely on any structure on $\mathcal O$, additional to the fact that it has a $\sigma$-algebra of course. This notion of variations relies uniquely on the structure of $\mmo$ as a normed vector space equipped with total variation.

A (straight) line is defined by two measures $µ,\nu \in \mmo$ as the set of points $\{µ + \lambda \nu, \lambda \in \R\}$, whereas the (straight) segment between $µ$ and $\nu$ in $\mmo$ is given by $\{(1-\theta)µ + \theta \nu, \theta \in [0,1]\}$. Of course, if $µ$ and $\nu$ are in $\mpo$, then so is the segment between them, thus providing us with a notion of variations on $\mpo$, and more generally on any convex subset of $\mmo$.

This notion of variation is compatible with the topology of total variation as, for any $µ,\nu \in \mmo$,
$$
\ba \R &\mapsto (\mmo, |\cdot|)\\
\lambda &\mapsto µ + \lambda \nu \ea\quad\text{ is continuous.}
$$

Note that if $\mo$ is a Polish space, $p\in [1,\infty)$, and $µ,\nu \in \mppo$, then $[0,1]\mapsto \mppo, \theta \mapsto (1-\theta)µ + \theta \nu$ is also continuous for $\mathcal W_p$.

Typical variations which are very different from straight lines are paths of the form $(\delta_{x(t)})_{t \in [0,1]}$, where $(x(t))_{t \in [0,1]}$ is a smooth path in $\mo$ (if it has a structure which allows us to talk about a smooth path). Recall that such paths are not continuous in total variation distance.%

\subsection{Action of functions}
Another natural type of variations on spaces of measures is given by the action of the semigroup $(\{f, f: \mathcal O \mapsto \mathcal O, f \text{ measurable}\},\circ)$ on $\mmo$ through the operation $f_\#µ$ for $f: \mathcal O\mapsto \mathcal O$ and $µ\in \mmo$.

Given a variation $(f_t)_{t \in [0,1]}$ on $(\mathcal C(\mo,\mo),\|\cdot\|_\infty)$ and $µ \in \mmo$, we can naturally consider the associated path $((f_t)_\#µ)_{t \in [0,1]}$ in $\mmo$. The previous map is not continuous in general with respect to total variation. Consider for instance the case of translation when $\mathcal O = \R$. The map $h \mapsto (\tau_h)_\#\delta_0$ is not continuous for the total variation distance. However, it behaves nicely with the weak topology. Indeed, if $\mo$ is a Polish space and $\phi \in \mathcal C_b(\mo)$ is uniformly continuous, then
$$
\lim_{t \to 0} \int_\mo \phi(f_t(x)) -\phi(f_0(x))µ(dx) = 0,
$$
since $\|f_t - f_0\|_\infty \to 0$ as $t \to 0$. The restriction to uniformly continuous functions is sufficient to obtain weak convergence, see for instance Theorem 3.A.5 in \citep{stroock2010probability}.\\

 Moreover, if we consider $µ \in \mpo$, still with $\mo$ a Polish space, then this type of action is quantifiable in terms of Wasserstein metrics. Indeed, for $p \in [1,\infty)$, $µ\in \mathcal P(\mo)$ and $f,g: \mo \mapsto \mo$, 
$$
\mathcal W_p(f_\#µ,g_\#µ) \leq \left(\mathbb E[d(f(X),g(X))^p]\right)^\frac1p \leq \|f-g\|_\infty,
$$
where $X \sim µ$.\\

There still remain variations out of reach for this formalism. Indeed, on $\mathcal P(\R)$, the path $((1-t)\delta_0 + \frac t2\delta_1 + \frac t2\delta_2)$ cannot be represented by means of action of functions as the image of $\delta_0$ by a map $f$ is $\delta_{f(0)}$. We shall focus later on such trajectories which split mass, but continue next by specifying certain variations given by actions of functions.

\subsection{Flows of ordinary differential equations and continuity equations}\label{subsec:ode}
A class of functions of particular interest is given by flows of ordinary differential equations (ODEs). Consider the case $\mathcal O = \T^d$ or $\mathcal O = \R^d$, a function $b: \R_+ \times \mo\mapsto \R^d$ and the ODE
\be\label{ode}
\dot X_t = b(t,X_t).
\ee
General smooth manifolds could also have been considered instead of $\mo$, but are omitted to keep the discussion notationally simple. If $b$ is continuous in time, and locally uniformly in time, Lipschitz in $x$, thanks to Cauchy-Lipschitz Theorem, \eqref{ode} defines a family of bi-Lipschitz homeomorphisms $(X_t^b)_{t \in \R}$ from $\mo$ into $\mo$, where $X^b_t(x)$ is the solution to \eqref{ode} at time $t$ with value $x$ at time $t = 0$. The family $(X_t^b(\cdot))_{t \in \R}$ acts on the space of measures through the continuity equation
\be\label{continuityequation}
\partial_t m + \text{div}(b \,m) = 0 \text{ in } (0,T)\times \mo,
\ee
meaning that if $m \in C([0,T],\mmo)$ (when $\mmo$ is equipped with the weak topology) is a weak solution to \eqref{continuityequation}, then, for all $t \geq 0$, $m_t = (X^\cdot_t)_\#m_0$.\\

In what follows, by definition, a \textit{weak solution} of a PDE of the form of \eqref{continuityequation} is an element $m \in \mathcal C([0,T], \mmo)$ such that, for any $t \in [0,T]$ and $\varphi \in \mathcal C^1_c([0,T]\times \mo)$, it holds that
$$
\int_0^t\int_\mo (-\partial_t \varphi(s,x) - \nabla_x \varphi(s,x) \cdot b(s,x))m_s(dx)ds = \int_{\mo}\varphi(0,x)m_0(dx) - \int_{\mo}\varphi(t,x)m_t(dx).
$$
We will mostly be interested in the case $m_0 \in \mpo$, which implies that $(X^\cdot_t)_\#m_0 \in \mpo$ for any time $t$. These variations on $\mpo$ enjoy natural regularity estimates with respect to Wasserstein metrics, just as it was the case in the previous section. For instance, it is easy to prove that if $m_0 \in \mppo$ for $p \in [1,\infty)$, then $(m_t)_{t \in [0,T]}$ is Lipschitz continuous for $\mathcal W_p$.

\subsection{Flows of stochastic differential equations and Fokker-Planck equations}
In the same spirit as the previous section, and keeping $\mo = \T^d$ or $\mo = \R^d$, we can associate to a stochastic differential equation (SDE) a natural evolution on $\mpo$. Consider a standard filtered probability space $(\Omega,\mathcal A, (\mathcal F_t)_{t \geq 0},\mathbb P)$, $(W_t)_{t \geq 0}$ a standard $k$ dimensional Brownian motion on it, and the SDE
\be\label{sde}
dX_t = b(t,X_t)\,dt + \sigma(t,X_t)\,dW_t
\ee
where $b:\R_+\times \mo \mapsto \R^d$ and $\sigma : \R_+\times \mo \mapsto M_{d\times k}(\R)$. If $(X_t)_{t \geq 0}$ is a weak solution to \eqref{sde} (i.e. a process on possibly another probability space such that for a $k$ dimensional Brownian motion, the integral representation of \eqref{sde} holds almost surely, together with some integrability conditions), then the evolution of its law $(m_t)_{t \geq 0}$ defines a natural path on $\mpo$. Under appropriate assumptions on $b$, $\sigma$ and the integrability of $X_0$, we can also characterize $(m_t)_{t \geq 0}$ as the unique weak solution to the following Fokker-Planck equation
\be\label{fp1}
\partial_t m - \frac12\sum_{i,j=1}^d\partial_{ij}((\sigma \sigma^T)_{ij}m) + \text{div}(b\,m) = 0 \text{ in } (0,\infty)\times \mo.
\ee
Weak solutions to \eqref{fp1} are defined by testing the previous equation against test functions as previously, but the test functions have to be taken with a $\mathcal C^2$ regularity in the space variable $x$ here.\\

Note that in this stochastic case, the regularity of $(m_t)_{t \geq 0}$ can still be measured in terms of Wasserstein metrics. Continuity can be obtained using quite standard stochastic analysis techniques. It is also a simple exercise to check that such variations are not Lipschitz continuous in general. This occurs for example in the case $\mo = \R$, $m_0 = \delta_0$, $b = 0$, $k =1$ and $\sigma =1$.
\begin{Rem}\label{rem:fpvar}
Other interpretations of the Fokker-Planck equation exist, in which it is derived by taking the limit $N \to \infty$ of the empirical measure of a family $(X^1_t,\dots,X^N_t)_{t \geq 0}$ of independent weak solutions to \eqref{sde}. Using this type of interpretation, by tuning properly the $\frac 1N$ in the empirical measure, we can obtain more general paths in $\mathcal M_+(\mo)$ instead of simply $\mpo$.
\end{Rem}

\subsection{Variations associated to couplings and geodesics of optimal transport}

In these last two sections, we consider variations which are not defined as images of functions but rather through some Lagrangian point of view. I restrict the following to probability measures over a convex Polish set $\mo$.
\begin{Rem}
The Lagrangian point of view is a general and quite ill-defined terminology which refers to the idea of looking at the trajectories of particles rather than at the trajectory of their density (or of the measure which describes their distribution). The latter is referred to as an Eulerian point of view. I come back to this distinction below.
\end{Rem} 
 
The most basic example of a Lagrangian variation consists in following a coupling. Given $\gamma \in \mathcal P(\mo^2)$ with marginals $µ$ and $\nu$, thanks to the convexity of $\mo$, we can associate to $\gamma$ a natural path between $µ$ and $\nu$. Indeed, for $x,y \in \mo$, we consider the path $((1-t)x + ty)_{t \in [0,1]}$. The natural Lagrangian variation associated to $\gamma$ is then the path $(((1-t)\pi_1 + t \pi_2)_\#\gamma )_{t \in [0,1]}$. At the level of random variables, this method is quite linear. Indeed, consider a couple of random variables $(X,Y)\sim \gamma$. Then, the path we are considering is simply the one given by $(\mathcal L( (1-t)X + t Y))_{t \in [0,1]}$.\\

\textbf{The case of optimal couplings}
When the coupling is optimal for some Wasserstein distance, we usually refer to it as a \textit{geodesic} between the two measures, in the associated Wasserstein space.
\begin{Def}
A constant speed geodesic (or simply geodesic in what follows) of $\mppo$ is a path $(µ_t)_{t \in [a,b]}$, with $a < b$ such that $$\forall t,s \in [a,b],\quad\mathcal W_p(µ_t,µ_s) = \frac{|t-s|}{b-a}\mathcal W_p(µ_a,µ_b).$$
\end{Def}
In particular, a geodesic is a path of minimal length between its extreme points. For instance, if $\gamma$ is an optimal coupling in $\mathcal W_p(µ,\nu)$ for $µ,\nu \in \mppo$, then it is a simple exercise to show that $µ_t := ((1-t)\pi_1 + t\pi_2)_\#\gamma$ is a geodesic in $\mppo$.\\

To consider more general trajectories than straight lines for the particles, it is natural to proceed as in the next case.

\subsection{Smooth Lagrangian flows and other continuity equations}\label{sec:lagrange}
I conclude this list of variations by considering extensions of the continuity equations presented above. This extension is mainly of interest for the optimal control problems studied in Part III. I consider here $\mo$ to be either $\T^d$ or $\R^d$.\\

There are two natural equivalent ways to consider general Lagrangian variations. 

First, we can consider a random variable $X_0$ on a probability space $(\Omega,\mathcal A,\mathbb P)$ with values in $\mo$, of law $µ$, as a way to keep track of the labels of the particles (the labels are here the elements of $\Omega$). Then the trajectories are represented by a stochastic process $(X_t)_{t \in [0,1]}$ which is almost surely continuous, with initial condition $X_0$, on the same probability space. In this interpretation, $X_t(\omega)$ is the position of the particle labeled $\omega$ at time $t$.  

Or, we can directly work with a probability measure $\eta \in \mathcal P(\mathcal C([0,1],\mo))$, which is none other than the law of the process $(X_t)_{t \in [0,1]}$. The first approach is clearly a lift of the second one, which can thus appear as more canonical.\\

 By definition, a \textit{Lagrangian path} $\eta$ is an element of $\mathcal P(\mathcal C([0,1],\mo))$. It is said to be smooth, $\mathcal C^k$,... if it is concentrated on such curves, or equivalently if $(X_t)_{t \in [0,1]}$ is smooth, $\mathcal C^k$,... almost surely.\\
 
 A Lagrangian path yields a variation in $\mpo$ by considering $((e_t)_\#\eta)_{t \in [0,1]}$, where I recall that $e_t$ is the evaluation mapping at time $t$. It is the canonical map from $\mathcal P(\mathcal C([0,1],\mo))$ into $\mathcal C([0,1],\mpo)$.\\

This Lagrangian perspective is to be thought of in contrast to the Eulerian one, which consists in looking at the continuity equation \eqref{continuityequation} with velocity $b$. When the associated vector field $b$ is smooth, there is a unique Lagrangian flow associated to it because of the uniqueness of the solution to the associated ODE \eqref{ode}.\\

When the velocity field is not sufficiently regular, the links between Eulerian and Lagrangian approaches become much more interesting. In such cases, the lack of uniqueness of solutions to \eqref{continuityequation} might raise some issues but a representation of solutions through Lagrangian curves is still possible. It is often referred to as Ambrosio's superposition principle and can be stated as follows.

\begin{Theorem}\label{ambrosio}
If $(m_t)_{t \in [0,1]}$ is a continuous curve in $\mpo$ for the weak topology and is a weak solution to \eqref{continuityequation} for a certain drift $b$ such that, for some $p\in (1,\infty)$ 
\be\label{cond:b}
\int_0^1\int_\mo |b(t,x)|^pm_t(dx)dt < \infty,
\ee
then there exists a process $(X_t)_{t \in [0,1]}$ such that $m_t = \mathcal L(X_t)$ for all $t$ and, almost surely in $\omega \in \Omega$
\be\label{eds}
X_t = X_0 + \int_0^tb(s,X_s)ds, \text{ for every } t \in [0,1].
\ee
Conversely, if $(X_t)_{t \in [0,1]}$ is an almost sure solution to \eqref{eds} and $\int_0^1\mathbb{E}[|b(s,X_s)|]ds < \infty$, then $(\mathcal{L}(X_t))_{t \in [0,1]}$ is a weak solution to \eqref{continuityequation} with drift $b$.
\end{Theorem}
Note the asymmetry of the assumptions with respect to the growth exponent $p$.\color{black}\\

The previous result is quite helpful, in particular when combined with the following result on \textit{absolutely continuous curves}. By definition, a curve $(x_t)_{t \in [0,1]}$ in a metric space $(X,d)$ is $p$-absolutely continuous if there exists a function $\lambda \in L^p([0,1],\R)$ such that 
$$
\forall 0 \leq t \leq s \leq 1, \quad d(x_t,x_s) \leq \int_t^s\lambda(u)du.
$$ 
\begin{Theorem}\label{thm:abs}
Let $p\in (1,\infty)$. Then $(m_t)_{t \in [0,1]}$ is an absolutely continuous curve in $\mppo$ if and only if there exists a measurable vector field $b : [0,1]\times \mo\mapsto \R^d$ such that \eqref{cond:b} holds and $(m_t)_{t \in [0,1]}$ is a weak solution to the continuity equation with velocity field $b$.
\end{Theorem}

The combination of the two previous results allows for a minimal representation of certain Lagrangian flows through a solution to \eqref{continuityequation} for a well chosen velocity $b$. The minimality is here with respect to the extra randomness carried in the Lagrangian formulation. Consider two processes $(X_t)_{t \in [0,1]}$ and $(Z_t)_{t \in [0,1]}$, such that almost surely, 
$$
X_t = \int_0^t Z_s ds + X_0.
$$
Assume that for some $p \in (1,\infty)$
$$
\int_0^1(\mathbb E[|Z_s|^p])ds < \infty.
$$
Then, the law $\mathcal L(X_t)$ is an absolutely continuous curve in $\mppo$. Thus, it can be represented by a continuity equation and the associated drift is formally given by $b(t,x) = \mathbb E[Z_t|X_t = x]$. But this can be represented by a process $(\tilde X_t)_{t \in [0,1]}$ solution to \eqref{eds}. By construction, $(\tilde X_t)_{t \in [0,1]}$ and $(X_t)_{t \in [0,1]}$ have the same marginal laws, i.e. $\mathcal L(X_t) = \mathcal L(\tilde X_t)$ for all $t$, but formally the former carries some minimality.

For control purposes (which shall be made clear in Part III), I also state another (rather trivial once seen as an extension of the previous) representation result, in which we do not have this minimality. Or in other words, the \emph{randomness} of the velocity is kept and somehow made precise.
\begin{Prop}\label{prop:ambrosiov2}
Let $(X_t,Z_t)_{t \in [0,1]}$ be an $\mo\times \R^d$ valued process such that, 
$$
\int_0^1\mathbb E[|Z_t|] dt < \infty
$$
 and, almost surely, for all $t \in [0,1]$
\be\label{cond16}
X_t = X_0 + \int_0^tZ_sds.
\ee
Then, there exists a measurable set $K \subset \R^N$, a measurable map $b: [0,1]\times \mo \times K \mapsto \R^d$ and a continuous path $(\tilde m_t)_{t \in [0,1]}$, valued in $\mathcal P(\mo\times K)$ such that 
\be\label{cond15}
((\pi_1,b(t,\pi_1,\pi_2))_\#\tilde m_t)_{t \in [0,1]} = (\mathcal L(X_t,Z_t))_{t \in [0,1]}
\ee
 and $(b,\tilde m)$ is a weak solution to the continuity equation
\be\label{ceK}
\partial_t \tilde m + \emph{div}_x(b\,\tilde m) = 0 \text{ in } (0,1)\times \mo\times K.
\ee
Conversely, given a measurable set $K\subset \R^N$, $p \in (1,\infty)$ and a solution $(b,\tilde m)$ to \eqref{ceK} satisfying 
$$
\int_0^1\int_{\mo\times K}|b(t,x,k)|^p\tilde m_t(dx,dk)dt < \infty,
$$
then, there exists a process $(X_t,Z_t)_{t \in [0,1]}$ such that \eqref{cond15} and \eqref{cond16} hold.
\end{Prop}
\begin{proof}
The proof is rather simple. On the first side, we denote by $(\Omega,\mathcal A,\mathbb P)$ the standard probability space on which the process $(X_t,Z_t)_{t \in [0,1]}$ is defined. Without loss of generality, we can assume that $(\Omega,\mathcal A,\mathbb P)= ([0,1]^N,\mathcal B, Leb)$. Taking $K = \Omega$, $b(t,x,k) = Z_t(k)$ for $(t,x,k) \in [0,1]\times \mo\times K$ and $\tilde m_t = \mathcal L(X_t,Id_\Omega)$ yields the solution, thanks to the second part of Theorem \ref{ambrosio}.

On the other side, we use the first part of Theorem \ref{ambrosio} to conclude.
\end{proof}

The paths considered in the previous result are natural if one thinks about controlled trajectories, i.e. trajectories in which the velocity is controlled. In game theory, this is equivalent to a mixed strategy, while in optimal control, it is often referred to as a relaxed control.
\begin{Rem}\label{rem:lc}
The main interest of the previous result is to characterize, when we are considering solutions to continuity equations, pairs $(b,\tilde m)$ which can be represented naturally through stochastic processes. This characterization is useful to prove equivalent formulations of the same optimal control problem.\\

 In particular, this makes apparent the need of a space $K$ to construct the randomization of speed upon. To highlight the difference between considering a process $(X_t,Z_t)$ or its law, observe that there are such laws for which it is not clear how to construct an associated suitable stochastic process. Take for instance $d = 1$ and a law given as $\mathcal L(X_t,Z_t) = m_t(dx)\frac{1}{\sqrt{2\pi}}e^{-\frac{z^2}{2}}dz$ for $(m_t)_{t \in [0,1]}$ a nice curve in $\mathcal P_p(\R)$. Then, Kolmogorov existence Theorem might yield a process $(X_t,Z_t)_{t \in [0,1]}$ such that $\mathcal L(X_t,Z_t)(dx,dz)= m_t(dx)\psi_t(x)(dz)$ for all $t > 0$, but such a process will fail to satisfy the relation
$$
X_t = X_0 + \int_0^tZ_sds,
$$
as this would require the construction of a process $(Z_t)_{t \in [0,1]}$ with too much independence.
\end{Rem}

\subsection*{Bibliographical comments}
ODEs, continuity equations and Fokker-Planck equations are now classical objects of mathematics and I will not give a history of their study. ODEs and continuity equations regained interest because of the DiPerna-Lions theory presented notably in \citep{dipernalions}. Recent developments in this line on second order equations can be found in the book of Le Bris and Lions \citep{lebris}. 
The link between absolutely continuous curves and the continuity equation has been key in the study of gradient flows in the book of Ambrosio, Gigli and Savaré \citep{ags}, notably following the seminal work of Otto in \citep{ottoporous}. Theorems \ref{ambrosio} and \ref{thm:abs} are respectively Theorems 8.2.1 and 8.3.1 in \citep{ags}, where detailed proofs are given. I proposed the formulation of Proposition \ref{prop:ambrosiov2} in \citep{bertucci2025tangent} to study optimal control problems in $\mpprd$. I also refer to the work of Jimenez, Marigonda and Quincampoix \citep{jimenez2023dynamical} for other related developments.

\newpage

\section{The flat or vertical derivative on spaces of measures}\label{sec:flat}
This section is concerned with a notion of derivative aligned with the vector space structure of $\mmo$, and more generally, with the inherited one of $\mpo$ as a convex set of a Banach space. What follows is a standard presentation of the concepts of Fréchet or Gateaux differentiability in Banach spaces, adapted to the present case. The objective is here to fix notation and ideas, as all the developments of this section are standard.\\

There are several ways to endow $\mmo$ (or $\mpo$) with a Banach like structure. In this section, I shall focus on the simplest way to do so, that is to use the norm in total variation $|\cdot|$ and look at the Banach space $(\mmo,|\cdot|)$. Other structures could have been investigated, and I refer to Section \ref{sec:W1} for such an example. Because I focus on the norm in total variation, in this section $\mo$ is simply a measurable space.\\

Several computations involving vertical derivatives are presented in Section \ref{sec:computations}.

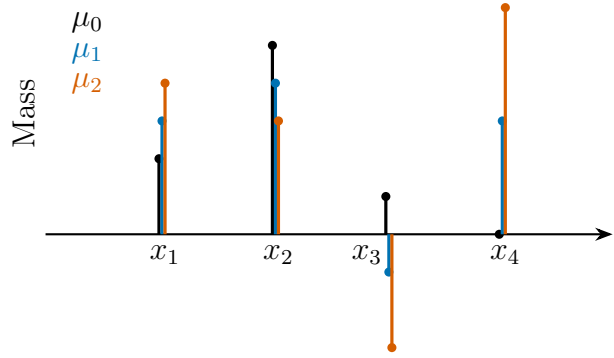
\begin{figure}[htbp]
  \centering
  \begin{tikzpicture}
    \definecolor{colorT1}{RGB}{0, 114, 178}   
    \definecolor{colorT2}{RGB}{213, 94, 0}    

    \draw[->, >=Stealth, thick] (0,0) -- (7.5,0);

    \draw[line width=1.2pt, black] (1.5,0) -- (1.5, 1.0);
    \fill[black] (1.5, 1.0) circle (0.06);
    \draw[line width=1.2pt, colorT1] (1.54,0) -- (1.54, 1.5);
    \fill[colorT1] (1.54, 1.5) circle (0.06);
    \draw[line width=1.2pt, colorT2] (1.58,0) -- (1.58, 2.0);
    \fill[colorT2] (1.58, 2.0) circle (0.06);

    \draw[line width=1.2pt, black] (3.0,0) -- (3.0, 2.5);
    \fill[black] (3.0, 2.5) circle (0.06);

    \draw[line width=1.2pt, colorT1] (3.04,0) -- (3.04, 2.0);
    \fill[colorT1] (3.04, 2.0) circle (0.06);

    \draw[line width=1.2pt, colorT2] (3.08,0) -- (3.08, 1.5);
    \fill[colorT2] (3.08, 1.5) circle (0.06);

    \draw[line width=1.2pt, black] (4.5,0) -- (4.5, 0.5);
    \fill[black] (4.5, 0.5) circle (0.06);

    \draw[line width=1.2pt, colorT1] (4.54,0) -- (4.54, -0.5);
    \fill[colorT1] (4.54, -0.5) circle (0.06);

    \draw[line width=1.2pt, colorT2] (4.58,0) -- (4.58, -1.5);
    \fill[colorT2] (4.58, -1.5) circle (0.06);

    \draw[line width=1.2pt, black] (6.0,0) -- (6.0, 0.0);
    \fill[black] (6.0, 0.0) circle (0.06);

    \draw[line width=1.2pt, colorT1] (6.04,0) -- (6.04, 1.5);
    \fill[colorT1] (6.04, 1.5) circle (0.06);

    \draw[line width=1.2pt, colorT2] (6.08,0) -- (6.08, 3.0);
    \fill[colorT2] (6.08, 3.0) circle (0.06);

    \node[fill=white, inner sep=1pt, anchor=north] at (1.58, -0.1) {$x_1$};
    \node[fill=white, inner sep=1pt, anchor=north] at (3.08, -0.1) {$x_2$};
    \node[fill=white, inner sep=1pt, anchor=north] at (4.25, -0.1) {$x_3$};
    \node[fill=white, inner sep=1pt, anchor=north] at (6.08, -0.1) {$x_4$};

    \node[anchor=south west, rotate=90] at (0, 1.0) {Mass};

    \node[anchor=west, text=black] at (0.2, 2.8) {$\mu_0$};
    \node[anchor=west, text=colorT1] at (0.2, 2.4) {$\mu_1$};
    \node[anchor=west, text=colorT2] at (0.2, 2.0) {$\mu_2$};

  \end{tikzpicture}
  \caption{Typical variation of interest in this section: $µ_t=µ + t\nu$, for some $µ,\nu$ which are here finitely supported.}
  \label{fig:signed_measure_variation}
\end{figure}

\subsection{On the normed vector space of measures}

Consider a function $U: \mathcal M(\mathcal O) \mapsto \R$. We would like, quite classically, to say that $U$ is (Fr\'echet) differentiable at $µ_0 \in \mmo$ if there exists a bounded linear operator $A: \mmo \mapsto \R$ such that for any $\nu \in \mmo$, 
$$
\lim_{|\nu| \to 0}\frac{U(µ_0 + \nu) - U(µ_0) - A(\nu)}{|\nu|} = 0.
$$
In the previous the boundedness of $A$ is with respect to the total variation norm. While working with this generality could certainly lead to interesting developments, I shall focus here on a more restrictive definition by imposing that $A$ can be represented by an element of $\mathcal B(\mo)$, the set of bounded measurable functions on $\mathcal O$. 
\begin{Def}\label{def:flat}
A function $U: \mmo \mapsto \R$ is said to be vertically, flat or V-differentiable at $µ_0 \in \mmo$ if there exists $\varphi \in \mathcal B(\mo)$ such that, for $\nu \in \mmo$
$$
\lim_{|\nu| \to 0}\frac{U(µ_0 + \nu) - U(µ_0) - \langle \varphi,\nu\rangle}{|\nu|} = 0.
$$
In such a case, we set $\varphi = \nabla_µ U(µ_0)$ or simply $\varphi = \nabla U(µ_0)$. The function $\nabla_µ U(µ_0)$ is called the flat or V-derivative of $U$ at $µ_0$.
\end{Def}
\begin{Rem}
Note that the V-derivative, if it exists, is a bounded function and not an element of $L^\infty(\mo,\R)$. Working in $L^\infty$ is extremely dangerous here as, for $\phi \in L^\infty(\mo,\R)$ and $µ\in \mmo$, the product $\langle \phi,µ\rangle$ is not well defined.
\end{Rem}
\begin{Def}\label{defC1}
A function $U: \mmo \mapsto \R$ is said to be $\mathcal C^1$ if it is V-differentiable everywhere and if the map $\nabla U: (\mmo,|\cdot|) \mapsto (\mathcal B(\mo),\|\cdot \|_\infty)$ is continuous.
\end{Def}
This Fr\'echet, flat or V-derivative enjoys natural properties on which I do not comment much, and only recall some of them. 

\begin{Prop}
Let $U: \mmo \mapsto \R$ be V-differentiable at $µ \in \mmo$. Then $U$ is continuous at $µ$ for $|\cdot |$.
\end{Prop}
\begin{Prop}
Given a map $m: [-1,1] \mapsto (\mmo,|\cdot|)$, differentiable at $t = 0$, and a map $U : \mmo \mapsto \R$, V-differentiable at $m(0)$, the composition $t \mapsto U(m(t))$ is differentiable at $t = 0$ and
\[
\frac{d}{dt}U(m(t))|_{t = 0} = \int_{\mathcal O}\nabla U(m(0))(x)m'(0)(dx).
\]
\end{Prop}

\subsection{Vertical derivatives on the space of probability measures} 
In several cases, we are given a function $U : \mpo \mapsto \R$. The notion of its V-derivative is then slightly more complex than what we had in the previous section. The only point is here to treat the fact that $\mpo$ is not a normed vector space, but only a convex subset of it, whose interior is empty. Once again in this section, we shall endow $\mpo$ with the total variations distance.

\subsubsection{Gateaux like derivative}
The most common approach in the literature has been to look at a sort of Gateaux derivative, namely by replacing usual variations in the vector space by convex combinations.
\begin{Def}\label{def:998}
A function $U: \mpo \mapsto \R$ is said to be Gateaux differentiable at $µ_0 \in \mpo$ if there exists $\varphi \in \mathcal B(\mo)$ such that, for any $\nu \in \mpo$
\be\label{Udif}
\lim_{\theta \to 0^+}\frac{U((1-\theta)µ_0 + \theta \nu) - U(µ_0) }{\theta} = \langle \varphi,\nu - µ_0\rangle.
\ee
\end{Def}

In this case, $\varphi$ is integrated against the difference $\nu-µ_0$, and not just $\nu$. This implies that if $\varphi$ satisfies \eqref{Udif} for any $\nu \in \mpo$, then the same holds for $\varphi + c$ for any $c \in \R$. Quite often, a convention is then arbitrarily chosen to define the differential of $U$.

\begin{Def}\label{def:gateaux}
Let $U: \mpo \mapsto \R$ be Gateaux differentiable at $µ_0 \in \mpo$. We denote by $\nabla^g_µ U(µ_0)$ or $\nabla^g U(µ_0)$ the unique element $\varphi$ satisfying \eqref{Udif} for any $\nu \in \mpo$, such that $\langle \varphi,µ_0\rangle = 0$. 
\end{Def}
\begin{Rem}
The previous is a correct definition since if two such elements $\varphi_1,\varphi_2$ exist, then, for all $\nu \in \mpo$, $\langle \varphi_1-\varphi_2,\nu-µ_0\rangle = 0$, which implies the equality up to a constant.
\end{Rem}
Note that this choice of normalization might have unpleasant consequences, in particular when linking some convexity properties of a function to the sign of its second order derivatives, see 2.7.3 in \citep{jakobsen2025master}. As usual Definition \ref{def:998} is much less restrictive than Fr\'echet derivative type definitions. For instance, it is well known that it does not imply continuity (even in finite dimension).

\subsubsection{Fr\'echet like derivative}
In order to establish chain rules, or more generally linear approximations of a map $U$ around a point of differentiability $µ$, it is much more natural to work with the following notion of (Fréchet) differentiability.
\begin{Def}
A function $U: \mpo \mapsto \R$ is Fr\'echet or V-differentiable at $µ_0 \in \mpo$ if there exists a bounded measurable $\varphi : \mathcal O \mapsto \R$ such that, with $\nu \in \mpo$,
$$
\lim_{|\nu - µ_0|\to 0}\frac{U(\nu) - U(µ_0) - \langle \varphi,\nu-µ_0\rangle}{|\nu - µ_0|} = 0.
$$
If this element $\varphi$ is such that $\langle \varphi,µ_0\rangle = 0$, we set $\varphi = \nabla_µ U(µ_0)$ or $\varphi = \nabla U(µ_0)$.
\end{Def}
I recall standard properties of such a derivative.
\begin{Prop}
Let $U: \mpo \mapsto \R$ be V-differentiable at $µ_0 \in \mpo$. Then,
\begin{enumerate}[(i)]
\item $U$ is Gateaux differentiable at $µ_0$ and the two derivatives coincide,
\item $U$ is continuous at $µ_0$ (for the total variation distance),
\item for a map $m: [-1,1] \mapsto (\mpo,|\cdot|)$, differentiable at $0$ such that $m(0) = µ_0$, $t \mapsto U(m(t))$ is differentiable at $0$ and
$$
\frac{d}{dt}U(m(t))|_{t = 0} = \int_\mo \nabla U(m(0))(x)m'(0)(dx).
$$
\end{enumerate}
\end{Prop}
\begin{Rem}
The notion of Fr\'echet derivative  allows one to  obtain Taylor-like expansions of $U(m(t))$ for more paths $m(t)$ than simply affine ones, which would have been the case for Gateaux derivative. However, note that such generality is still quite restrictive as paths of the form $(\delta_{x(t)})_{t\in [-1,1]}$ for $(x(t))_{t \in [-1,1]}$ a $C^1$ curve in $\mathcal O$ are still prohibited.
\end{Rem}

The following result gives a sufficient condition on the Gateaux derivative of a function to state that it is equal to its Fr\'echet derivative, as is always the case in Banach spaces.

\begin{Prop}
Let $U: \mpo \mapsto \R$ be Gateaux differentiable everywhere and such that 
$$
\ba
(\mpo,|\cdot|) &\mapsto (\mathcal B(\mo),\|\cdot\|_\infty)\\ 
µ\quad &\mapsto \quad \nabla^g U(µ)
\ea
$$
is continuous. Then $U$ is $V$-differentiable everywhere.
\end{Prop}

\subsection{Functions valued in a normed vector space}
The previous presentation is still meaningful when working with functions $U : \mmo \mapsto E$, for $(E,\|\cdot\|_E)$ a separable Banach space. The case of functions defined on $\mpo$ is also quite similar. In such a context, the corresponding notion of Fréchet differential is given by
\begin{Def}
A function $U: \mmo \mapsto E$ is said to be Fr\'echet, vertically or V-differentiable at $µ_0 \in \mmo$ if there exists a strongly measurable bounded function $\varphi : \mathcal O \mapsto E$ such that
$$
\lim_{|\nu| \to 0}\frac{\|U(µ_0 + \nu) - U(µ_0) - \langle \varphi,\nu\rangle\|_E}{|\nu|} = 0.
$$
In such a case, we set $\varphi = \nabla_µ U(µ_0)$ or simply $\varphi = \nabla U(µ_0)$. The function $\nabla_µ U(µ_0)$ is called the flat or Fr\'echet derivative of $U$ at $µ_0$.
\end{Def}
In the previous, the integral in $\langle \varphi,\nu\rangle$ is defined in the Bochner sense and requires the notion of strong measurability. I do not enter into such developments but refer to II.2 in \citep{diestel} instead. Then, what has been previously seen can be adapted to this case. An important case is when $E$ is a function space over a space $X$. In this case, we write $U(µ)(x)$ for the evaluation at $x \in X$ of $U(µ)$ for $µ \in \mmo$. The derivative is then denoted through $
\nabla_µ U(µ)(x)(y)$ for $µ \in \mmo, x \in X, y \in \mo$. The order between $x$ and $y$ clearly matters (they do not necessarily live in the same space, even if in most applications in MFG, $X = \mo$). With this convention, we have for instance for all $µ \in \mmo, x \in X, y \in \mo$,
$$
\left(\lim_{\eps \to 0} \frac{U(µ + \eps \delta_y) - U(µ)}{\eps} \right)(x) = \nabla_µ U(µ)(x)(y).
$$
This might seem confusing but it is the convention used in most of the literature.

\subsection*{Bibliographical comments}
The previous standard developments can be found in numerous textbooks, see for instance the one of Brezis \citep{brezis}. It is probable that mathematicians have considered this notion of differential in the special case of measures or probability measures for quite some time, and it is very difficult to give a starting point. An influential early work is the approach used at the formal level in the so-called Otto calculus \citep{ottoporous} as an intermediate computational step. It has also been extensively studied to address MFG master equations, see the course of Lions \citep{lions20072008}, the book of Cardaliaguet, Delarue, Lasry and Lions \citep{cdll} or the book of Carmona and Delarue \citep{carmona2017probabilistic} for instance. 
\newpage

\section{Horizontal derivatives}\label{sec:horiz}
\subsection{Preliminary discussion}
In the previous section, I only used the fact that $\mathcal O$ is a measurable space. In most applications, we are interested in cases in which $\mathcal O$ has much more structure, e.g. $\R^d,\T^d$, a domain of $\R^d$ or more generally smooth manifolds. Thus, when trying to use topological properties of the underlying space, we are led to work with the topology of weak convergence on the space of measures.

When considering a more geometrical notion of derivative, aligned with the Wasserstein distances, it is much simpler to work on $\mpo$ rather than on $\mmo$, and the latter case is not presented in this book. Of course the normalization to $1$ is quite arbitrary and it is the non-negativity and the fixed mass which are of interest.

This section is organized as follows: first I introduce the main concept of derivative which is well-suited for all the variations considered above which are not linear combinations. I will then spend some time justifying this notion, making various links with more standard notions, and then present computations.

I end this section with the important case of a derivative which is somehow between the vertical one of the previous section and the horizontal one that I am about to present.\\

In all that follows in this section, $\mo$ is the closure of a smooth convex domain of $\R^d$, equipped with the Euclidean distance, and $p \in (1,\infty)$.

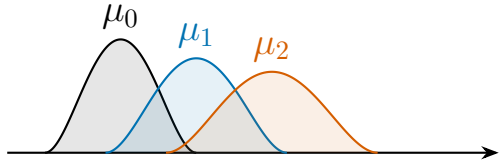
\begin{figure}[htbp]
  \centering
  \begin{tikzpicture}[
      declare function={
        mu(\x) = 1.5 * (1 - \x*\x) * (1 - \x*\x);
      }
    ]

    \definecolor{colorT0}{RGB}{0, 0, 0}       
    \definecolor{colorT1}{RGB}{0, 114, 178}   
    \definecolor{colorT2}{RGB}{213, 94, 0}    

    \draw[->, >=Stealth, thick] (-1.5, 0) -- (5.0, 0);

    \draw[thick, colorT0, fill=colorT0, fill opacity=0.1, domain=-1:1, samples=100]
      plot ({\x}, {mu(\x)});
    \node[above, text=colorT0, font=\large] at (0, 1.5) {$\mu_0$};

    \draw[thick, colorT1, fill=colorT1, fill opacity=0.1, domain=-1:1, samples=100]
      plot ({\x*1.2 + 1}, {mu(\x) / 1.2});
    \node[above, text=colorT1, font=\large] at (1, 1.25) {$\mu_1$};

    \draw[thick, colorT2, fill=colorT2, fill opacity=0.1, domain=-1:1, samples=100]
      plot ({\x*1.4 + 2}, {mu(\x) / 1.4});
    \node[above, text=colorT2, font=\large] at (2, 1.07) {$\mu_2$};

  \end{tikzpicture}
  \caption{Typical variation of interest in this section: evolution of a measure $\mu_t := (Id + t\phi)_\#\mu$ starting from an initial compactly supported density $\mu_0$. The vector field $\phi(x)$ induces both a translation and a dilation of the support.}
  \label{fig:pushforward_evolution}
\end{figure}

\subsection{Main definition}
The main notion is the following.
\begin{Def}\label{def:geoD}
A function $U: \mathcal P_p(\mathcal O) \mapsto \R$ is horizontally or H-differentiable at $µ \in \mathcal P_p(\mathcal O)$, if there exists $\phi \in L^{p'}((\mathcal O,µ),\R^d)$ and a modulus of continuity $\omega$ such that, for any $\nu \in \mathcal P_p(\mathcal O), \gamma \in \Pi(µ,\nu)$
\be\label{def:DU}
\left|U(\nu) - U(µ) - \int_{\mo^2}\phi(x)\cdot(y-x)\gamma(dx,dy)\right| \leq \omega(C_p(\gamma))C_p(\gamma).
\ee
We set $\phi= D_µU(µ)$ or simply $DU(µ)$ sometimes and call it the $H$-derivative. We say that $U$ is H-differentiable if it is H-differentiable everywhere.
\end{Def}

Several remarks are in order after this definition:
\begin{itemize}
\item As we shall see, at most one function $\phi$ satisfies the requirements.
\item The integral term is well defined since, using H\"older's inequality, we find that
$$
\ba
 \int_{\mathcal O^2}\phi(x)\cdot(y-x)\gamma(dx,dy)&\leq \left( \int_{\mathcal O^2}|\phi(x)|^{p'} \gamma(dx,dy)\right)^\frac{1}{p'}\left( \int_{\mathcal O^2}|x-y|^{p} \gamma(dx,dy)\right)^\frac{1}{p}\\
 &=\|\phi\|_{L^{p'}_µ}C_p(\gamma).
 \ea
$$
\item The set $L^{p'}_µ(\mathcal O,\R^d)$ might be quite uninteresting sometimes, making the present notion difficult to use (e. g. if $µ$ is a Dirac mass, then $L^{p'}_µ(\mathcal O,\R^d)$ is isomorphic to $\R^d$ once it has been quotiented by the equality $µ$ almost everywhere).
\item I somehow drop the $p$ in the notation. However confusing it may seem, the notation is sufficiently complex as is, and so we shall have to be careful with it !
\end{itemize}
As expected, differentiability implies continuity.
\begin{Prop}
If $U$ is H-differentiable at $µ$, then it is continuous at $µ$ (for the $\mathcal W_p$ distance).
\end{Prop}
\begin{proof}
Let $(µ_n)_{n \geq 0}$ be a sequence in $\mppo$ such that $\mathcal W_p(µ_n,µ) \to 0$ as $n \to \infty$. Then, for all $n \geq 0$, consider an optimal coupling $\gamma_n$ for $\mathcal W_p(µ,µ_n)$. Evaluating \eqref{def:DU} on $\gamma_n$ yields the result, since 
$$
\left|\int_{\mo^2}D_µU(µ)(x)\cdot(y-x)\gamma_n(dx,dy) \right| \leq \|\phi\|_{L^{p'}_µ}\mathcal W_p(µ_n,µ) \underset{n \to \infty}{\to} 0.
$$
\end{proof}
The following justifies that the notation $D_µU(µ_0)$ is meaningful.
\begin{Prop}\label{prop:deruniqbound}
Two functions $\phi_1$ and $\phi_2$ satisfying \eqref{def:DU} are equal in $ L^{p'}_µ(\mathcal O,\R^d)$.
\end{Prop}
\begin{proof}
\textit{Case $\mo = \R^d$:} Consider two such functions. Take $\psi \in \mathcal C_c(\R^d,\R^d)$. Consider the natural deterministic coupling $\gamma_h$ between $µ$ and $(Id + h \psi)_{\#}µ$ for $h \ne 0$. Using \eqref{def:DU} with $\nu = (Id + h \psi)_{\#}µ$, dividing by $h$ and letting $h \to 0$, we obtain that 
$$
\int_{\R^d}\phi_1(x)\cdot\psi(x)µ(dx)= \int_{\R^d}\phi_2(x)\cdot\psi(x)µ(dx).
$$
Since the equality holds for any $\psi\in \mathcal C_c(\R^d,\R^d)$, the result follows by density.

\textit{General case:} In general, we need to adapt the previous argument since we cannot guarantee that, for $µ \in \mppo, \psi \in \mathcal C_c^\infty(\mo,\R^d)$, $(Id + h \psi)_\#µ \in \mppo$. Instead, consider the set 
$$
K := \{ \psi \in \mathcal C^\infty_c(\R^d,\R^d) | \forall t \geq 0, X^\psi_t(\mo) \subset \mo\},
$$
where $X^\psi_t$ denotes the flow of the ODE $\dot x = \psi(x)$. In other words, the flow of the ODE associated to $\psi$ leaves $\mo$ invariant. Arguing as above, and replacing $(Id + h \psi)_{\#}µ$ by $(X^\psi_h)_\#µ$ we find that if $\phi_1,\phi_2$ are two elements as in the definition, then, for any $\psi \in K$
\be\label{eq:854}
\int_{\mo}\phi_1(x)\cdot\psi(x)µ(dx)= \int_{\mo}\phi_2(x)\cdot\psi(x)µ(dx).
\ee
It turns out that the previous is sufficient to conclude that $\phi_1= \phi_2$. The arguments which follow are quite classical so I only sketch them. First, the proposition would hold if we were able to take any $\psi \in \mathcal C^\infty_c(\R^d,\R^d)$. Hence, it is purely a result about smooth functions. Second, if we show that it is true for the case $\mo = \{ x_1\geq 0\}$, then, by the smoothness of the domain, we can recover the required case, by localization and by using diffeomorphisms which will send $\partial \mo$ (locally) into a hyperplane. Then, it suffices to remark that for $\mo = \{x_1\geq 0\}$, \eqref{eq:854} is equivalent to having any smooth test functions. Indeed, in this case, $K$ is characterized as the set of functions $\psi \in  \mathcal C^\infty_c(\R^d,\R^d)$ such that $\psi_1((0,x_2,\dots,x_d)) \geq 0$, and \eqref{eq:854} against any $\psi$ or any $\psi \in K$ is equivalent, namely since $K-K$ spans a sufficiently large space.

\end{proof}

The notion of horizontal derivative  allows one to  measure appropriately the Lipschitz regularity of functions as the next result shows.

\begin{Prop}
Let $U: \mathcal P_p(\R^d) \mapsto \R$ be an H-differentiable function at $µ_0$. Then, 
$$
\limsup_{\nu, \gamma \in \Pi(µ_0,\nu) C_p(\gamma)\to 0}\frac{|U(\nu)-U(µ_0)|}{C_p(\gamma)} = \|D_µ U(µ_0)\|_{L^{p'}_{µ_0}}.
$$
\end{Prop}
\begin{proof}
The inequality 
$$
\limsup_{C_p(\gamma)\to 0}\frac{|U(\nu)-U(µ_0)|}{C_p(\gamma)} \leq \|D_µ U(µ_0)\|_{L^{p'}_{µ_0}}
$$
simply follows from H\"older's inequality. Let $\psi \in \mathcal C_c(\R^d,\R^d)$. Evaluating the difference $U((X^\psi_t)_{\#}µ_0) - U(µ_0)$, we arrive at
$$
\limsup_{C_p(\gamma) \to 0}\frac{|U(\nu)-U(µ_0)|}{C_p(\gamma)} \geq \frac{\int_{\R^d}D_µU(µ_0)(x)\cdot\psi(x)µ_0(dx)}{\|\psi\|_{L^p_{µ_0}}}.
$$
Arguing once again by density as previously yields the result.
\end{proof}

\subsection{Interpretation of the derivative}
Above, I defined the differentiability of $U: \mathcal P_p(\mathcal O) \mapsto \R$ at $µ$ as the existence of a first order expansion of $U$ of the type
$$
U(\nu) = U(µ) + \int_{\mathcal O^2}D_µU(µ)(x)\cdot(y-x)\gamma(dx,dy) + o(C_p(\gamma))
$$
for $\gamma \in \Pi(µ,\nu)$. I here give some intuition about this definition.\\

\textbf{The direction of the increment:} In this point of view, we evaluate the derivative of $U$ in the direction of the coupling $\gamma$, which is a Lagrangian direction. Indeed, recall that taking a coupling consists in giving mass to pairs of elements: the pair $(x,y)$ has a mass $\gamma(dx,dy)$. The underlying Lagrangian evolution follows the directions $y-x$ with a weight of $\gamma(dx,dy)$. At the level of random variables, if we consider $X \sim µ$ and $Y \sim \nu$, we are simply considering the path $((1-t)X + t Y)_{t \in [0,1]}$. In this case, the direction of the increment is $Y-X$, and it is this term that we recover in the $(y-x)\gamma(dx,dy)$. Quite naturally, the size of the displacement associated to a coupling $\gamma$ is quantified by $C_p(\gamma)$.

The more abstract definition of Section \ref{sec:manifold} might also help to understand the term $(y-x)$.\\

\textbf{How we evaluate along this increment:} There is a quite non-obvious fact in the way we evaluate along such Lagrangian evolutions, and the next section illustrates it. It can be phrased as follows: we are considering Lagrangian paths, but the derivative is somehow Eulerian. This can be observed in \eqref{def:geoD} as the fact that $D_µU(µ)$ only depends on $x$. A priori, this is not obvious. I refer to Theorem \ref{thm:structurelocal} for reasons why, or to Section \ref{sec:superdiff} for the different situation of sub/super differentials, in which it is not the case.

\subsection{Links with derivatives of functions of random variables}\label{sec:lift}
I focus here on the case $\mo = \R^d$. Recall that $\mathbb L^p:= L^p(\Omega, \R^d)$ for $(\Omega, \mathcal A, \mathbb P)$ a standard probability space. Let $U: \mpprd \mapsto \R$ and consider its lift $\mathcal U: \mathbb L^p \mapsto \R$ defined by
$$
\forall X \in \mathbb L^p, \quad \mathcal U (X) = U(\mathcal L(X)).
$$
The derivative $D_µU$ is intrinsically linked to the Fréchet derivative $\nabla \mathcal U(X)$ (in the usual Banach sense) of $\mathcal U$ at any $X$ of law $µ$. We describe this link through the next two results.
\begin{Prop}\label{prop:gg}
Assume that $U$ is H-differentiable at $µ \in \mpprd$, then, for any $X \sim µ$, $\mathcal U$ is differentiable at $X$ and $\nabla \mathcal U(X) = D_µU(µ)(X)$ (in $\mathbb L^p$).
\end{Prop}
\begin{proof}
Evaluating \eqref{def:DU} for any $\gamma = \mathcal L(X,Y)$ for $X\sim µ$, $Y \in \mathbb L^p$ yields the result, since
\be\label{eq:identif}
\int_{\mo^2}D_µU(µ)(x)\cdot(y-x)\gamma(dx,dy) = \mathbb E[D_µU(µ)(X)\cdot(Y-X)],
\ee
and $C_p(\gamma) = \|X-Y\|_p$.
\end{proof}
\begin{Theorem}\label{thm:structurelocal}
Assume that $\mathcal U$ is Fréchet differentiable at $X \in \mathbb L^p$, then $U$ is H-differentiable at $\mathcal L(X)$ and $\nabla \mathcal U(X) = D_µU(\mathcal L(X))(X)$.
\end{Theorem}
\begin{proof}
Arguing as in the previous proof, the result holds if we are able to show that $\nabla \mathcal U(X) =g(X)$ for some deterministic measurable function $g: \R^d\mapsto \R^d$. Let $Z = \nabla \mathcal U(X)$. We define for any $x \in \R^d$, $g(x)$ as the median of $Z$ conditioned on $X=x$, component by component. That is, $g$ is given by the relation $g_i(x) = \inf \{z \in \R | \mathbb P(Z_i \in (-\infty,z]|X = x) \geq \frac 12\}$, which indeed defines a measurable mapping. Note that here, and in the rest of the proof, the conditioning on $\{X=x\}$ is obtained through an implicit use of the disintegration theorem on the couple $(Z,X)$ to obtain the conditional law of $Z$ with respect to $X$.

Assume first that there exists $\zeta$, a uniform random variable on $[0,1]$, defined on the same probability space as $X$, which is independent of $X$.  The rest of the argument follows a simple idea. We want to construct increments along which to test the differentiability of $\mathcal U$, which shall yield the result. In order to do so, let $i \in \{1; \dots ;d\}$ and define $v_i^\pm(x)$ by
$$
\ba
v_i^\pm (x) := \inf\bigg\{&u \in [0,1] | \mathbb P(Z_i = g_i(x), \zeta \leq u|X=x) \\
&\geq \big(\mathbb P(Z_i < g_i(x)|X = x)-\mathbb P(Z_i > g_i(x)|X=x)\big)_\pm\bigg\}.
\ea
$$
Consider now the random variables (our increments to come) $\xi^+_i := \mathbb 1_{\{Z_i > g_i(X)\}} + \mathbb 1_{\{Z_i = g_i(X), \zeta \leq v_i^+(X)\}}$ and $\xi^-_i := \mathbb 1_{\{Z_i < g_i(X)\}} + \mathbb 1_{\{Z_i = g_i(X), \zeta \leq v_i^-(X)\}}$. We now show that, for any $\eps > 0$, $X + \eps\xi^+_i$ and $X+\eps \xi^-_i$ have the same law. We condition the two laws on $\{X = x\}$ and remark that we are facing Bernoulli like random variables. The parameter of the first is given by 
$$
\ba
\mathbb P& (Z_i > g_i(X) | X = x ) + \mathbb P(Z_i = g_i(X), \zeta \leq v^+_i(X) | X = x)\\
&= \mathbb P (Z_i > g_i(X) | X = x )+  (\mathbb P(Z_i < g_i(X)|X=x) - \mathbb P (Z_i > g_i(X) | X = x ))_+\\
& = \max(\mathbb P (Z_i > g_i(X) | X = x ),\mathbb P(Z_i < g_i(X)|X=x)),
\ea
$$
where the first equality comes from the definition of $v^+_i$ and the fact that, almost surely, we can argue as if, conditioned on $\{X = x\}$, $\zeta$ is a uniform random variable. By symmetry, it follows that the announced equality of laws holds, and thus that 
$$
\mathbb E[Z_i \xi^+_i] = \lim_{ \eps \to 0} \frac{\mathcal U(X + \eps \xi^+_ie_i) - \mathcal U(X)}{\eps} = \lim_{ \eps \to 0} \frac{\mathcal U(X + \eps \xi^-_ie_i) - \mathcal U(X)}{\eps} = \mathbb E[Z_i \xi^-_i].
$$
Conditioning on $\{ X = x\}$, we remark that $\mathbb E[g_i(X) \xi^+_i] = \mathbb E [g_i(X)\xi^-_i]$, and thus we deduce from the previous equality that 
$$
0 = \mathbb E[(Z_i - g_i(X))(\xi^+_i - \xi^-_i)] = \mathbb E[|Z_i - g_i(X)|].
$$
Hence the result is proved if we are able to assume we can always consider an independent uniform random variable on $[0,1]$. Lemma \ref{lemma:structureofgradient} just below allows us to do it by placing ourselves at possibly another random variable of law $\mathcal L(X)$. Hence, the result is proved.
\end{proof}
\begin{Lemma}\label{lemma:structureofgradient}
Let $p \in (1,\infty)$, $\mathcal U: \mathbb L^p\mapsto \R$ be Fréchet differentiable at $X \in \mathbb L^p$ and verify, for all $Y,Y' \in \mathbb L^p$ such that $\mathcal L(Y) = \mathcal L(Y')$, $\mathcal U(Y) = \mathcal U(Y')$. Then, $\mathcal U$ is also Fréchet differentiable at any $X'$ such that $\mathcal L(X') = \mathcal L(X)$. Furthermore, $\mathcal L(X,\nabla \mathcal U(X)) = \mathcal L(X',\nabla \mathcal U(X'))$.
\end{Lemma}
\begin{proof}
Let $X$ and $X'$ be as in the statement. For any $\eps > 0$, consider a map $\tau_\eps$ given by Theorem \ref{thm:proba} such that $\|X'\circ \tau_\eps - X\|_\infty \leq \eps$. Let us now compute, for $o(\cdot)$ the remainder appearing in the definition of the differentiability of $\mathcal U$ at $X$, and any $Y \in \mathbb L^p, \eps > 0$,
\be\label{eq:1111}
\ba
\mathcal U(X' + Y) &= \mathcal U(X'\circ \tau_\eps + Y\circ\tau_\eps -X + X)\\
& = \mathcal U(X) + \mathbb E [\nabla \mathcal U(X)\cdot (X'\circ \tau_\eps + Y\circ\tau_\eps -X )] + o(\|X'\circ \tau_\eps + Y\circ\tau_\eps -X \|_p)\\
&= \mathcal U(X') + \mathbb E[\nabla \mathcal U(X)\circ \tau_\eps^{-1} \cdot Y]  +\mathbb E [\nabla \mathcal U(X)\cdot (X'\circ \tau_\eps -X)]\\
& \quad \quad+ o(\|X'\circ \tau_\eps + Y\circ\tau_\eps -X \|_p).
\ea
\ee
Write $Z_\eps = \nabla \mathcal U(X)\circ \tau_\eps^{-1}$. Fix $\delta > 0$ and use the previous equality for $\eps, \eps' > 0$ and $\|Y\|_p \leq \delta$, then 
$$
\mathbb E [(Z_\eps - Z_{\eps'})\cdot Y] \leq C(\eps + \eps') + o(\eps + \delta) + o(\eps' + \delta),
$$
which leads to 
$$
\delta \|Z_\eps - Z_{\eps'}\|_{p'} \leq C(\eps + \eps') + o(\eps + \delta) + o(\eps' + \delta).
$$
Therefore, $(Z_\eps)_{\eps > 0}$ is a Cauchy sequence, set $Z_*$ its limit in $\mathbb L^{p'}$. Passing to the limit $\eps \to 0$ in \eqref{eq:1111} implies that $Z_* = \nabla \mathcal U(X')$. To conclude, remark now that for all $\eps > 0$, $(X \circ \tau_\eps^{-1}, Z_\eps)$ has the law of $(X,\nabla \mathcal U(X))$, but the former converges in $\mathbb L^p \times \mathbb L^{p'}$ to $(X',\nabla \mathcal U(X'))$, hence the result follows.
\end{proof}

\begin{Rem}
The previous approach, and in particular Theorem \ref{thm:structurelocal}, can be extended to more general cases than $\mo = \R^d$, as we shall see in the Section \ref{sec:superdiff}.
\end{Rem}

\begin{Rem}
In fact, the derivative $D_µU(µ)$ has more structure than simply being an element of $L^{p'}_µ$, as it should almost be a gradient. In general, as soon as $U$ is differentiable at $µ$,
$$
D_µ U(µ) \in \overline {\nabla_x C^{\infty}(\R^d,\R)}^{L^{p'}_µ}.
$$
However, such a structure will not be helpful in what follows, see Gangbo and Tudorascu \citep{gangbotudorascu} for more details on this.
\end{Rem}

\subsection{Links with the V-derivative}
Certain smooth functions are expected to be both V-differentiable and H-differentiable, i.e. differentiable in the senses of both Definitions \ref{def:flat} and \ref{def:geoD}. For such functions $U$, there is a natural link between $\nabla_µ U$ and $D_µ U$ that I now present.\\

The following two results state that if $U$ is differentiable in a sense, and its derivative sufficiently smooth, then it is also differentiable in the other sense. In both cases, a notion of regularity has to be used for the derivative and I refer to Section \ref{sec:C1} for more on this topic. The main result of this section is the formula \eqref{formula}.
\begin{Theorem}
Let $\mo$ be the closure of a bounded smooth convex domain of $\R^d$. Let $U : \mpo \mapsto \R$ be V-differentiable and assume that for all $µ\in \mpo$, $x \mapsto \nabla_µ U(µ)(x) \in \mathcal C^1(\mo)$ and that $$
\ba
&(\mo\times \mpo)&\mapsto &\R^d,\\
&(x,µ) &\mapsto &\nabla_x\nabla_µ U(µ)(x)
\ea
$$ is a continuous function for the product topology of the usual one on $\mo$ and the weak topology on $\mpo$. Then $U$ is also H-differentiable (for any $p\in (1,\infty)$) and 
\be\label{formula}
\boxed{D_µU(µ) = \nabla_x \nabla_µ U(µ), \quad µ \text{ almost everywhere}.}
\ee
\end{Theorem}
\begin{proof}
Let $µ,\nu \in \mpo$. From the regularity of $U$, we can compute the derivative $\frac{d}{d\theta} U((1-\theta)µ + \theta \nu)$ and obtain that
$$
U(\nu) - U(µ) = \int_0^1\int_\mo\nabla_µU((1-\theta)µ + \theta\nu)(x)(\nu-\mu)(dx)d\theta.
$$
Let $\gamma \in \Pi(µ,\nu)$, and remark that the previous can be written
$$
\ba
U(\nu) - U(µ) &= \int_0^1\int_{\mo^2}\nabla_µU((1-\theta)µ + \theta\nu)(y) - \nabla_µU((1-\theta)µ + \theta\nu)(x)\gamma(dx,dy)d\theta\\
&=: \int_0^1\int_{\mo^2}\nabla_x\nabla_µU((1-\theta)µ+ \theta\nu)(x)\cdot(y-x)\gamma (dx,dy)d\theta+\\
& \quad +R(µ,\nu),
\ea
$$
where $R(µ,\nu)$ is just defined as the rest of the expression. Denote by $\omega$ a uniform modulus of continuity of the family $(\nabla_x\nabla_µ U(µ))_{µ \in \mpo}$. Such a modulus exists because $\mpo$ is compact for the weak topology. This allows  allows one to to estimate
$$
|R(µ,\nu)| \leq \int_0^1\int_{\mo^2}\int_0^{|y-x|}\omega(u)\,du\,\gamma(dx,dy)d\theta.
$$
Using the concavity of $\omega$ and Jensen's inequality, we deduce that
$$
\int_0^{|y-x|}\omega(u)\,du \leq |y-x| \omega\left(\frac{|y-x|}{2}\right).
$$
Choosing $(X,Y)\sim \gamma$, we obtain $|R(µ,\nu)| \leq \|X-Y\|_p\|\omega(\frac{|X-Y|}{2})\|_{p'}$, which is indeed a $o(\|X-Y\|_p)$ by dominated convergence (Recall that $\mo$ is bounded so that $X-Y$ is bounded almost surely). Hence, it only remains to show that 
$$
 \int_0^1\int_{\mo^2}(\nabla_x\nabla_µU((1-\theta)µ+ \theta\nu)(x)- \nabla_x\nabla_µU(µ)(x))\cdot(y-x)\gamma (dx,dy)d\theta
$$
is also a $o(C_p(\gamma))$. Thanks to the regularity of $\nabla_x \nabla_µ U$, it is also the case, hence the result is proved.
\end{proof}
\begin{Rem}
The convexity of the domain is not at all a mandatory assumption and can be removed easily. It just allows an easier use of the Taylor expansion in the estimate of $R(µ,\nu)$. The boundedness of $\mpo$ can also be removed, but then precise assumptions on the growth of the derivative have to be made.
\end{Rem}

I do not prove the result in the other direction, but simply state it in the way it can be found in the literature.
\begin{Theorem}
In the case $p =2$, $\mo = \R^d$, consider a function $U: \mpt \mapsto \R$ such that its lift $\mathcal U$ is $\mathcal C^1$ in $\mathbb H$, that $\nabla_X \mathcal U$ is Lipschitz continuous in $\mathbb H$ and that for all $µ \in \mpt$, $D_µ U \in L^2((\R^d,µ),\R^d)$ can be extended into a continuous function on $\R^d$. Then, there exists a continuous map $\Phi: \mpt \times \R^d \mapsto \R$ such that for any $µ,\nu\in \mpt$
$$
U(\nu) - U(µ) = \int_0^1\int_{\R^d}\Phi((1-\theta)µ + \theta\nu,x)(\nu-\mu)(dx)d\theta.
$$ 
\end{Theorem}
This statement may seem arbitrarily stated and I now comment on it. The statement on the regularity of $U$ is made through its lift, because defining a notion of $\mathcal C^1$ map on $\mpt$ is not obvious and I refer to Section \ref{sec:C1} for more details on this question. We are also facing the fact that $D_µU(µ)$ is only defined on the support of $µ$ and not everywhere. Finally, the previous is not sufficient to prove that $U$ is V-differentiable as such a notion was defined with respect to the total variation. For instance, $µ \mapsto \int_{\R^d}|x|^2µ(dx)$ satisfies the requirements of the previous Theorem, but it is not vertically differentiable, since it is not continuous for the total variation. Indeed, consider for instance the case $\mo = \R$, $µ = \delta _0$, $(µ_n)_{n \geq 1}$ defined by $µ_n:= (1-\frac1n)\delta_0 + \frac 1n \delta_{\text{exp}(n)}$, which does converge to $\delta_0$ in total variation but not for $\mathcal W_2$.
\begin{Rem}
The previous result can be extended to more general $p \in (1,\infty)$ and domains $\mo$ without much difficulties.
\end{Rem}

\subsection{The flat $\mathcal W_1$ case}\label{sec:W1}
At this point, I recall that the description of the vertical (or flat) derivative of Section \ref{sec:flat} was deliberately made provocative by the use of the total variation norm, and other distances could have been considered, with of course some appropriate changes. In this section, I present such a possible point of view, in the case of derivatives of functions $U: \mathcal P_1(\R^d) \mapsto \R$.\\

Denote by $Lip(\R^d)$ the space of Lipschitz functions on $\R^d$, and recall the Kantorovich duality formula valid for any $µ,\nu \in \mathcal P_1(\mo)$
$$
\mathcal W_1(µ,\nu) = \sup_{\phi \in Lip(\R^d) \\ Lip(\phi)\leq 1} \int_{\R^d}\phi \,d(µ-\nu).
$$
The previous duality formula hints that $\mathcal P_1(\R^d)$ is a closed convex subset of a Banach space. Instead of choosing this Banach space to be a version of $W^{-1,1}(\R^d)$, defined usually as the topological dual of $Lip(\R^d)$, I introduce the following.\\

Let $E$ be the Banach space obtained as the completion of the span of Dirac masses on $\R^d$ for the norm
$$
\|µ\|_E := \sup_{\phi} \int_{\R^d}\phi \,dµ,
$$
where the supremum is taken over all $\phi \in Lip(\R^d)$ such that $Lip(\phi)\leq 1$ and $\sup_{x \in \R^d} \frac{|\phi(x)|}{1+ |x|} \leq 1$. The following properties of $E$ hold.
\begin{Prop}
The topological dual of $E$ is isomorphic to $Lip(\R^d)$ endowed with the norm $$\|\phi\|:= Lip(\phi) + \sup_{x \in \R^d} \frac{|\phi(x)|}{1+ |x|}.$$
\end{Prop}
\begin{proof}
Let $L\in E'$. Denote $\phi(x) = L(\delta_x)$. Observe that for $x,y \in \R^d$, 
$$
|\phi(x) - \phi(y)| \leq \|L\|_{E'} \mathcal W_1(\delta_x,\delta_y) = \|L\|_{E'}|x-y|.
$$
Hence, $\phi$ is Lipschitz continuous and it is now easy to verify that we can identify $E'$ with a subspace of $Lip(\R^d)$. Conversely, every Lipschitz function defines an element of $E'$ by construction, thus the result holds. 
\end{proof}

\begin{Prop}
The Wasserstein space $\mathcal P_1(\R^d)$ is a closed convex subset of $E$.
\end{Prop}
\begin{proof}
Let $µ \in \mathcal P_1(\R^d)$. By considering a sequence of discrete measures converging toward $µ$ for $\mathcal W_1$, we obtain that $µ \in E$. It is closed as the distance of $E$ yields the same topology as the $\mathcal W_1$ norm and the Wasserstein space is complete.
\end{proof}

The following definition is an analogue of defining the Fréchet derivative for a function defined on a closed convex subset of the Banach space $E$.

\begin{Def}
A function $U: \mathcal P_1(\R^d) \mapsto \R$ is said to be $\mathcal W_1$ differentiable at $µ_0 \in \mathcal P_1(\R^d)$ if there exists $\varphi \in Lip(\R^d)$ such that
$$
\lim_{\mathcal W_1(\nu,µ_0) \to 0}\frac{U(\nu) - U(µ_0) - \langle \varphi,\nu-µ_0\rangle}{\mathcal W_1(\nu,µ_0)} = 0.
$$
If $\varphi$ is such that $\langle \varphi,µ_0\rangle = 0$, then we set $\nabla^E_µ U(µ_0) = \varphi$.
\end{Def}
As it was the case for the V-derivative, one can check that two functions $\varphi_1,\varphi_2$ satisfying the previous requirement differ by at most a constant, and thus that notation $\nabla^E_µ U$ is also meaningful here. Furthermore, standard properties of Fréchet derivatives also hold in this case.\\

It is tempting to think that the previous notion of derivative  allows one to  unify the horizontal and vertical points of view in the case $\mathcal P_1(\R^d)$. Indeed, if $U$ is $\mathcal W_1$ differentiable at $µ_0$, then because $\nabla^E_µ U(µ_0)$ is Lipschitz, it is almost everywhere differentiable on $\R^d$. Considering its almost everywhere gradient $\nabla_x \nabla_µU(µ_0)$, we are tempted to say that we recover, because of formula \eqref{formula}, the H-derivative of $U$ at $µ_0$. It turns out that it is not the case, because this almost everywhere differentiability of Lipschitz functions is not enough. Indeed, consider the case $U(µ) = \int_{\R^d}|x|µ(dx)$. This function $U$ is clearly $\mathcal W_1$ differentiable everywhere with derivative $x \mapsto |x| \in Lip(\R^d)$. However, it is not horizontally differentiable at $\delta_0$ for instance. 

\subsection{Links with a more geometric notion of derivative}\label{sec:geo}
When defining the horizontal derivative of a function $U: \mppo \mapsto \R$, we have allowed to consider variations along all couplings. This may seem quite demanding, and a more "geometric" point of view would have consisted in asking information only for couplings which are optimal for $\mathcal W_p$. This may seem natural since, after all, we want to go from one measure to the other, not necessarily by arbitrary path. This leads to the definition
\begin{Def}\label{def:geoDD}
A function $U: \mathcal P_p(\mathcal O) \mapsto \R$ is geometrically differentiable at $µ \in \mathcal P_p(\mathcal O)$, if there exists $\phi \in L^{p'}((\mathcal O,µ),\R^d)$ and a modulus of continuity $\omega$ such that, for any $\nu \in \mathcal P_p(\mathcal O)$, $\gamma$ optimal coupling in $\mathcal W_p(µ,\nu)$
$$
\left|U(\nu) - U(µ) - \int_{\mo^2}\phi(x)\cdot(y-x)\gamma(dx,dy)\right| \leq \omega(\mathcal{W}_p(\nu,µ))\mathcal{W}_p(\nu,µ).
$$
\end{Def}
Of course, Definition \ref{def:geoD} is more demanding than Definition \ref{def:geoDD}, so horizontal differentiability implies geometrical differentiability. A striking result is in fact the following, which is due to Gangbo and Tudorascu in \citep{gangbotudorascu}, in the case $p =2, \mo = \R^d$.
\begin{Theorem}\label{thm:geoDD}
Let $U: \mpt \mapsto \R$ be geometrically differentiable at $µ$. Then $U$ is horizontally differentiable at $µ$.
\end{Theorem} 
I omit the long and technical proof of this result. Also, the restrictions on $p$ and $\mo$ can be relaxed.

\subsection{The case of probability measures on manifolds}\label{sec:manifold}
When we are concerned with probability measures over manifolds, and not simply over subsets of $\R^d$, some care has to be taken when manipulating the horizontal differential. I present briefly this setting here, both to show this difference, but also as I believe it is a nice way to understand in more details Definition \ref{def:geoD}, and in particular the term $(y-x)$ in the integral. Furthermore, I focus on the case of a compact manifold without boundary, mainly to lighten the notation, but it plays little role in the following. Recall that the vertical point of view does not see the geometry of the space $\mo$, so it does not change anything here.\\

Let $M$ be a smooth compact Riemannian manifold without boundary. I adopt the usual notation: for $x \in M$, $T_xM$ is the tangent space to $M$ at $x$, which is equipped with the norm $|\cdot|_x$. The tangent bundle is denoted by $TM$ and I write $\phi: M \mapsto TM$, when, for any $x \in M, \phi(x) \in T_xM$. The exponential map is denoted by $\exp_x(v)$, it is a map $TM \mapsto M$ and it yields the value at time $t=1$ of the unique geodesic $\theta : [0,1]\mapsto M$ such that $\theta(0) = x, \theta'(0) = v$. When considering a probability measure $\eta \in \mathcal P(TM)$, we denote its size by 
$$
C_M(\eta) := \sqrt{\int_M\int_{T_xM}|v|^2_x\eta(dx,dv)}.
$$

\begin{Def}\label{def:dM}
Let $U: \mathcal P(M) \mapsto \R$. It is said to be horizontally differentiable at $µ \in \mathcal P(M)$ if there exists a measurable $\phi: M \mapsto TM$ such that $$\int_{M}|\phi(x)|^2_xµ(dx) < \infty,$$ and a modulus of continuity $\omega$ such that, for any $\eta \in \mathcal P(TM)$ such that $(\pi_1)_\#\eta = µ$,
\be\label{def:DUM}
\left|U(\exp_\#\eta) - U(µ) - \int_{TM}\phi(x)\cdot v\,\eta(dx,dv)\right| \leq \omega(C_M(\eta))C_M(\eta).
\ee
We still write $\phi= D_µU(µ)$. 
\end{Def}

Similar developments on this notion of derivative can be carried out and are omitted here. I insist on the Definition \ref{def:dM} and its link with Definition \ref{def:geoD}. Somehow, the term $(y-x)$ in the integral in Definition \ref{def:geoD} is not only the difference between the paired points in the coupling $\gamma$, but also the initial speed of the path joining them. Hence, instead of writing everything in terms of the coupling $\gamma \in \mathcal P(\R^d\times \R^d)$ joining the two measures, I could have used (perfectly equivalently) a coupling $\eta \in \mathcal P(\R^d\times \R^d)$ defined by $\eta = (\pi_1,\pi_2-\pi_1)_\#\gamma$. On manifolds, we express the Lagrangian paths in terms of pairs (initial position, speed), while on $\R^d$ we can equivalently do it in terms of pairs (initial point, arriving points). 

Finally, the integrability requirement on $\phi$ is imposed so that the integral in \eqref{def:DUM} is well defined. If we consider other costs that $C_M$, then we can change the integrability requirement on $\phi$ as well, namely by using the $(L^p,L^{p'})$ duality.

\subsection{Some useful computations of vertical and horizontal derivatives}\label{sec:computations}
Here, I present some computations around the notions of derivatives introduced above. As most of what follows can be computed quite easily, some proofs are left as exercises to the reader interested in improving her or his skills in calculus.

\subsubsection{Around polynomials and their derivatives}
Let $\mo$ be the closure of a smooth bounded domain of $\R^d$ (more general domains could be considered). Recall that a monomial on $\mmo$ is a function $U: \mmo \mapsto \R$ of the form
\be\label{def:poly}
U(µ) := \prod_{i =1}^k\int_{\mo}\phi_i(x)µ(dx),
\ee
for some continuous real functions $(\phi_i)_{1 \leq i \leq k}$ on $\mo$.
We can compute the differentials of such functions. 
\begin{Prop}
Let $U$ be a monomial given by \eqref{def:poly}. Then $U$ is V-differentiable and, for all $µ \in \mmo, x \in \mo,$
$$
\nabla_µ U(µ)(x) = \sum_{i = 1 }^k \phi_i(x)\prod_{j \ne i}\int_{\mo}\phi_j(y)µ(dy).
$$
Furthermore, if the $\phi_i$ are differentiable and their first order derivatives are bounded, then the restriction of $U$ to $\mpo$ is also horizontally differentiable (for any $p$) and for all $µ \in \mpo$
$$
D_µU(µ)(x) = \sum_{i = 1 }^k \nabla_x\phi_i(x)\prod_{j \ne i}\int_{\mo}\phi_j(y)µ(dy).
$$
\end{Prop}
\begin{Rem}
The continuity of the $\phi_i$ is not needed to obtain the $V$-differentiability.
\end{Rem}
Building on the previous, we can also consider cases in which the function inside the integral depends on $µ$ itself, as is the case for $U: \mmo \mapsto \R$ defined by
\be\label{def:nonlinearpoly}
U(µ) := \int_{\mo}\phi(µ,x)µ(dx),
\ee
for some $\phi: \mmo\times \mo \mapsto \R$.
\begin{Prop}
Let $U: \mmo \mapsto \R$ be defined by \eqref{def:nonlinearpoly}. If $\phi$ is uniformly continuous on $\mmo\times \mo$ equipped with the product topology of the total variation distance and the Euclidean one, and that $\phi(\cdot,x)$ is V-differentiable at $µ_0 \in \mmo$ for any $x \in \mo$, with some uniformity in the sense that there exists a modulus of continuity $\omega$ such that, for all $x \in \mo$
$$
\left| \phi(\nu,x) - \phi(µ_0,x) - \int_\mo \nabla_µ \phi(µ_0,x)(y)(\nu-µ_0)(dy) \right| \leq \omega (|µ_0 - \nu|)|µ_0-\nu|,
$$
 then $U$ is V-differentiable at $µ_0$ and, for any $x \in \mo$,
$$
\nabla U(µ_0)(x) = \phi(µ_0,x) + \int_\mo \nabla_µ\phi(µ_0,y)(x)µ_0(dy).
$$
\end{Prop}
\begin{Rem}
The same types of result hold if we work with the $\mathcal W_1$ norm instead of the total variation, using then the definition of Section \ref{sec:W1}.
\end{Rem}
Many times, it is not for such smooth functionals that one wants to compute derivative, but a similar structure is apparent. Consider for instance the so-called \textit{entropy} $E: \mpo \mapsto \R$ defined by
\be\label{def:entropy}
E(µ) = \begin{cases} \int_\mo \log(µ(x))µ(x)dx \text{ if } µ << Leb,\\ + \infty \text{ else,}\end{cases}
\ee
where $µ << Leb$ stands for the fact that $µ$ is absolutely continuous with respect to the Lebesgue measure, with density still denoted by $µ$. Formally, we can then compute the vertical derivative
$$
\nabla_µ E(µ)(x) = \log(µ(x)) + 1,
$$
where, extremely formally, I have used $\nabla_µ \log(µ) = µ^{-1}$, at the level of the density. Pursuing in this heuristics, we then obtain the horizontal derivative as 
$$
D_µE(µ)(x) = \nabla_x(\log(µ(x))).
$$
This  allows one to  express the \textit{Fisher information} $\mathcal I(µ)$ of $µ$ through the entropy thanks to the formula
\be\label{fisher}
\mathcal I(µ) := \int_\mo \frac{|\nabla_xµ(x)|^2}{µ(x)}dx = \|D_µE(µ)\|^2_{L^2_µ}.
\ee

\textbf{Interaction kernels:}\\

As a slight extension of polynomials, we can also consider the case in which the $\phi_i$ are functions of several variables. For instance let $U: \mpo \mapsto \R$ be formally given by
$$
U(µ) := \int_{\mo^2}\psi(x,y)µ(dx)µ(dy),
$$
then it follows that if $\psi$ is continuous, $U$ is $V$-differentiable and, for all $µ \in \mpo$, $x \in \mo$
$$
\nabla_µ U(µ)(x) = \int_{\mo}\psi(x,y)µ(dy) + \int_{\mo}\psi(y,x)µ(dy).
$$
These functions typically arise as the potential energy or entropy of some mean field systems, namely when $\psi(x,y) = \phi(x-y)$ and $\phi$ is itself the pairwise potential of the interaction. For instance in the case of the mean field system associated to the so-called Dyson Brownian motion, the functional
$$
E(µ) := -\frac 12\int_{\R^2}\log(|x-y|)µ(dx)µ(dy)
$$
appears. Such models are sometimes called log-gas. We have the relation $E(µ) > - \infty$ as soon as $µ$ has a logarithmic moment, but $E(µ) = + \infty$ if $µ$ is too concentrated, in particular if it is of the form $µ = \nu + \eps \delta_{x_0}$ for any $\nu \in \mathcal M_+(\mo)$, $\eps > 0$ and $x_0 \in \R$. Since such measures are dense in any $\mathcal P_p(\R)$, $E$ is clearly not differentiable. However, we can still make the formal computation
$$
-D_µE(µ)(x) = \int_{\R}\frac{1}{x-y}µ(dy) =: H[µ](x),
$$
defining the middle term through principal values, as in the theory of distributions. In particular, when we are concerned with $\mathcal P_2(\R)$, it is quite tempting to say that $E$ is horizontally differentiable when $H[µ] \in L^2_µ(\R)$ and obtain a formula similar to \eqref{fisher}.

\subsubsection{Derivatives along continuity equations and Fokker-Planck equations}
I focus here on the case $\mo = \R^d$ to simplify notation, and my objective is to compute the limit
$$
\lim_{t \to 0^+} \frac{U(m_t)-U(m_0)}{t},
$$
for $U: \mpo \mapsto \R$ and $(m_t)_{t \geq 0}$ weak solution to the PDE
\be\label{ce5}
\partial_t m - \sigma\Delta m + \text{div}(b \,m) = 0 \text{ in } (0,\infty)\times \R^d,
\ee
where $b : \R^d \mapsto \R^d$ is a given Lipschitz function and $\sigma \in \mathbb R_+$. I will give several ways to compute this limit, depending on whether $U$ is V-differentiable or H-differentiable. Because the flow generated by \eqref{ce5} is not regular in either total variation or Wasserstein metrics, additional requirements are needed on the derivatives of $U$ in order to be able to compute the limit.
\begin{Prop}\label{prop:georges}
Let $U: \mprd \mapsto \R$ be a V-differentiable map such that
\begin{itemize}
\item  for all $µ$, $\nabla_µ U(µ): \R^d \mapsto \R$ is a $\mathcal C^2$ map, such that, uniformly in $µ$, $\nabla_µ U(µ)(x)$ and $\nabla_x \nabla_µ U(µ)(x)$ converge toward $0$ as $|x| \to \infty$,
\item  the map $\mprd \mapsto \mathcal C^2(\R^d), µ \mapsto \nabla_µ U(µ)$ is continuous, when $\mprd$ is equipped with the topology of weak convergence, and $\mathcal C^2_b(\R^d)$ is equipped with its usual norm.
\end{itemize}
Then,
$$
\ba
\lim_{t \to 0^+}\frac{U(m_t) - U(m_0)}{t} =& \int_{\R^d} \nabla_x\nabla_µ U(m_0)(x)\cdot b(x)m_0(dx) \\
&+\sigma \int_{\R^d} \Delta_x \nabla_µU(m_0)(x)m_0(dx).
\ea
$$
\end{Prop}
\begin{proof}
Let $t > 0$. Since $U$ is V-differentiable, we can compute $\frac{d}{ds}U((1-s)m_0 + sm_t)$ and obtain
$$
U(m_t) - U(m_0) = \int_0^1\int_{\R^d}\nabla_µU((1-s)m_0 + sm_t)(x)(m_t-m_0)(dx)ds.
$$
Using the fact that $(m_t)_{t \in [0,1]}$ is a weak solution to \eqref{ce5}, we can express the right hand side as follows
$$
\ba
U(m_t) - U(m_0) =\int_0^1\int_0^t&\int_{\R^d}\{\nabla_x \nabla_µU((1-s)m_0 + sm_t)\cdot b(x)\\
& \quad \quad + \sigma \Delta_x\nabla_µU((1-s)m_0 + s m_t)(x)\}m_{s'}(dx)ds'ds.
\ea
$$
The previous holds because, thanks to the decay of $\nabla_µ U(µ)(x)$ when $|x| \to \infty$, we can indeed use it as a test function in the weak formulation of \eqref{ce5}. Dividing by $t$, letting $t\to 0^+$ and using the regularity of $U$ yields the result. 
\end{proof}
\begin{Rem}
Because the flow of \eqref{ce5} is not continuous for the total variation distance, the continuity requirement on $\nabla_µ U$ has to be made with respect to the weak topology. If $m_0$ is sufficiently regular, then we can replace this assumption by continuity with respect to total variation. A similar fact is used below when working with the horizontal derivative.
\end{Rem}

To present the results for the horizontal derivative, I split the computation of the limit in the cases $b = 0$ and $\sigma = 0$. I also focus on the case of $\mpt$.
\begin{Prop}
Assume $\sigma = 0$ and $U: \mpt \mapsto \R$ is H-differentiable at $m_0 \in \mpt$. Then
$$
\lim_{t \to 0^+}\frac{U(m_t) - U(m_0)}{t} = \int_{\R^d}D_µ U(m_0)(x)\cdot b(x)\, m_0(dx).
$$
\end{Prop}
\begin{proof}
Consider a random variable $Y_0$ of law $m_0$.  For $t > 0$, define $X^b_t$ the flow of the ODE. Consider the coupling $\gamma_t := \mathcal L(Y_0,X^b_t(Y_0))\in \Pi(m_0,m_t)$. Evaluating the H-differentiability of $U$ in $m_0$ yields, for some modulus of continuity $\omega$
\be\label{eq:1348}
\left|U(m_t) - U(m_0) - \int_{\R^d}D_µU(m_0)(x)\cdot(X^b_t(x) - x)m_0(dx) \right| \leq C_2(\gamma_t)\omega(C_2(\gamma_t)).
\ee
Remark that, using the Lipschitz regularity of $b$, for some $C$ depending on $b$ and on the second moment of $m_0$
$$
\int_{\R^d}|X^b_t(x) - x|^2m_0(dx) = \int_{\R^d}\left|\int_0^tb(X^b_s(x))ds\right|^2m_0(dx) \leq C t^2.
$$
This  allows one to  control the remainder in \eqref{eq:1348}. Furthermore,
$$
\int_{\R^d}\left|\frac{X^b_t(x) - x}{t} - b(x)\right|^2m_0(dx) = \int_{\R^d}\left| \int_0^t\frac{b(X^b_s(x)) - b(x)}{t}ds\right|^2m_0(dx).
$$
The latter goes to $0$ when $t \to 0$, which implies that
$$
\lim_{t \to 0} \int_{\R^d}D_µU(m_0)(x)\cdot\frac{X^b_t(x) - x}{t}m_0(dx) = \int_{\R^d}b(x)\cdot D_µ U(m_0)(x)m_0(dx).
$$
Hence the result is proved.
\end{proof}
\begin{Prop}
Assume $b \equiv 0$, $\sigma = 1$, that $m_0$ has finite Fisher information, and that $U: \mpt \mapsto \R$ is H-differentiable at $m_0$. Then
$$
\lim_{t \to 0^+}\frac{U(m_t) - U(m_0)}{t} = \int_{\R^d} D_µ U(m_0)(x)\cdot \frac{\nabla_x m_0(x)}{m_0(x)}m_0(dx).
$$
\end{Prop}
\begin{proof}
The main idea of the proof consists in writing the heat equation
$$
\partial_t m = \text{div}_x\left( \frac{\nabla_x m}{m} m  \right) \text{ in } (0,\infty)\times \R^d,
$$
and using the formalism of the previous proposition. We now need several classical estimates on the heat equation that I will not prove but simply refer to. The way we define H-differentiability makes quite transparent the estimates we need. Indeed, we already have for couplings $\gamma$ such that $(\pi_1)_\#\gamma = m_0$,
$$
\left|U((\pi_2)_\#\gamma)- U(m_0) - \int_{\R^d}D_µU(m_0)(x)\cdot(y-x)\gamma(dx,dy)\right| \leq \omega(C_2(\gamma))C_2(\gamma)
$$
We want to take $\gamma (dx,dy) = \mathcal L(X_0,X_t)$ where $(X_t)_{t \geq 0}$ is an $\mathbb H$ valued process given as the solution to
$$
\frac{d X_t}{dt} = \frac{\nabla_x m_t}{m_t}(X_t) \text{ for } t > 0,
$$
and $X_0 \sim m_0$. Such a process exists as a consequence of Theorem \ref{ambrosio}, which holds because the decay of the Fisher information along the heat flow yields the required integrability condition. Observe that for any $t > 0$, $X_t \sim m_t$. In order to obtain the result, it then remains to show that 
\be\label{est1a}
C_2(\gamma) = (\mathbb E[|X_t - X_0|^2])^\frac 12 \text{ is of order } t,
\ee
and that
\be\label{est2a}
\lim_{ t \to 0}\mathbb E \left[ \left| \frac{X_t - X_0}{t} - \frac{\nabla_x m_0}{m_0}(X_0) \right|^2\right] = 0.
\ee
Estimate \eqref{est1a} follows from the continuity of the Fisher information along the heat flow, when starting from an initial condition with finite Fisher information, see for instance Lemma 5.5 in \citep{daudinseeger}. 

To obtain \eqref{est2a}, start with the computation
$$
\ba
\mathbb E \left[ \left| \frac{X_t - X_0}{t} - \frac{\nabla_x m_0}{m_0}(X_0) \right|^2\right] &= \mathbb E \left[ \left| \frac1t\int_0^t  \frac{\nabla_x m_s}{m_s}(X_s)- \frac{\nabla_x m_0}{m_0}(X_0) ds\right|^2\right]\\
 &\leq \frac1t\int_0^t\mathbb E \left[ \left|  \frac{\nabla_x m_s}{m_s}(X_s)- \frac{\nabla_x m_0}{m_0}(X_0) \right|^2\right] ds.
\ea
$$
So the estimate holds if we are able to prove that $(\frac{\nabla_x m_s}{m_s}(X_s))_{s> 0}$ converges toward $\frac{\nabla_x m_0}{m_0}(X_0)$ in $\mathbb H$. First, the convergence of the norm holds, because of the continuity of the Fisher information. Now, it suffices to identify the weak limit of this sequence, and standard PDE arguments (e.g. multiply by smooth test functions and integrate) show that this limit can only be $\frac{\nabla_x m_0}{m_0}(X_0)$. Hence, the result is proved.
\end{proof}

\subsubsection{Discretization of the space of measures or measures on discrete spaces}\label{sec:disccomp}
Let $(µ_i)_{i \in I}$ be a finite collection of measures in $\mmo$ that we can imagine as a discretization of $\mmo$ for instance. Consider the function $v: \R^I \mapsto \R$ defined by
\be\label{defvNflat}
v(z) = U\left(\sum_{i \in I} z_i µ_i\right)
\ee
for a function $U : \mmo \mapsto \R$.
\begin{Prop}
If $U$ is V-differentiable, then $v$ is differentiable and, for any $i \in I$, $z \in \R^I$,
$$
\partial_i v (z) = \int_{\mo}\nabla_µU\left(\sum_{j \in I} z_j µ_j\right)(x)µ_i(dx).
$$
\end{Prop}
The proof is a simple exercise.\\

There are two useful applications of this result. The first one is the discretization of $\mo$ through a finite partition $(A_i)_{i \in I}$ of measurable subsets, and the setting of $µ_i = \mathbb 1_{A_i}$. In such a case, $z_i$ counts the (algebraic) mass put on $A_i$, and it should appear natural that it is the vertical derivative that appears here. The second one is when there exists a vector $y \in (\R^d)^N$ such that $µ_i = \delta_{y^i}$ and a similar interpretation holds.\\

Note that when $\mo$ is finite, given by $\{e_1,\dots,e_N\}$, then setting $I = N$ and $µ_i = \delta_{e_i}$ yields an equality between $U$ and $v$, and between the V-derivative of $U$ and the gradient of $v$ as well.

\subsubsection{Particles approximations}
Recall the approximation result of Theorem \ref{thm:glivenko}. With such results in mind, it is natural, given $U : \mathcal P(\mo)\mapsto \R$, $N > 1$, to consider its $N$-particle projection $u: \mo^N \mapsto \R$ defined by
\be\label{def:uproj}
\forall x \in \mo^N, \quad u(x) = U\left( \frac1N\sum_{i =1}^N \delta_{x_i}\right).
\ee
\begin{Prop}
Let $p \in (1,\infty)$ and assume $U: \mppo \mapsto \R$ is H-differentiable, then $u$ defined by \eqref{def:uproj} is differentiable and 
$$
\forall 1\leq i \leq N, x \in \mo^N, \quad \nabla_{x_i} u(x) = \frac1ND_µU\left( \frac1N\sum_{i =1}^N \delta_{x_i}\right)(x_i).
$$
\end{Prop}
The proof is again a simple exercise once one remarks that to evaluate the difference $u(y) - u(x)$, one needs to use the coupling $\frac 1N \sum_{i= 1}^N \delta_{(x_i,y_i)}$.

\subsubsection{Wasserstein distances}
In this section, I give a simple result of differentiability of Wasserstein distances. To lighten notation, I restrict my attention to the case $\mo = \R^d$ and $p=2$, even though this does not play a strong role in the following.
\begin{Prop}
Let $µ_0,\nu \in \mpt$ be such that the unique optimal coupling in $\mathcal W_2(µ_0,\nu)$ is given by $(Id,T)_\#µ_0$ for some measurable map $T$. Then, $\Phi: µ \mapsto \mathcal W_2^2(µ,\nu)$ is H-differentiable at $µ_0$ and 
$$
D_µ \Phi(µ_0)(x) = 2(x - T(x)).
$$
\end{Prop}
\begin{proof}
Let $X \in \mathbb H$ be of law $µ_0$ and consider a sequence $(\xi_n)_{n \geq 0}$ converging toward $0$ in $\mathbb H$. Finally, for $n \geq 0$, denote by $Y_n$ an element of $\mathbb H$ of law $\nu$ such that $\mathbb E[|Y_n - (X+ \xi_n)|^2] \leq  \Phi(\mathcal L(X+ \xi_n)) + \eps_n$, for $\eps_n=\|\xi_n\|^2_2$. By definition of $\mathcal W_2^2$:
$$
\Phi(\mathcal L(X+\xi_n)) \leq \mathbb E[|X+\xi_n - T(X)|^2] = \Phi(µ_0) + 2\mathbb E[(X-T(X))\cdot\xi_n] + \eps_n,
$$
\be\label{estPhi}
\Phi(µ_0) \leq \mathbb E[|X-Y_n|^2] = \mathbb E[|X + \xi_n - Y_n|^2] -2 \mathbb E[ (X-T(X))\cdot\xi_n] - \eps_n + 2 \mathbb E[(Y_n-T(X))\cdot\xi_n].
\ee
Rearranging leads to
$$
|\Phi(\mathcal L(X+\xi_n))- \Phi(µ_0) - 2 \mathbb E[(X-T(X))\cdot\xi_n]| \leq \|\xi_n\|_2(2\|\xi_n\|_2 + \|Y_n -T(X)\|_2) + \eps_n.
$$
Hence the result follows, provided that $\lim_{n \to \infty}\|Y_n -T(X)\|_2= 0$. Passing to the limit in \eqref{estPhi} yields that $((X,Y_n))_{n \geq 0}$ converges in law toward $(X,T(X))$. Let $\delta > 0 $ and denote by $f: \R^d \mapsto \R^d$ a continuous bounded function such that $\|f(X) - T(X) \|_2 \leq \delta$. Now remark that
$$
\ba
\mathbb E[|Y_n - T(X)|^2] &\leq  2\mathbb E[|Y_n - f(X)|^2] + 2 \mathbb E[|f(X) - T(X)|^2]\\
&\leq 2 \mathbb E[|T(X)|^2] - 4 \mathbb E[Y_n\cdot f(X)] + 2 \mathbb E [|f(X)|^2] + 2 \delta^2.
\ea
$$
Using that $(Y_n)_{n \geq 0}$ is bounded in $\mathbb H$ and $f$ is bounded, we obtain that $(Y_n\cdot f(X))_{n \geq 0}$ is uniformly integrable and thus we can pass to the limit $n \to \infty$ in the previous to obtain
$$
\lim_{n \to \infty} \mathbb E[|Y_n - T(X)|^2] \leq 2 \mathbb E[|f(X) - T(X)|^2] + 2 \delta^2 \leq 4 \delta^2.
$$
Hence the result follows.
\end{proof}
\begin{Rem}\label{rem:nondiff}
The previous result is in fact an equivalence, as the H-differentiability of $\Phi$ implies the uniqueness of an optimal coupling, which has to be deterministic. I refer to Alfonsi and Jourdain \citep{alfonsi} Theorem 3.2 for a proof of this reciprocal. I encourage the reader to verify as an exercise, when $\mo = \R$ and $\rho$ is the normal distribution, that $W_2^2(\cdot,\rho)$ is not differentiable at $\delta_0$.
\end{Rem}

\subsection*{Bibliographical comments}
The starting point of the development on horizontal derivatives is often taken as the introduction of a formal calculus by Otto in \citep{ottoporous}. The study of gradient flows in Wasserstein spaces grew rapidly, and with it more precise developments of the horizontal derivative (sometimes called Wasserstein gradient). A sub-differential calculus can be found in the book of Ambrosio, Gigli and Savaré \citep{ags}. The latter is much more in the spirit of the geometric derivative presented in Section \ref{sec:geo}. 

A breakthrough came through the lift of Lions (Section \ref{sec:lift}), which made possible to define a notion of derivative of a function $U : \mpt \mapsto \R$ through the notion of Fréchet differentiability of its lift $\mathcal U: \mathbb H \mapsto \R$. This made possible to show the structure of the gradient of the lift as a consequence of the dependence on the law only, Lemma \ref{lemma:structureofgradient}. These developments, mainly motivated by the study of MFG, were presented by Lions in his lectures \citep{lions20072008,lions20082009} and can be found in the notes of Cardaliaguet \citep{cardaliaguet2010notes}. More precise developments can be found in Chapter 5 of the book of Carmona and Delarue \citep{carmona2017probabilistic}, especially on the link between the vertical and the horizontal notions. As mentioned above, the equivalence between the restriction or not to optimal couplings was established by Gangbo and Tudorascu in \citep{gangbotudorascu}. Theorem \ref{thm:structurelocal} is a local refinement of previous results which were stated for functions which were differentiable at more than one single point. It can be found in  \citep{gangbotudorascu} or in Alfonsi and Jourdain \citep{alfonsi}. The more restrictive results which were generalized by it can be found in \citep{lions20082009,cardaliaguet2010notes,carmona2017probabilistic}.

See also the thesis of Ceccherini Silberstein \citep{giacomothesis} for general notions of derivatives on Wasserstein spaces on smooth manifolds.

\newpage

\section{Convexity, monotonicity and sub-differential calculus} 
This section is concerned with two notions of convexity, each linked with a different geometry on spaces of measures. One is the usual one and is related to the notion of vertical differentiability, while the other one is concerned with Lagrangian displacements and thus naturally linked to horizontal differentiability. \emph{En passant}, I will also give associated definitions and properties of monotone mappings, namely because of the fundamental role that monotonicity plays in the study of mean field game master equations.

After exploring notions of differentiability in the two previous sections, I will introduce here the weaker notion of sub-differentiability. This notion is particularly helpful in optimization and it will play a strong role in the study of Hamilton-Jacobi-Bellman equations in Part III. The main link between convexity and sub-differentiability is through the classical result: \emph{a convex function is everywhere sub-differentiable in the interior of its domain.} This result is true in Banach spaces, or more generally in topological vector spaces in which the geometric version of Hahn-Banach separation Theorem holds.

\subsection{Flat convexity and monotonicity}
 In this section, $\mo$ is a measurable set.\\

The notion of flat convexity, as the one of flat differentiability, is directly inherited from the structure of vector space of $\mmo$. It can be stated as follows.
\begin{Def}
\begin{itemize}
\item A set $X \subset \mmo$ is flat convex, or simply convex, if, for all $µ,\nu \in X$, $t \in [0,1]$,
$
(1-t) µ + t \nu \in X.
$
\item Let $X\subset \mmo$ be a convex set. The function $U: X \mapsto \R\cup\{+ \infty\}$ is convex if its epigraph $\{(v,µ) \in \R\times X | v \geq U(µ)\}$ is convex for the product convexity on $\R\times X$ issued from the flat convexity, or equivalently, if  $µ,\nu \in X$, $t \in [0,1]$
$$
U((1-t)µ + t \nu) \leq (1-t)U(µ) + t U(\nu).
$$
\item A function $U$ is concave if $-U$ is convex.
\end{itemize}
\end{Def}

Standard properties of convex functions and sets hold in this case, and thus I state without proof the usual result.
\begin{Prop}
\begin{itemize}
\item An intersection of convex sets is convex.
\item For any bounded measurable $f: \mo \mapsto \R$, the map $µ \mapsto \langle f, µ\rangle$ is convex.
\item  Let $X\subset \mmo$ be a convex set, then a point-wise supremum of convex functions defined on $X$ is convex.
\end{itemize}
\end{Prop}

Furthermore, this notion is also compatible with the discretization evoked in Section \ref{sec:disccomp}.
\begin{Prop}
Let $U: \mmo \mapsto \R$ be flat convex and consider, for a finite family of measures $(µ_i)_{i \in I}$, the function $v$ defined in \eqref{defvNflat}. Then $v$ is convex on $\R^I$.
\end{Prop}
The immediate proof is left as an exercise.\\

As in finite dimensional cases, we can define an associated notion of monotonicity.
\begin{Def}
Let $A\subset \mmo$ be a convex set. A set valued function $F: A \rightrightarrows (\mmo)'$ is said to be
\begin{itemize}
\item (flat) monotone if for all $µ,\nu \in A$, $f\in F(µ), g\in F(\nu)$
$$
\langle f-g,µ-\nu \rangle \geq 0,
$$
\item strictly (flat) monotone if for all $µ,\nu \in A$, with $µ \ne \nu$, $f\in F(µ), g\in F(\nu)$
$$
\langle f-g,µ-\nu \rangle >0.
$$
\end{itemize}
\end{Def}
The previous applies naturally when $F$ is a single valued mapping. The next classical result is stated without proof.
\begin{Prop}
Let $U: \mmo \mapsto \R$ be a V-differentiable function. Then $\nabla_µ U: \mmo \mapsto (\mmo)'$ is flat monotone if and only if $U$ is convex.
\end{Prop}

\subsection{Flat super-differentials}
Let $\mo$ be a measurable space. Since $(\mmo,|\cdot|)$ is a Banach space, a super-differential calculus can be defined on it. Here, I recall such a classical notion. To obtain results such as the sub-differentiability of convex functions, I will not be able to restrict my attention to elements in $(\mmo)'$ which can be represented as bounded measurable functions, as I did for vertical derivatives. 
\begin{Def}
Let $U: \mmo \mapsto \R$ be a usc function. We say that $\ell\in (\mmo)'$ belongs to the flat super-differential of $U$ at $µ$ and we write $\ell \in \partial_{flat}^+U(µ)$ when there exists a modulus of continuity $\omega(\cdot)$ such that for all $\nu \in \mmo$
$$
U(\nu) \leq U(µ) + \ell(\nu-µ) + \omega(|\nu-µ|)|\nu-µ|.
$$
By definition, $\partial^+_{flat}U(µ)$ is the set of such $\ell$.

For a lsc function $V: \mmo \mapsto \R$, the set $\partial^-_{flat}V(µ)$ is equal to $-\partial^+(-V)(µ)$.
\end{Def}
\begin{Rem}
Following Section \ref{sec:W1}, a similar notion can be developed if we are interested in $(E,\|\cdot\|_{W^{-1,1}})$, which has the advantage of having a more convenient dual.
\end{Rem}
These super-differentials enjoy several standard properties such as the ones listed in the following proposition.
\begin{Prop}
Let $U: \mmo \mapsto \R$ be a usc function.
\begin{enumerate}[(i)]
\item If $U$ attains a maximum at $µ\in \mmo$, then $0 \in \partial_{flat}^+U(µ)$.
\item If $U$ is V-differentiable at $µ \in \mmo$, then $\partial_{flat}^+U(µ) = \partial_{flat}^-U(µ) = \nabla U(µ)$.
\item If $U$ is concave, then for all $µ \in \mmo$, $\partial_{flat}^+U(µ)\ne \emptyset$. Furthermore, $µ \mapsto -\partial_{flat}^+U(µ)$ is a monotone operator. 
\end{enumerate}
\end{Prop}
\begin{proof}
Points $(i)$ and $(ii)$ are trivial. We do not have the usual other implication in $(ii)$ since I restricted the notion of differentiability to elements which are bounded measurable functions. Point $(iii)$ follows from the Hahn-Banach Theorem. Indeed, since $U$ is concave, the epigraph of $-U$ is convex, and because $U$ is defined everywhere and usc, the interior of this epigraph is non-empty. Separating $(-U(µ),µ)$ from this interior yields that $\partial_{flat}^+U(µ)\ne \emptyset$. The second part of $(iii)$ follows from an immediate computation.
\end{proof}
\begin{Rem}\label{rem:loctoglob}
As usual in Banach spaces, when we are given an element in the sub-differential of a flat convex function, we know that we can choose the associated modulus of continuity to be $0$, given that we are in the interior of its domain of definition. Thus for flat convex functions, sub-differentiability is in fact a global property. Note that the previous proof directly yields such a global element.
\end{Rem}

\subsection{Coupling convexity}
In this section, $\mo$ is the closure of a smooth convex domain of $\R^d$. I introduce here a not so usual notion of convexity on $\mpo$, for which I will quickly make links with more standard notions.
\subsubsection{Basic definitions and examples}
\begin{Def}
\begin{itemize}
\item A set $A \subset \mpo$ is coupling convex if, for all $µ, \nu \in A$, for all $t \in [0,1], \gamma \in \Pi(µ,\nu)$, it holds that $((1-t)\pi_1 + t\pi_2)_\#\gamma \in A$.
\item Let $A\subset \mpo$ be a coupling convex set. A function $U: A \mapsto \R\cup\{+ \infty\}$ is coupling convex if its epigraph $\{(v,µ) \in \R\times A | v \geq U(µ)\}$ is convex for the product convexity on $\R\times A$ issued from the coupling convexity on $A$, or equivalently, if for all $µ,\nu \in A$, $\gamma \in \Pi(µ,\nu)$, $t \in [0,1]$
$$
U(((1-t)\pi_1 + t\pi_2)_\#\gamma) \leq (1-t)U(µ) + t U(\nu).
$$
\item A function $U$ is concave if $-U$ is convex.
\end{itemize}
\end{Def}
This notion of convexity enjoys some non-usual properties listed in the next example.
\begin{Ex}
\begin{itemize}
\item The set $\{µ\}$ is coupling convex if and only if $µ$ is a Dirac mass.
\item As a consequence: let $µ_0 \in \mppo$ and let $U: \mppo \mapsto \R\cup \{+ \infty\}$ defined as $0$ on $µ_0$ and $+\infty$ everywhere else. Then $U$ is convex along couplings if and only if $µ_0$ is a Dirac mass.
\item Assume $\mo = [-1,1]$. The function $µ \mapsto \int_\mo -|x|^2 µ(dx)$ is linear (seen as the restriction of a linear function on $\mmo$), but not coupling convex. For instance consider $µ = \delta_{-1}$ and $\nu = \delta_1$.
\end{itemize}
\end{Ex}
At this point, it would be natural to be worried about the existence of non trivial coupling convex functions. The next result aims at being reassuring in this direction, as it states that many functions are indeed coupling convex.
\begin{Prop}\label{prop:couplconv1} The following assertions hold.
\begin{itemize}
\item If $\phi : \mo \mapsto \R$ is a convex function bounded from below, then $µ \mapsto \langle \phi,µ\rangle$ is coupling convex.
\item The pointwise supremum of coupling convex functions is coupling convex.
\end{itemize}
\end{Prop}
\begin{proof}
Let $\phi : \mo \mapsto \R$ is a convex function bounded from below, and $µ,\nu \in \mpo$, $\gamma \in \Pi(µ,\nu), t \in [0,1]$. Then
$$
\ba
\int_{\mo^2} \phi((1-t)x + t y)\gamma(dx,dy) &\leq \int_{\mo^2}(1-t)\phi(x) + t \phi(y)\gamma(dx,dy)\\
&= (1-t) \int_\mo\phi(x)µ(dx) + t \int_\mo \phi(y)\nu(dy),
\ea
$$
hence the first point follows, even if the previous integrals are equal to $+ \infty$. Concerning the second point, the proof is the usual one: Let $(U_i)_{i \in I}$ be a family of coupling convex functions and define for all $µ \in \mpo$, $V(µ) := \sup_{i \in I} U_i(µ)$. Then for $µ,\nu \in \mpo$, $\gamma \in \Pi(µ,\nu)$ and $t\in [0,1]$
$$
\ba
V(((1-t)\pi_1 + t\pi_2)_\#\gamma) &= \sup_{i \in I}\{U_i(((1-t)\pi_1 + t\pi_2)_\#\gamma)\}\\
&\leq \sup_{i \in I}\{(1-t)U_i(µ)+ tU_i(\nu)\} \leq (1-t)V(µ) + t V(\nu).
\ea
$$
\end{proof}
At this point, the following result should be easy to prove for the interested reader.
\begin{Prop}\label{prop:liftconv}
Let  $p \in [1,\infty)$, $U : \mpprd \mapsto \R\cup\{+\infty\}$ and $\mathcal U$ its lift on $\mathbb L^p$. Then $U$ is coupling convex if and only if $\mathcal U$ is convex (on $\mathbb L^p$).
\end{Prop}

In addition, coupling convexity behaves well with particle approximation of measures as the next simple result shows.
\begin{Prop}\label{prop:easyccproj}
Let $U : \mpo \mapsto \R$ be coupling convex. Then, for any $N \geq 1$, the function $u: \mo^N \mapsto \R$ defined by \eqref{def:uproj} is convex on $\mo^N$.
\end{Prop}
\begin{proof}
It suffices to remark that for $x,y \in \mo^N, t \in [0,1]$, $\gamma= \frac1N \sum_{i=1}^N\delta_{(x_i,y_i)}$ is a coupling between $\frac1N \sum_{i=1}^N \delta_{x_i}$ and $\frac1N \sum_{i=1}^N\delta_{y_i}$ and that $u((1-t)x + t y) = U(((1-t)\pi_1 + t \pi_2)_\#\gamma)$.
\end{proof}

\subsubsection{Regularity of coupling convex functions}
Maybe a striking argument to link the previous notion of convexity with the \emph{geometry} of Wasserstein spaces is the fact that convexity along couplings implies regularity in Wasserstein distances. I present two results here, the first one in the case $\mo = \R^d$ and the second one which deals with a more compact case. I do not consider here functions which take the value $+\infty$, even if the case in which $\mo$ is compact can be seen as a case in which $U(µ) = + \infty$ if $µ \notin \mpo$. 
\begin{Prop}
Let $p \in [1,\infty)$, $U: \mpprd \mapsto \R$ be coupling convex and locally bounded (i.e. bounded on bounded sets). Then $U$ is locally Lipschitz continuous, i.e. for any $R > 0$, there exists $C_R > 0$ such that for any $µ, \nu \in \mpprd$ with $M_p(µ),M_p(\nu) \leq R^p$, 
$$
|U(µ)-U(\nu)|\leq C_R \mathcal{W}_p(µ,\nu).
$$
\end{Prop}
\begin{proof}
The proof follows the usual one in vector spaces. For $R> 0$ define $B_R := \{µ \in \mpprd| M_p(µ) \leq R^p\}$. Let $R > 0$ and $µ,\nu \in B_R, µ \ne \nu$, with $\mathcal{W}_p(µ,\nu) < 1$. Then, for any $t \in (0,1)$ there exists $\tilde µ \in \mpprd$ and $\gamma \in \Pi(\nu,\tilde µ)$ such that $µ = ((1-t) \pi_1 +  t\pi_2)_\#\gamma$. Indeed, for this, we consider $\tilde \gamma \in \Pi(µ,\nu)$ and define $\tilde µ := (t^{-1}(\pi_1 -\pi_2) + \pi_2)_\#\tilde \gamma$. Then, by coupling convexity
$$
U(µ)  \leq (1-t)U(\nu) + tU(\tilde µ).
$$
Setting $t = \mathcal{W}_p(µ,\nu)$, we obtain 
$$
U(µ) - U(\nu) \leq (U(\tilde µ) - U(\nu))\mathcal{W}_p(µ,\nu).
$$
Remark now that if we choose $\tilde \gamma$ as an optimal coupling for $\mathcal{W}_p$, then $\tilde µ \in B_{R+1}$. Hence, we obtain
$$
U(µ) - U(\nu) \leq 2\sup_{µ' \in B_{R+1}}\{|U(µ')|\} \mathcal{W}_p(µ,\nu).
$$
Arguing by symmetry, and using the triangle inequality for $µ,\nu$ such that $\mathcal W_p(µ,\nu) \geq 1$, the result is proved by setting $C_R := 2\sup_{µ' \in B_{R+1}}\{|U(µ')|\}$.
\end{proof}

When we are in the case, $\mo \ne \R^d$, one has to be quite careful with the boundary of $\mo$. Consider for instance the case $\mo = [0,1]$ and $U(µ) = -\int_0^1\sqrt{x}µ(dx)$, which is coupling convex, bounded, but not Lipschitz continuous.

\begin{Prop}\label{prop:lipwinf}
Let $\mo$ be the closure of a convex bounded domain of $\R^d$. Let $U : \mpo \mapsto \R$ be a bounded, coupling convex function. Then, for any $\eps > 0$, there exists $C_\eps > 0$ such that for any $µ,\nu \in \mpo$ with support at distance at least $\eps$ of $\partial \mo$, 
$$
|U(µ) - U(\nu)| \leq C_\eps \mathcal{W}_\infty(µ,\nu).
$$
\end{Prop}
\begin{proof}
The proof follows a similar argument. Let $\eps > 0$ and consider $B_\eps:= \{µ \in \mpo | \mathbb P_µ(d(x,\partial \mo) \leq \eps ) = 0\}$. Let $µ,\nu \in B_\eps, µ\ne \nu$. Let $(X,Y)$ be a coupling of $µ,\nu$. Introduce $t= \|X-Y\|_\infty > 0$. Define
$$
X' = Y + \frac{\eps}{2}\frac{X-Y}{t}.
$$
Observe that
$$
 \left\|\frac\eps2\frac{X-Y}{t}\right\|_\infty \leq \frac\eps2.
$$
Hence, $\mathcal L(X') \in B_{\frac \eps2}$ and the rest of the proof follows as above.
\end{proof}

\begin{Rem}
This type of Lipschitz regularity might seem quite disappointing: why not obtain a stronger Lipschitz estimate, in $\mathcal{W}_1$ norm for instance ? It so happens that there are some monsters in the set of coupling convex functions over $\mpo$ as the next example shows. It is somehow coherent with the fact that, if we interpret $\mpo$ as a subset of some $\mpprd$, then, it contains no open ball in $\mathcal{W}_p$ distances, only open balls in $\mathcal{W}_\infty$.
\end{Rem}

The following result highlights the possible irregularity that coupling convex functions might have.
\begin{Prop}
Let $\mo$ be the closure of a bounded convex domain of $\R^d$ and $\ell$ be a linear (continuous) form on $L^\infty(\Omega,\R^d)$. Define $U: \mpo\mapsto \R$ with
$$
U(µ) := \sup_{X \sim µ}\ell(X).
$$
Then $U$ is coupling convex on $\mpo$.
\end{Prop}
\begin{proof}
First of all, note that $U$ is well defined, since, for all $µ \in \mpo$, $X\sim µ$, $\|X\|_\infty \leq K$ where $K$ depends only on $\mo$. Hence, for $C > 0$ such that
$$
\forall X \in L^\infty, \quad \ell(X)\leq C \|X\|_\infty,
$$
we have $U(µ) \leq CK$. Now, let $µ,\nu \in \mpo$, $\gamma \in \Pi(µ,\nu)$ and $t \in [0,1]$. Let $\eps > 0$ and consider $Z_\eps$ such that $Z_\eps \sim ((1-t)\pi_1 + t \pi_2)_\#\gamma$ and 
$$
U(((1-t)\pi_1 + t \pi_2)_\#\gamma) \leq \ell(Z_\eps) + \frac{\eps}{2}.
$$
Consider now a couple of random variables $(X_\eps,Y_\eps)\sim \gamma$ such that $\|Z_\eps - (1-t)X_\eps - t Y_\eps\|_\infty \leq \frac{\eps}{2C}$. It then follows that
$$
\ba
U(((1-t)\pi_1 + t \pi_2)_\#\gamma)& \leq \ell((1-t)X_\eps + t Y_\eps) + C\|Z_\eps - (1-t)X_\eps - t Y_\eps\|_\infty + \frac\eps2\\
& \leq (1-t)\ell(X_\eps) + t \ell(Y_\eps) + \eps\\
& \leq (1-t) U(µ) + t U(\nu) + \eps.
\ea
$$
Since $\eps > 0$ is arbitrary, we conclude that $U$ is indeed coupling convex.
\end{proof}

\subsection{Horizontal super-differentials}\label{sec:superdiff}
In this section, $\mo$ is still the closure of a smooth convex domain of $\R^d$, and I focus on the case of super-differentials of functions on $\mppo$ for $1 < p < \infty$. For most of what follows, the convexity assumption on $\mo$ can be relaxed quite easily, but then the formalism of manifolds introduced in Section \ref{sec:manifold} is more natural. The following notion of super-differential is much more involved and I start by explaining briefly why a naive definition is not appropriate, even if the main motivation shall come later on, when trying to prove comparison of viscosity solutions of HJB equations on $\mppo$, in Part III.

\subsubsection{A first attempt}\label{sec:badsuper}
Let $U: \mppo \mapsto \R$ be a usc function. A first natural attempt to define horizontal super-differentials (which I will simply call super-differential in this section) is to take the definition of the horizontal differential and only ask for an inequality, i.e. we could say that $\phi$ is in the super-differential of $U$ at $µ$ if there exists a modulus of continuity $\omega$ such that for any $\nu \in \mppo$, $\gamma \in \Pi(µ,\nu)$
\be\label{eq:badsuper}
U(\nu) \leq U(µ) + \int_{\mo^2}\phi(x)\cdot(y-x)\gamma(dx,dy) + C_p(\gamma)\omega(C_p(\gamma)).
\ee
This point of view is too restrictive. Even if it is clear that such functions $\phi \in L^{p'}_µ$ are clearly elements of the super-differential, there are other ones that we shall need later on, especially in Part III.  

The best arguments in favor of seeking a more general notion are to come later on. At this time, we can only argue using the interpretation of the lift of Section \ref{sec:lift}. Theorem \ref{thm:structurelocal} states that if $\mathcal U: \mathbb H$ depends only on the law of its argument and is differentiable at some $X \in \mathbb H \mapsto \R$, then $\nabla \mathcal U(X)$ is a deterministic function of $X$. Such a statement has no reason to hold at the level of super-differential: if $Z \in \partial^+\mathcal U(X)$, then $Z$ is not necessarily a deterministic function of $X$. Thus we shall keep track of this "extra" randomness, and the objective of what follows is to explain how we can do it.

\subsubsection{Extending the super-differentials}
In this section, let $p \in (1,\infty)$.
In order not to miss essential elements of super-differentials, I shall use the following definition, which can be seen as an extension of \eqref{eq:badsuper}, which we could have first tried to use. In the following, for any $\Gamma \in \mathcal P(\mo\times \R^d\times \mo)$, I use the notation $\tilde C_p(\Gamma) = C_p((\pi_1,\pi_3)_\#\Gamma)$.
\begin{Def}
Let $U: \mppo\mapsto \R$ be a usc function. A measurable function $\psi : \mo \mapsto \mathcal P(\R^d)$ is an element of the super-differential of $U$ at some $µ \in \mppo$ if 
\begin{itemize}
\item $$
\int_\mo \int_{\R^d}|z|^{p'}\psi_x(dz)µ(dx) < \infty.
$$
\item there exists a modulus of continuity $\omega$ such that for any $\nu \in \mppo$, $\Gamma \in \Pi(µ(dx)\psi_x(dz),\nu(dy))$
\be\label{eq:superdiff}
U(\nu) \leq U(µ) + \int_{\mo\times \R^d\times \mo}z\cdot(y-x)\Gamma(dx,dz,dy) + \tilde C_p(\Gamma)\omega(\tilde C_p(\Gamma)).
\ee
\end{itemize}
In this case, we say that $U$ is super-differentiable at $µ$ and write $\psi \in \partial^+U(µ)$, while $\partial^+U(µ)$ is the set of such elements. The sub-differential of an lsc $V: \mppo \mapsto \R$ is defined as $\partial^-V(µ) = - \partial^+(-V)(µ):= \{ \psi: \mo \mapsto \mathcal P(\R^d) | x \mapsto (-Id)_\# \psi_x \in \partial^+(-V)(µ)\}$.
\end{Def}
\begin{Rem}
Thanks to the disintegration Theorem, considering $\psi$ as above is equivalent to considering $\gamma \in \mathcal P_{pp'}(\mo\times \R^d)$ such that $\gamma(dx,dy) = µ(dx)\psi_x(dy)$. Hence, if $\psi \in \partial^+U(µ)$, I sometimes write $\gamma \in \partial^+U(µ)$ for the associated coupling $\gamma$.
\end{Rem}

The previous definition allows one to mix, or randomize, the value of the element in the super-differential over $x \in \mo$, namely in such a way that the law of this randomization is $\psi_x(dz)\in \mathcal P(\R^d)$. The presence of this (quite annoying notation wise) coupling $\Gamma$ highlights that we allow general interdependence between the randomization of this element in the super-differential and the direction of displacement $(y-x)$ which is drawn according to $(\pi_1,\pi_3)_\#\Gamma$. 

If this definition seems quite monstrous to the reader (which it probably should at first sight), she or he can find comfort in the fact that, if we are looking at $\mo = \R^d$, it is simply the translation of its meaning at the level of the lift in the space of random variables as the next result shows.

\begin{Prop}\label{prop:superdifflift}
Let $p \in (1,\infty)$, $U: \mpprd \mapsto \R$ a usc function and consider $\mathcal U: \mathbb L^p \mapsto \R$ defined by $\mathcal U(X) = U(\mathcal L(X))$. Then
\begin{enumerate}
\item For all $X \in \mathbb L^p$, $Z \in \partial^+ \mathcal U(X)$, we have $\psi \in \partial^+U(µ)$ where $\mathcal L(X,Z) = µ(dx)\psi_x(dz)$.
\item For all $µ \in \mpprd, \psi \in \partial^+ U(µ)$, for any $(X,Z) \in \mathbb L^p \times \mathbb L^{p'}$ of law $µ(dx)\psi_x(dz)$, $Z \in \partial^+\mathcal U(X)$.
\end{enumerate}
\end{Prop}
I used the notation $\partial^+ \mathcal U(X)$ for the usual Fréchet super-differential of $\mathcal U$ at $X \in \mathbb L^p$, which is thus a subset of $\mathbb L^{p'}$.
\begin{proof}
\begin{enumerate}
\item 
Let $Z \in \partial^+\mathcal U(X)$, then there exists a modulus of continuity $\omega$ such that for all $Y \in \mathbb L^p$
\be\label{eq:superlift}
\mathcal U(Y) \leq \mathcal U(X) + \mathbb E[Z\cdot(Y-X)] + \omega(\|Y-X\|_p)\|Y-X\|_p.
\ee
Furthermore, for any $(X',Z')\sim \mathcal L((X,Z))$, we consider $\tau^n$ a bi-measurable mapping given by Theorem \ref{thm:proba} such that $\|(X,Z)-(X'\circ \tau^n,Z'\circ\tau^n)\|_\infty \leq n^{-1}$. For any $Y \in \mathbb L^p$
$$
\ba
&\mathcal U(Y) - \mathcal U(X') - \mathbb E[Z'\cdot(Y-X')]  -  \omega(\|Y-X'\|_p)\|Y-X'\|_p\\
& = \mathcal U(Y\circ\tau^n) - \mathcal U(X')- \mathbb E[Z'\circ\tau^n\cdot(Y\circ\tau^n-X'\circ\tau^n)]  - \omega(\|Y-X'\|_p)\|Y-X'\|_p\\
&\leq \mathbb E [Z\cdot(Y\circ\tau^n - X)] - \mathbb E[Z'\circ\tau^n\cdot(Y\circ\tau^n-X'\circ\tau^n)]\\
&\quad + \omega(\|Y\circ\tau^n-X\|_p)\|Y\circ\tau^n-X\|_p -\omega(\|Y-X'\|_p)\|Y-X'\|_p\\
&\quad\underset{n \to \infty}{\longrightarrow} 0.
\ea
$$
 Hence, for any such $(X',Z')$, $Z' \in \partial^+\mathcal U(X')$. Thus, denote $\gamma = \mathcal L(X,Z)$ and $\Gamma \in \Pi(\gamma,µ')$ for $µ' \in \mpprd$. Consider a random variable $(X'',Z'',Y'')$ of law $\Gamma$. With what we just proved, $Z'' \in \partial^+\mathcal U(X'')$ (with modulus $\omega$), which can be stated precisely as
$$
U(\mu') \leq U(µ) + \int_{(\R^d)^3}z\cdot(y-x)\Gamma(dx,dz,dy) + \tilde C_p(\Gamma)\omega(\tilde C_p(\Gamma)).
$$

\item The other side is more immediate. Let $\psi \in \partial^+U(µ)$ and denote $\omega$ a modulus of continuity such that for all $\Gamma \in \mathcal P((\R^d)^3)$ with $(\pi_1,\pi_2)_\#\Gamma = µ(dx)\psi_x(dz)$
$$
U((\pi_3)_\#\Gamma) \leq U(µ) + \int_{(\R^d)^3}z\cdot(y-x)\Gamma(dx,dz,dy) + \tilde C_p(\Gamma)\omega(\tilde C_p(\Gamma)).
$$
Let $(X,Z) \in \mathbb L^p\times \mathbb L^{p'}$ of law $µ(dx)\psi_x(dz)$ and consider any $Y \in \mathbb L^p$. Evaluating the previous inequality on $\Gamma = \mathcal L(X,Z,Y)$ yields the result.
\end{enumerate}
\end{proof}

This notion of super-differential indeed extends the previous one (of Section \ref{sec:badsuper}) in the following sense.
\begin{Prop}
Let $U: \mppo \mapsto \R$ be a usc function and consider $µ \in \mppo$, $\phi \in L^{p'}_µ$ such that for any $\nu \in \mppo$, $\gamma \in \Pi(µ,\nu)$, \eqref{eq:badsuper} holds and define $\psi: x \mapsto \delta_{\phi(x)}$. Then $\psi$ satisfies \eqref{eq:superdiff}.
\end{Prop}
The proof is trivial so I omit it. Furthermore, there is also a sort of reverse of the previous as the next result shows.
\begin{Prop}\label{prop:barycentric}
Let $U: \mpt \mapsto \R$ be a usc function and consider $µ \in \mpt, \psi \in \partial^+U(µ)$, then $x \mapsto \delta_{\int_{\R^d}z\psi_x(dz)}$ satisfies \eqref{eq:badsuper}. 
\end{Prop}
I stated this result only in the case $p =2, \mo = \R^d$, because it is only proven in the literature in this setting, even if it is quite clear that it can be extended to more general settings. The proof has been done by Gangbo and Tudorascu \citep{gangbotudorascu}.\\

In the previous, given $\psi : \R^d \mapsto \mathcal P(\R^d)$, the element $x \mapsto \int_{\R^d}z\psi_x(dz)$ is often called the \textit{barycentric projection} of $\psi$.\\

Standard results concerning super-differentials are still valid for this horizontal super-differentials. The following presentation is restricted to the case $\mo = \R^d$ to avoid some technicalities in the proof.

\begin{Prop}\label{prop:propsuper}
Let $p \in (1,\infty)$ and $U : \mpprd \mapsto \R$.
\begin{enumerate}[(i)]
\item If $V: \mpprd \mapsto \R$ is H-differentiable at $µ$, then 
$$
\psi \in \partial^+U(µ) \Leftrightarrow x \mapsto (Id + DV(µ)(x))_\#\psi_x(dz) \in \partial^+(U+V)(µ).
$$
\item Let $V: \mpprd \mapsto \R$. Then, for all $\psi\in \partial^+U(µ), \phi\in \partial^+V(µ)$, $\Theta \in \mathcal P((\R^d)^3)$ such that $(\pi_1,\pi_2)_\#\Theta = µ(dx)\psi_x(dz)$ and $(\pi_1,\pi_3)_\#\Theta = µ(dx)\phi_x(dz)$, it holds that $(\pi_1,\pi_2 + \pi_3)_\#\Theta \in \partial^+(V+U)(µ)$.
\item If $\psi \in \partial^+U(µ)$, then for all $\lambda \geq 0$, $(x \mapsto (\lambda Id)_\#\psi_x(dz)) \in \partial^+(\lambda U)(µ)$.
\item Let $\psi,\phi \in \partial^+U(µ)$, then, for any $\Theta \in \mathcal P((\R^d)^3)$ such that $(\pi_1,\pi_2)_\#\Theta = µ(dx)\psi_x(dz)$ and $(\pi_1,\pi_3)_\#\Theta = µ(dx)\phi_x(dz)$, it holds that for all $t \in [0,1]$, $(\pi_1,(1-t)\pi_2 + t\pi_3)_\#\Theta \in \partial^+U(µ)$.
\end{enumerate}
\end{Prop}
\begin{Rem}
Those results are standard if we simply lift the problem to $\mathbb L^p$. I provide a proof at the level of $\mpprd$ to hint how it can be generalized to other cases than $\mo = \R^d$.
\end{Rem}
\begin{Rem}
The last property is the equivalent of the convexity of the sub/super-differentials in this coupling convexity.
\end{Rem}
\begin{proof}
\begin{enumerate}[(i)]
\item Let $V: \mpprd \mapsto \R$ be horizontally differentiable. It then follows that there exists a modulus of continuity $\omega$ such that for any $\nu \in \mpprd, \gamma \in \Pi(µ,\nu)$
$$
\left|V(\nu) - V(µ) - \int_{(\R^d)^2}DV(µ)(x)\cdot(y-x)\gamma(dx,dy)\right| \leq C_p(\gamma)\omega( C_p(\gamma)).
$$
 Consider $\psi \in \partial^+U(µ)$. Once again, without loss of generality, we can assume that the modulus of continuity used to express the super-differentiability of $U$ at $µ$ is $\omega$. Since $\psi \in \partial^+U(µ)$, for any $\nu \in \mpprd$, $\Gamma \in \Pi(µ(dx)\psi_x(dz),\nu(dy))$
$$
U(\nu) \leq U(µ) + \int_{ (\R^d)^3}z\cdot(y-x)\Gamma(dx,dz,dy) + \tilde C_p(\Gamma)\omega(\tilde C_p(\Gamma)).
$$
Choosing $\gamma = (\pi_1,\pi_3)_\#\Gamma$ and using the two previous inequalities yields 
$$
\ba
(U+V)(\nu)  \leq (U+V)(µ) + \int_{( \R^d)^3}&(z+ DV(µ)(x))\cdot(y-x)\Gamma(dx,dz,dy)\\
& + 2\tilde C_p(\Gamma)\omega(\tilde C_p(\Gamma)).
\ea
$$
The reverse inequality works in a similar manner. Once again, we assume that the super-differentiability of $(U+V)$ at $µ$ can be expressed with the same modulus of continuity $\omega$. Then, the same computation as above yields that if $\psi \in \partial^+(U+V)(µ)$, then $x \mapsto (Id - DV(µ)(x))_\#\psi_x(dz) \in \partial^+U(µ)$.
\item Let $\psi ,\phi, \Theta$ be as in the statement. For any $\nu \in \mpprd$, $\bar \Gamma \in \Pi(\Theta,\nu)$, we can express the super-differentiability of $U$ and $V$ on respectively $(\pi_1,\pi_2,\pi_4)_\#\bar \Gamma$ and $(\pi_1,\pi_3,\pi_4)_\#\bar \Gamma$. Assuming that the two moduli of continuity are $\omega$ (which we can always do without loss of generality) and summing the two inequalities leads to
$$
\ba
(U+V)(\nu)  \leq (U+V)(µ) + \int_{(\R^d)^4}&(z+z')\cdot(y-x)\bar\Gamma(dx,dz,dz',dy)\\
& + 2\tilde C_p(\Gamma^U)\omega(\tilde C_p(\Gamma^U)).
\ea
$$
It is then almost possible to conclude that $µ(dx)\psi^*_x(dz):=(\pi_1,\pi_2+\pi_3)_\#\bar\Gamma \in \partial^+(V+U)(µ)$, as it only remains to remark that all couplings $\Gamma \in \Pi(µ(dx)\psi^*_x(dz),\nu)$ can be obtained in this way, which is true. 
\item This point is trivial.
\item The proof follows from the previous points used on the maps $(1-t)U$ and $tU$ for $t \in [0,1]$.
\end{enumerate}
\end{proof}

I now turn to a result whose statement should look familiar to people used to super-differentials, but which turns out to involve more sophisticated arguments than in the usual case. This time, I present several cases of base spaces $\mo$. 

\begin{Prop}\label{prop:alfonsiO}
Let $p \in (1,\infty)$ and $U: \mpprd\mapsto \R$ be (horizontally) sub and super-differentiable at $µ \in \mpprd$. Then, $U$ is H-differentiable at $µ$ and $\partial^-U(µ) = \partial^+U(µ) = \{ x \mapsto \delta_{D_µU(µ)(x)}\}$.
\end{Prop}
\begin{proof}
The first part of the proof consists in showing that if $\partial^+U(µ)\ne \emptyset$ and $\partial^-U(µ)\ne \emptyset$, then the two sets contain in fact a unique element which is the same. Remark that this is a reformulation of the fact that the lift of $U$ is Fr\'echet differentiable at $µ$. Concluding using Theorem \ref{thm:structurelocal}, we obtain the result.

 Take $\psi^+ \in \partial^+ U(µ)$ and $\psi^-\in \partial^-U(µ)$. Note that we can choose, without loss of generality, the two moduli of continuity $\omega$ associated to $\psi^+$ and $\psi^-$ to be the same. Then, for any $\nu \in \mpprd$, consider two elements $\Gamma^\pm \in \Pi(µ(dx)\psi_x^\pm(dz),\nu)$. By definition of sub/super differentiability, 
 $$
U(\nu) \leq U(µ) + \int_{ (\R^d)^3}z\cdot(y-x)\Gamma^+(dx,dz,dy) + \tilde C_p(\Gamma^+)\omega(\tilde C_p(\Gamma^+)),
$$
$$
U(\nu) \geq U(µ) + \int_{(\R^d)^3}z\cdot(y-x)\Gamma^-(dx,dz,dy) - \tilde C_p(\Gamma^-)\omega(\tilde C_p(\Gamma^-)).
$$
Remark now that we can always choose $\Gamma^\pm$ such that $(\pi_1,\pi_3)_\#\Gamma^+ = (\pi_1,\pi_3)_\#\Gamma^- =: \gamma$. Furthermore, by the Gluing Lemma, we can obtain a probability measure $\Theta \in \mathcal P((\R^d)^4)$ such that  $(\pi_1,\pi_2,\pi_4)_\# \Theta = \Gamma^+$ and $(\pi_1,\pi_3,\pi_4)_\# \Theta = \Gamma^-$. Then, from the two previous inequalities, we obtain
$$
 \int_{ (\R^d)^4}(z^- - z^+)\cdot(y-x)\Theta(dx,dz^+,dz^-,dy) \leq 2  C_p(\gamma)\omega( C_p(\gamma)).
$$
We now know along which coupling $\gamma$ we should test the super/sub-differentiability of $U$. Indeed, setting $\gamma':= (\pi_1,\pi_1 + \eps(\pi_3 - \pi_2)|\pi_3 -\pi_2|^{p'-2})_\#\Theta$ for some $\eps > 0$, we obtain, using the last inequality,
\be\label{eq:1232}
\eps\int_{ (\R^d)^4}|z^- - z^+|^{p'}\Theta(dx,dz^+,dz^-,dy) \leq 2  C_p(\gamma')\omega( C_p(\gamma')).
\ee
It is indeed possible to use the previous coupling since: i) the choice is not circular as $\gamma'$ only depends on $(\pi_1,\pi_2,\pi_3)_\#\Theta$ which can be fixed a priori, ii)
$$
C^p_p(\gamma') = \eps^p\int_{(\R^d)^4} |z^--z^+|^{p(p'-1)}\Theta(dx,dz^+,dz^-,dy) =: \eps^pK^p < \infty,
$$
where the finiteness of $K$ comes from the growth assumption in the definition of sub/super-differentials. Using $p(p'-1) = p'$, we then obtain from \eqref{eq:1232} and the previous definition of $K$
$$
\eps K^p \leq 2\eps K \omega( \eps K).
$$
Since the previous holds for any $\eps > 0$, it follows that $K = 0$ and thus the first part of the proof is obtained. Thanks to Theorem \ref{thm:structurelocal}, as mentioned at the beginning of the proof, the result thus holds. 

\end{proof}

The previous result can be extended to more general domains, as the next result shows.

\begin{Prop}
Let $\mo$ be the closure of a smooth convex bounded domain. Let $U: \mpo\mapsto \R$ and $µ \in \mpo$ whose support is at distance $\delta > 0$ of $\partial \mo$. If $U$ is both (horizontally) sub and super-differentiable at $µ$ (in the $\mathcal P_2(\mo)$ sense), then, $U$ is H-differentiable at $µ$ and both $\partial^+U(µ)$ and $\partial^-U(µ)$ are equal to $\{x \mapsto \delta_{D_µU(µ)(x)}\}$.
\end{Prop}
\begin{proof}
The argument of the first step of the previous proof can be easily adapted to this setting. Indeed, instead of going in the direction given by $\eps (z^--z^+)|z^--z^+|^{p'-2}$, we simply evaluate both the sub and super-differentiability in the direction $\eps (z^--z^+)|z^--z^+|^{-1}$ (which is taken as $0$ when $z^- =z^+$), and obtain that
$$
\int_{(\R^d)^4} |z^--z^+|\Theta(dx,dz^+,dz^-,dy) = 0.
$$
Thus we deduce that for some $\psi$, $\partial^+ U(µ) = \partial^-U(µ) = \{\psi\}$. In particular, for any $\nu \in \mpo$, $\Gamma \in \Pi(µ(dx)\psi_x(dz),\nu)$, 
\be\label{ineq2281}
\left|U(\nu) - U(µ) - \int_{\mo\times \R^d\times \mo}z\cdot (y-x) \Gamma(dx,dz,dy)\right|  \leq \tilde C_2(\Gamma)\omega(\tilde C_2(\Gamma)).
\ee
Therefore, it only remains to prove that $\psi$ is deterministic to obtain the $H$-differentiability of $U$ at $µ$. To obtain the result, it suffices to remark that the proof of Theorem \ref{thm:structurelocal} also holds in this case. Indeed, let $E$ be the space of $\mo$ valued random variables (on a standard probability space) and $\mathcal U: E \mapsto \R$, the lift of $U$ defined by $\mathcal U(X) = U(\mathcal L(X))$. Inequality \eqref{ineq2281} implies that for any $X,Y \in E, Z \in \mathbb H$ such that $\mathcal L(X,Z,Y) = \Gamma$, we have
$$
\left|\mathcal U(Y) - \mathcal U(X) - \mathbb E[Z\cdot(Y-X)]\right|  \leq \|X-Y\|_2\omega(\|X-Y\|_2).
$$
Just as in the proof of Theorem \ref{thm:structurelocal}, we can define a mapping $g$ and variations $\xi^\pm$. The only part where the space $\mo$ plays a role is in the admissibility of the variations $\eps \xi^\pm$ we constructed. Because these variations are bounded and $µ$ has support away from $\partial \mo$, we know that for $\eps$ small enough they are admissible. Hence, the same consequence holds and $U$ is indeed H-differentiable at $µ$.
\end{proof}

\subsubsection{Sub-differentiability of coupling convex functions}
The aim of this section is to present what the analogue of the result: "a convex function is sub-differentiable on the interior of its domain" might be in this context of coupling convex functions. As the interior of the domain of a coupling convex function is technical to manipulate, I will present only some partial results here.\\

I start by giving a result in the direction of Remark \ref{rem:loctoglob}: for coupling convex functions, sub-differentiability is a global property.
\begin{Prop}\label{prop:loctoglob}
Let $\mo$ be a convex closed set, $p \in (1,\infty)$ and $U: \mppo\mapsto \R$ be coupling convex. Consider $\psi \in \partial^-U(µ)$. Then, for any $\nu \in \mppo, \Gamma \in \Pi(µ(dx)\psi_x(dz),\nu)$,
$$
U(\nu) \geq U(µ) + \int_{\mo\times\R^d\times \mo}z\cdot(y-x)\Gamma(dx,dz,dy).
$$
\end{Prop}
\begin{proof}
Once again, even in this unusual setting, the proof follows usual arguments. Let $t \in (0,1)$. Using the coupling convexity of $U$ we have
$$
\frac{U(((1-t)\pi_1 + t \pi_3)_\#\Gamma) - U(µ)}{t} \leq U(\nu) - U(µ).
$$
Denoting by $\omega$ the modulus of continuity appearing in $\psi \in \partial^-U(µ)$, we obtain, by testing the sub-differentiability along $(\pi_1,\pi_2,(1-t)\pi_1 + t \pi_3)_\#\Gamma=: \Gamma_t$,
$$
 \int_{\mo\times\R^d\times \mo}z\cdot(y-x)\Gamma(dx,dz,dy) - t^{-1}\omega(\tilde C_p(\Gamma_t))\tilde C_p(\Gamma_t) \leq U(\nu) - U(µ).
$$
Passing to the limit $t \to 0$ yields the result.
\end{proof}

I now prove that coupling convex functions are horizontally sub-differentiable, first in the simpler case of functions on $\mpprd$, and then on $\mpo$.
\begin{Theorem}
Let $p \in (1,\infty)$ and $U: \mpprd \mapsto \R$ be coupling convex and locally bounded. Then $U$ is (horizontally) sub-differentiable everywhere, furthermore, at any point $µ$, the modulus of continuity $\omega$ in the expression of the sub-differentiability can be replaced by $0$. 
\end{Theorem}
\begin{proof}
The proof follows a lifting argument. We consider the lift of $U$ as $\mathcal U: \mathbb L^p\mapsto \R, X \mapsto U(\mathcal L(X))$. The function $\mathcal U$ is convex in $\mathbb L^p$ (in the usual sense) and locally bounded. Hence, it is sub-differentiable everywhere (Proposition 5.2 in \cite{ekeland1999convex}), that is: for any $X \in \mathbb L^p$, there exists $Z \in \mathbb L^{p'}$ such that for all $Y \in \mathbb L^p$
$$
U(\mathcal L(Y)) \geq U(\mathcal L(X)) + \mathbb E[Z\cdot(Y-X)].
$$
The result then follows from Proposition \ref{prop:superdifflift}.
\end{proof}

As it is the case in finite dimension, results such as the one above are not true up to the boundary of the domain of definition of the convex function. For instance $x \mapsto -\sqrt{x}$ is not sub-differentiable everywhere on $[0,\infty)$. In addition, we also have to deal with the fact that, as we saw in Proposition \ref{prop:lipwinf}, we have a priori to work with $\mathcal{W}_\infty$ sub-differentiability, which is not pleasing. Hence, I first present a result under an additional regularity property and then come back to a sort of general sub-differentiability. I insist upon the fact that the latter is given purely for the sake of completeness, I am not aware of any use of such a result at the moment.
\begin{Prop}
Assume that $\mo$ is the closure of a bounded convex domain of $\R^d$. Let $U: \mpo \mapsto \R$ be coupling convex and $\mathcal{W}_p$ Lipschitz continuous for some $1< p < \infty$. Then, for every $µ \in \mpo$ whose support is at distance at least $\eps > 0$ of $\partial \mo$, $U$ is (horizontally) sub-differentiable at $µ$. More precisely, at every such $µ$, there exists a measurable $\psi : \mo \mapsto \mathcal P(\R^d)$ such that
$$
\int_\mo\int_{\R^d}|z|^{p'}\psi_x(dz)µ(dx) < \infty,
$$
and for any $\nu \in \mpo, \Gamma \in \Pi(µ(dx)\psi_x(dz),\nu)$
$$
U(\nu) \geq U(µ) + \int_{\mo\times \R^d\times \mo}z\cdot(y-x)\Gamma(dx,dz,dy).
$$
\end{Prop}
\begin{proof}
Let $U, µ$ be as in the statement, for some $\eps > 0$. For $N > 1$, define $U^N: \mo^N \mapsto \R$ by 
$$
U^N(x_1,\dots,x_N) = U\left(\frac1N \sum_{i =1}^N \delta_{x_i}\right).
$$
Recall that from Proposition \ref{prop:easyccproj}, $U^N$ is convex on $\mo^N$. Since it is finite everywhere, for any $(x_1,\dots,x_N)$ in the interior of $\mo^N$, there exists $(y_1,\dots,y_N) \in (\R^d)^N$ such that for all $(x'_1,\dots,x'_N) \in \mo^N$
\be\label{eq:gg}
U^N(x'_1,\dots,x'_N) \geq U^N(x_1,\dots,x_N) + \frac1N\sum_{i =1}^N y_i\cdot(x'_i - x_i).
\ee
In particular, if for all $1 \leq i \leq N$, $d(x_i,\partial \mo) > \eps$, choosing 
$$
x_i'= x_i + \eps \frac{y_i}{|y_i|},
$$
yields
$$
\frac{1}{N}\sum_{i=1}^N|y_i| \leq \frac{U^N(x'_1,\dots,x'_N)-U^N(x_1,\dots,x_N)}{\eps}.
$$
We now use the $\mathcal{W}_p$ Lipschitz regularity of $U$ to obtain a stronger estimate. Here, we choose alternatively $x'_i = x_i + \delta y_i|y_i|^{p'-2}$, for $\delta > 0$ sufficiently small so that $x' \in (\mo)^N$, we obtain from \eqref{eq:gg} that
$$
\ba
\delta \frac1N\sum_{i=1}^N|y_i|^{p'} &\leq  C \mathcal{W}_p\left(\frac1N\sum_{j =1}^N\delta_{x_j},\frac1N\sum_{j =1}^N\delta_{x'_j}\right)\\
& \leq \delta C \left(\frac{1}{N}\sum_{i =1}^N|y_i|^{p'}\right)^\frac1p.
\ea
$$
where $C$ is the Lipschitz constant of $U$.
In particular, for some (other) $C > 0$ independent of $N$ and $x$,
\be\label{boundy}
\frac1N\sum_{i=1}^N|y_i|^{p'} \leq  C.
\ee
We now consider, thanks to Theorem \ref{thm:glivenko} a sequence $(x_{i,N})_{1 \leq i \leq N}$ valued in $\mo$ such that for all $1 \leq i \leq N$, $d(x_{i,N},\partial \mo) > \eps$ and $\lim_{N \to \infty}\frac1N\sum_{i =1}^N \delta_{x_{i,N}} = µ$. For any $N > 1$, we denote by $(y_{1,N},\dots,y_{N,N})$ an element in the sub-differential of $U^N$ at $(x_{1,N},\dots,x_{N,N})$. Observe that the sequence $\frac1N\sum_{i=1}^N \delta_{(x_{i,N},y_{i,N})}$ is bounded in $\mathcal P_{p'}(\mo\times \R^d)$ because of \eqref{boundy}. Hence, extracting a subsequence if necessary, it converges in $\mathcal P_1(\mo\times \R^d)$ to some $\gamma \in \mathcal P_{p'}(\mo\times \R^d)$. Consider now $\Gamma \in \Pi(\gamma,\nu)$ for some $\nu \in \mpo$, as well as a sequence $(x'_{i,N})_{1\leq i \leq N}$ such that $\lim_{N \to \infty}\frac1N\sum_{i=1}^N \delta_{(x_{i,N},y_{i,N},x'_{i,N})} = \Gamma$ (for the $\mathcal W_1$ distance). See Lemma \ref{lemma:existbs} for the existence of such a sequence. For all $N > 1$,
$$
U^N(x'_1,\dots,x'_N) \geq U^N(x_1,\dots,x_N) + \frac1N\sum_{i =1}^N y_i\cdot(x'_i - x_i).
$$
Passing to the limit in the previous expression using the continuity of $U$, we deduce that $\gamma \in \partial^-U(µ)$.
\end{proof}

\begin{Lemma}\label{lemma:existbs}
Let $\gamma \in \mathcal P_1(\R^d\times \R^k)$ of first marginal $µ$. Let $(x_{i,N})_{1 \leq i \leq N}$ be a sequence such that $N^{-1}\sum_{i=1}^N \delta_{x_{i,N}}$ converges in $\mathcal P_1(\R^d)$ toward $µ$. Then, there exists a sequence $(y_{i,N})_{1 \leq i \leq N}$ such that, extracting a subsequence if necessary, $N^{-1}\sum_{i=1}^N \delta_{(x_{i,N},y_{i,N})}$ converges in $\mathcal P_1(\R^d\times \R^k)$ toward $\gamma$.
\end{Lemma}
\begin{proof}
Let $\eps > 0$. Because of Theorem \ref{thm:glivenko}, there exists a measure $\gamma_\eps$ of the form $K^{-1}\sum_{i =1}^K\delta_{(\tilde x_i,\tilde y_i)}$ such that $\mathcal W_1(\gamma_\eps,\gamma) \leq \eps$. Write $µ_\eps = (\pi_1)_\#\gamma_\eps$. Let $N$ be given as $qK$ for $q \geq 1$ and consider an optimal pairing between $(x_{i,N})_{1 \leq i \leq N}$ and the points in $(\tilde x_i)_{ 1 \leq i \leq K}$ repeated $q$ times each. Then, for any $i \leq N$, set $y^\eps_{i,N}$ as $\tilde y_j$ if $x_{i,N}$ is paired with an $\tilde x_j$. It then follows that
$$
\mathcal W_1\left( \frac 1N \sum_{i =1}^N\delta_{(x_{i,N},y^\eps_{i,N})}, \gamma_\eps \right) \leq \mathcal W_1\left( \frac 1N \sum_{i =1}^N\delta_{x_{i,N}}, µ_\eps \right)\underset{N \to \infty}{\longrightarrow}\mathcal{W}_1(µ,µ_\eps) \leq \eps.
$$
Since the convergence above is uniform in $\eps$, we conclude the proof by a diagonal argument to obtain the required sequence, not defined on all integers but only on an increasing sequence, but this does not raise any issue.
\end{proof}

I now state this quite impractical result for completeness. 
\begin{Prop}
Assume that $\mo$ is compact (in addition to being convex closed already). Let $U: \mpo \mapsto \R$ be coupling convex. Then, for every $µ \in \mpo$ whose support is at distance at least $\eps > 0$ of $\partial \mo$, there exists $\ell \in (L^\infty(\Omega,\R^d))^*$ such that for any $\nu \in \mpo$, $\gamma \in \Pi(µ,\nu)$ and $(X,Y) \sim \gamma$
$$
U(\nu) \geq U(µ) + \ell(Y-X).
$$
\end{Prop}
\begin{proof}
Consider the lift $\mathcal U$ of $U$ in $E_{\mo}:= \{ X \in L^\infty(\Omega,\R^d), \mathbb P(X  \in\mo) =1\}$, which is a convex subset of $L^\infty(\Omega,\R^d)$. That is, $\mathcal U (X) = U(\mathcal L(X))$ if $X \in E_\mo$, and $\mathcal U(X) = + \infty$ else. Take $µ$ as in the statement. As we saw above, any $X \sim µ$ lies in an open ball of $E_\mo$ for the $\|\cdot\|_\infty$ norm. Since $\mathcal U$ is convex on $E_\mo$ and $X$ lies in the interior of the domain of $\mathcal U$, we deduce that $\mathcal U$ is sub-differentiable at that point. We denote by $\ell \in (L^\infty)^*$ this element in the sub-differential. The result follows.
\end{proof}

\begin{Rem}
This remark is not directly related to the previous Proposition but had to be made somewhere. In this non-usual geometry, Dirac masses are not necessarily the extreme points of $\mpo$, these are rather the measures which put mass on $\partial \mo$, contrary to what happens in the \emph{flat} geometry where Dirac masses are the extreme points of $\mpo$. 
\end{Rem}

\subsubsection{Coupling monotonicity}
I focus on the case $\mo = \R^d$ to simplify notations. Let $p \in (1,\infty)$ and $µ \in \mpprd$. Consider the set 
$$
\mathcal P_{p'}^µ := \{\gamma \in \mathcal P_{pp'}(\R^d\times\R^d), (\pi_1)_\#\gamma = µ\}.
$$
Recall that thanks to the disintegration Theorem, sub and super-differentials of functions $U: \mpprd \mapsto \R$ at $µ$ can be seen as subsets of $\mathcal P^µ_{p'}$. As a convention, I no write te $A: \mpprd \rightrightarrows \cup_{µ \in \mpprd} \mathcal P_{p'}^µ$ when for all $µ \in \mpprd$, $A(µ) \subset \mathcal P_{p'}^µ$.

 The notion of monotonicity associated to coupling convexity is defined as follows.
\begin{Def}\label{def:couplingmon}
Let $A: \mpprd \rightrightarrows \cup_{µ \in \mpprd} \mathcal P_{p'}^µ$. It is said to be coupling monotone if for all $µ,\nu \in \mpprd, \gamma_µ \in A(µ), \gamma_\nu \in A(\nu)$ and $\Theta \in \Pi(\gamma_µ,\gamma_\nu)$,
$$
\int_{\R^{4d}}(z_1-z_2)\cdot(y-x)\Theta(dx,dz_1,dy,dz_2) \geq 0.
$$
\end{Def}
The next result links coupling monotonicity and coupling convexity.
\begin{Prop}
Let $U : \mpprd \mapsto \R$ be coupling convex and consider $A:\mpprd \rightrightarrows \cup_{µ \in \mpprd} \mathcal P_{p'}^µ$ defined for $µ \in \mpprd$ by
$$
A(µ) = \{\gamma \in \mathcal P_{pp'}(\R^d\times \R^d) | \gamma \in \partial^- U(µ)\}.
$$
Then $A$ is coupling monotone.
\end{Prop}
\begin{proof}
The proof is standard. Let $µ,\nu,\gamma_µ,\gamma_\nu,\Theta$ be as in Definition \ref{def:couplingmon}. Then using $\gamma_µ \in \partial^-U(µ)$ for the first inequality and $\gamma_\nu \in \partial^-U(\nu)$ for the second one (without modulus of continuity thanks to Proposition \ref{prop:loctoglob}),
$$
\int_{\R^{4d}}z_1\cdot(y-x)\Theta(dx,dz_1,dy,dz_2) \leq U(\nu) - U(µ) \leq\int_{\R^{4d}}z_2\cdot(y-x)\Theta(dx,dz_1,dy,dz_2),
$$
which yields the result.
\end{proof}

Just as it was the case for coupling convexity, coupling monotonicity can also be understood through the lifting procedure. Consider an operator $A:\mpprd \rightrightarrows \cup_{µ \in \mpprd} \mathcal P_{p'}^µ$. Define its lift $\tilde{\mathcal A}: \mathbb L^p \rightrightarrows \mathbb L^{p'}$ through
$$
\tilde{\mathcal A} = \{Y \in \mathbb L ^{p'} | \mathcal L(X,Y) \in A(\mathcal L(X))\}.
$$
The following holds trivially, thanks to Theorem \ref{thm:existrv}.
\begin{Prop}
The operator $A$ is coupling monotone if and only if $\tilde{\mathcal A}$ is monotone (in the usual sense in Banach spaces).
\end{Prop}

\subsubsection{Examples and computations}
Below are some examples of super-differentials of functions.
\begin{Prop}
Let $1< p < \infty$, $\phi: \R^d \mapsto \R$ be a concave function such that, for some $C > 0$, for any $x \in \R^d$, $|\phi(x)| \leq C (1 + |x|^p)$. Then $U:µ \mapsto \langle \phi,µ\rangle$ is (horizontally) super-differentiable everywhere and, for any $µ \in \mpprd$, $\gamma \in \mathcal P(\R^d\times \R^d)$ supported in the graph of $\partial^+\phi$ such that $(\pi_1)_\#\gamma = µ$, we have $\gamma \in \partial^+U(µ)$.
\end{Prop}
\begin{proof}
Let $\gamma$ be a coupling as in the statement, take $\nu \in \mprd$ and consider $\Gamma \in \Pi(\gamma,\nu)$. Because it is concave, $\phi$ is everywhere super-differentiable and, for all $x,y \in \R^d$ and $z \in \partial^+\phi(x)$
$$
\phi(y) \leq \phi(x) + z\cdot(y-x).
$$
If we are able to integrate the previous against $\Gamma$, and that $\gamma \in \mathcal P_{pp'}(\R^d \times \R^d)$ then the result holds. Setting in the previous, for $\delta > 0$, $y = x - \delta z |z|^{p'-2}$, we obtain
$$
\delta|z|^{p'}\leq \phi(x) - \phi(x -  z |z|^{p'-2}) \leq \phi(x) + C ( 1+ |x -\delta z |z|^{p'-2}|^p) \leq C(1 + |x|^p + \delta^p|z|^{p'}).
$$
Taking $\delta$ sufficiently small yields the result.
\end{proof}

\begin{Prop}\label{prop:superdiffWp}
Let $\mo$ be the closure of a convex domain of $\R^d$ and $p \in (1,\infty)$. Let $\nu\in \mppo$ and consider $U: \mppo \mapsto \R$ defined by $U(µ):= \mathcal W^p_p(µ,\nu)$. Then $U$ is super-differentiable everywhere and, for any $µ \in \mppo$, $\gamma$ optimal coupling for $\mathcal{W}_p(µ,\nu)$, $(\pi_1,p(|\pi_1-\pi_2|^{p-2}(\pi_1-\pi_2)))_\#\gamma \in \partial^+U(µ)$.
\end{Prop}
\begin{proof}
The proof is given through the probabilistic representation. Let $\gamma$ be an optimal coupling for $\mathcal{W}_p(µ,\nu)$, $µ' \in \mppo$ and consider any triple of random variables $(X,Y,X')$ such that $(X,Y)\sim \gamma$ and $X'\sim µ'$. Then, by optimality of $\gamma$
$$
U(µ') - U(µ) \leq \mathbb E[|X'-Y|^p- |X-Y|^p].
$$
Consider first the case $1< p \leq 2$. By convexity of $x \mapsto |x|^p$
$$
\ba
|X'-Y|^p- |X-Y|^p&\leq -p\langle |X'-Y|^{p-2}(X'-Y),X-X'\rangle\\
& = p\langle|X-Y|^{p-2}(X-Y) -|X'-Y|^{p-2}(X'-Y),X-X'\rangle\\
&\quad \quad + p\langle |X-Y|^{p-2}(X-Y),X'-X\rangle
\ea
$$
Recalling Lemma 4.4 in \cite{dibenedetto2012degenerate}, there exists $\lambda_p > 0$ such that for any $a,b \in \R^d$, $\langle|a|^{p-2}a -|b|^{p-2}b,a-b\rangle \leq \lambda_p|a-b|^p$. Hence,
$$
U(µ') - U(µ) \leq \mathbb E[p |X-Y|^{p-2}(X-Y)\cdot(X'-X)] + p\lambda_p C_p^p(\mathcal L(X,X')).
$$
Hence the result holds in this case. In the case, $p > 2$, we use the $\mathcal C^2$ smoothness of $x \mapsto |x|^p$ to obtain
$$
\ba
|X'-Y|^p- |X-Y|^p&\leq  p\langle |X-Y|^{p-2}(X-Y),X'-X\rangle\\
& \quad \quad + p(p-1)\max(|X-Y|^{p-2},|X'-Y|^{p-2})|X'-X|^2.
\ea
$$
Observe now that by H\"older's inequality
$$
\mathbb E[|X-Y|^{p-2}|X'-X|^2] \leq \left( \mathbb E[|X-Y|^p]\right)^\frac{p-2}{p}\left( \mathbb E [|X'-X|^p]\right)^\frac{2}p
$$
The second term of the product on the right hand side being equal to $C_p(\mathcal L(X',X))^2$, the result once again follows.
\end{proof}

\begin{Prop}\label{prop:superI}
Let $1 < p < \infty$, $\mo$ be the closure of a convex domain of $\R^d$, $\nu \in \mathcal P_{p'}(\mo)$ and define $U(µ) := \mathcal I_p(µ,\nu)$. Then, for any $µ \in \mppo$, $\gamma$ optimal coupling in $\mathcal I_p(µ,\nu)$, we have $(\pi_1,-\pi_2)_\#\gamma \in \partial^+U(µ)$. Furthermore, $U$ is coupling concave.
\end{Prop}
\begin{proof}
Let $(X,Y)\sim \gamma$, and take $X' \sim µ'$ for any $µ' \in \mppo$, then
$$
\mathcal I_p(µ',\nu) - \mathcal I_p(µ,\nu) \leq - \mathbb E[(X'-X)\cdot Y].
$$
Hence the super-differentiability follows. The concavity follows from the formula, for any $X \sim µ$
$$
\mathcal I_p(µ,\nu) = \inf_{Y \sim \nu}- \mathbb E[X \cdot Y],
$$
and the use of the second part of Proposition \ref{prop:couplconv1} and an immediate analogue of Proposition \ref{prop:liftconv}.
\end{proof}

\subsection{Sub-differential calculus on the torus}\label{sec:torus}
In this section, I take some time to make explicit extensions of several results above to the case of the flat torus $\T^d$, since it will be of interest in Part III. Morally, sub-differentials have the same properties as the ones we just saw, the only issue is that the term $y-x$ in the definition of the sub/super-differential makes no longer sense on the torus. Up to the addition of an element in $\mathbb Z^d$, $y-x$ is the speed of a geodesic joining $x$ to $y$ in time $1$. However, there may be multiple geodesics in this situation and thus some care is needed. Observe that in this flat case, the exponential map is simply defined as $\exp(x,v) = x +v$ modulo $\mathbb Z^d$.

I start by giving the definitions of sub/super-differential and by recalling that, since $\T^d$ is a proper manifold, the Definition \ref{def:dM} is the appropriate notion of H-differential in this context. In this case, for $\eta \in \mathcal P_2(\T^d\times \R^d)$, I define
$$
C_{\T^d}(\eta) = \sqrt{\int_{\T^d\times \R^d}|v|^2\eta(dx,dv)}.
$$
\begin{Def}
Let $U: \mptd \mapsto \R$. A measurable function $\psi : \T^d \mapsto \mathcal P(\R^d)$ is an element of the super-differential of $U$ at $µ \in \mptd$ if 
\begin{itemize}
\item $$ \int_{\T^d}\int_{\R^d}|z|^2\psi_x(dz)µ(dx) < \infty,$$
\item there exists a modulus of continuity $\omega$ such that for any $\eta \in \mathcal P_2(\T^d\times \R^d)$ such that $(\pi_1)_\#\eta = µ$, for any $\Gamma \in \mathcal P(\T^d\times \R^d\times \R^d)$ such that $(\pi_1,\pi_2)_\#\Gamma = µ(dx)\psi_x(dz)$ and $(\pi_1,\pi_3)_\#\Gamma = \eta$, 
$$
U(\exp_\# \eta) \leq U(µ) + \int_{\T^d \times \R^d\times \R^d} z\cdot v\, \Gamma(dx,dz,dv) + C_{\T^d}(\eta) \omega(C_{\T^d}(\eta)).
$$
\end{itemize}
The set of such elements is denoted $\partial^+U(µ)$. We sometimes identify $\psi$ with $µ(dx)\psi_x(dz)$. The sub-differential of a function $V$ is defined as $\partial^-V(µ) = \{x \mapsto (-Id)_\#\psi_x | \psi \in \partial^+(-V)(µ)\}$.
\end{Def}
\begin{Rem}
Because we are working on $\mathcal P(\T^d)$, the choice of $C_{\T^d}$ and of the integrability condition might seem surprising, as we could be looking for more general $(L^p,L^{p'})$ duality. This is in fact the case and the current choice has the advantage to be convenient.
\end{Rem}

Since notions of convexity on $\T^d$ are more involved, I do not enter in the link between sub-differentiability and convexity in this context. Instead, I start by showing that $\mathcal W_2^2$ is still super-differentiable everywhere in this setting.
\begin{Prop}\label{prop:superdifftorus}
For any $\nu, µ_0 \in \mptd$, the map $µ \mapsto \mathcal W_2^2(µ,\nu)$ is super-differentiable at $µ_0$ and if $\gamma$ is an optimal coupling in $\mathcal W_2^2(µ_0,\nu)$, then for any measurable map $\zeta : \T^d \times \T^d \mapsto \R^d$ which associates to $(x,y) \in (\T^d)^2$ the velocity of a constant speed geodesic joining $y$ to $x$ in time $1$, then
$$
(\pi_1,2\zeta)_\#\gamma \in \partial^+\mathcal W_2^2(\cdot,\nu)(µ_0).
$$
\end{Prop}
\begin{proof}
Observe that for any $x,y\in \T^d$, $v \in \R^d$
\be\label{ineqtorus}
\ba
d(\exp(x,v),y)^2 &= \inf_{z \in \mathbb Z^d}|x + v -y + z|^2\\
&\leq |v + x-y - z_{x,y}|^2\\
&= |x-y+z_{x,y}|^2 +|v|^2 + 2 v\cdot(x-y-z_{x,y})\\
&= d(x,y)^2 + |v|^2 + 2 v\cdot(x-y-z_{x,y}),
\ea
\ee
where we use the notation $z_{x,y}$ to denote a point of optimum in $\inf_{z \in \mathbb Z^d}|y-x + z|^2$. 

Let $\eta \in \mathcal P_2(\T^d \times \R^d)$ such that $(\pi_1)_\# \eta = µ_0$. Let $\Gamma \in \mathcal P(\T^d\times \R^d\times \T^d)$ such that $(\pi_1,\pi_2)_\#\Gamma = \eta$ and $(\pi_1,\pi_3)_\#\Gamma = \gamma$, for $\gamma$ an optimal coupling in $\mathcal W_2^2(µ_0,\nu)$. Integrating \eqref{ineqtorus} against $\Gamma$ yields the result since any $\zeta$ as in the statement has to satisfy for any $x,y \in \T^d$, $\zeta(x,y) = x-y - z_{x,y}$ for some $z_{x,y}$. (Note that $x,y \mapsto z_{x,y} := \zeta(x,y) -x+y$ is measurable because $\zeta$ is.)
\end{proof}
\begin{Rem}
Let me note that the set of $\zeta$ as in the statement of the theorem is not empty and that, for instance, $\zeta(x,y)$ can be defined as the unique element in $[-\frac12,\frac12)^d$ such that, modulo the addition of an element in $\mathbb Z^d$, $\zeta(x,y) = x-y$.
\end{Rem}

I now turn to the analogue of Proposition \ref{prop:alfonsiO}, which has the virtue of insisting upon the fact that the proof of Theorem \ref{thm:structurelocal}, originally due to Alfonsi and Jourdain in \citep{alfonsi}, is much more general than $\mpprd$.
\begin{Prop}
Let $U: \mptd \mapsto \R$ and $µ \in \mptd$ such that $U$ is both sub and super-differentiable at $µ$. Then $U$ is $H$-differentiable at $µ$ and $\partial^+U(µ) = \partial^-U(µ) = \{x \mapsto \delta_{D_µU(µ)(x)}\}$.
\end{Prop}
\begin{proof}
I argue as in Proposition \ref{prop:alfonsiO} and leave the extension of the first step of the proof to this case as an exercise (i.e. to show that $\partial^+U(µ) = \partial^-U(µ) = \{\psi\}$). It thus remains to show that this map $\psi$ is deterministic. Once again, we can construct exactly the same perturbations $\xi^\pm$ as we did in the proof of Theorem \ref{thm:structurelocal}, around a certain random variable $X \sim µ$. The crucial point on $\T^d$ is now to understand the sum $X + \eps \xi^\pm$. For this we use the exponential map and consider $X^\pm = \exp(X,\eps \xi^\pm)$. To conclude the proof as we did in Theorem \ref{thm:structurelocal}, it then suffices to show that $\mathcal L(X^+) = \mathcal L(X^-)$.  But since $\T^d$ is a smooth manifold whose injectivity radius is bounded from below uniformly, the previous holds for $\eps$ small enough, since $|\xi^\pm| \leq 1$. In other words, for $\eps$ small enough, the variations $\eps \xi^\pm$ do not cross the cut-locus of $\T^d$. The rest of the proof then follows.
\end{proof}
Finally, I state without proof (which is a good exercise for the interested reader), the following.
\begin{Prop}
Let $U : \mptd \mapsto \R$.
\begin{enumerate}[(i)]
\item If $V: \mptd \mapsto \R$ is (horizontally) differentiable at $µ$, then 
$$
\psi \in \partial^+U(µ) \Leftrightarrow x \mapsto (Id + DV(µ)(x))_\#\psi_x(dz) \in \partial^+(U+V)(µ).
$$
\item Let $V: \mptd \mapsto \R$. Then, for all $\psi\in \partial^+U(µ), \phi\in \partial^+V(µ)$, $\Theta \in \mathcal P(\T^d\times \R^d\times \R^d)$ such that $(\pi_1,\pi_2)_\#\Theta = µ(dx)\psi_x(dz)$ and $(\pi_1,\pi_3)_\#\Theta = µ(dx)\phi_x(dz)$, it holds that $(\pi_1,\pi_2 + \pi_3)_\#\Theta \in \partial^+(V+U)(µ)$.
\item If $\psi \in \partial^+U(µ)$, then for all $\lambda \geq 0$, $(x \mapsto (\lambda Id)_\#\psi_x(dz)) \in \partial^+(\lambda U)(µ)$.
\item Let $\psi,\phi \in \partial^+U(µ)$, then, for any $\Theta \in \mathcal P(\T^d \times\R^d\times \R^d)$ such that $(\pi_1,\pi_2)_\#\Theta = µ(dx)\psi_x(dz)$ and $(\pi_1,\pi_3)_\#\Theta = µ(dx)\phi_x(dz)$, it holds that for all $t \in [0,1]$, $(\pi_1,(1-t)\pi_2 + t\pi_3)_\#\Theta \in \partial^+U(µ)$.
\end{enumerate}
\end{Prop}

\subsection{Geodesic convexity}
I conclude this section on convexity by mentioning the notion of geodesic convexity, which plays an important role, namely in the proof of several functional inequalities or in the study of \textit{gradient flows}, but not in this book. I consider here $\mo$ to be the closure of a smooth convex domain of $\R^d$ and start with the following.

\begin{Def}
\begin{itemize}
\item A set $A \subset \mppo$ is geodesic convex if it contains all geodesics (in $\mppo$) between any pairs of its elements. 
\item A function $U: \mppo \mapsto \R\cup\{+ \infty\}$ is geodesic convex if, for any $µ,\nu \in \mppo$, $\gamma$ optimal in $\mathcal{W}_p(µ,\nu)$ and $t \in [0,1]$
$$
U(((1-t)\pi_1+ t \pi_2)_\#\gamma) \leq (1-t)U(µ) + t U(\nu).
$$
\end{itemize}
\end{Def}
The notion of geodesic convexity is more intuitive, for instance, we recover here that any singleton is geodesic convex. Clearly, coupling convexity is a stronger notion than geodesic convexity. Interestingly, for continuous functions which are finite everywhere, the two notions are equivalent, provided that $d$ the dimension of the base space, is at least $2$.
\begin{Prop}\label{prop:geoimplconv}
Let $d \geq 2$, $U: \mppo \mapsto \R$ be a continuous geodesic convex function for $1\leq p < \infty$. Then $U$ is coupling convex.
\end{Prop}
\begin{proof}
I argue by approximation, proving first that it is the case for coupling between finitely supported measures and then argue by density.
Let $N> 1$ be an integer and consider $C_N$ the set of discrete measures $µ$ of the form
$$
µ_x := \sum_{i =1}^N\frac1N\delta_{x_i},
$$
for $x=(x_i)_{1 \leq i \leq N} \in (\mo)^N$, which are not supposed to be distinct. The two main observations on which this proof relies are: i) for paths of the form $(µ_{x_t})_{t \in [0,1]}$ where there are no collisions in $(x_t)_{t \in [0,1]}$, there exists a finite open cover $(I_k)$ of $[0,1]$ such that $(µ_{x_t})_{t \in [0,1]}$ restricted to each of the $I_k$ is a geodesic; ii) around any segment $x \in (\mo)^N$, there are collision-less segments in $(\mo)^N$. I prove those two points first, and then explain how to use them to obtain the result.

\textit{Step 1: Convexity along collisionless segment}\\
A segment in $(\mo)^N$ is a path of the form $((1-t)x + t y)_{t \in [0,1]}$ for $x,y \in (\mo)^N$. It is called collisionless if the cardinality of the support of $µ_{x_t}$ is a constant in time. If the constant is equal to $N$ for all time, then the terminology is quite transparent, but more generally, it allows for two (or more) particles to stick together during the whole time interval $[0,1]$. It is now a simple observation, relying on the continuity of the path in each of the coordinates, that around such a configuration at any time $t_0 \in (0,1)$, there exists some $\eps > 0$ such that $(µ_{x_t})_{t \in [t_0-\eps,t_0+\eps]}$ is a geodesic, with an $\eps$ which depends only on the infimum between the positive distances between the particles. This follows from a computation that I leave to the interested reader.
From the previous fact, we indeed obtain the open cover of $[0,1]$, on which the map $t \mapsto U(µ_{x_t})$ is convex in any subintervals. Thus the map $t \mapsto U(µ_{x_t})$ is in fact convex on $[0,1]$.

\textit{Step 2: Density of collisionless segments}\\
Let $x,y \in (\mo)^N$ and consider the segment $((1-t)x + ty)_{t \in [0,1]}$. Clearly, such a segment has collisions in general, for instance consider the case $N = 2, \mo = \R^2$ and $x = ((0,0),(1,0)), y = ((1,1),(0,1))$, which has a collision at time $t=\frac12$. However, such a system can be perturbed slightly so that no collision happens, for instance by replacing $y$ by $((1,1-\eps),(0,1))$. I insist upon the fact that a collision is not the crossing of the trajectories of the particles, but the fact that, when the trajectories cross, the particles are at this crossing point at the same time. I do not reproduce the quite long computations of such a fact in general. I hope the reader is convinced by my small example and refer her or him to Proposition 6.4 in \citep{cavagnari2023lagrangian} for a complete proof.

\textit{Step 3: Conclusion}\\
Let $µ,\nu \in \mppo$ and $\gamma \in \Pi(µ,\nu)$. I want to prove that $U(((1-t)\pi_1 + t\pi_2)_\#\gamma) \leq (1-t)U(µ) + t U(\nu)$. Consider a sequence of empirical measures of the form $N^{-1}\sum_{i = 1}^N \delta_{(x^N_i,y^N_i)}$ converging toward $\gamma$ for $\mathcal W_p$. Consider the associated segment $((1-t)x^N+ty^N)_{t \in [0,1]}$. From \textit{Step 2}, there exists $(\tilde x^N,\tilde y^N)$, at $\|\cdot\|_\infty$ distance of $(x^N,y^N)$ arbitrary close to $0$, such that $((1-t)\tilde x^N+t\tilde y^N)_{t \in [0,1]}$ is collisionless. Applying \textit{Step 1}, for any $t \in [0,1]$, $U(µ_{(1-t)\tilde x^N + t \tilde y^N}) \leq (1-t) U(µ_{\tilde x^N}) + t U(µ_{\tilde y^N})$. Passing to the limit in the approximation of \textit{Step 2}, it follows that $U(µ_{(1-t) x^N + t y^N}) \leq (1-t) U(µ_{ x^N}) + t U(µ_{ y^N})$ because of the continuity of $U$. Passing to the limit $N \to \infty$ yields the required result.

\end{proof}
Note the importance of the dimension criteria $d \geq 2$ in \textit{Step 2} of the proof. In fact, the result does not hold in dimension of the base space $d =1$, a counterexample can be found in \citep{parker2024some}.

From the previous result, the interest of geodesic convex functions (compared to coupling convexity) only arises when looking at more singular functions. A typical example of such a function is the entropy $E: \mpt \mapsto \R\cup\{+ \infty\}$ defined in \eqref{def:entropy}.
\begin{Prop}\label{prop:entropy}
The entropy $E$ is geodesic convex on $\mpt$ but not coupling convex.
\end{Prop}
\begin{proof}
I only sketch the proof in $\mathcal P_2(\R)$, and refer to \cite{mccann} for more details. Let $µ,\nu \in \mathcal P_2(\R)$ and consider $(µ_t)_{t \in [0,1]}$ a geodesic joining the two. If $µ$ has no density, then, there is nothing to prove. If it is not the case, we know that there exists a non-decreasing map $T: \R \mapsto \R$ such that $T_\#µ =\nu$. I argue as if $T$ were differentiable. Furthermore, $µ_t = ((1-t)Id + t T)_\#µ $. For all $t \in [0,1)$, $µ_t$ has a density, and we can compute
$$
E(µ_t) = \int_\Rµ(x) \log\left( \frac{µ(x)}{(1-t) + tT'(x)} \right)dx
$$
Using the concavity of the $\log$, we thus obtain that $E$ is geodesic convex.\\

To show that it is not coupling convex, consider for instance the case in which $µ$ is the uniform distribution on $[0,1]$ and compute the entropy of $(\frac{\pi_1+\pi_2}{2})_\#µ(dx)µ(dy)$.
\end{proof}
It turns out that many evolutionary PDEs can be written as gradient flows (with respect to some Wasserstein metrics) of geodesic convex functionals, such as the entropy $E$. In the case of $E$, I leave as an exercise to verify that this PDE is the heat equation, which can be intuited by combining the computation for the H-derivative of $E$ of Section \ref{sec:computations} and formula $\Delta m = \text{div}(\frac{\nabla m}{m}m)$. The geodesic convexity of the potential of the gradient flow translates into various useful quantitative stability properties of the associated PDE.

\subsection*{Bibliographical comments}
The notion of convexity, or the study of convex functions on normed vector space is the subject of many excellent textbooks, see for instance the book of Ekeland and Temam \citep{ekeland1999convex}. The notion of flat monotonicity shall play a fundamental role in the study of MFG master equations below.

The study of coupling convexity for functions of probability measures has been much less studied. There have been recent developments, namely using coupling monotonicity, to provide complete studies of evolution equations in Wasserstein spaces, see for instance the work of Cavagnari, Savar\'e and Sodini \citep{savare1,savare2,cavagnari2023lagrangian}. In particular Proposition \ref{prop:geoimplconv} is borrowed from \citep{cavagnari2023lagrangian}. In these works, the terminology used is the one of totally convex functions for the convexity, and the operators for which I defined a notion of coupling monotonicity are called multi-valued probability vector field. I refer to the work of Pinzi and Savaré \citep{pinzi} for developments on the link between coupling convexity and general horizontal sub-differentiability. The notion of coupling convexity was studied by Lions, namely for its associated notion of monotonicity \citep{lions2014seminar}. See also the respective books of Villani and Santambrogio \citep{villani,santambrogio} or the work of Parker \citep{parker2024some} for developments on geodesic convexity. In fact, the notion of geodesic convexity (which is classical notion in the study of Riemannian manifolds) has attracted most of the attention on Wasserstein spaces, since the seminal work of McCann \citep{mccann}. The main reason for this is the importance of geodesic convex (and not coupling convex) functions which appear in the study of gradient flows in Wasserstein spaces, such as the entropy for instance. An associated notion of monotonicity has also been studied, namely for its role in MFG, see for instance the works of Gangbo, Meszaros, Mou and Zhang \citep{gangbo,gangbo2} for instance.

I choose to focus the presentation on coupling convex functions, since, when interested in regular, say continuous, functions, this notion seems more natural. Moreover, the associated notion of monotonicity has also been proved of interest in works of Meynard \citep{meynard2,meynard3}. I believe the current presentation is quite new, in particular its links with the notion of horizontal sub-differentials. The current horizontal sub-differential calculus is presented in more details in Chapter 10 of the book of Ambrosio, Gigli and Savaré \citep{ags}. It plays a fundamental role in the study of HJB equation on Wasserstein spaces presented in Part III. Proposition \ref{prop:barycentric} is a consequence of Theorem 3.17 in \citep{gangbotudorascu}. Proposition \ref{prop:entropy} is proven in \citep{mccann}.

\newpage
\section{$\mathcal C^{1,\alpha}$ regularity of a function of a probability measure}\label{sec:C1}
I present here $\mathcal C^{1,\alpha}$ spaces of regularity for functions of measures. There will be two definitions, one associated with the vertical point of view, and one with the horizontal one. As a convention, $\mathcal C ^{1,1}$ stands for the Lipschitz regularity of the derivative and $\mathcal C^{0,1}$ for the Lipschitz regularity of the function itself.

 This section is quite independent of the rest of this book, but serves as a good illustration of differences between vertical and horizontal derivatives. I will not address the question of second order derivatives.

\subsection{H\"older regularity of the flat differential}
Let $\mo$ be a measurable space and $U : \mmo \mapsto \R$ be a flat differentiable function. For any $µ \in \mmo$, recall that the V-derivative $\nabla U(µ)$ belongs to $ \mathcal B(\mo)$, the space of bounded measurable functions, which is a normed vector space equipped with the $\|\cdot\|_\infty$ norm. Hence, continuity or H\"older regularity of $\nabla U$ is understood quite naturally through the following definition.
\begin{Def}
Let $U: \mmo \mapsto \R$ be V-differentiable and consider the map
$$
\ba
\nabla U: (\mmo,|\cdot|) &\mapsto (\mathcal B(\mo),\|\cdot\|_\infty),\\
µ\quad\quad&\mapsto\quad \nabla U(µ).
\ea
$$
The function $U$ is said to be:
\begin{itemize}
\item V-$\mathcal C^{1,0}$ if $\nabla U$ is continuous.
\item V-$\mathcal C^{1,\alpha}$ for $\alpha \in (0,1]$ if $\nabla U$ is $\mathcal C^{0,\alpha}$ H\"older continuous. 
\end{itemize}
\end{Def}

Standard developments can then be made and I do not go into further details here.

\subsection{H\"older regularity of the horizontal differential}
The natural aspect of the previous section completely disappears as soon as we are trying to use this approach on the horizontal differential. Indeed, let $\mo$ be the closure of a smooth convex domain of $\R^d$ and consider a H-differentiable map $U : \mppo \mapsto \R$ for $1< p < \infty$. Then, trying to give a notion of regularity of $DU$, we have to compare, as $\nu$ and $µ$ get close to one another, $DU(µ)$ and $DU(\nu)$. However, those two elements lie in respectively $L^{p'}_µ(\R^d)$ and $L^{p'}_\nu(\R^d)$, which makes such a comparison far from obvious...

\subsubsection{A notion of parallel transport}
Such a problem is classical in geometry. Indeed, for $\phi$ a smooth function over a Riemannian manifold $M$, at any $x \in M$, the differential of $\phi$ at $x$ belongs to $T_xM$. Moreover, for $x \ne y$, comparing elements in $T_xM$ and $T_yM$ is not obvious. The standard tool is to use what is called \textit{parallel transport}. The parallel transport of vector $v \in T_xM$ along a smooth curve $\theta:[0,1]\mapsto M$ is a smooth curve valued in $TM$, the tangent bundle of $M$. This curve has to enjoy several properties, notably that it preserves angles (hence the terminology parallel transport). I will not need to use the concept of parallel transport in a general manifold in detail, and refer to classical textbooks on differential geometry for more on this concept. Observe that parallel transport is trivial in $\R^d$, it is simply a translation. In particular, the parallel transport along a curve in $\R^d$ simply depends on the starting and ending points of the curve.

The aim of the present section is to explain how we can use the classical parallel transport as a heuristic for "parallel transporting" an element $\phi \in L^{p'}_µ$ into something more fitting for a comparison with an element $\psi \in L^{p'}_\nu$. As we shall see, we will not be able to transport $\phi$ to an element of $L^{p'}_\nu$ in general.

To realize such a parallel transport, a curve between $µ$ and $\nu$ needs to be prescribed. It is the curve along which the transport shall be made. Since we are concerned with the horizontal derivative, it should be natural to look for a Lagrangian curve. Let $µ \in \mppo$, $\eta \in \mathcal P(\mathcal C([0,1],\mo))$ such that $(e_0)_\#\eta = µ$ and $\phi \in L^{p'}_µ(\R^d)$. Given a path $\theta$ in the support of $\eta$, the objective is to parallel transport $\phi(\theta(0))$ along the curve $\theta$. Since at this level, parallel transport is trivial (we are in $\R^d$), this procedure simply consists in assigning $\phi(\theta(0))$ to $\theta(1)$. Hence, the path $\theta$ plays a role only through its starting and ending point. Thus, we can simply work with the coupling $\gamma =(e_0,e_1)_\#\eta \in \mathcal P_p(\mo\times \mo)$ instead of the whole $\eta$. In cases in which the base space is not $\R^d$, much more care is needed. The previous discussion leads to the following.

\begin{Def}\label{def:para}
Let $µ,\nu \in \mppo$. A parallel transport of $\phi \in L^{p'}_µ(\mo,\R^d)$ along the coupling $\gamma \in \Pi(µ,\nu)$ is a probability measure $\Gamma \in \Pi(µ(dx)\delta_{\phi(x)}(dz),\nu)$ such that $(\pi_1,\pi_3)_\#\Gamma = \gamma$. Furthermore, we say that $\phi$ is parallel transported (by $\Gamma$, along $\gamma$) onto $\phi'$ where $(\pi_3,\pi_2)_\#\Gamma = \nu(dy)\phi'_y(dz)$.
\end{Def}

This definition is quite transparent at the level of the lift in the case $p =2, \mo = \R^d$. Indeed, if we take $X \in \mathbb H$ of law $µ$ and interpret $Z=\phi(X) \in \mathbb H$ as an element of the tangent space to $\mathbb H$ at $X$, then its parallel transport toward $Y \in \mathbb H$ of law $\nu$ is still given by $Z$. Thus starting from $(X,Z)$, we end up in $(Y,Z)$. At the level of $\mpt$, all this information is given by the coupling $\Gamma$ as in the previous definition.\\

More generally, we have the following existence and uniqueness result.
\begin{Prop}
Let $µ,\nu \in \mppo, \phi \in L^{p'}_µ(\mo,\R^d)$ and $\gamma \in \Pi(µ,\nu)$. Then, there exists a unique parallel transport of $\phi$ along $\gamma$.
\end{Prop}
\begin{proof}
By using Lemma \ref{lemma:gluing} on $µ(dx)\delta_{\phi(x)}(dz)$ and $\gamma(dx,dy)$ whose first marginal is $µ$, we obtain the existence. Uniqueness is in fact also a consequence of Lemma \ref{lemma:gluing}. Indeed, the coupling $µ(dx)\delta_{\phi(x)}(dz)$ is deterministic.
\end{proof}
This notion of parallel transport will help us to estimate the continuity of $D_µU$.

\begin{Rem}
Definition \ref{def:para} is concerned with the parallel transport of a deterministic coupling in $\mathcal P_{pp'}(\mo\times \R^d)$. We could have defined parallel transport for general couplings, but I do not present it because it is of no interest here.
\end{Rem}

\subsubsection{Regularity of the differential}
Let $µ,\nu \in \mppo$, $\gamma \in \Pi(µ,\nu)$ and $\phi \in L^{p'}_µ$. Our original aim was to parallel transport $\phi$ along $\gamma$, so that it can be seen as an element of $L^{p'}_\nu$. However it turns out that it is not the case. By doing such a procedure, we only end up with a probability measure $\xi \in \mathcal P_{pp'}(\mo\times \R^d)$ whose first marginal is $\nu$. But there is no reason that $\xi$ is a deterministic coupling of its marginals, as the next example shows.
\begin{Ex}
Let $\mo = [-\frac12,\frac12]$, $µ$ be the Lebesgue measure and $\phi = Id$. Then parallel transporting $\phi$ along $\gamma(dx,dy)=dx\,\delta_0(dy)$, we end up with $\delta_0(dy)\mathbb1_{z \in \mo}dz \in \mathcal P_2([-\frac12,\frac12]\times \R)$, which is not deterministic.
\end{Ex}

Hence, to compare the parallel transport of $\phi$ along the coupling $\gamma$ with elements in $L^{p'}_\nu$, we cannot use the norm of $L^{p'}_\nu$. Instead, I shall use the following quantity. Let $\mu \in \mppo$, $\xi_1,\xi_2 \in \mathcal P_{p'}^µ:=\{\xi \in \mathcal P_{pp'}(\mo\times \R^d)|(\pi_1)_\#\xi = µ\}$, and define
\be\label{eqretour}
d_{µ,p'}(\xi_1,\xi_2):= \left(\inf_\gamma \int_{\mo\times \R^d\times \R^d}|z-z'|^{p'}\gamma(dx,dz,dz')\right)^{\frac{1}{p'}},
\ee
where the infimum is taken over $\gamma \in \mathcal P(\mo\times \R^d\times \R^d)$ such that $(\pi_1,\pi_2)_\#\gamma = \xi_1$ and $(\pi_1,\pi_3)_\#\gamma = \xi_2$.
\begin{Prop}
The function $d_{µ,p'}$ defines a metric on $\mathcal P_{p'}^µ$.
\end{Prop}
The proof is left as an exercise to the reader. When $\xi_1$ and $\xi_2$ are deterministic, i.e. when there are maps $\phi_1,\phi_2$ such that $\xi_i = (Id,\phi_i)_\#µ$, we recover the $L^{p'}_µ$ distance:
$$
d_{µ,p'}(\xi_1,\xi_2) = \|\phi_1-\phi_2\|_{L^{p'}_\mu}.
$$
If one of the two is deterministic, say $\xi_1$ is given by a map $\phi_1$, then
$$
d_{µ,p'}(\xi_1,\xi_2) = \left(\int_{\mo\times \R^d}|\phi_1(x) - z'|^{p'}\xi_2(dx,dz')\right)^{\frac{1}{p'}}.
$$
In this case, I abuse notation and denote $d_{µ,p'}(\xi_1,\xi_2)= d_{µ,p'}(\phi_1,\xi_2)$. In particular, for $µ,\nu \in \mppo$, $\gamma \in \Pi(µ,\nu)$, $\phi \in L^{p'}_µ$, $\psi \in L^{p'}_\nu$ and $\mathcal T^\gamma(\phi)$ the parallel transport of $\phi$ along $\gamma$,
$$
\left(d_{\nu,p'}(\mathcal T^\gamma(\phi),\psi)\right)^{p'} = \int_{\mo\times \mo}|\phi(x) - \psi(y)|^{p'}\gamma(dx,dy).
$$

I can now introduce the main definition of this section.
\begin{Def}
Let $p \in (1,\infty)$ and $U: \mppo \mapsto \R$ be an H-differentiable map. 
\begin{itemize}
\item The function $U$ is H-$\mathcal C^{1,0}$ if for any $µ \in \mppo$, any sequences $(\nu_n)_{n \geq 0}$, $(\gamma_n)_{n \geq 0}$ such that $\gamma_n \in \Pi(µ,\nu_n)$ and $\lim_{n \to \infty}C_p(\gamma_n) = 0$
$$
\lim_{n \to \infty} \int_{\mo\times \mo} |D_µU(µ)(x)-D_µU(\nu_n)(y)|^{p'}\gamma_n(dx,dy) = 0.
$$
\item The function $U$ is H-$\mathcal C^{1,\alpha}$ for $\alpha \in (0,1)$ if there exists $K >0$ such that for any $µ,\nu \in \mppo$, $\gamma \in \Pi(µ,\nu)$ with $C_p(\gamma) \leq 1$ 
$$
\left( \int_{\mo\times \mo}|D_µU(µ)(x)-D_µU(\nu)(y)|^{p'}\gamma(dx,dy)\right)^\frac{1}{p'} \leq K (C_p(\gamma))^\alpha.
$$
\item The function $U$ is H-$\mathcal C^{1,1}$ if there exists $K >0$ such that for any $µ,\nu \in \mppo$, $\gamma \in \Pi(µ,\nu)$
$$
\left( \int_{\mo\times \mo}|D_µU(µ)(x)-D_µU(\nu)(y)|^{p'}\gamma(dx,dy)\right)^\frac{1}{p'} \leq K C_p(\gamma).
$$
\end{itemize}
\end{Def}
Note that, as usual, H\"older spaces are defined locally and Lipschitz regularity is defined globally.

Once again, the lift can help to understand the heavy notation.
\begin{Prop}\label{prop:c1a}
Let $p \in (1,\infty)$, $U : \mathcal P_p(\R^d) \mapsto \R$ be a H-differentiable map and consider its lift $\mathcal U : \mathbb L^p \mapsto \R$. Then, for $\alpha \in [0,1]$, $U$ is H-$\mathcal C^{1,\alpha}$ if and only if $\mathcal U$ is $ \mathcal C^{1,\alpha}$ in $\mathbb L^p$(in the usual sense).
\end{Prop}
\begin{proof}
We already know that $\mathcal U$ is differentiable everywhere in $\mathbb L^p$, thanks for instance to Proposition \ref{prop:gg}. Let $X,Y \in \mathbb L^p$, and define $\gamma= \mathcal L(X,Y)$. Compute
$$
\ba
\mathbb E[|\nabla \mathcal U(X) - \nabla \mathcal U(Y)|^{p'}] &= \mathbb E \left[ \left| D_µU(\mathcal L(X))(X) - D_µU(\mathcal L(Y))(Y)  \right|^{p'} \right]\\
&= \int_{\R^d\times \R^d}\left| D_µU(\mathcal L(X))(x) - D_µU(\mathcal L(Y))(y) \right|^{p'} \gamma(dx,dy),
\ea
$$ 
from which the result easily follows.
\end{proof}

\begin{Rem}
As it was the case for the notion of derivative, it is an interesting question to understand if we recover the same notion by restricting our attention to optimal couplings (see Section \ref{sec:geo}). I am not aware of such a result at the moment.
\end{Rem}

\subsubsection{Examples of $\mathcal C^{1,\alpha}$ functions}
The examples below highlight how to verify on some cases $\mathcal C^{1,\alpha}$ regularity.
\begin{Prop}
Let $p \geq 2$ and $f\in \mathcal C^{1,1}(\mo,\R)$, then $U: µ \mapsto \langle f,µ\rangle$ is H-$\mathcal C^{1,1}$. 
\end{Prop}
\begin{proof}
Let $µ,\nu \in \mppo$ and $\gamma \in \Pi(µ,\nu)$. First, $U(µ)$ is well defined since $f$ has at most quadratic growth. Then 
$$
\ba
\int_{\mo\times \mo}|D_µU(µ)(x)-D_µU(\nu)(y)|^{p'}\gamma(dx,dy) &= \int_{\mo\times \mo}|\nabla_xf(x)-\nabla_xf(y)|^{p'}\gamma(dx,dy)\\
&\leq \|\nabla_x f\|^{p'}_{Lip}(C_{p'}(\gamma))^{p'} \leq \|\nabla_x f\|^{p'}_{Lip}(C_{p}(\gamma))^{p'},
\ea
$$
since $p \geq p'$.
\end{proof}

A more involved function is the following.
\begin{Prop}\label{prop:supconvW}
Let $\nu \in \mpt$ and $\delta \in (0,1)$, then define $U: \mpt \mapsto \R$ by
$$
U(µ) = \sup_{µ' \in \mpt} \left\{\mathcal{W}_2^2(µ',\nu) - \frac{1}{\delta}\mathcal{W}_2^2(µ',µ)\right\}.
$$
Then $U $ is H-$ \mathcal C^{1,1}$.
\end{Prop}
\begin{Rem}
The function $U$ in the proposition is in fact an approximation of $\mathcal W_2^2$. It is an instructive exercise to verify that the convergence of $U$ toward $\mathcal W_2^2$ holds (locally uniformly in $µ$ for instance), when $\delta \to 0$.
\end{Rem}
\begin{proof}
I present the proof using the lift. Let $\mathcal U: \mathbb H \mapsto \R$ be defined by $\mathcal U(X) = U(\mathcal L(X))$. For $X \in \mathbb H$, define
$$
\Phi(X) := \inf_{Y \sim \nu} \mathbb E[|X-Y|^2].
$$
By invoking Theorem \ref{thm:proba}, $\Phi(X) = \mathcal{W}_2^2(\mathcal L(X),\nu)$. Note that the function $X \mapsto \Phi(X) - \mathbb E[|X|^2]$ is concave, as the infimum of linear functions in $X$. The main argument of the proof consists in looking at $\mathcal U$ (or at $U$) as a sup convolution of $\Phi$.\\
\textit{Step 1: semi-concavity of $\mathcal U$}\\
For any $X \in \mathbb H$,
$$
\mathcal U(X) = \sup_{X'\in \mathbb H} \left\{ \inf_{Y \sim \nu} \{\mathbb E[|X'-Y|^2]\} - \frac1\delta \mathbb E[|X'-X|^2]\right\}.
$$
In order to establish the semi-concavity, we define $\Psi(X,X') := \inf_{Y \sim \nu} \{\mathbb E[|X'-Y|^2]\} - \frac1\delta \mathbb E[|X'-X|^2] - C\mathbb E[|X|^2]$ for $C > 0$ to be chosen later on. Observe that $ \Psi(X,X') = g(X') + (1-\delta^{-1})\|X'\|^2 +2\delta^{-1}\mathbb E[X'\cdot X] - (C+\delta^{-1})\|X\|^2$ for $g$ a concave function. Hence if 
$$
\begin{pmatrix}
-(C + \delta^{-1}) & \delta^{-1}\\ \delta^{-1} & ( 1-\delta^{-1})
\end{pmatrix} \leq 0,
$$
then $\Psi$ is concave. It is the case for $C \geq \frac{1}{1-\delta}$. We now use a simple argument to show that $\mathcal U(X) - C\|X\|^2= \sup_{X'} \Psi(X,X')$ is also concave. Let $X_1,X_2 \in \mathbb H$, and consider $Y_1,Y_2 \in \mathbb H$ which are $\eps$-optimal in resp. $\sup_{X'} \Psi(X_i,X')$ for $i=1,2$. Then, for $t \in [0,1]$, it follows that
$$
\ba
\sup_{X'} \Psi((1-t)X_1 + tX_2,X')&\geq \Psi((1-t)X_1 + tX_2,(1-t)Y_1 + tY_2)\\
& \geq (1-t)\Psi(X_1,Y_1) + t\Psi(X_2,Y_2)\\
& \geq (1-t)\sup_{X'} \Psi(X_1,X') + t \sup_{X'} \Psi(X_2,X') - \eps.
\ea
$$
Hence, $X \mapsto \mathcal U(X) - C\|X\|^2$ is indeed concave.\\
\textit{Step 2: semi-convexity of $\mathcal U$}\\
A direct computation yields that $X \mapsto \mathcal U(X) + \frac1\delta \|X\|^2$ is convex.\\
\textit{Step 3: semi-concavity and semi-convexity implies $\mathcal C^{1,1}$}\\
From the semi-concavity, we deduce that the super-differential of $\mathcal U$ is non-empty everywhere, and from the semi-convexity, that it is also the case for the sub-differential. Hence, it follows that $\mathcal U$ is differentiable everywhere. Furthermore, from the semi-convexity, we obtain that, for $X,Y\in \mathbb H$
$$
\mathbb E[(\nabla\mathcal U(X) - \nabla\mathcal U(Y))\cdot (X-Y)] \geq -\frac2\delta\|X-Y\|^2.
$$
The similar inequality from the semi-concavity finally yields that $\mathcal U$ is $\mathcal C^{1,1}$. Thanks to Theorem \ref{thm:structurelocal}, we deduce that $U$ is H-differentiable everywhere, and the result now holds thanks to Proposition \ref{prop:c1a}.
\end{proof}
\subsection*{Bibliographical notes}
Once again, at the level of the flat derivative, the notion of $\mathcal C^{1,\alpha}$ regularity is quite common and I do not comment more on it here.\\

At the level of the horizontal derivative, the spaces $\mathcal C^{1,\alpha}(\mppo)$ have been the subject of very few works, even though the space $\mathcal C^{1,1}(\mpt)$ plays an important role in several results of Chapter 5 in the book of Carmona and Delarue \citep{carmona2017probabilistic}.
It turns out that for certain mean field optimal control problems (see Section \ref{sec:optcontr}), when some singular terms are present, $\mathcal C^{1,1}$ regularity is helpful to prove several properties of the associated value function. We noticed this idea with Lions and Souganidis in \citep{contrdyson} and proved with Lions that we can create $\mathcal C^{1,1}$ approximations of functions using sup-convolution (as in Proposition \ref{prop:supconvW}) in \citep{bertucci2024approximation}. This idea is an extension of the famous Lasry-Lions regularization in Hilbert spaces \citep{lasrylions}.

The regularity spaces and the notion of parallel transport presented here are the subject of \citep{bertucci2025tangent} which is concerned with more general statements and settings. The notion of parallel transport in Wasserstein spaces has been studied before, notably by Ambrosio and Gigli \citep{ambrosiogigli}, who studied a more Eulerian notion. Gigli also introduced the quantity defined in \eqref{eqretour} in \citep{giglitesi,gigli2}. More generally, the notion of tangent space to the Wasserstein space is often considered by analogy with Riemannian geometry and can lead to several interesting developments. I refer to Gigli \citep{giglitesi}, Ambrosio, Gigli and Savaré \citep{ags}, Aussedat \citep{aussedat} and my work \citep{bertucci2025tangent} for more on this topic.

\newpage

\section{Optimality conditions and perturbed optimization}\label{sec:optim}
In this section I treat two topics which are linked to optimization problems over spaces of measures. The first topic is optimality conditions, and the second one concerns the existence and uniqueness of optimizers and how it can be achieved by perturbing slightly the problem.

\subsection{Optimality conditions}
Optimality conditions shall take different forms, depending on the geometry with which we look at the problem.
\subsubsection{The vertical point of view}
Optimality conditions in terms of the vertical derivative are quite standard. I only present a simple result as an illustration. Let $\mo$ be a measurable space, $E$ be either $\mmo$ or $\mpo$, $F: E \mapsto \R$ and consider
\be\label{optim1}
\inf_{µ\in E} F(µ).
\ee
When $E= \mmo$, $µ_0$ is optimal in \eqref{optim1} and $F$ is V-differentiable at $µ_0$, we of course have $\nabla_µ F(µ_0)= 0$. When dealing with $E = \mpo$, because of the many corners and boundaries of $\mpo$, such optimality conditions only yield inequalities in general. Recall also that in this case, I choose a normalization for $\nabla_µ F$.
\begin{Prop}\label{prop:vertopti}
Let $µ_0$ be a solution to \eqref{optim1} with $E =\mpo$ and assume that $F$ is (vertically) Gateaux differentiable at $µ_0$, then $\nabla_µ F(µ_0)\geq 0$ and $\nabla_µ F(µ_0) = 0$ $µ_0$ almost everywhere.
\end{Prop}
\begin{proof}
For any $y \in \mo$, by optimality of $µ_0$, for any $t \in (0,1)$
$$
\frac{F((1-t)µ_0 + t \delta_y) - F(µ_0)}{t} \geq 0.
$$
Passing to the limit $t \to 0$, we deduce that $\int_\mo\nabla_µF(µ_0)(x)(\delta_y - µ_0)(dx) \geq 0$. Using $\langle \nabla_µF(µ_0),µ_0\rangle = 0$, it follows that $\nabla_µ F(µ_0)(y) \geq 0$ for all $y \in \mo$, which proves the first part of the claim. Using once again the normalization $\langle \nabla_µF(µ_0),µ_0\rangle = 0$, we finally obtain that $\nabla_µ F(µ_0) = 0$ $µ_0$-almost everywhere.
\end{proof}
Second order conditions could also be considered but they are not in the scope of this book. 

\subsubsection{The horizontal point of view}
Let $\mo$ be the closure of a smooth domain of $\R^d$, $p\in(1,\infty)$ and $F : \mppo \mapsto \R$. Consider 
\be\label{optim2}
\inf_{µ\in \mppo} F(µ).
\ee
In this problem again, we shall face the fact that, if the minimizer is an extreme point of $\mppo$, only inequalities will be available. As mentioned above, in this geometry associated to horizontal derivatives, the extreme points are the measures which put mass close to the boundary of $\mo$. This can be seen through the following results.

\begin{Prop}\label{prop:uyt}
Let $\mo = \R^d$ and $µ_0$ be optimal for \eqref{optim2}. If $F$ is H-differentiable at $µ_0$, then $D_µF(µ_0) = 0$, in $L^{p'}_{µ_0}$.
\end{Prop}
\begin{proof}
Let $\phi \in \mathcal C^{1}_c(\R^d,\R^d)$. Recalling the optimality of $µ_0$, for all $t > 0$
$$
\frac{F((Id + t \phi)_\#µ_0) - F(µ_0)}{t} \geq 0.
$$
Passing to the limit $t \to 0$ and recalling that $C_p((Id,Id + t \phi)_\#µ_0) = t \|\phi\|_{L^p_{µ_0}}$, it follows that
$$
\int_{\R^d}D_µF(µ_0)(x)\cdot \phi(x) \,µ_0(dx) \geq 0.
$$
The result follows by density of $\mathcal C^{1}_c(\R^d,\R^d)$ in $L^p_{µ_0}(\R^d,\R^d)$.
\end{proof}

When $\mo$ has a boundary, the previous argument still works if we restrict our attention to more precise variations than simply ones of the form $(Id+t\phi)_\#µ_0$. Namely, we need to focus, just as we did in Proposition \ref{prop:deruniqbound}, on vector fields in the set $K:= \{\phi \in \mathcal C^{1}_c(\mo,\R^d) | \forall t \geq 0, X_t^\phi(\mathcal O) \subset \mo\}$, where $X^\phi_t$ is the flow of the ODE $\dot x = \phi(x)$.

\begin{Prop}
Let $µ_0$ be optimal for \eqref{optim2}, and $F$ be H-differentiable at $µ_0$. Then for any $\phi \in K$,
\be\label{eq:1697}
\int_{\mo}D_µF(µ_0)(x)\cdot\phi(x)\,µ_0(dx) \geq 0.
\ee
\end{Prop}
\begin{proof}
The proof follows the same argument as the previous one, except for the fact that $(Id + t \phi)_\#µ_0$ is replaced by $(X^\phi_t)_\#µ_0$, for $(X^\phi_t(\cdot))_{t \geq 0}$ the flow of the ODE $\dot x = \phi(x)$.
\end{proof}
Note that if such an optimal $µ_0$ does not put mass close to $\partial \mo$, say at a distance less than $\eps > 0$, then the conclusion of Proposition \ref{prop:uyt} holds since then $K$ is dense in $L^p_{µ_0}(\mo,\R^d)$. On the other hand, considering situations such as $\mo = \R_+$, $F(µ) = \langle Id, µ\rangle$, it is clear that one cannot hope for more than \eqref{eq:1697}.

\subsection{Perturbed optimization}\label{sec:perturbed}
As we shall see when studying HJB equations in non finite dimensional spaces, results in perturbed optimization are quite handy. If we were working on a Banach space, the typical form of results we need is the following, often called Stegall's Lemma.
\begin{Lemma}\label{lemma:stegall}
Let $E$ be a reflexive Banach space and $f : E \mapsto \R\cup\{+\infty\}$ be a lsc function bounded from below which satisfies
$$
\lim_{|x|_E \to \infty} \frac{f(x)}{|x|} = + \infty.
$$
Then, for all $\eps > 0$, there exists $y \in E'$, $\|y\|_{E'}\leq \eps$ such that 
\begin{itemize}
\item The lsc function $x \mapsto f(x) + \langle y,x\rangle$ has a minimum at some $x_* \in E$.
\item For any $(x_n)_{n\geq 0}$ such that $\lim_{n \to \infty}f(x_n) + \langle y,x_n\rangle = f(x_*) + \langle y,x_*\rangle$, $(x_n)_{n \geq 0}$ converges toward $x_*$.
\end{itemize}
\end{Lemma}
A minimum point such as the one in the previous Lemma is called strongly exposed. The previous result has two main aspects: first it creates points of minimum, second it creates them such that they are exposed. In compact cases, creating points of minima of lsc functions is not a problem, however, making them strict (which is equivalent to exposed in compact cases) is useful to obtain various stability properties of PDEs, as we shall see in Part II.

In infinite dimension, of course even the existence of a point of minimum is not guaranteed. Consider for instance the case of a separable Hilbert space $H$ with basis $(e_n)_{n \geq 1}$, and the set $S\subset H$ defined by $S= \{(1 + \frac1n)e_n| n \geq 1\}$ which is closed. The function distance to $S$ is continuous on the closed unit ball, but never attains its minimum which is $0$.

I also insist upon the fact that the perturbation used in the previous lemma is smooth (linear here) which is important in several applications, notably when trying to derive the optimality conditions at the point of minimum of the perturbed function.\\

 In the remainder of this section, I will establish two variants of Stegall's Lemma. The first one is aligned with the vertical derivative, and is concerned with a compact case where we want to extend the existence of points of minima to the existence of a strict minimum for a perturbation. The second one will be in the spirit of horizontal derivatives and deals with non-compact cases, thus when existence of minima is not clear without perturbations.

\subsubsection{Vertical variants in compact cases}
Let $E$ be either $\mpo$ or $\mathcal M_1(\mo)$, for $\mo$ the closure of a smooth bounded domain of $\R^d$, and $F: E \mapsto \R$ be lsc for the weak topology. Since $E$ is compact for the weak topology, $F$ admits a minimum simply by lower semi-continuity. The following variant of Stegall's Lemma  allows one to  obtain strict / strongly exposed points of minimum. Recall that $H^k(\mo)$ is the usual Sobolev space.
\begin{Prop}\label{prop:stegallvert}
For any $\eps > 0$, $k > \frac d2$, there exists $\phi \in H^k(\mo)$ with $\|\phi\|_{H^k} \leq \eps$ such that $µ \mapsto F(µ) + \langle \phi,µ\rangle$ admits a unique minimizer on $E$, which is thus strongly exposed for the weak topology.
\end{Prop}
\begin{proof}
We use two important properties of $H^k(\mo)$. First, from standard Sobolev embeddings, $H^k(\mo)\subset \mathcal C(\mo)$, hence, $E \subset (H^k(\mo))'$. Second, it is a Hilbert space. Define the correspondence $A: H^k(\mo) \rightrightarrows E$ by
$$
A(\phi) := \text{argmin}_{µ \in E} F(µ) + \langle µ,\phi\rangle.
$$
For every $\phi \in H^k(\mo)$, $A(\phi)\ne \emptyset$ since $µ \mapsto F(µ) + \langle \phi,µ\rangle$ is lower semi continuous and $E$ is compact, both for the weak topology. Furthermore, $-A$ is cyclical monotone, i.e. for every $n \geq 2$, $\phi_1,\dots,\phi_n \in H^k(\mo)$ and  $µ_1 \in A(\phi_1),\dots,µ_n \in A(\phi_n)$, 
$$
\sum_{i=1}^n\langle\phi_i - \phi_{i-1},µ_i\rangle \leq 0,
$$
where $\phi_0 := \phi_n$. Indeed, noting $µ_{n+1} = µ_1$, we have
$$
\sum_{i=1}^n\langle\phi_i - \phi_{i-1},µ_i\rangle = \sum_{i=1}^n\langle \phi_i,µ_i - µ_{i+1}\rangle.
$$
By optimality of $µ_i$, we know that $\langle \phi_i,µ_i - µ_{i+1}\rangle \leq F(µ_{i+1})- F(µ_i)$. Hence, the cyclical monotonicity of $-A$ follows because $\sum_{i=1}^nF(µ_{i+1})- F(µ_i) = 0$. 

As a consequence of cyclical monotonicity, $-A: H^k(\mo) \rightrightarrows (H^k(\mo))'$ is included in the sub-differential of a convex function over $H^k(\mo)$, see Theorem 2.5 in \cite{brezis}. This convex function has full domain since $A$ has full domain, and it is lsc as it can be defined as the supremum of affine functions. Since $H^k(\mo)$ is separable, Mazur's Theorem (Theorem 1.20 in \citep{phelps}) states that it is Gateaux differentiable on a dense set, hence its sub-differential is reduced to a singleton on each of those points of Gateaux-differentiability (Proposition 1.8 in \citep{phelps}). Hence, on those points, $A$ is single valued, so there is a unique point of minimum of the associated function. Since we are working on a compact set, it follows that this unique point of minimum is also strongly exposed. Indeed, if a strict minimum $µ^*$ was not exposed, there would exist a minimizing sequence living outside a small ball centered around $µ^*$ (for any $\mathcal W_1$ for instance). By compactness and lsc, this would imply the existence of a point of minimum other than $µ^*$ which is not possible.
\end{proof}
\begin{Rem}
This proof makes apparent a link between properties of differentiability of convex functions on a dense set and the existence of strongly exposed points of minimum of perturbed functions, one at the level of the dual of the other. This is at the core of the theory of perturbed optimization. Spaces in which convex functions are differentiable at many points are called Asplund spaces, and spaces in which one can create exposed points of almost minimum are said to satisfy the Radon-Nikodym property, and there is an important duality between the two, see for instance the book of Phelps \citep{phelps} or the one of Diestel and Uhl \citep{diestel}.
\end{Rem}
I shall often use the previous result with a $k$ sufficiently large so that the perturbation $\phi$ is in fact a $\mathcal C^2$ function.

\subsubsection{A horizontal variant}
The previous variant only dealt with the uniqueness of points of minimum as their existence was clear. In non-compact cases, it is not so clear how to make the previous proof work, as one cannot guarantee immediately that $A$ is not empty valued on sufficiently many points. I do not address the difficult question of the direct extension of the previous result to $\mpprd$, as it would require many additional developments. I rather present another variant, which is more in line with the horizontal regularity of functions, lifting approach or coupling convexity. I restrict myself to the case $\mo = \R^d$.

 \begin{Prop}\label{variational:principle}
 Let $p \in (1,\infty)$ and $U: \mpprd \mapsto \R$ be lsc, bounded from below which satisfies
$$
 \lim_{M_p(µ) \to \infty} \frac{U(µ)}{(M_p(µ))^\frac1p} = + \infty.
 $$
 Then, for any $\eps > 0$, there exists $\nu \in \mathcal P_{p'}(\R^d)$ with $M_{p'}(\nu) \leq \eps$ such that
 $$
 \ba
 \mpprd &\mapsto \quad \quad \R\\
 µ \quad&\mapsto U(µ) + \mathcal I_p(µ,\nu)
 \ea
 $$
 has a strongly exposed minimum at some $µ^*\in \mpprd$.
 \end{Prop}
 \begin{proof}
Consider the lift $\mathcal U: \mathbb L^p\mapsto \R$ of $U$. Because of Proposition \ref{prop:samecont}, it is also lsc, bounded from below and coercive on $\mathbb L^p$. Thanks to Lemma \ref{lemma:stegall}, for any $\eps > 0$, there exists $Y \in \mathbb L^{p'}$, such that $\|Y\|^{p'}_{p'} \leq \eps$ and $X \mapsto \mathcal U (X) + \mathbb E[Y\cdot X]$ has a strongly exposed minimum at $X^* \in \mathbb L^p$. We now want to pull back the perturbation $X \mapsto \mathbb E[Y\cdot X]$ at the level of probability measures, in order to obtain our perturbation. Of course this map is not law invariant, but it turns out that $µ \mapsto \inf_{X \sim µ}\mathbb E[Y\cdot X]$ is. Furthermore, using once again Theorem \ref{thm:proba}, we obtain that $\inf_{X \sim µ}\mathbb E[Y\cdot X]= \mathcal I_p(µ,\mathcal L(-Y))$. 

It now remains to prove that the map $ µ \mapsto U(µ) +\mathcal I_p(µ,\mathcal L(-Y))$ has a strongly exposed minimum at $µ^*:= \mathcal L(X^*)$. Remark first that $\mathcal L(X^*,-Y)$ is optimal for $\mathcal I_p(µ^*,\mathcal L(-Y))$ since $\mathcal U$ is law invariant. It then follows that for any $X \in \mathbb L^p$ such that $\mathcal L(X) \ne \mathcal L(X^*)$,
$$
\ba
U(\mathcal L(X)) + \mathcal I_p(\mathcal L(X),\mathcal L(-Y)) &= \inf_{X' \sim \mathcal L(X)} \mathcal U(X') + \mathbb E [Y\cdot X']\\
&> U(µ^*) + \mathcal I_p(µ^*,\mathcal L(-Y)).
\ea
$$
The strict inequality comes from the fact that $X^*$ is a point of strongly exposed minimum. Thus, $µ^*$ is a point of strict minimum. Finally, let $(µ_n)_{n \geq 0}$ be such that $U(µ_n) +\mathcal I_p(µ_n,\mathcal L(-Y)) \to_{n \to \infty} U(µ^*) +\mathcal I_p(µ^*,\mathcal L(-Y))$. Then, using Theorem \ref{thm:proba}, we can find for all $n \geq 0$, $(X_n,-Y_n)$ whose law is optimal for $\mathcal I_p(µ_n,\mathcal L(-Y))$, such that $\|Y_n -Y\|_\infty \leq (n+1)^{-1}$. It then follows that $\mathcal U(X_n) + \mathbb E[Y\cdot X_n] \to_{n \to \infty}\mathcal U(X^*) + \mathbb E[Y\cdot X^*]$. But since $X^*$ is strongly exposed, the previous implies $\|X^*-X_n\|_p \to_{n \to \infty} 0$, and thus that $\mathcal W_p(µ_n,µ^*) \to_{n \to \infty} 0$.
 \end{proof}
 
 The downside of this result is that the perturbation obtained is no longer smooth. However, thanks to Proposition \ref{prop:superI}, it is super-differentiable, and it will be sufficient in the applications at hand below. In fact we shall use either one of two extensions presented above. The first one is a corollary to Proposition \ref{variational:principle} whose proof is left as an exercise, and the second one requires a bit of care. 
\begin{Cor}\label{cor1s}
Let $p \in (1,\infty)$, $U: \mpprd \times \mpprd \mapsto \R$ be lsc, bounded from below and satisfying 
$$
\lim_{M_p(µ) + M_p(\nu)\to \infty}\frac{U(µ,\nu)}{(M_p(µ))^\frac1p + (M_p(\nu))^\frac1p} = +\infty.
$$
 For every $\eps > 0$, there exist $\xi_1,\xi_2 \in \mathcal P_{p'}(\R^d)$, such that $M_{p'}(\xi_1),M_{p'}(\xi_2)  \leq \eps$ and $(µ,\nu) \mapsto U(µ,\nu) + \mathcal I_p(µ,\xi_1) + \mathcal I_p(\nu,\xi_2)$ has a strongly exposed point of minimum.
\end{Cor}
\begin{proof}
The proof is the exact analogue of Proposition \ref{variational:principle} applied on the lift $\mathcal U(X,Y) = U(\mathcal L(X),\mathcal L(Y))$ and is left as an exercise.
\end{proof}

In some applications, notably to obtain stability properties of viscosity solutions in Part III, we would like to use a variant of Proposition \ref{variational:principle} on only a certain marginal of the measure argument. In order to state this properly, we need the following notation.

Let $p,q \in (1,\infty)$, $\mo$ be a smooth manifold, take $\nu \in \mppo$, set $\mathcal P^\nu_q(\mo\times \R^d) := \{\gamma \in \mathcal P(\mo \times \R^d) | (\pi_1)_\#\gamma = \nu, (\pi_2)_\#\gamma \in \mathcal P_q(\R^d)\}$. Let $d_{\nu,p}: (\mathcal P^\nu_p(\mo\times \R^d))^2 \mapsto \R_+$ be defined as in \eqref{eqretour}, but for the change $p' \leftarrow p$, which does not alter the fact that $d_{\nu,p}$ is a metric on $\mathcal P^\nu_p(\mo\times \R^d)$. We also define $\mathcal I_p^\nu: \mathcal P^\nu_p(\mo\times \R^d) \times \mathcal P^\nu_{p'}(\mo\times \R^d)\mapsto \R$ by
$$
I_p^\nu(\gamma,\eta) = \inf \left\{ -\int_{\mo\times \R^d\times \R^d}x\cdot y \Gamma(dz,dx,dy)\right\},
$$
where the infimum is taken over $\Gamma \in \mathcal P(\mo \times \R^d\times \R^d)$ such that $(\pi_1,\pi_2)_\#\Gamma = \gamma$ and $(\pi_1,\pi_3)_\#\Gamma = \eta$. We then have the following.

\begin{Prop}\label{cor2s}
Let $p \in (1,\infty)$, $\nu\in \mppo$ for $\mo$ a smooth manifold, and $U: \mathcal P^\nu_p(\mo\times \R^d)\mapsto \R$ be lsc (for $d_{\nu,p}$), bounded from below and such that 
$$
\lim_{M_p((\pi_2)_\#\gamma) \to \infty}\frac{U(\gamma)}{(M_p((\pi_2)_\#\gamma))^\frac1p} = + \infty.
$$
 Then, for every $\eps > 0$, there exists $\eta \in \mathcal P^\nu_{p'}(\mo\times \R^d)$ such that $M_{p'}((\pi_2)_\#\eta) \leq \eps$ and $\gamma \mapsto U(\gamma) + \mathcal I^\nu_p(\gamma,\eta)$ has a strongly exposed point of minimum (for $d_{\nu,p}$) on $\mathcal P^\nu_p(\mo\times \R^d)$.
\end{Prop}
\begin{proof}
The proof follows the one of Proposition \ref{variational:principle} once a more precise probabilistic setting has been fixed. In this proof, we consider a probability space $\Omega$ which is given as the product of two standard probability space $\Omega_1$ and $\Omega_2$. We consider a random variable $Z: \Omega \mapsto \mo, \omega \mapsto Z(\omega)$ of law $\nu$ which depends only on $\pi_1(\omega)$. This can be done thanks to Theorem \ref{thm:existrv} since $\Omega_1$ is standard. Note that, because $\Omega_2$ is standard as well, for any $\gamma \in \mathcal P_p^\nu(\mo\times \R^d)$, there exists $X$ defined on $\Omega$ such that $\mathcal L(Z,X) = \gamma$. This is not a trivial consequence of Theorem \ref{thm:existrv} but still only a simple exercise. We now define a lift $\mathcal U(X) = U(\mathcal L(Z,X))$. Note that $\mathcal U$ is not law invariant, but that we can still applies Stegall's Lemma to it. Thus we find $Y \in \mathbb L^{p'}$ such that $X \mapsto U(\mathcal L(Z,X)) + \mathbb E[X\cdot Y]$ has a strongly exposed minimum at $X^*$. Write $\gamma^* = \mathcal L(Z,X^*)$. We can now argue similarly to the proof of Proposition \ref{variational:principle} to conclude.
\end{proof}

\subsection*{Bibliographical comments}
The optimality conditions presented above are classical. In the vertical case, they can be found in my paper \citep{bertucci2023monotone} for instance.\\

Topics of perturbed optimization have arisen both in optimization problems and in questions of geometry of Banach spaces. Stegall's Lemma (Lemma \ref{lemma:stegall}) was first proven by Stegall in \citep{stegall1978optimization} and a self contained presentation can be found in the book of Phelps \citep{phelps}, in particular, see Theorem 5.15 for a statement in the more general case of Banach spaces with the Radon-Nikodym property (RNP). See also Bourgain \citep{bourgain} and Diestel and Uhl \citep{diestel} for more on the topic of the RNP.

I presented first Proposition \ref{prop:stegallvert} in \citep{bertucci2023monotone} to study MFG master equations, with a proof very similar to Rockafellar's on cyclically monotone operators. We established Proposition \ref{variational:principle} with Lions in \citep{bertucci2026stegall} to study HJB equations on non-compact spaces of measures.

\newpage

\section*{Bridge to Parts II and III}
In this Part I, I presented a sort of dichotomy in the study of functions of a measure argument, which was essentially based on the type of variation we want to consider for this measure argument.\\

Considering arguments in $\mpo$ and restricting ourselves to transport or more generally to variations associated with displacement on the underlying space $\mo$, the notion of horizontal derivative naturally arises, and with it notions of convexity and monotonicity for instance. In Part III, I will study problems of optimally controlling such variations, as is done in optimal transport for instance. Hence, for all of Part III, the notion of H-derivative will play a central role. More precisely, it will be the associated notion of sub/super-differentiability which will be of importance, as the PDE of interest does not necessarily have H-differentiable solutions. I insist here upon the fact that it will be the underlying variations of interest which will "select" the H-derivative.\\

When dealing with mean field games in Part II, the general structure of the underlying variations is much more general. In this setting, the measure argument $µ$ describes the distribution of players in the state space $\mo$. Without entering into too much detail for the moment, it is quite reasonable that players anticipate more general variations than simply ones given by transport. For instance, the mass of players taking part in the game could naturally vary. This leads us to the setting of the vertical derivative, and associated notions of monotonicity and convexity. Hence, only the V-derivative will appear in Part II. Once again, solutions to the PDE of interest will not be smooth in general, and very few properties will be needed, mostly sub/super-differentials and optimality conditions. Furthermore, it is the vertical derivative associated to the $\mathcal W_1$ distance that will be of interest, as in Section \ref{sec:W1}. In particular, it will be compatible with the weak topology.

\newpage

\part{Mean Field Game master equations}
\section{Introduction to master equations}
\subsection{A (very) short background on mean field games}
Mean field game (MFG) theory was originally motivated by the mathematical study of dynamic nonatomic games. Its development is mainly based on the seminal work of Lasry and Lions \citep{lasry2006jeux,lasry2006jeux2,lasry2007mean,lions20062007}. Of course, earlier studies of such games precede this work as for example, the economic models by Krusell and Smith \citep{krusell1998income}, Scheinkman and Weiss \citep{scheinkman1986borrowing}, or the work in engineering by Huang, Caines, and Malhamé \citep{huang2003individual,huang2006large,huang2007large}.

Using a terminology from economics or physics, the main objective of MFG is to characterize macroscopic quantities by aggregating some microscopic effects. In most interesting (mathematically at least) situations, the equation characterizing the macro quantity is forward, while the one characterizing microscopic effects is backward, with respect to time. The forward evolution comes from the usual computation of the evolution of certain quantities, when the backward structure comes from the usual accounting for anticipation: decision making through dynamic programming or computation of the value of an asset as the sum of its future expected gains for instance.

I believe it is fundamental to insist upon the fact that, despite its name, MFG theory is not only concerned with games, but with more general micro-macro configurations or forward-backward structures. This important remark opens many models, which are not games, to the powerful mathematical developments made in the theory of MFG, see for instance \citep{bitcoin,edmond,clemence}. It will also allow me to consider freely classes of equations, without the need to justify the existence of a proper game from which they originate.

In any case, the aggregate of the microscopic scale usually takes the form of a non-negative measure $m$ over some set $\mo$ which describes the possible states of the microscopic quantities we wish to aggregate. Hence, MFG naturally leads to some analysis questions on spaces of measures. In some cases, it will only take the form of an evolution equation for a measure argument, and in other ones, it leads to the study of PDEs on spaces of measures.\\

Due to its strong modeling power, this theory is currently the subject of extensive research by many mathematicians and is being applied in many fields such as engineering, physics, finance, economics, and more generally in social sciences.\\

The most challenging mathematical object in MFG theory is arguably the master equation (ME), which is usually a nonlinear PDE posed on a space of measures\footnote{Note that spaces of measures can be finite dimensional when the base space is finite.}. When there is uniqueness of the underlying micro-macro equilibria, the solution of the ME allows for an almost complete description of the model. The study of this equation is the topic of this part. For a general presentation of MFG theory and the role of the master equation, I refer to the book by Carmona and Delarue \citep{carmona2017probabilistic,carmona2017probabilistic2}, and to the book of Cardaliaguet, Delarue, Lasry and Lions \citep{cdll} for early developments on this equation.

\subsection{A formal derivation of a master equation}\label{sec:derive}
In this section, I give a formal and quite abstract derivation of the ME. As mentioned above, MFGs usually involve the coupling of two scales: a microscopic one and a macroscopic one. In this abstract setting, it is not possible to provide a meaning to these scales. Hence, for the moment, I shall simply consider two features of the model, one is computed backwardly in time and the other one is computed forwardly in time. Let me call $v$ the backward quantity and $y$ the forward one. Here, I assume that $v$ and $y$ are valued into respectively $F$ and $E$, two topological vector spaces. 

In most cases, the coupling between those two quantities appears as follows. Given $T> 0$ and a trajectory $(y_t)_{0\leq t \leq T}$, one can compute an associated trajectory $(v_t)_{0 \leq t \leq T}$ through the backward equation
\be\label{eq:backward}
\ba
\frac{d}{dt}v_t + B[y_t,v_t] = 0,\\
v_T = G[y_T],
\ea
\ee
where $B: E\times F \mapsto F$ and $G: E \mapsto F$ are given operators. The notation $[ \cdot]$ is simply here to insist upon the fact that the dependence can be non-local, notably in the important cases in which $E$ and $F$ are functional spaces.

Similarly, given a time $T > 0$ and a trajectory $(v_t)_{0\leq t \leq T}$, one can compute the forward evolution of the quantity labeled $y$ through
\be\label{eq:forward}
\ba
\frac{d}{dt}y_t + A[y_t,v_t] = 0,\\
y_0 \text{ given},
\ea
\ee
where $A: E\times F \mapsto E$ is another given operator.\\

Numerous questions can be asked on the connections between the two previous equations. The most natural ones are questions of existence and uniqueness of a solution $(v_s,y_s)_{t\leq s \leq T}$ for given $t\leq T$ and $y_t$ to the system \eqref{eq:backward}-\eqref{eq:forward} (on the time interval $[t,T]$). If an underlying game is present, such a solution corresponds to a Nash equilibrium of the game. Many developments around existence, uniqueness or stability of such systems are possible. I pass over them here, and just insist upon the fact that the main mathematical difficulty lies in the forward-backward nature of the coupling.\\

Another point of view consists in trying to condense all the information of the system \eqref{eq:backward}-\eqref{eq:forward} into a PDE whose unique solution $U$ is the function which maps $(t,y_t) \in (-\infty,T]\times E$ to $v_t \in F$, where $(v_s,y_s)_{t\leq s \leq T}$ is the unique solution to \eqref{eq:backward}-\eqref{eq:forward} on $[t,T]$. This (ambitious) program requires first the uniqueness of the solutions of \eqref{eq:backward}-\eqref{eq:forward}, and I shall come back to this important requirement later on. Assuming that this uniqueness result holds here, we find that, given a solution $(v_s,y_s)_{t\leq s \leq T}$ to \eqref{eq:backward}-\eqref{eq:forward} and such a map $U$, the following holds true
\be\label{eq:626}
\forall s \in [t,T], \quad U(s,y_s) = v_s.
\ee
Indeed, because of the uniqueness result, given a solution $(v_s,y_s)_{t \leq s \leq T}$ to \eqref{eq:backward}-\eqref{eq:forward}, for any $s \in [t,T]$, $(v_{s'},y_{s'})_{s \leq s' \leq T}$ is the unique solution to \eqref{eq:backward}-\eqref{eq:forward} on $[s,T]$ starting from $y_s$ at time $s$. Taking the time derivative of \eqref{eq:626} we obtain
$$
\forall s \in (t,T), \quad \partial_s U(s,y_s) - \langle  A[y_s,U(s,y_s)],D_y\rangle U(s,y_s) = - B[y_s,U(s,y_s)],
$$
where $\langle \xi,D_y\rangle U(s,y)$ just stands in a quite formal way, for the fact that we are taking the derivative of $U$, with respect to the variable $y$, in the direction given by $\xi$. Thus, using the definition of $U$, we finally arrive at the \textit{master equation}
\be\label{firstme}
-\partial_t U(t,y) + \langle A[y,U],D_y \rangle U = B[y,U] \text{ in } (- \infty,T)\times E .
\ee
Equation \eqref{firstme} has the typical form of a MFG master equation. In this general framework, the master equation is a non-linear transport equation. The forward-backward system \eqref{eq:backward}-\eqref{eq:forward} represents the characteristics of the ME \eqref{firstme}. In this setting, solving the system or the ME is equivalent. I just explained how to solve the ME given solutions to the systems for any initial conditions (in a uniqueness regime). Conversely, given a solution $U$ to \eqref{firstme}, it is a simple exercise to verify that the unique solution $(v_t,y_t)_{0 \leq t \leq T}$ to \eqref{eq:backward}-\eqref{eq:forward} on $[0,T]$ starting from $y_0$ at time $0$ satisfies
$$
\frac{d}{dt} y_t = -A[y_t,U(t,y_t)].
$$
Equation \eqref{firstme} is naturally associated to the terminal condition
$$
U(T,y) = G(y) \text{ in } E.
$$
Of course, it is quite often natural that the equation is only posed on a subset of $E$, but I do not enter into this question in this introductory part. 

\subsection{The various master equations}\label{sec:various}
In several situations, the exact form of  \eqref{firstme} is perturbed, for instance by terms modelling stochasticity. Below are some comments on the variations the previous structure can possess.

\subsubsection{Different spaces $E$ and $F$}
When $E,F$ are finite dimensional, equation \eqref{firstme} is self explanatory. The notation $\langle \cdot, \cdot \rangle$ then stands for the Euclidean scalar product and $D_y$ for the gradient operator in the direction $y$. Those types of MEs are called finite dimensional.\\

Exactly the same interpretation holds when $E=F$ is a separable Hilbert space. In this case, we talk about Hilbertian MEs.\\

In some cases, as we shall see in Section \ref{subsec:example}, $F$ will be some functional sub-space of $\mathcal C( \mo)$ and $E$ a subset of $\mmpo$. In this case, terms of the form $\langle  v,D_y\rangle U(s,y)$ have to be understood as a derivative with respect to a measure argument, and I refer to Part I for more details on this. I shall call such equations MEs on sets of measures.

\subsubsection{Terms modelling noises}
In certain important cases, the forward equation \eqref{eq:forward} is in fact of a stochastic nature. In those situations, the backward equation \eqref{eq:backward} is quite ill-adapted. Indeed, to solve it, one has to be given a path $(y_t)_{t \geq 0}$. But if this path is stochastic, typically adapted to a natural filtration $(\mathcal F_t)_{t \geq 0}$, then simply plugging it in \eqref{eq:backward} leads to a solution to this equation which is not at all adapted to $(\mathcal F_t)_{t \geq 0}$. Hence, this backward equation should be replaced by some backward stochastic equation, in the sense of backward SDEs introduced by Pardoux and Peng \citep{pardoux}. The study of the system \eqref{eq:backward}-\eqref{eq:forward} then becomes much more involved whereas the study of the associated master equation usually does not increase too much in difficulty. I said usually since, depending on the nature of the stochasticity, the terms arising from it can be quite different as we shall see. In these stochastic cases, the derivation of the master equation stays the same, it is still defined through the now stochastic forward-backward system. The only change is the need to use an appropriate variant of Ito's Lemma instead of the chain rule when differentiating the relation \eqref{eq:626}.\\

\textbf{Terms modelling jumps}\\

Random jumps in the dynamics of the forward variable $y$ can be modelled by considering the following extension of \eqref{firstme}
$$
-\partial_t U + \langle A[y,U(t,y)],D_y \rangle U + \lambda (U - \mathcal T^* U(t,\mathcal T y))= B[y,U] \text{ in } (- \infty,T)\times E,
$$
where $\lambda \geq 0$, $\mathcal T: E\mapsto E$ is a linear operator, and $\mathcal T^*: F \mapsto F$ is a form of an adjoint of $\mathcal T$. Such terms appear if, at random times given by a Poisson process of intensity $\lambda$, the variable $y$ is transformed through the map $\mathcal T$, i.e. the underlying variable jumps from $y$ to $\mathcal T y$. I do not recall the derivation of such terms and refer to the second section of our paper with Lasry and Lions \citep{bertucci2019some}.

Such terms are not so difficult to treat and provide in general good examples to test whether techniques on MEs without noise terms have a chance to be generalized to more complex ones. Typical examples of such operators $\mathcal T$ are translation.\\

\textbf{Additional stochastic parameter}\\

In other cases, randomness can be carried out through an additional variable $\theta \in \R^k$ for instance. In such models, the operators $A$ and $B$ depend on the value $\theta$ of the parameters as well, and this value satisfies an SDE of the form
$$
d \theta_t= b(\theta_t)dt + \sigma(\theta_t)dW_t,
$$
where $(W_t)_{t \geq 0}$ is a standard $k$ dimensional Brownian motion on $\R^k$, and $b:\R^k \mapsto \R^k, \sigma : \R^k \mapsto \mathcal M_{k}(\R)$ are given Lipschitz functions. In such cases, the solution to the ME has to depend on the value of the parameters as well. We are then led to consider MEs of the form
$$
\ba
-\partial_t U(t,y,\theta) + \langle A[y,U,\theta],D_y \rangle U - &b(\theta)\nabla_\theta U -\frac 12Tr(\sigma\sigma^T(\theta)D^2_\theta U)\\
&= B[y,U,\theta] \text{ in } (- \infty,T)\times E\times \R^k.
\ea
$$

Typical examples of such situations include cases in which $\theta$ represents an environmental factor, such as the climate, or a price which affects all players.\\

\textbf{The variable $y$ is perturbed by a diffusion}\\

Finally, cases have been considered in which the forward variable $y$ is directly affected by a diffusion, that is cases in which \eqref{eq:forward} is replaced by an SDE or stochastic PDE for instance. These kinds of models need to be treated with a lot of care. Indeed, if \eqref{eq:forward} is stochastic, the variants of Ito's Lemma needed to differentiate \eqref{eq:626} might be non-trivial and involve second order derivatives on $E$. General forms of such equations are then simply MEs of the form of \eqref{firstme} with second order terms in $y$.\\

I insist upon the fact that such models, despite the mathematical attractiveness that follows from the analysis of these difficult second order terms, can often be avoided and be recasted simply through an additional stochastic parameter. 

\subsection{Creation of singularities and how to prevent them}\label{sec:shocks}
As mentioned in Section \ref{sec:derive}, the derivation of the ME only makes sense in a uniqueness regime of the forward-backward system \eqref{eq:backward}-\eqref{eq:forward}. Here, I explain on a simple example what happens when we are not in such a regime, and give a formal definition of the monotone regime, that will be detailed in several examples later on.

\subsubsection{Singularities on the simplest master equation}
The simplest non-linear ME one can think of might be the Burgers' equation in one dimension
\be\label{burgers}
\partial_t u+ u\, \partial_x u = 0, \text{ in } (0,\infty)\times \R,
\ee
with initial condition $u_0: \R \mapsto \R$ given. (The time has simply been reversed to study a forward equation instead of a backward one.) The solution is here a function $u : [0,\infty)\times \R \mapsto \R$. One can check that there is an explicit formula for the solution to \eqref{burgers} which is given by
\be\label{formula:burgers}
\forall t \geq 0, x \in \R, \quad u(t,x) = u_0( (Id + t u_0)^{-1}(x)),
\ee
where $(Id + t u_0)^{-1}$ refers to the inverse function of $Id + tu_0$, and where $Id:x \mapsto x$ is the identity.\\

The previous formula turns out to be extremely insightful on the type of singularity that can happen for the ME. Indeed, this formula only makes sense when $Id + t u_0$ is invertible. If $u_0$ is $C$-Lipschitz, we know that such a condition will be true for all time $t < C^{-1}$. It is also a simple exercise to verify that the solution given by formula \eqref{formula:burgers} is the unique reasonable solution to \eqref{burgers}.

Starting from a $C$-Lipschitz initial condition $u_0$, if we reach a time $t^*$ such that $Id + t^* u_0$ is no longer invertible, it happens that $u(t,\cdot)$ will develop a singularity as $t \to t^*$ and that no continuous solution exists after $t^*$. In other words, a discontinuity appears, even if we started from a Lipschitz initial condition, and it is quite difficult to make sense of the equation after that time.

\begin{Rem}\label{rem:multid}
In this one dimensional setting, one can select a "good" discontinuous solution after $t^*$, namely by using the notion of entropic solution to the associated scalar conservation law. However, if we consider the multi-dimensional equation given by
$$
\partial_t u+ (u \cdot\nabla_x) u = 0, \text{ in } (0,\infty)\times \R^d,
$$
where now $u_0: \R^d \mapsto \R^d$ and $u: [0,\infty)\times \R^d \mapsto \R^d$, the situation is much more complex. The formula \eqref{formula:burgers} is still valid and the study of the solution before the appearance of discontinuities is possible. But what happens after that is still only poorly understood...
\end{Rem}

\subsubsection{Monotone regimes to prevent the creation of singularities}
The formula \eqref{formula:burgers} makes quite apparent the requirements needed on $u_0$ to allow its use for all time. Indeed, if $u_0$ is non-decreasing, then $Id + t u_0$ is invertible for all time $t \geq 0$, and the equation is somehow well posed, even if $u_0$ is not (Lipschitz) continuous for instance. This fact generalizes easily to the multi-dimensional setting of Remark \ref{rem:multid} if we replace the non-decreasing property with the notion of monotonicity.

In addition to preventing the creation of singularities, this non-decreasing assumption can be interpreted in terms of the associated forward-backward system, precisely as a uniqueness regime. Indeed, the forward-backward system associated to \eqref{burgers} is given by
$$
\begin{cases}
\ba
\dot v_t = 0,&\quad v_0 = u_0(y_0),\\
\dot y_t = v_t,&\quad y_T \text{ given}.
\ea
\end{cases}$$
Taking $(v^1,y^1)$ and $(v^2,y^2)$ two solutions of the system, we see that
$$
\ba
-(u_0(y^1_0)-u_0(y^2_0))(y^1_0 - y^2_0) = \int_0^T\frac{d}{dt}\left[(v^1_t-v^2_t)(y^1_t-y^2_t)\right]dt = \int_0^T(v^1_t-v^2_t)^2dt \geq 0.
\ea
$$
Hence, because $u_0$ is non-decreasing, we obtain that $v^1_t = v^2_t$ for almost all $t$, hence the fact that $(v^1,y^1)=(v^2,y^2)$ since both are constant.\\

A simple generalization of the previous argument leads to the following formal definition of a monotone regime for \eqref{firstme}. 
\begin{Def}
The monotone regime for \eqref{firstme} is the conjunction of the following properties:
\begin{itemize}
\item The space $F \subset E'$.
\item The map $G: E \mapsto E'$ is monotone.
\item The map $(A,B): E\times E'\mapsto E\times E'$ is monotone, that is for all $(x^1,y^1),(x^2,y^2) \in E\times E'$
$$
\langle A[x^1,y^1] - A[x^2,y^2],y^1-y^2\rangle_{E\times E'} + \langle B[x^1,y^1] - B[x^2,y^2],x^1-x^2\rangle_{E'\times E} \geq 0.
$$
\end{itemize}
\end{Def}
Of course, several extensions can be considered, for instance by considering maps $A,B,G$ which are not defined everywhere. To obtain uniqueness of the associated forward-backward system, some strict monotonicity is needed. Furthermore, it is an open question to understand to what extent such monotone regimes are limiting assumptions, or in other words if there exist other general classes of assumptions which yield uniqueness of the forward-backward system.

 Note that on the simple Burgers equation mentioned above, the formula of the explicit solution is defined for all times if and only if we are in the monotone regime, as asking that $Id + t u_0$ is invertible for all $t \geq 0$ is equivalent to demanding that $u_0$ is non-decreasing. Hence, on this simple equation, it is not possible to avoid this monotone regime to consider regular solutions on arbitrary long time intervals.

\subsection{The most studied example}\label{subsec:example}
Here I present the original example of a MFG considered by Lasry and Lions \citep{lasry2007mean}. As I shall show, considering slight variations on this model, we can arrive at the various MEs described in section \ref{sec:various}.\\

\subsubsection{Definition of the mean field game}
The following description of the game is only formal. We consider a continuum of players of mass $1$. Each player is described by a parameter $x \in \R^d$. We call $\R^d$ the state space, and $x$ the state of the player. The players are initially distributed according to $µ \in \mprd$. All players control the evolution of their state. The evolution of the state of a player is given by the process $(X_t)_{t \geq 0}$, solution to the SDE
\be\label{stateplayer}
dX_t = \alpha_t dt + \sqrt{2\sigma}dW_t,
\ee
where $\sigma \geq 0$ is a parameter, $(W_t)_{t \geq 0}$ is a standard Brownian motion on a standard probability space $(\Omega, \mathcal A, (\mathcal F_t)_{t \geq 0},\mathbb P)$ and $(\alpha_t)_{t \geq 0}$ is the control of the player considered. The control has to be an adapted $\R^d$ valued process, in particular, a player cannot anticipate on the future of the Brownian motion. Each player chooses its control $(\alpha_t)_{t \geq 0}$ on its own and two distinct players are facing independent realizations of the Brownian motion. A player whose state evolves through $(X_t)_{t \geq 0}$ and whose control is $(\alpha_t)_{t \geq 0}$ faces the cost
\be\label{eq:costex}
\mathbb E \left[ \int_0^TL(X_t,\alpha_t,m_t)dt + G(X_T,m_T)  \right],
\ee
where $L: \R^d \times \R^d \times \mprd \mapsto \R_+$ and $G: \R^d \times \mprd \mapsto \R_+$ are given cost functions, $T > 0$ is the duration of the game and $(m_t)_{t \geq 0}$ is the $\mprd$ valued process which describes the evolution of the distribution of players. That is, for a measurable $A \subset \R^d$ and $t \in [0,T]$, $m_t(A)$ is the proportion of players which are in $A$ at time $t$. The presence of the expectation in \eqref{eq:costex} stands for the fact that the players are risk neutral with respect to the randomness present in the evolution of their state.

\subsubsection{Equilibria of the MFG}
I insist upon the fact that \eqref{eq:costex} does not simply represent a cost that a given player has to minimize in $(\alpha_t)_{t \geq 0}$. It is only the case if $(m_t)_{t \in [0,T]}$ is known, which is clearly not the case in general. If it is the case, then the players can find their optimal response to $(m_t)_{t \geq 0}$ by means of dynamic programming. Namely they can introduce the value function $u:[0,T]\times \R^d \mapsto \R$ defined by
$$
u(t,x) := \inf_{(\alpha_s)_{s \in [t,T]}} \mathbb E \left[ \int_t^TL(X_s,\alpha_s,m_s)ds + G(X_T,m_T) \bigg| X_t = x \right],
$$
where the infimum is taken over admissible controls and \eqref{stateplayer} holds in $(t,T)$. They then solve the associated HJB equation which is given by
\be\label{eq:hjb1}
\ba
-\partial_t u - \sigma \Delta u + H(x,\nabla_x u,m_t) = 0 \text{ in } (0,T)\times \R^d,\\
u(T,x) = G(x,m_T) \text{ in } \R^d,
\ea
\ee
where $H$, the Hamiltonian of the problem, is given by $H(x,p,m) = \sup_{\alpha \in \R^d}\{-L(x,\alpha,m) - \alpha\cdot p\}$. Given a sufficiently smooth solution $u$ to \eqref{eq:hjb1} (which thus requires a path $(m_t)_{t \geq 0}$ to be given as well), a player can choose the control 
$$
\alpha^*_t := -D_pH(X_t,\nabla_x u(t,X_t),m_t)
$$
to minimize \eqref{eq:costex}.\\

As just mentioned, players have to decide which path $(m_t)_{t \geq 0}$ they are going to anticipate. Particularly interesting anticipations are called equilibria\footnote{There is a slight abuse of terminology here, as equilibria usually refer to the actions of the players, and not to their distribution.}, which are defined in the following way: a path $(m_t)_{t \geq 0}$ is an equilibrium of the MFG if, given that all players anticipate $(m_t)_{t \geq 0}$, there is a way for the players to play optimally according to this anticipation, which is such that the evolution of their distribution in the state space is actually going to be $(m_t)_{t \in [0,T]}$.\\

I leave as an exercise to the interested reader to verify, using the interpretation of Fokker-Planck equation of Remark \ref{rem:fpvar}, that we can read in the coupling of the following system the notion of equilibrium.
\be\label{mfg:sys}
\begin{cases}
-\partial_t u - \sigma \Delta u + H(x,\nabla_x u,m_t) = 0 \text{ in } (0,T)\times \R^d,\\
\partial_t m - \sigma \Delta m - \text{div}(D_pH(x,\nabla_x u,m)m) = 0 \text{ in } (0,T)\times \R^d,\\
m|_{t = 0} = µ, \quad u(T,x) = G(x,m_T) \text{ in } \R^d.
\end{cases}
\ee
All equilibria $(m_t)_{t \geq 0}$ are associated to a solution to \eqref{mfg:sys} and vice-versa.\\

I do not give a precise statement about the associated monotone regime for the moment, and refer to Hypothesis \ref{hyp:mon} below for the case in which $H(x,p,m) = \tilde H(x,p) - f(m)(x)$.

\subsubsection{The associated master equation}
Before deriving the associated master equation, I need to make precise what type of derivative on the space of measures I am going to use. As mentioned at the end of Part I, it is natural in this MFG context that the underlying measure evolves more generally than simply by following flows of ODEs or SDEs. For instance, players could jump in the state space, or they could leave or enter the game. This strongly suggests to work with the vertical derivative. This might seem strange, as for the present game, such behaviour is not authorized, and it will remain the case in the rest of this part, because I choose to concentrate the study on a simple master equation. However, I strongly encourage the reader to extend the following study to a case in which the Fokker-Planck equation in \eqref{mfg:sys} contains a term in $\lambda m$ for $\lambda$ a constant, to realize that in such a case, the use of the V-derivative is necessary.\\

Following the derivation of the ME given in Section \ref{sec:derive}, the ME associated to the forward-backward system \eqref{mfg:sys} is
\be\label{me}
\ba
-\partial_t U - &\left\langle \nabla_µ U(t,µ,x,\cdot),\text{div}(D_pH(\cdot,\nabla_x U(t,µ,\cdot),µ)µ) + \sigma \Delta µ\right\rangle\\
&\quad - \sigma \Delta_x U + H(x,\nabla_x U,µ) = 0 \quad \text{in } (0,T)\times \mprd\times \R^d.
\ea
\ee
It is naturally associated to the terminal condition 
$$
U(T,µ,x) = G(µ)(x) \text{ in } \mprd \times \R^d.
$$
In \eqref{me}, the notation $\langle \cdot,\cdot \rangle$ stands for the extension of the $L^2(\R^d)$ scalar product. Sticking with the terminology of Section \ref{sec:derive}, the forward variable is here $µ$ or $m$ and $U(t,µ)$ is a function on the state space. However, sometimes we use the natural identification between $[0,T]\times \mprd \mapsto (\R^d \mapsto \R)$ and $[0,T]\times \mprd \times \R^d \mapsto \R$, such as is done for $G$ just above.\\

The equation \eqref{me} is the master equation associated to the MFG I just presented. Its solution $U$ should be referred to as the value of the MFG. Indeed, when it is well defined, $U(t,µ,x)$ is the (unique) cost a player in state $x \in \R^d$, at time $t \in [0,T]$, facing the distribution of other players $µ \in \mprd$ is going to pay, under the fundamental condition that players are going to follow the (unique) equilibrium. This is yet another opportunity to insist upon the fact that the solution to the ME is well defined only under a uniqueness assumption.\\

Note that by considering a finite discretization of $\R^d$, we can construct a finite dimensional ME which is an approximation of \eqref{me}. Indeed, it suffices to consider standard finite difference schemes as approximation to the forward-backward system \eqref{mfg:sys}, and then derive the associated ME. I do not give all the details of such an approximation as it is: i) standard in view of the literature on finite difference schemes, ii) quite heavy in terms of notation that we shall not need in the following.

\subsection{The Hilbertian framework}
There is another ME that we can derive directly from this MFG in the case $\sigma = 0$, which is called the Hilbertian ME. It consists in a much more Lagrangian approach to the description of the distribution of players. In terms of modeling, this approach consists in replacing the probability measure $µ\in \mprd$ by a random variable $X: \Omega \mapsto \R^d$, where $(\Omega,\mathcal A,\mathbb P)$ is a standard probability space, which is the space of labels. One can for instance think about $\Omega$ as the space of names of players. Hence, a distribution of players is here simply $X: \Omega \mapsto \R^d$ which associates to a player (through its label) its position. In this framework, a strategy profile (i.e. a strategy for all players) can be summarized by a map $\alpha: [0,T]\times \Omega \mapsto \R^d$ and it yields the following evolution from an initial distribution of players $X_0: \Omega \mapsto \R^d$
$$
dX_t = \alpha_t dt.
$$
Hence the Hilbertian analogue of \eqref{mfg:sys} is given by 
$$
\begin{cases}
-\partial_t u(t,x)  + H(x,\nabla_x u,\mathcal L(X_t)) = 0 \text{ in } (0,T)\times \R^d,\\
dX_t = -D_pH(X_t,\nabla_x u(t,X_t),\mathcal L(X_t))dt  \text{ in } (0,T)\times \Omega,\\
X_0 \text{ given}, \quad u(T,x) = G(x,\mathcal L(X_T)) \text{ in } \R^d.
\end{cases}
$$
Instead of directly considering the associated master equation, we differentiate with respect to $x$ the HJB equation above and set $v = \nabla_x u$ to arrive at the forward-backward system
$$
\begin{cases}
-\partial_t v(t,x)  + D_xH(x,v,\mathcal L(X_t)) + D_pH(x,v,\mathcal L(X_t))D_x v = 0 \text{ in } (0,T)\times \R^d,\\
dX_t = -D_pH(X_t,v(t,X_t),\mathcal L(X_t))dt  \text{ in } (0,T)\times \Omega,\\
X_0 \text{ given}, \quad v(T,x) = D_xG(x,\mathcal L(X_T)) \text{ in } \R^d.
\end{cases}
$$
Finally, if we considered directly the ME associated to the previous forward-backward system, we would end up with a solution as a function $ \mathcal V: [0,T]\times (\Omega \mapsto \R^d) \mapsto (\R^d \mapsto \R^d)$. Instead of such a function, we shall consider, to simplify notation, a function $V: [0,T]\times (\Omega \mapsto \R^d) \mapsto (\Omega \mapsto \R^d)$ defined through the following relation
\be\label{rel:lift}
V(t,X)(\omega) := \mathcal V(t,X)(X(\omega)).
\ee
Finally, restricting ourselves to the subset $L^2(\Omega,\R^d) \subset (\Omega \mapsto \R^d)$, we are able to write the ME
\be\label{me:hilbertian}
-\partial_t V + D_xH(X,V,\mathcal L(X)) + \langle D_pH(X,V,\mathcal L(X)), \nabla_X \rangle V(t,X) = 0 \text{ in } (0,T)\times L^2(\Omega,\R^d),
\ee
where $\langle \cdot,\cdot \rangle$ here denotes the canonical scalar product in $\mathbb H :=L^2(\Omega,\R^d)$ and $\nabla_X$ the gradient operator in $\mathbb H$.\\

Several comments are here in order:
\begin{itemize}
\item The derivation of \eqref{me:hilbertian} might seem mysterious if one does not remark that differentiating \eqref{rel:lift}, say formally in $\mathbb H$, leads to 
$$
\nabla_X V(t,X)(\omega) = \nabla_x \mathcal V(t,X)(X(\omega)) + \nabla_X\mathcal V(t,X)(X(\omega))(\omega).
$$
\item The restriction to $L^2(\Omega,\R^d)$ is purely arbitrary at this formal level, and only has the advantage of using the Hilbertian structure. In practice, this kind of restriction has to be justified by the fact that at the initial time, $\mathbb E[|X_0|^2] < \infty$, as well as growth conditions on $H$ to ensure that the dynamics will stay inside $L^2(\Omega,\R^d)$. Similar equations in $L^p(\Omega,\R^d)$ space for $p \ne 2$ could also have been considered.
\item However probabilistic this framework may seem, note that the underlying model is purely deterministic, and the language of random variables is used to arrive at a compact formulation.
\end{itemize}

One of the main interests of this Hilbertian reformulation is that it makes apparent another monotone structure in the MFG, different from the one described in Hypothesis \ref{hyp:mon} below. Indeed, in this case, the natural monotone structure is given by the following.
\begin{Def}
The monotone regime associated to the Hilbertian setting consists of the following two assumptions:
\begin{itemize}
\item The map $X \mapsto D_xG(\mathcal L(X))(X)$ is monotone as a map from $\mathbb H $ to $ \mathbb H$.
\item The map $(X,Y) \mapsto (-D_xH(X,Y,\mathcal L(X)),D_pH(X,Y,\mathcal L(X)))$ is monotone as a map from $\mathbb H\times \mathbb H$ to $ \mathbb H \times \mathbb H$.
\end{itemize}
\end{Def}

\begin{Rem}
The fact that this Hilbertian framework not only works but can produce similar results of existence and uniqueness as the previous setting, is an excellent justification of the fact that, contrary to what can be read in several papers or books, indistinguishability is absolutely not a feature of MFG. Indeed, in the Hilbertian case, players can be distinguished by their label and the study can be carried out in a similar fashion.
\end{Rem}

\begin{Rem}
The Hilbertian framework is well suited for an analysis of the associated master equation through the use of the H-derivative, but I do not expand more on it in this simple presentation.
\end{Rem}

\subsection{A master equation without a game}
In contrast with the previous example which can be looked at as the most typical MFG, I now insist on a model which yields a ME, without any proper underlying game, but simply because of the existence of an underlying forward-backward structure.\\

Consider a market in which a good (or product) is exchanged. We assume that there is a constant stream of demand $D$ for the product, when expressed in money. The product is supplied by small producers. In order to produce the product, and then to sell it on the market, they need to buy a unit of production capital, which will produce a constant unit stream of product. The products are then sold on the market at a certain price. For the sake of simplicity, we assume that the market will only be open during the time interval $[0,T]$. At time $t \in [0,T]$, we denote by $K_t> 0$ the number of active units of production, by $V_t$ the value that a unit of production will generate in $[t,T]$, and by $P_t$ the price at which the product is sold.\\

I assume perfect competition between the producers. Because the amount of cash available to buy the product is constant, the price is given by $P_t = D/K_t$. Given a trajectory of number of production units $(K_t)_{t \in [0,T]}$, the value $V_t$ is computed backwardly by summing over the future revenue it will generate as follows
$$
V_t = \int_t^T\frac{D}{K_s} ds.
$$
On the other hand, I make the assumption that units of production are bought and sold depending on the revenue they generate, so that, at time $t$, if machines (units of production) are generating a value $V_t$, then we have
$$
\frac{d}{dt}K_t = \lambda(V_t - c),
$$
for $\lambda,c >0$ two constants. Hence, we have a forward backward structure, simply by differentiating the integral definition of $V_t$, which can be summarized in
$$
\begin{cases}
\frac{d}{dt}V_t = -\frac{D}{K_t},\\
\frac{d}{dt}K_t = \lambda(V_t - c),\\
K_0 \text{ given}, \quad V_T = 0.
\end{cases}
$$
There is a ME associated to this forward-backward system, which is satisfied by the value $U(t,K)$ a unit of production will generate, given that the current number of active unit of production is $K$ and that we are at time $t$. It is the following

$$
\ba
-\partial_t &U -\lambda(U - c) \partial_K U =  \frac{D}{K} \text{ in } (0,T)\times (0,\infty),\\
&U(T,K) = 0 \text{ in } (0,\infty).
\ea
$$
It is worth noting that this ME is quite easy to understand, as it does not involve analysis in infinite dimension.

Finally, no proper game has been defined between the producers in the market, only rules, but this ME still makes sense.

\begin{Rem}
In this case, we are in a monotone regime, in fact an anti-monotone one (that is $-U$ solves a ME in monotone regime). Indeed, the terminal value is constant thus monotone, and $\lambda >0,D>0$ yield the correct monotonicity properties.
\end{Rem}

\subsection{Bibliographical comments}
For a nice introduction to MFG, I refer the reader to Cardaliaguet and Porretta \citep{cardaliaguet2021introduction}. In addition to the references to which I pointed at the beginning of this section, I also refer the reader interested in MFG systems to Achdou and Capuzzo-Dolcetta \citep{achdouitalo} for numerical methods and toward Porretta \citep{porretta} for weak notions of solutions to MFG systems. I will not comment further on the general literature on MFG but rather focus on the one associated to the ME. Other economic models were derived in Gabaix et al \citep{moll1} and in Achdou et al \citep{moll2} using MFGs. The ME was derived by Lasry and Lions while discussing MFG in Chicago with Lucas, who informed them of the Krusell-Smith model \citep{krusell1998income} which has a stochastic component in the variation of macroscopic quantities (the so-called common noise). It then became apparent that MEs are essential tools in macro-economics and links were made with previous models which were in fact finite dimensional MEs such as Scheinkman and Weiss \citep{scheinkman1986borrowing} or Conze, Lasry and Scheinkman \citep{conzelasry}. Recently, numerical computations to provide a solution to the Krusell-Smith model were provided by Achdou, Lasry and Lions \citep{achdoukrusell}.

The first mathematical developments on MEs are due to PL Lions, in finite dimensional cases, and were presented in his lectures at Collège de France \citep{lions20112012}. In 2015, he presented an unpublished work on the Hilbertian case, extending his analysis in finite dimension. In the meantime, an extensive study of the ME through the formula \eqref{eq:626} was done by Chassagneux-Crisan-Delarue \citep{chassagneux}, which was later extended in a sort of complete study in the book of Cardaliaguet-Delarue-Lasry-Lions \citep{cdll}, which was also continued by Jakobsen and Rutkowski \citep{jakobsen2025master} in a slightly different setting. Accounts of these studies are also presented in the book of Carmona-Delarue \citep{carmona2017probabilistic2}. More recent and efficient studies have been carried out since, and they will be presented or mentioned below.

The monotone regime was identified in the first paper on MFG by Lasry and Lions \citep{lasry2007mean} as the natural condition to ensure uniqueness of solutions to the MFG system, and thus of the underlying Nash equilibria. They also found variants of these assumptions, for instance in finite dimensional frameworks, in cases in which the Hamiltonian $H$ depends on both $\nabla_x u$ and $m$, as well as in the Hilbertian setting. This latter setting was also used in Gangbo and Meszaros \citep{gangbo} and Gangbo, Meszaros, Mou and Zhang \citep{gangbo2}, independently, under the name displacement monotonicity.

The last model of a master equation without a game originates from the modelling done in either Achdou et al \citep{edmond} or in Bertucci, Bertucci, Lasry and Lions \citep{bitcoin}.

\newpage

\section{Uniqueness and stability of the solution to the master equation}
As mentioned in Section \ref{sec:shocks}, MEs will develop singularities in general and their study in the presence of singularities is not very realistic at the moment. In this section, we shall see that, under monotonicity assumptions, we can define a notion of solution to the ME which enjoys stability properties, as well as uniqueness. In the next part, when dealing with questions of existence, we shall see that monotonicity assumptions are also important to prevent the creation of singularities.\\

The rest of this section is divided into two parts. First, I shall be concerned with finite dimensional MEs, then with MEs on spaces of probability measures.

\subsection{A first result of well-posedness}
Consider a bounded domain $\Omega$ of $\R^d$ whose boundary is $\mathcal C^2$, functions $U_0: \bar \Omega \mapsto \R^d, F,G : \bar \Omega\times \R^d \mapsto \R^d$ defined on the closure $\bar \Omega$ of $\Omega$ and the following ME
\be\label{me:finite}
\ba
\partial_t U(t,q) + \langle F(q,U(t,q)),\nabla_q\rangle U(t,q) = G(q,U(t,q)) \text{ in } (0,\infty)\times \bar \Omega,\\
U|_{t=0} = U_0.
\ea
\ee
The unknown $U$ in the previous equation is valued in $\R^d$. In \eqref{me:finite} and in all that follows which is concerned with this finite dimensional ME, $\langle \cdot,\cdot \rangle$ denotes the Euclidean scalar product. Note also that the previous equation has to hold on the boundary of $\bar \Omega$. I apologize for the use of $\Omega$ which has nothing to do with a probability space in this context of finite dimensional MEs.

Assumptions on $F,G$ and $U_0$ shall be made later on. As I do not want to enter a rather complex discussion about what happens at $\partial \Omega$, I shall assume that the following condition holds.
\begin{hyp}\label{hyp:F}
For all $q \in \partial \Omega, V \in \R^d$, 
$$
\langle F(q,V),\eta(q)\rangle \geq 0,
$$
where $\eta(q)$ is the unit outward normal vector to $\partial \Omega$ at the point $q$.
\end{hyp}
This assumption is purely technical and should not worry the reader too much. In the usual infinite dimensional master equation, the analogous Hypothesis \ref{hyp:F} is verified. It simply states that any potential underlying trajectory of the forward variable stays inside $\bar \Omega$. Hypothesis \ref{hyp:F} shall always be valid and I will not recall it in every result.\\

In this finite dimensional setting, the monotone regime of interest is summarized in the following assumption.
\begin{hyp}\label{hyp:monf}
The mapping $(G,F): \bar\Omega \times \R^d \mapsto \R^{2d}$ is monotone, and so is $U_0: \bar\Omega\mapsto \R^d$.
\end{hyp}

The next result is probably the most important one on the monotone regime. Its proof highlights how monotonicity implies uniqueness of solutions to \eqref{me:finite}.
\begin{Prop}\label{prop:propagationN}
Assume that Hypothesis \ref{hyp:monf} holds. Then, there exists at most one smooth solution $U: [0,\infty)\times \bar \Omega \mapsto \R^d$ of \eqref{me:finite}. Furthermore, such a solution is a monotone map in $q$, for all $t \geq 0$.
\end{Prop}
\begin{proof}
Consider $U$ and $V$ two such solutions and define $W: \R_+ \times \bar \Omega^2 \mapsto \R$ by
$$
W(t,q,p) = \langle U(t,q) - V(t,p),q-p\rangle.
$$

\textbf{Step 1: $W$ is non-negative.} Note that for any $t \geq 0, q,p \in \Omega$, we have the relations
$$
\nabla_q W(t,q,p) = D_qU(t,q)\cdot(q-p) + U(t,q) - V(t,p),
$$
$$
\nabla_p W(t,q,p) = D_pV(t,p)\cdot (p-q) + V(t,p) -U(t,q).
$$
Using the fact that both $U$ and $V$ solve \eqref{me:finite}, we thus obtain that
$$
\ba
\partial_t& W + F(q,U)\cdot \nabla_q W + F(p,V) \cdot \nabla_p W = \\
&= \langle G(q,U) - G(p,V),q-p\rangle + \langle F(q,U)-F(p,V),U-V\rangle \text{ in } (0,\infty) \times \Omega^2.
\ea
$$
Observe that, since $U_0$ is monotone, $W|_{t = 0} \geq 0$. From the previous PDE, by using the fact that $(G,F)$ is monotone, it follows
$$
\partial_t W + F(q,U)\cdot \nabla_q W + F(p,V) \cdot \nabla_p W \geq 0.
$$
We want to show that $W \geq 0$. Assume that it is not the case. Thus, there exists $\lambda > 0$ such that it is also not the case for $W + \lambda t$. Consider $T > 0$ and a point of minimum $(t_*,q_*,p_*)$ of $W +\lambda t$ on $[0,T]\times \bar\Omega^2$. Clearly $t_* > 0$. Note that if $t_* < T$, then $\partial_t W(t_*,q_*,p_*) = -\lambda$ and if $t_* = T$ then $\partial_t W(t_*,q_*,p_*) \leq -\lambda$. If $q_*,p_* \in \Omega$, then we obtain that $\nabla_q W(t_*,q_*,p_*) = \nabla_p W(t_*,q_*,p_*) = 0$ and thus that
$$
- \lambda \geq 0,
$$
which is a contradiction. Hence, it remains to treat the case where $q_*$ or $p_*\in \partial \Omega$. It suffices to remark that in this case, Hypothesis \ref{hyp:F} implies that $\langle F(q_*,U(t_*,q_*)), \nabla_q W(t_*,q_*,p_*) \rangle \leq 0$ if $q_* \in \partial \Omega$. Hence, we still obtain that $-\lambda \geq 0$ which is a contradiction. Thus $W \geq 0$.\\

\textbf{Step 2: $U= V$.} Consider a point $(t,q) \in [0,\infty)\times \Omega$. For $\epsilon > 0$ sufficiently small, the ball of center $q$ and of radius $\eps$ lies in $\Omega$. Thus, for any $v\in \R^d, \|v \| \leq 1$, 
$$
\langle U(t,q + \eps v) - V(t,q),\eps v\rangle \geq 0.
$$
This implies that $U= V$ on $\R_+ \times \Omega$, since they are continuous, the equality also holds on $\bar \Omega$. The monotonicity of $U$ then follows by non-negativity of $W$ since $V= U$.
\end{proof}
As we shall see in Section \ref{sec:meexist}, the monotone regime is in fact not essential in the previous result, as uniqueness can be recovered directly from the regularity. Nonetheless, the proof of this result will allow us to understand a key feature of the stability of the solutions of MEs in monotone regimes.

\subsection{Monotone solutions to finite dimensional master equations}

In the previous proof, the uniqueness of solutions to \eqref{me:finite} was obtained by proving the positivity of $W$ through a comparison principle-type argument. The main idea I want to follow here is to define a concept of solution with the minimum regularity requirements, so that the previous proof can still be carried out.

The use of this argument formally required three steps:
\begin{itemize}
\item Using the relation that for all $i \in \{1,\ldots,d\}$,
\be \label{relationwn}
\partial_{q_i} W(t,q,p) = \langle D_{q_i}U(t,q), q - p \rangle + U^i(t,q) - V^i(t,p),
\ee
as well as the equation satisfied by $U$ (and its counterpart for $V$) to derive an equation for $W$.
\item Considering a point of minimum of $W$.
\item Using \eqref{relationwn} at this point of minimum to derive a contradiction.
\end{itemize}

Therefore, to maintain this strategy, we only need information on $U$ at the point of minimum of $W$. Since there is no way to know a priori what those points are going to be, we need to require information on $U$ at points which are minima of functions of the form $(t,q) \mapsto \langle U(t,q) - V, q - p \rangle$ for $V \in \mathbb{R}^d$, $p \in \bar \Omega$. Moreover, at such a point of minimum, the equation \eqref{me:finite} satisfied by $U$ is only used to obtain information about the variations of $\langle U(t,q), q - p \rangle$.\\

To illustrate the above, let us consider the following case: let $V \in \mathbb{R}^d$, $p \in \bar \Omega$, and $(t_*,q_*)$ be a minimum point of $(t,q) \mapsto \langle U(t,q) - V, q - p \rangle$. Using the fact that $U$ is a solution to \eqref{me:finite} on $\R_+ \times \bar \Omega$, we obtain
\be
\partial_t \langle U(t_*,q_*), q_* - p \rangle + F(q_*, U(t_*,q_*)) \cdot D_q U(t_*,q_*) \cdot (q_* - p) = \langle G(q_*, U(t_*,q_*)), q_* - p \rangle.
\ee

Assume that $q_*$ is in the interior of $\Omega$. Then we have
\be
0 = U(t_*,q_*) - V + \langle D_q U(t_*,q_*), q_* - p \rangle.
\ee

This is clearly the analogue of \eqref{relationwn}. We then deduce that
\be\label{eq1234}
\partial_t \langle U(t_*,q_*), q_* - p \rangle = \langle G(q_*, U(t_*,q_*)), q_* - p \rangle + \langle F(q_*, U(t_*,q_*)), U(t_*,q_*) - V \rangle.
\ee
On the other hand, if $q_* \in \partial \Omega$, thanks to Hypothesis \ref{hyp:F}, we still have
\be \label{eq1}
\partial_t \langle U(t_*,q_*), q_* - p \rangle \geq \langle G(q_*, U(t_*,q_*)), q_* - p \rangle + \langle F(q_*, U(t_*,q_*)), U(t_*,q_*) - V \rangle.
\ee
This is exactly the type of information that will be encoded in the definition of so-called monotone solutions, or M-solutions.\\

A fundamental observation is now in order: the relations \eqref{eq1234} and \eqref{eq1} above do not involve any derivatives of $U$ with respect to $q$! By invoking the theory of viscosity solutions, initiated by Crandall and Lions \citep{crandall1983viscosity}, we will handle the time derivative term using a test function\footnote{More details on viscosity solutions are given in Part III of this book.}. To this end, observe that if $\vartheta$ is a smooth function of $t$ and $(t_*, q_*)$ is a minimum point of $(t,q) \mapsto \langle U(t,q) - V, q - p \rangle - \vartheta(t)$, then we can replace the first term in the above inequality to obtain
\be
\vartheta'(t_*) \geq \langle G(q_*, U(t_*, q_*)), q_* - p \rangle + \langle F(q_*, U(t_*, q_*)), U(t_*, q_*) - V \rangle.
\ee

We have now seen all the main motivations behind the following definition.

\begin{Def} \label{def:monn}
A continuous function $U : \mathbb{R}_+ \times \bar\Omega \mapsto \mathbb{R}^d$ is called a \emph{M-solution} to \eqref{me:finite} if for all $T > 0$, $V \in \mathbb{R}^d$, $p \in \bar \Omega$, $\vartheta \in \mathcal{C}^1(\mathbb{R})$, and $(t_*, q_*) \in (0,T] \times \bar\Omega$, point of strict minimum of the function $(t,q) \mapsto \langle U(t,q) - V, q - p \rangle - \vartheta(t)$ on $[0,T] \times \bar\Omega$, we have
\be
\vartheta'(t_*) \geq \langle G(q_*, U(t_*, q_*)), q_* - p \rangle + \langle F(q_*, U(t_*, q_*)), U(t_*, q_*) - V \rangle.
\ee
\end{Def}
I called such solutions monotone solutions when I introduced them. In this book addressing broader topics, I use this terminology M-solutions to avoid confusion.
\begin{Rem}
I only require this relation to hold at points of strict minimum for two main reasons: first, because it is sufficient, as we will see; and second, because in some cases it is not clear whether it would hold at all minimum points. It is for instance the case for MEs of MFGs of optimal stopping, as I explained in \citep{bertucci2021monotone}. It will also allow us to obtain strong stability properties for such solutions.
\end{Rem}

\subsection{Uniqueness and stability in finite dimension}

Given the computations that preceded Definition \ref{def:monn}, we immediately have the following result.

\begin{Prop}
Let $U:[0,\infty)\times \bar \Omega \mapsto \R^d$ be a classical solution to \eqref{me:finite}. Then $U$ is also an M-solution to \eqref{me:finite}.
\end{Prop}

Given the way we constructed Definition \ref{def:monn}, the following uniqueness result should not come as a surprise. 

\begin{Theorem}\label{thm:uniqmonf}
Under Hypothesis \ref{hyp:monf}, there exists at most one M-solution to \eqref{me:finite} with initial condition $U_0$. Moreover, if such a solution $U$ exists, then for all $t \geq 0$, $U(t,\cdot)$ is a monotone map.
\end{Theorem}

\begin{proof}
Let $U$ and $V$ be two M-solutions. We want to show that $W$ defined by $W(t,q,p):= \langle U(t,q)-V(t,p),q-p\rangle$ is non-negative. Assume that there exists $(T,\tilde q,\tilde p) \in (0,\infty) \times \bar\Omega^2$ such that
$$
W(T,\tilde q,\tilde p) < 0.
$$

Then there exist $\kappa, \epsilon > 0$ such that for all $\alpha > 0$,
$$
\inf_{t,s \in [0,T], p,q} \langle U(t,q) - V(s,p), q - p \rangle + \alpha (t - s)^2 + \epsilon (t + s) \leq -\kappa < 0.
$$

Using a trivial variant of Stegall's Lemma in finite dimension (see Section \ref{sec:perturbed}), for any $\eta > 0$, there exist $a,b \in \mathbb{R}^d$, $\xi_1, \xi_2 \in \mathbb{R}$ such that $|a| + |b| + |\xi_1| + |\xi_2| \leq \eta$ and
$$
(t,s,q,p) \mapsto \langle U(t,q) - V(s,p), q - p \rangle + \alpha(t - s)^2 + \epsilon (t + s) + \langle a,q \rangle + \langle b,p \rangle + \xi_1 t + \xi_2 s
$$
has a strict minimum on $[0,T]^2 \times \bar \Omega^2$ at $(t_*, s_*, q_*, p_*)$. Assume first that $t_*, s_* > 0$. Since $U$ is an M-solution to \eqref{me:finite}, we obtain
$$
\ba
-\epsilon - 2\alpha(t_* - s_*) - \xi_1 \geq &\langle G(q_*, U(t_*, q_*)), q_* - p_* \rangle\\
& + \langle F(q_*, U(t_*, q_*)), U(t_*, q_*) - V(s_*, p_*) + a \rangle.
\ea
$$

The analogous inequality for $V$ gives
$$
\ba
-\epsilon - 2\alpha(s_* - t_*) - \xi_2 \geq &\langle G(p_*, V(s_*, p_*)), p_* - q_* \rangle\\
& + \langle F(p_*, V(s_*, p_*)), V(s_*, p_*) - U(t_*, q_*) + b \rangle.
\ea
$$

Adding the two inequalities yields
$$
\begin{aligned}
-2\epsilon \geq & \langle G(q_*, U(t_*, q_*)) - G(p_*, V(s_*, p_*)), q_* - p_* \rangle \\
& + \langle F(q_*, U(t_*, q_*)) - F(p_*, V(s_*, p_*)), U(t_*, q_*) - V(s_*, p_*) \rangle + O(\eta).
\end{aligned}
$$

Using the monotonicity of $(G,F)$, we obtain $\epsilon \leq 0$ (for small enough $\eta$ compared to $\eps$), a contradiction. Thus it simply remains to treat the case $t_* = 0$ or $s_* = 0$.

Let us consider the case $t_* = 0$ (the other is symmetric). From the definition of $(t_*, s_*, q_*, p_*)$ we get
\be\label{eq:star}
\langle U_0(q_*) - V(s_*, p_*), q_* - p_* \rangle + \alpha s_*^2 + \epsilon s_* + \langle a,q_* \rangle + \langle b,p_* \rangle + \xi_2 s_* \leq -\frac{\kappa}{2}.
\ee

Choosing $\eta$ sufficiently small implies in particular that $s_*$ is of order $(\sqrt{\alpha})^{-1}$. Since $V$ is continuous and satisfies the initial condition, we also have that
$$
\liminf_{s \to 0} \langle U_0(q_*) - V(s, p_*), q_* - p_* \rangle \geq 0.
$$

This yields a contradiction in \eqref{eq:star} for large enough $\alpha$ because of the monotonicity of $U_0$ since, along a subsequence, because of the compactness of $\bar \Omega$, $p^*\to \bar p$ for some $\bar p \in \bar \Omega$ and thus  $V(s_*,p_*) \to U_0(\bar p)$.\\

Hence we have shown $W \geq 0$. Now, for any $t \geq 0$, $q \in \Omega $, and $\epsilon > 0$ small enough so that $q \pm \epsilon v \in \Omega$ for all $|v| \leq 1$, we get
$$
\langle U(t,q) - V(t,q - \epsilon v), \epsilon v \rangle \geq 0.
$$
Dividing by $\epsilon$ and letting $\epsilon \to 0$, we obtain
$$
\langle U(t,q) - V(t,q), v \rangle \geq 0.
$$
Therefore, $U(t,q) = V(t,q)$, and the equality also holds on $[0,T]\times\bar \Omega$ by continuity. Setting $U = V$, we have also shown the monotonicity of $U$ in $q$.
\end{proof}

The previous uniqueness result justifies the notion of M-solution introduced in Definition \ref{def:monn}. Another good argument in favor of this notion is the stability that they enjoy, which is highlighted in the next result.

\begin{Prop}\label{prop:stabfinite}
Let $(G_n, F_n)_{n \geq 0}$ be a sequence of functions converging locally uniformly to $(G,F): \bar\Omega \times \mathbb{R}^d \mapsto \mathbb{R}^{2d}$. Let $(U_n)_{n \geq 0}$ be a sequence of M-solutions to \eqref{me:finite} (associated to $(G_n, F_n)$), such that $U_n \to_{n \to \infty} U$ locally uniformly. Then $U$ is an M-solution to \eqref{me:finite}.
\end{Prop}

\begin{proof}
Let $T > 0$, $V \in \mathbb{R}^d$, $p \in \bar \Omega$, $\vartheta$ a smooth real-valued function, and $(t_*, q_*) \in (0,T] \times \bar \Omega$ a point of strict minimum of $(t,q) \mapsto \langle U(t,q) - V, q - p \rangle - \vartheta(t)$. For all $n \geq 0$, once again by Stegall's Lemma, there exist sequences $(a_n)_{n \geq 0}$ and $(\xi_n)_{n \geq 0}$ with $|a_n| + |\xi_n| \leq 1/n$, such that the function
$$
(t,q) \mapsto \langle U_n(t,q) - V - a_n, q - p \rangle - \vartheta(t) - \xi_n t
$$
has a strict minimum on $[0,T] \times \bar\Omega$ at $(t_n, q_n)$. Since $(t_*, q_*)$ is a strict minimum, $(t_n, q_n) \to (t_*, q_*)$ as $n \to \infty$. Because $t_n > 0$ and $U_n$ is an M-solution, we have:
$$
\vartheta'(t_n) + \xi_n \geq \langle G_n(q_n, U_n(t_n, q_n)), q_n - p \rangle + \langle F_n(q_n, U_n(t_n, q_n)), U_n(t_n, q_n) - V - a_n \rangle.
$$
Passing to the limit yields the result.
\end{proof}

\begin{Rem}
Considering points of strict minimum is very useful in this proof. Indeed, it is not clear that general (non-strict) minima can be approximated by a sequence of minima, whereas this is guaranteed for strict minima. This is mainly due to the fact that we constrained the shape of the function we are minimizing. In particular, we do not allow terms such as $|q-q^*|^2$.
\end{Rem}

\subsection{Monotone solutions of master equations on spaces of measures}

In this section, I want to introduce an analogue notion of M-solution for MEs on spaces of probability measures such as 
\be\label{me1}
\ba
\partial_t U - &\left\langle \nabla_µ U(t,µ,x,\cdot),\text{div}(D_pH(\cdot,\nabla_x U(t,µ,\cdot))µ) + \sigma \Delta µ\right\rangle\\
&\quad - \sigma \Delta_x U + H(x,\nabla_x U) = f(µ)(x) \quad \text{in } (0,\infty)\times \motd\times \T^d,\\
U|_{t = 0} &= U_0.
\ea
\ee
where $\langle \cdot,\cdot \rangle$ denotes the extension of the $L^2(\T^d,\R)$ scalar product. The data of the problem are $H: \T^d\times \R^d \mapsto \R$, $f: \motd \mapsto \mathcal C(\T^d,\R)$, $U_0: \motd \mapsto \mathcal C(\T^d,\R)$, and $\sigma> 0$. Here, as in Section \ref{sec:flat}, $\nabla_µ$ denotes the vertical (or flat) derivative with respect to the measure argument $µ$. Even if this use of $\nabla_µ$ is only formal, because typical solutions are not necessarily V-differentiable, the point of view here is more aligned with the one of Section \ref{sec:W1}. Namely, $\motd$ will be equipped with the natural extension of the $\mathcal W_1$ distance and not with the distance in total variation. Hence, define for $µ,\nu \in \motd$,
\be\label{defW1mo}
\mathcal W_1(µ,\nu) = \sup \int_{\T^d}\phi(x)(µ-\nu)(dx),
\ee
where the supremum is taken over $\phi : \T^d \mapsto \R$ such that $Lip(\phi) \leq 1$ and $\|\phi\|_\infty\leq 1$. Note the notation $\mathcal W_1$ is slightly abused here, since the previous does not exactly coincide with the $1$-Wasserstein distance. Indeed, the diameter of $\T^d$ is greater than $1$ in general. This should not be a worry for the careful reader, and the not careful one can continue to ignore it.\\

 Note that compared to \eqref{me}, there are three main differences.
\begin{enumerate}
\item First, the state space is not $\R^d$ but $\T^d$. This implies that any set of the form $[0,T]\times\motd\times\T^d$ is compact. This will spare us the discussion on boundary conditions or behaviour of the solution at infinity. The latter are interesting questions, but they are not treated here to focus on the monotonicity and its implications.
\item I have here assumed that the dependences in $\nabla_x U$ and $µ$ are separated in $H$. This will ease the presentation of the associated monotone regime.
\item The ME is set on $\motd$ rather than on $\mptd$. This is mainly a choice of modelling: I want to insist upon the fact that the restriction to $\mptd$ is useless and that we can consider situations in which fewer players can be present. This will allow me to avoid some uninteresting technicalities.
\end{enumerate}

The notion of M-solution to \eqref{me1} shall be similar to the one for \eqref{me:finite}. The main idea is to adapt the definition from the previous section, drawing the following analogy:
\be
\ba
G^i(q,U) \leftarrow f(µ)(x) - H(x,\nabla_x U) + \sigma\Delta_x U,\\
F^i(q,U) \leftarrow -\text{div}(D_pH(x,\nabla_x U)µ) - \sigma \Delta µ,
\ea
\ee
where $i$ corresponds to $x$ and $q$ corresponds to $µ$. The main difference is that, whereas in the previous section we considered solutions $U(t,q)$ valued in a finite dimensional space, $U(t,µ)$ will now be valued in some infinite dimensional subspace of $(\T^d \mapsto \R)$. Hence, the regularity of $U(t,µ)$ with respect to the state variable $x$ is a fundamental question in the analysis of \eqref{me1}, notably because our new operators $F$ and $G$ are differential (in $x$) operators. To keep the present discussion quite simple, I will simply restrict my attention to solutions $U$ such that for all $t,µ$, $U(t,µ) \in \mathcal C^2(\T^d,\R)$, and will comment later on possible extensions in this direction.\\

As I did in the finite dimensional case, the starting point on M-solutions shall be to consider a proof of uniqueness / propagation of monotonicity. Such uniqueness is a priori only valid in the monotone setting, which will be here defined by
\begin{hyp}\label{hyp:mon}
The Hamiltonian $H$ is convex in its second argument. The mappings $f,U_0: \motd \mapsto \mathcal C (\T^d,\R)$ are monotone, that is for all $µ,\nu \in \motd$
$$
\langle f(µ) - f(\nu),µ-\nu\rangle \geq 0,
$$ 
and the same holds for $U_0$.
\end{hyp}
The following result holds.
\begin{Prop}\label{prop:monc}
Assume Hypothesis \ref{hyp:mon} holds. Then, there exists at most one smooth solution to \eqref{me1}.
\end{Prop}
Smooth refers here to the fact that the solution is V-differentiable everywhere, and that its V-derivative is a smooth function of all its arguments.
\begin{proof}
The proof is quite similar to the one of Proposition \ref{prop:propagationN}. Consider $U$ and $V$ two solutions, and define $W$ by
$$
W(t,µ,\nu) = \langle U(t,µ) - V(t,\nu),µ- \nu\rangle.
$$
As a consequence of the identity
$$
\nabla_µW(t,µ,\nu) = U(t,µ) - V(t,\nu) + \langle \nabla_µU(t,µ),µ-\nu\rangle,
$$
and the similar one for $\nabla_\nu$, it is easy to check that $W$ is a solution to 
$$
\ba
\partial_t W - &\left\langle \nabla_µ W(t,µ,\nu,\cdot),\text{div}(D_pH(\cdot,\nabla_x U(t,µ,\cdot))µ) + \sigma \Delta µ\right\rangle\\
- &\left\langle \nabla_{\nu} W(t,µ,\nu,\cdot),\text{div}(D_pH(\cdot,\nabla_x V(t,\nu,\cdot))\nu) + \sigma \Delta \nu\right\rangle\\
 =\quad& \langle f(µ) - f(\nu),µ - \nu\rangle + \langle \sigma \Delta U - H(\cdot,\nabla_x U(t,µ,\cdot)),µ- \nu\rangle\\
&\quad \quad +  \langle \sigma \Delta V - H(\cdot,\nabla_x V(t,\nu,\cdot)),\nu- µ\rangle\\
- &\left\langle U(t,µ,\cdot) - V(t,\nu,\cdot),\text{div}(D_pH(\cdot,\nabla_x U(t,µ,\cdot))µ) + \sigma \Delta µ\right\rangle\\
- &\left\langle V(t,\nu,\cdot) - U(t,µ,\cdot),\text{div}(D_pH(\cdot,\nabla_x V(t,\nu,\cdot))\nu) + \sigma \Delta \nu\right\rangle.
\ea
$$
Remark that the terms in $\sigma$ cancel out. Using the convexity of $H$, we obtain
$$
\ba
\langle \nabla_x(U(t,µ,\cdot) - V(t,\nu,\cdot))\cdot D_pH(\cdot,\nabla_x U(t,µ,\cdot)),µ\rangle\\
 + \langle H(\cdot,\nabla_x V(t,\nu,\cdot))- H(\cdot,\nabla_x U(t,µ,\cdot)),µ\rangle  \geq 0.
 \ea
$$ Using the monotonicity of $f$ we obtain 
\be\label{eq:ppmaxW}
\ba
\partial_t W - &\left\langle \nabla_µ W(t,µ,\nu,\cdot),\text{div}(D_pH(\cdot,\nabla_x U(t,µ,\cdot))µ) + \sigma \Delta µ\right\rangle\\
- &\left\langle \nabla_{\nu} W(t,µ,\nu,\cdot),\text{div}(D_pH(\cdot,\nabla_x V(t,\nu,\cdot))\nu) + \sigma \Delta \nu\right\rangle\\
 &\geq 0.
\ea
\ee
Since $W|_{t = 0} \geq 0$, we can conclude, as in the finite dimensional case, that $W \geq 0$, namely by invoking a slight generalization of the first order conditions of Proposition \ref{prop:vertopti} to $\motd$ and the invariance by the associated flows of $\motd$. This can also be obtained directly by considering a point of minimum of $W(t,\cdot,\cdot)$, and using the computations of Proposition \ref{prop:georges}. \\

To conclude, it now remains to show that $W \geq 0 \Rightarrow U = V$, since the non-negativity of $W$ already proved the propagation of monotonicity. Take $t \geq 0$ and $µ \in \motd$ such that $|µ|< 1$. For $\eps \in (0,1-|µ|)$, we can compute for any $\nu \in \mptd$

$$
\eps\langle U(t,µ + \eps \nu) - V(t,µ),\nu\rangle \geq 0.
$$
Hence, by dividing by $\eps$ and by using the continuity of $U$, we deduce that $U(t,µ,\cdot) \geq V(t,µ,\cdot)$ in $\T^d$. By symmetry we deduce that $U(t,µ) = V(t,µ)$ and by continuity that the equality holds for all $µ \in \motd$.

\end{proof}

Based on the same reasoning as previously, we are thus led to consider the following definition.

\begin{Def}\label{def:mon}
A continuous function $U : \mathbb{R}_+  \times \M_1(\mathbb{T}^d) \times \mathbb{T}^d\mapsto \mathbb{R}$ is an M-solution to \eqref{me1} if:
\begin{itemize}
\item For every $t > 0$ and $µ \in \motd$, the map $x \mapsto U(t,µ,x)$ belongs to $\mathcal{C}^2(\mathbb{T}^d, \mathbb{R})$.
\item The map $(t,µ,x) \mapsto \nabla_x U(t,µ,x)$ is bounded over $[0,T]\times \motd\times \T^d$ for any $T > 0$. 
\item For every $T > 0$, $\phi \in \mathcal{C}^2(\mathbb{T}^d, \mathbb{R})$, $\nu \in \mathcal{M}_1(\mathbb{T}^d)$, $\vartheta \in \mathcal{C}^1(\mathbb{R})$, and $(t_*, µ_*) \in (0,T] \times \motd$ a point of strict minimum of the map
\[
(t,µ) \mapsto \langle U(t,µ,\cdot) - \phi, µ - \nu \rangle - \vartheta(t)
\]
on $[0,T] \times \motd$, we have:
\be\label{eq505}
\ba
\vartheta'(t_*) \geq& \langle f(µ_*) - H(\cdot,\nabla_x U(t_*,µ_*)) + \sigma \Delta_x U(t_*,µ_*,\cdot), µ_* - \nu \rangle\\
&+ \langle U(t_*,µ_*,\cdot) - \phi, -\emph{div}(D_p H(\cdot,\nabla_x U(t_*,µ_*,\cdot))µ_*) - \sigma \Delta µ_* \rangle.
\ea
\ee 
\end{itemize}
\end{Def}

\begin{Rem}
This definition of a solution to the master equation only requires continuity with respect to the measure variable $µ$, the bound on $\nabla_x U$ is discussed below.
\end{Rem}

\subsection{Uniqueness and stability on spaces of measures}
As we could expect, the uniqueness of M-solutions to \eqref{me1} in the monotone regime holds, as the next result shows.
\begin{Theorem}
Assume Hypothesis \ref{hyp:mon}. Then there exists at most one M-solution to \eqref{me1} and if such a $U$ exists, it is such that $U(t,\cdot)$ is monotone for all times $t\geq 0$.
\end{Theorem}
\begin{proof}
Consider two M-solutions $U$ and $V$ and define the function $W$ by
$$
W(t,s,µ,\nu) = \langle U(t,µ) - V(s,\nu),µ-\nu\rangle.
$$
 I want to show that for all time $t \geq 0$, $W(t,t,\cdot,\cdot) \geq 0$. Arguing by contradiction, I assume that there exists $T > 0, \eps,\kappa, \bar\eta > 0$ such that 
 \be\label{eq:1231}
 \ba
 \sup_{\tiny{\ba\alpha > 0, \,\|\phi_1\|_{\mathcal C^2} +\|\phi_2\|_{\mathcal C^2}\\+|\xi_1| + |\xi_2| \leq \bar\eta\ea}}&\left\{\inf_{\tiny{\ba t,s \in [0,T],\\ µ,\nu \in \motd\ea}} W(t,s,µ,\nu) + \alpha (t-s)^2+\eps(t+s) + \xi_1t + \xi_2s + \langle \phi_1,µ\rangle + \langle \phi_2,\nu\rangle\right\}\\
 & \leq -\kappa.
 \ea
 \ee
Because for any set of parameters, the function minimized in \eqref{eq:1231} is continuous and the minimization is done over a compact set, the infimum in \eqref{eq:1231} is always reached at some point $(t_*,s_*,µ_*,\nu_*)$ (I omit the dependence on $\eps,\alpha$ etc...) From Proposition \ref{prop:stegallvert}, we know that for any $\eta \in (0,\bar \eta)$, there exists $\xi_1,\xi_2, \phi_1$ and $\phi_2$ such that $ \|\phi_1\|_{C^2} +\|\phi_2\|_{C^2}+|\xi_1| + |\xi_2| \leq \eta$ and $(t_*,s_*,µ_*,\nu_*)$ is a point of strict minimum.

Assume first that $t_*,s_* > 0$. Using that $U$ and $V$ are M-solutions to \eqref{me1}, we obtain
$$
 \ba
2(s_*-t_*)\alpha -\eps \geq& \langle f(µ_*) - H(\cdot,\nabla_x U(t_*,µ_*)) + \sigma \Delta_x U(t_*,µ_*,\cdot), µ_* - \nu_* \rangle + \xi_1\\
&+ \langle U(t_*,µ_*) - V(s_*,\nu_*), -\text{div}(D_p H(\cdot,\nabla_x U(t_*,µ_*,\cdot))µ_*) - \sigma \Delta µ_* \rangle\\
& +  \langle \phi_1, -\text{div}(D_p H(\cdot,\nabla_x U(t_*,µ_*,\cdot))µ_*) - \sigma \Delta µ_* \rangle,
\ea
$$
and
$$
 \ba
2(t_*-s_*)\alpha - \eps \geq& \langle f(\nu_*) - H(\cdot,\nabla_x V(s_*,\nu_*)) + \sigma \Delta_x V(s_*,\nu_*,\cdot), \nu_* - \mu_* \rangle + \xi_2\\
&+ \langle V(s_*,\nu_*) - U(t_*,\mu_*), -\text{div}(D_p H(\cdot,\nabla_x V(s_*,\nu_*,\cdot))\nu_*) - \sigma \Delta \nu_* \rangle\\
& +  \langle \phi_2, -\text{div}(D_p H(\cdot,\nabla_x V(s_*,\nu_*,\cdot))\nu_*) - \sigma \Delta \nu_* \rangle.
\ea
$$
Summing the two relations, using Hypothesis \ref{hyp:mon} and the uniform bound on $\nabla_x U,\nabla_x V$ leads to 
$$
-2 \eps + C \eta \geq 0. 
$$
Since $\eta >0$ can be arbitrarily small compared to $\eps > 0$, we arrive at a contradiction. Hence, either $s_*= 0$ or $t_* = 0$. Arguing as in Theorem \ref{thm:uniqmonf}, the limit $\alpha \to \infty$ leads to a contradiction as well since $W|_{t = s =0}$ is non-negative.\\

Hence, $W (t,t,\cdot,\cdot) \geq 0$ and the rest of the proof follows from the end of the proof of Proposition \ref{prop:monc}.
\end{proof}
\begin{Rem}
Note that the bound in $\nabla_x U$ is fundamental to estimate the remaining terms in $\phi_1,\phi_2$ in the proof.  Such perturbations are useful to prove stability results as the one below, as well as to study some more involved problems like optimal stopping. However, note that to simply establish uniqueness of solutions to \eqref{me1}, we could have define M-solutions by asking for the required relation to hold at any point of minimum, thus avoiding the use of the perturbations, and the uniform estimate on $\nabla_x U$ as well.
\end{Rem}

There is also a stability of M-solutions in this infinite dimensional setting. It can be formulated as follows.
\begin{Prop}\label{prop:stabc}
Let $(H_n)_{n \geq 0}$ be a sequence of functions in $\mathcal C^1(\T^d\times \R^d)$ which converges locally uniformly (in $ \mathcal C^1(\T^d\times \R^d)$) toward some function $H$ and a sequence $(f_n)_{n \geq 0}$ in $\mathcal C(\motd,\mathcal C(\T^d,\R))$ which converges locally uniformly toward some function $f$. Consider $(U_n)_{n \geq 0}$ a sequence of M-solutions to \eqref{me1} on $[0,T]\times \motd\times \T^d$ (associated with $H_n$ and $f_n$) and assume that $(U_n,\nabla_x U_n,\Delta_x U_n)_{n \geq 0}$ converges uniformly toward $(U,\nabla_x U,\Delta_x U)$ for some function $U: [0,T]\times \motd \mapsto \mathcal C^2(\T^d,\R)$. Then, $U$ is an M-solution to \eqref{me1} associated with $H$ and $f$.
\end{Prop}
The proof, which is a simple adaptation of Proposition \ref{prop:stabfinite} to this new setting, is left as an exercise.\\

\textbf{Comments on the regularity in the $x$ variable}\\
In all this section, the regularity of the function $U$ with respect to the spatial variable $x$ has been treated in a quite simplistic manner. I assume sufficient regularity (namely $\mathcal C ^2$ in $x$) so that all the techniques used in the finite dimensional case could be applied here. This is a strong requirement and it has as a consequence that existence of such solutions is quite hard to establish. A natural question is to understand whether or not these types of regularity assumptions can be avoided. 

A first answer is that when those uniformly elliptic terms are present in the equation (the $\sigma \Delta$), second order regularity in $x$ typically holds. This type of answer is not entirely satisfactory as it leaves open cases in which those second order terms are not uniformly elliptic. Another natural idea is to try to change slightly the definition of M-solution to restrict our attention to probability measures with more regularity, so that the derivatives of $U$ could be understood in a weaker sense. Such ideas have been investigated by Cardaliaguet and Souganidis \citep{sougacarda} in the case $\sigma = 0$, and partially by Cecchin and me \citep{bertucci2022mean} in the case $\sigma > 0$.

It is fair to say that the problem of lowering the regularity requirement in $x$ of the value function $U$ of MFG MEs is mostly open at the moment.

\subsection{Bibliographical comments}
The proof of Proposition \ref{prop:propagationN} is due to PL Lions and was presented in \citep{lions20112012}. I introduced the notion of monotone solutions (M-solutions) in \citep{bertucci2021monotone,bertucci2023monotone}, to treat respectively finite dimensional MEs and ones on sets of measures. This notion was also used by Cardaliaguet and Souganidis \citep{sougacarda} to tackle the case $\sigma = 0$, see also \citep{bertucci2022monotone} for cases on Hilbert spaces. In \citep{bertucci2022mean}, we used with Cecchin the notion of monotone solution to show convergence of the finite dimensional approximations of MEs in continuous state space, notably when the state space is discretized.
Other known results of uniqueness or stability of MEs relied on the regularity of the solution, and will be mentioned at the end of the next section.

\newpage
\section{Existence of solutions to master equations}\label{sec:meexist}
To study questions of existence of solutions to MEs, I shall take a quite similar route to the one of the uniqueness section above. Based on another proof of uniqueness, namely one involving regularity rather than monotonicity, I will introduce the notion of Lipschitz solution, and suitable properties that they enjoy.\\

In this section, I want to provide an existence theory for solutions to MEs \emph{before} the appearance of shocks. The notion of solution considered here shall be different from the one of M-solutions. This will lead to developments that could seem at first sight orthogonal to the ones of uniqueness of the previous section. I shall clarify the links between the two notions at the end of this section.\\

I will start by giving the heuristics behind Lipschitz solutions in the finite dimensional setting. Then, I will present the notion of Lipschitz solution in finite dimension, as well as the main results behind it. I will then study in less detail the analogue notion in infinite dimension. Comments on the links between Lipschitz solutions and how to prove that shocks do not appear in monotone regimes will be the last part of this section before the concluding bibliographical comments.

\subsection{Regularity implies uniqueness for master equations}
In this section, I elaborate on the comment I made after Proposition \ref{prop:propagationN}, that when the solution to the ME is smooth, we do not need the monotone regime to claim uniqueness. Let me insist upon the fact that these kinds of general statements are quite standard in the study of non-linear PDEs, and I simply present them because their proofs are instructive.

\begin{Prop}\label{prop:azert}
Assume that $(G,F): \bar\Omega\times \R^d \mapsto \R^{2d}$ is locally Lipschitz continuous. Let $U: [0,T]\times \bar\Omega \mapsto \R^d$ be a smooth solution to \eqref{me:finite}. Then, there exists no other smooth solution to \eqref{me:finite} with initial condition $U|_{t = 0}$.
\end{Prop}
\begin{proof}
Consider $V$ a smooth solution to \eqref{me:finite} with initial condition $U|_{t = 0}$. Define $W(t,q) = U(t,q) - V(t,q)$. Observe that $W$ solves
\be\label{eq:W5}
\partial_t W + F(q,V)\cdot \nabla_q W = G(q,U) - G(q,V) +(F(q,V) - F(q,U))\cdot \nabla_q U \text{ in } (0,T)\times \bar\Omega.
\ee
Take $t \geq 0$, and consider $i,q$ such that $W^i(t,q) = \max_{1 \leq j \leq d, p \in \bar\Omega}\{W^j(t,p)\}$, which always exist since $\bar \Omega$ is compact. Evaluating \eqref{eq:W5} at $i,q$ and recalling the assumption I made on the behaviour of $F$ at the boundary of $\Omega$, we obtain that
$$
\partial_t W^i(t,q) \leq C\|W(t,\cdot)\|_\infty (1+\|D_qU(t,\cdot)\|_\infty),
$$
where $C$ is a Lipschitz constant of $(G,F)$ on $\bar \Omega \times B(0,\|U\|_\infty + \|V\|_\infty)$. Since a similar computation could be made at points of minimum of $W$, we have in fact proven that
$$
\frac{d}{dt}\|W(t,\cdot)\|_\infty \leq C\|W(t,\cdot)\|_\infty (1 +\|D_qU(t,\cdot)\|_\infty),
$$
because such an estimate holds pointwise at any point of extremum of $W(t,\cdot)$. Since $U$ is smooth and $[0,T]\times \bar \Omega$ is compact, we know that there exists $K> 0$ such that $\sup_{t \leq T}\|D_qU(t,\cdot)\|_\infty \leq K-1$. Hence, using Gronwall's Lemma, we find that
\be\label{eq:1307}
\|W(t,\cdot)\|_\infty \leq \|W(0,\cdot)\|_\infty e^{CKt},
\ee
from which the result follows.
\end{proof}
The key insight from the previous result is the following. In the end, we are able to derive a quantitative propagation of error \eqref{eq:1307} which relies uniquely on the regularity of the data (through the constant $C$), which is fair to assume, and on the Lipschitz regularity of \emph{one of the two} solutions (the constant $K$). Hence, we can hope that Lipschitz regularity is a pivotal one and that it might be enough to define a concept of solution which yields existence easily.

\subsection{Lipschitz solutions in finite dimension}
In this section, I show how to define a concept of solution to the ME \eqref{me:finite} whose sole required regularity is Lipschitz continuity with respect to the variable $q\in \bar \Omega$, as it was the crucial quantity in the proof of Proposition \ref{prop:azert}. I will freely use the notation $\|D_qU\|_\infty$ to denote the Lipschitz constant of a map $U: \bar \Omega \mapsto \R^d$.

First, note that \eqref{me:finite} is of the form
\be\label{menlin}
\partial_t U(t,q) + (V(t,q)\cdot \nabla_q)U(t,q) = B(t,q) \text{ in } (0,\infty)\times \bar \Omega,
\ee
with $V$ and $B$ given by $V(t,q) = F(q,U(t,q))$ and $B(t,q) = G(q,U(t,q))$. Observe that for equations of the type \eqref{menlin}, representation formulas are available. Indeed, given (sufficiently regular) $V$ and $B$, the unique solution to \eqref{menlin} with initial condition $U_0$ is given by
\be\label{defUlip}
U(t,q) = \int_0^tB(t-s,q(s))ds + U_0(q(t)),
\ee
where $(q(s))_{s \in [0,t]}$ is the solution to
\be\label{charlip}
\frac{d}{ds}q(s) = -V(t-s,q(s)),
\ee
with initial condition $q(0) = q$.\footnote{There is a slight abuse of notation here since the trajectory $(q_s)_{s \in [0,t]}$ obviously depends on $t$.} Let us note that this characteristic formula is well-defined as soon as $V$ is Lipschitz in its second variable, uniformly with respect to its first. Indeed, in such a case, the usual Cauchy-Lipschitz theory applies nicely to \eqref{charlip}. Denote by $\Phi(V,B)$ the function given by \eqref{defUlip}-\eqref{charlip}. We shall see that this mapping $\Phi$ is well defined on Lipschitz functions, as it is claimed above. Furthermore, the following shows how to define a solution to \eqref{me:finite} based on this operator $\Phi$. 
\begin{Def}
Let $T > 0$. A function $U: [0,T]\times \bar\Omega \mapsto \R^d$ is a Lipschitz solution to \eqref{me:finite} on $[0,T]$ if
\begin{itemize}
\item There exists $C> 0$ such that for all $t\in [0,T], q \in \bar\Omega, \|D_qU(t,\cdot)\|_\infty\leq C$.
\item $U = \Phi(F(\cdot,U(\cdot,\cdot)), G(\cdot,U(\cdot,\cdot)))$.
\end{itemize}
\end{Def}
\begin{Rem}
In the previous definition and in what follows, I abuse slightly notation to denote the Lipschitz constant of $U(t,\cdot)$ with respect to $q$ by $\|D_qU(t,\cdot)\|_\infty$.
\end{Rem}

As expected, this notion of solution is consistent with the classical one.
\begin{Prop}
Let $U: [0,T]\times \bar\Omega \mapsto \R^d$ be a smooth solution to \eqref{me:finite}. Then $U$ is a Lipschitz solution to \eqref{me:finite}.
\end{Prop}
The proof is trivial and left as an exercise to the reader. The main result on Lipschitz solution is the following.
\begin{Prop}\label{prop:lipn}
Assume that $U_0$ is Lipschitz, that $F$ and $G$ are locally Lipschitz. Then there exists $T_c \in (0,\infty]$ and $U:[0,T_c)\times \bar \Omega \mapsto \R^d$ such that: for any $T < T_c$, $U$ is a Lipschitz solution to \eqref{me:finite} on $[0,T]$, for any $V$ Lipschitz solution to \eqref{me:finite} on $[0,T']$, we have $T' < T_c$ and $V = U|_{t \leq T'}$. Moreover, either $T_c = \infty$, or $\limsup_{t \to T_c}\|D_qU(t,\cdot)\|_\infty = \infty$.
\end{Prop}
\begin{proof}
Let $T > 0$. For $C> 0$ define the set
$$
E_C := \{ U : [0,T]\times \bar\Omega \mapsto \R^d,\|U\|_{\infty} \leq C, \|D_qU\|_{\infty} \leq C\}.
$$
Let $\Psi(U):= \Phi(F(\cdot,U(\cdot,\cdot)), G(\cdot,U(\cdot,\cdot)))$. We first show that if $C$ is large enough and $T$ small enough, then $\Psi(E_C) \subset E_C$. Let $U \in E_C$, then 
$$
\|\Psi(U)\|_\infty \leq TK(1 + C) + K,
$$
where $K$ denotes a constant which bounds $\|D_qG\|_\infty$ and $\|U_0\|_\infty$. Hence, we see that if $C > K$, then for $T$ sufficiently small, $\|\Psi(U)\|_\infty \leq C$. From now on, we can argue as if $F$ and $G$ are globally Lipschitz continuous. For $t \in [0,T]$ and $q,\tilde q\in \bar\Omega$, compute
$$
\ba
|\Psi(U)(t,q) - \Psi(U)(t,\tilde q)| &= \bigg| \int_0^tG(q(s), U(t-s,q(s))) - G(\tilde q(s),U(t-s,\tilde q(s)))ds\\
&\quad \quad \quad \quad \quad + U_0(q(t))- U_0(\tilde q(t))  \bigg|\\
& \leq \bigg(t(\|D_qG\|_{\infty} + \|D_uG\|_{\infty}\|D_qU\|_{\infty}) + |D_qU_0|_{\infty}\bigg) \sup_{0\leq s \leq t}|q(s) - \tilde q(s)|,
\ea
$$
where $(q(s))_{s \in [0,T]}$ and $(\tilde q(s))_{s \in [0,T]}$ are the solutions to
$$
\frac{d}{ds}X(s) = - F(X(s),U(t-s,X(s)))
$$
with respective initial conditions $q$ and $\tilde q$. The chain rule yields
$$
\frac{d}{ds}|q(s) - \tilde q(s)| \leq \|D_qF\|_{\infty}|q(s)-\tilde q(s)| + \|D_uF\|_{\infty}\|D_qU\|_{\infty}|q(s)- \tilde q(s)|.
$$
By Gronwall's Lemma, we get
$$
\sup_{0\leq s \leq t}|q(s) - \tilde q(s)| \leq e^{t(\|D_qF\|_{\infty} +\|D_uF\|_{\infty}\|D_qU\|_{\infty})}|q-\tilde q|.
$$
Thus 
$$
\|D_q\Psi(U)\|_\infty \leq (TK_1(1+C) + K_1)e^{TK_1(1+C)},
$$
where $K_1$ is a constant greater than the Lipschitz constants of $G,F$ and $U_0$. It then follows easily that $\Psi(E_C)\subset E_C$ for $T$ small enough, provided $C$ is larger than $K_1$.\\

Now we want to show that $\Psi$ is a contraction in the supremum norm. Let $U,V \in E_C$. We compute for $t\in [0,T], q \in \Omega$
$$
\ba
|\Psi(U)(t,q) - \Psi(V)(t,q)|&= \bigg| \int_0^tG(q(s), U(t-s,q(s))) - G(\tilde q(s),V(t-s,\tilde q(s)))ds\\
&\quad \quad \quad \quad + U_0(q(t))- U_0(\tilde q(t))  \bigg|\\
& \leq \bigg(t(\|D_qG\|_{\infty} + \|D_uG\|_{\infty}\|D_qU\|_{\infty}) + \|D_qU_0\|_{\infty}\bigg) \sup_{0\leq s \leq t}|q(s) - \tilde q(s)|\\
& \quad+ t\|D_uG\|_{\infty}\|U-V\|_{\infty},
\ea
$$
where $(q(s))_{s \in [0,T]}$ solves
$$
\frac{d}{ds}X(s) = - F(X(s),U(t-s,X(s)))
$$
with initial condition $q$ and $(\tilde q(s))_{s \in [0,T]}$ solves the same equation with $U$ replaced by $V$, also with initial condition $q$. A similar computation as in the first part of the proof yields
$$
\frac{d}{ds}|q(s) - \tilde q(s)| \leq ( \|D_q F\|_{\infty} + \|D_uF\|_{\infty}C)|q(s)-\tilde q(s)| + \|D_uF\|_{\infty}\|U-V\|_{\infty}.
$$
By Gronwall's Lemma again,
$$
|q(s)-\tilde q(s)| \leq (e^{(\|D_q F\|_{\infty} + \|D_uF\|_{\infty}C)s} -1)\frac{\|D_uF\|_{\infty}\|U-V\|_{\infty}}{\|D_q F\|_{\infty} + \|D_uF\|_{\infty}C}.
$$
It follows that, possibly choosing $T$ smaller, $\Psi$ is a contraction in the supremum norm. By a classical Picard iteration argument, we can then show that for all $U\in E_C$, the sequence $(\Psi^k(U))_{k \geq 0}$ converges, provided $C$ is large enough and $T$ small enough. This limit is unique since $\Psi$ is a contraction.\\

It remains to show the existence of a maximal time of existence. It is quite straightforward as it simply follows from the remark that the constant $C$ and $T$ in the previous depend only on $\|U_0\|_\infty + \|D_q U_0\|_\infty$ (in addition to the dependence on $F$, $G$ and $\Omega$ of course...). Hence, we can repeat the previous argument by replacing $U_0$ by $U(T,\cdot)$, and then actually do it an infinite number of times to create an increasing sequence of times $(T_n)_{n \geq 0}$, on which a unique Lipschitz solution exists. Then, only one of two things happens. Either the sequence of times $T_n$ that we construct in such a way is unbounded, or it is bounded. In the second case, it must be that $\|U(T_n,\cdot)\|_\infty + \|D_q U(T_n,\cdot)\|_\infty\to \infty$ as $n \to \infty$, and the divergence has to come from the derivative, since otherwise we can easily show that $U$ has to remain bounded.
\end{proof}

This result shows that Lipschitz regularity in $q$ implies uniqueness of solutions to \eqref{me:finite}. This uniqueness has to be taken with care since it says nothing about uniqueness in weaker classes of notions of solutions. Indeed, a priori, the existence of a Lipschitz solution $U$ prevents the existence of another Lipschitz solution on the same time interval. For instance a weaker solution $V$ (e.g. a monotone solution) to \eqref{me:finite} could co-exist with $U$ and still $U \ne V$. The next result argues in the opposite direction. It is a form of weak-strong uniqueness result in the sense that the existence of a Lipschitz solution prevents the existence of certain weak solutions, namely ones that can be approximated uniformly by smooth functions.

\begin{Prop}
Consider a Lipschitz solution $U$ to \eqref{me:finite} on the time interval $[0,T_c)$. Consider for $T \in (0,T_c)$, $\epsilon > 0$, a classical solution $V_{\epsilon}$ of
\be
\begin{aligned}
\left|\partial_t V_{\epsilon} + \langle F(q,V_{\epsilon}),\nabla_q\rangle V_{\epsilon} - G(q,V_{\epsilon})\right| \leq \epsilon \text{ in } (0,T)\times \bar\Omega,\\
|V_{\epsilon}(0,q) - U_0(q)| \leq \epsilon \text{ in } \bar \Omega.
\end{aligned}
\ee
Then the following holds for some constant $C$ depending only on $T, U, F$ and $G$
\be\label{eq:ffx}
\sup_{t \leq T, q \in \bar \Omega}|U(t,q) - V_{\epsilon}(t,q)| \leq C \epsilon.
\ee
\end{Prop}
\begin{proof}
Let us consider $t > 0$ and $i \in \{1,...,d\}, q^* \in \bar \Omega$ such that 
$$
U^i(t,q^*) - V_\eps^i(t,q^*) = \max_{j,q}\{U^j(t,q) - V_\eps^j(t,q)\}
$$
 Consider $\kappa > 0$ and the solution $q(\cdot)$ of the ODE
\be\label{ode12}
\frac{d}{ds}q(s) =- F(q(s),U(t+ \kappa -s,q(s))),
\ee
with initial condition $q^*$ and $\tilde q(\cdot)$ the solution to the ODE
\be\label{ode22}
\frac{d}{ds}\tilde q(s) = -F(\tilde q(s),V_\eps(t+ \kappa -s,\tilde q(s))),
\ee
also with initial condition $q^*$. Using the definition of Lipschitz solutions, the time derivative of $U^i(t,q^*) - V^i(t,q^*)$ can be computed along the ODEs \eqref{ode12} and \eqref{ode22} as follows.
\be\label{eq2er}
\begin{aligned}
U^i&(t+ \kappa, q^*) - V_\eps^i(t+ \kappa, q^*) - U^i(t,q^*) + V_\eps^i(t,q^*)=\\
&=  U^i(t+ \kappa, q^*) - U^i(t,q(\kappa))  - V_\eps^i(t+ \kappa, q^*) + V_\eps^i(t,\tilde q(\kappa))\\
& \quad + U^i(t,q(\kappa))- U^i(t,q^*) + V_\eps^i(t,q^*) - V_\eps^i(t,\tilde q(\kappa))\\
& \leq \int_0^{\kappa} G^i(q(s),U(t+\kappa -s,q(s))) - G^i(\tilde q(s),V_\eps(t + \kappa -s, \tilde q(s))) + \epsilon ds\\
&+ U^i(t,q(\kappa)) - U^i(t,\tilde q(\kappa)).
\end{aligned}
\ee
Since $U$ is uniformly Lipschitz, we deduce that
$$
|U^i(t,q(\kappa)) - U^i(t,\tilde q(\kappa))| \leq C |q(\kappa) - \tilde q(\kappa)|.
$$
From \eqref{ode12} and \eqref{ode22} it follows that
$$
|U^i(t,q(\kappa)) - U^i(t,\tilde q(\kappa))| \leq C \kappa \|U-V_\eps\|_{\infty}.
$$
Hence, dividing by $\kappa$ and taking the limit $\kappa \to 0$ in \eqref{eq2er}, we obtain
$$
\frac{d}{dt}(U^i(t, q^*) - V_\eps^i(t, q^*)) \leq C \|U(t,\cdot) - V_\eps(t,\cdot) \|_{\infty} + \epsilon.
$$
Thus,
$$
\frac{d}{dt}\|U(t,\cdot) - V_\eps(t,\cdot)\|_{\infty} \leq C  \|U(t,\cdot) - V_\eps(t,\cdot) \|_{\infty} + \epsilon.
$$
Hence, from Gronwall's Lemma
\be\label{lasteqprop}
\|U(t,\cdot) - V_\eps(t,\cdot)\|_{\infty} + \epsilon \leq 2\epsilon e^{Ct},
\ee
from which the proposition immediately follows.
\end{proof}
\begin{Rem}
From the previous result we can indeed deduce that the Lipschitz solution attracts the limits of approximations of \eqref{me:finite}, since the constant $C$ in \eqref{lasteqprop} or \eqref{eq:ffx} does not depend on $V_{\epsilon}$.
\end{Rem}

\subsection{Lipschitz solutions in infinite dimension}
In this section, I present how to extend the notion of Lipschitz solution to the infinite dimensional equation \eqref{me1}. I will go over this generalization relatively briefly as the main ideas are the same ones as in the finite dimensional case. The main difficulty in generalizing the previous case lies, as always, in the fact that \eqref{me1} involves derivatives of $U$ with respect to $x$, or in other words that the infinite dimensional analogues of $F$ and $G$ are differential operators.\\

In this case, we will be interested in the gradient of the value function defined by $W := \nabla_x U$. Formally, it is a solution to
\be\label{lipW}
\ba
&\partial_t W - \langle \nabla_µ W(t,µ,x,\cdot), \text{div}(D_pH(\cdot,W(t,µ,\cdot)) µ) + \sigma \Delta µ\rangle-\sigma\Delta_x W \\
& + D_xH(x,W) + D_xW \cdot D_pH(x,W) = \nabla_x f(µ)(x) \quad \text{in } (0,\infty) \times \motd\times \T^d.
\ea
\ee
This equation only involves linear terms in the derivatives of $W$ with respect to $x$. One can therefore give an analogue definition of a Lipschitz solution in this case, once again based on a fixed point procedure on the underlying linearized transport equation. Of course, we need here to specify the distance with which $\motd$ will be equipped, so that we can talk about Lipschitz solutions. In this section on Lipschitz solutions, all statements shall be made on the space $(\motd,\mathcal W_1)$, see \eqref{defW1mo}.
\begin{Def}
Let $T > 0$. A function $W: [0,T]\times \motd \times \T^d\mapsto \R^d$ is a Lipschitz solution to \eqref{lipW} with initial condition $\nabla_x U_0$ if 
\begin{itemize}
\item $W$ is Lipschitz in $(µ,x)$, uniformly in $t\in [0,T]$.
\item For all $(t,µ,x)\in [0,T]\times \motd\times \T^d$, 
$$
\ba
W(t,µ,x) = \mathbb{E}\bigg[\int_0^t&\nabla_x f(m(s))(X_s) - D_x H(X_s,W(t-s,m(s),X_s))ds\\
& + \nabla_x U_0(m(t),X_t)\bigg],
\ea
$$
where $(X_s,m(s))_{s \in [0,t]}$ is the strong solution to
$$
\begin{cases}
dX_s = -D_pH(X_s,W(t-s,m(s),X_s))ds + \sqrt{2\sigma}dB_s \quad \text{for } s \in (0,t)\\
\partial_s m -\sigma \Delta m -\text{div}(D_pH(x,W(t-s,m(s),x))m(s)) = 0 \quad \text{in } (0,t) \times \T^d,
\end{cases}
$$
with initial condition $(µ,x)$, and where $(B_s)_{s \geq 0}$ is a $d$ dimensional Brownian motion.
\end{itemize}
\end{Def}
\begin{Rem}
Note that because of the nature of the equation, we have to take into account the evolution of the individual state variable $x$.
\end{Rem}

The following result is the exact analogue of Proposition \ref{prop:lipn}.
\begin{Prop}
Assume that $\nabla_x U_0$ and  $\nabla_x f$ are Lipschitz in $(µ,x)$ and that $D_x H$ and $D_p H$ are Lipschitz continuous in $(x,p)$. Then there exists $T_c \in (0,\infty]$ and $W: [0,T_c)\times \motd \times \T^d\mapsto \R^d$ such that: for any $T < T_c$, $W$ is a Lipschitz solution to \eqref{lipW} on $[0,T]$, for any $T$ and Lipschitz solution $\tilde W$ to \eqref{lipW} on $[0,T]$, we have that $T < T_c$ and $\tilde W = W$ in $[0,T]$. Moreover, either $T_c = +\infty$, or the $\limsup$ of the Lipschitz constant of $W(t,\cdot,\cdot)$ goes to $+ \infty$ as $t \to T_c$.
\end{Prop}
\begin{proof}
The proof is similar to the one of Proposition \ref{prop:lipn}. In fact, the argument follows exactly if we use Lemma \ref{lemma:lip} below. I leave the details of this adaptation and simply state and prove Lemma \ref{lemma:lip}.
\end{proof}
\begin{Lemma}\label{lemma:lip}
Consider $T > 0$ and $f^\alpha,f^\beta :[0,T]\times \motd \times \T^d \mapsto \R^d$ such that
$$
K_j:=\sup_{t \in [0,T]}\sup_{x\ne y, µ \ne \nu}\frac{|f^j(t,µ,x)-f^j(t,\nu,y)|}{|x-y| +\mathcal W_1(µ,\nu)} < \infty,
$$
for $j = \alpha,\beta$. We have the two following estimates.
\begin{enumerate}
\item There exists $C > 0$ depending on $K_\alpha$ such that for any $µ^1,µ^2 \in \motd, x^1,x^2 \in \T^d$, if we denote by $((X^i_t,m^i(t))_{t \in [0,T]})_{i = 1,2}$ the solution to
$$
\begin{cases}
dX^i_t = f^\alpha(t,m^i(t),X^i_t)dt + \sqrt{2\sigma}dB_t,\\
\partial_t m^i -\sigma \Delta m^i +\emph{div}(f^\alpha(t,m^i(t))m^i(t)) = 0 \quad \text{in } (0,T) \times \T^d,\\
m^i(0) = µ^i, \quad X^i_0 = x^i,
\end{cases}
$$
it holds that almost surely
$$
\sup_{t \in [0,T]} |X^1_t-X^2_t| + \mathcal W_1(m^1(t),m^2(t)) \leq C(|x^1-x^2| + \mathcal W_1(µ^1,µ^2)).
$$
\item There exists $C$ such that, for any $µ\in \motd,x \in \T^d$, and $((X^j_t,m^j(t))_{t \in [0,T]})_{j = \alpha,\beta}$ solution to 
$$
\begin{cases}
dX^j_t = f^j(t,m^j(t),X^j_t)dt + \sqrt{2\sigma}dB_t,\\
\partial_t m^j -\sigma \Delta m^j +\emph{div}(f^j(t,m^j(t))m^j(t)) = 0 \quad \text{in } (0,T) \times \T^d,\\
X^j_0 = x, \quad m^j(0) = µ,
\end{cases}
$$
the following holds
$$\sup_{t \in [0,T]} |X^\alpha_t-X^\beta_t| + \mathcal W_1(m^\alpha(t),m^\beta(t)) \leq C\|f^\alpha - f^\beta\|_\infty.$$
\end{enumerate}
\end{Lemma}
\begin{proof}
I do not justify the existence of solutions to the systems as i) it is classical, ii) it can easily be proven with the type of estimate we are about to deal with.\\

Let me start with the first estimate. Take a couple $((X^i_t,m^i(t))_{t \in [0,T]})_{i = 1,2}$ as in the statement of the lemma. From the assumptions we made, we obtain easily that $(t,x) \mapsto f^\alpha(x,m^i(t))$ is uniformly in space $\frac12$ H\"older continuous in time, and uniformly in time Lipschitz continuous in space. Take $t \in (0,T]$ and $g: \T^d \mapsto \R$ a Lipschitz function. Consider the solution $\varphi$ to 
\be\label{def:varphi}
\ba
-\partial_s \varphi - \sigma \Delta \varphi - f^\alpha(x,m^1(s))\cdot \nabla_x \varphi = 0 \text{ in } (0,t)\times \T^d,\\
\varphi|_{s = t} = g.
\ea
\ee
From the regularity of $f^\alpha$ and $g$, we also deduce that $\varphi \in \mathcal C^{1,2}([0,t]\times \T^d,\R)$ and that $\nabla_x \varphi$ is bounded by a constant $C$ which depends only on $\|\nabla_x g\|_\infty$, $K^\alpha$, $T$ and $d$. Hence, we can use $\varphi$ as a test function in the weak formulation of the PDE satisfied by $m^1-m^2$, between $0$ and $t$, to obtain
$$
\ba
\int_{\T^d}g(x)(m^1(t)-m^2(t))(dx) = &\int_0^t\int_{\T^d}(f^\alpha(m^1(s))(x)-f^\alpha(m^2(s))(x))\cdot \nabla_x \varphi(s,x)m^2(s)(dx)\\
& + \int_{\T^d} \varphi(0,x)(µ^1-µ^2)(dx).
\ea
$$ 
We thus obtain, since $g$ was arbitrary, that
$$
\mathcal W_1(m^1(t),m^2(t)) \leq C\left( \int_0^t \mathcal W_1(m^1(s),m^2(s))ds +\mathcal W_1(µ^1,µ^2)\right).
$$
We then conclude by Gronwall's Lemma to obtain the estimate on $\mathcal W_1(m^1(t),m^2(t))$. The estimate on $|X^1- X^2|$ then easily follows from Gronwall's Lemma as well.\\

We now turn to the second estimate. The strategy of proof is similar. Consider $g$, a real Lipschitz function over $\T^d$ and $\varphi$ defined by \eqref{def:varphi}, where $m^1$ is replaced by $m^\alpha$. It then follows that

$$
\ba
\left|\int_{\T^d}g(x)(m^\alpha(t)-m^\beta(t))(dx) \right|&= \bigg|\int_0^t\int_{\T^d}(f^\beta(m^\beta(s))(x)-f^\alpha(m^\alpha(s))(x))\cdot \nabla_x \varphi(s,x)m^\beta(s)(dx)\bigg|\\
& \leq \|\nabla_x \phi\|_\infty\int_0^t \|f^\beta(m^\beta(s))-f^\alpha(m^\alpha(s))\|_\infty ds\\
& \leq Ct \|f^\beta - f^\alpha\|_\infty + CK_\alpha\int_0^t\mathcal W_1(m^\alpha(s),m^\beta(s))ds.
\ea
$$
Taking the supremum over such $g$, the left hand side becomes $\mathcal W_1(m^\alpha(t),m^\beta(t))$. Hence, we can conclude once again thanks to Gronwall's Lemma, exactly as we did at the end of the proof of Proposition \ref{prop:lipn}.
\end{proof}

Having defined a concept of Lipschitz solution to \eqref{lipW}, we can now reconstruct a solution to \eqref{me1} as the following definition shows.

\begin{Def}
Let $T> 0$ and  $U:[0,T]\times \motd\times \T^d\mapsto \R$. We say that $U$ is a Lipschitz solution to \eqref{me1} on $[0,T]$ if
\begin{itemize}
\item $U$ is differentiable in $x$. We then denote $W = \nabla_x U$.
\item $W$ is a Lipschitz solution to \eqref{lipW} on $[0,T]$.
\item For all $(t,µ,x) \in [0,T]\times \motd\times \T^d$
$$
\ba
U(t,µ,x) = \mathbb{E}\bigg[\int_0^t &f(m(s))(X_s) + L(X_s,-D_pH(X_s,W(t-s,m(s),X_s)))\,ds\\
& + U_0(m(t),X_t)\bigg],
\ea
$$
where $(X_s,m(s))_{s \in [0,T]}$ is the strong solution to
$$
\begin{cases}
dX_s = -D_pH(X_s,W(t-s,m(s),X_s))\,ds + \sqrt{2\sigma}\,dB_s \quad \text{for } s \in (0,t)\\
\partial_s m -\sigma \Delta m -\emph{div}(D_pH(x,W(t-s,m(s),x))m(s)) = 0 \quad \text{in } (0,t) \times \T^d,
\end{cases}
$$
with initial condition $(x,µ)$, and where $(B_s)_{s \geq 0}$ is a $d$-dimensional Brownian motion.
\end{itemize} 
\end{Def}
Similar results of well-posedness at the level of $U$ can be established but I do not focus on these types of questions here.

\subsection{Filling in the blanks between existence and uniqueness}
Up to now, we saw two different notions of solutions to the ME. One is based on an a priori short time horizon result of well-posedness, and another one based on the fact that, in the monotone regime, we can hope for global well-posedness (through the propagation of monotonicity and the uniqueness of M-solutions). In this section, I present links between these concepts, namely that Lipschitz solutions are M-solutions, and that propagation of monotonicity can lead to propagation of regularity under reasonable assumptions.\\

Once again, I will mainly focus on the case of the finite dimensional ME \eqref{me:finite}. There are, of course, several important and meaningful difficulties when trying to generalize the following in infinite dimension, but they exceed the scope of this book and references are given below.

\subsubsection{From Lipschitz solutions to M-solutions}
For master equations in finite dimension, Lipschitz solutions are also M-solutions as the next result shows.
\begin{Prop}
Assume that the assumptions of Proposition \ref{prop:lipn} hold. Let $U:[0,T]\times \bar\Omega \mapsto \R^d$ be a Lipschitz solution to \eqref{me:finite} on $[0,T]$, then $U$ is an M-solution to \eqref{me:finite} on $[0,T]$.
\end{Prop}
\begin{proof}
First remark that $U$ is continuous in time, thus continuous. It follows from the definition of a Lipschitz solution as well as from the regularity of $F$ and $G$ (namely local Lipschitz regularity). It is an easy fact to verify and the interested reader can estimate the difference between $U(t,q)$ and $U(s,p)$ using the definition of $U$ and the stability of the ODE quite simply.\\

Let $V \in \R^d, \tilde{q} \in \bar\Omega,$ a smooth real function $\vartheta$ and $(t_0,q_0)\in (0,T]\times\bar\Omega $ a point of strict minimum of the function $(t,q) \mapsto \langle U(t,q)- V,q- \tilde{q}\rangle - \vartheta(t)$ on $[0,T]\times \bar\Omega$.
For all $t \in [0,T], q \in \Omega$,
\be\label{eq64}
\langle U(t,q)- V,q- \tilde{q}\rangle - \vartheta(t) \geq \langle U(t_0,q_0)- V,q_0- \tilde{q}\rangle - \vartheta(t_0).
\ee
Introduce $(q(s))_{s\in[0,t_0]}$ the solution to
$$
\frac{d}{ds}q(s) = -F(q(s),U(t_0-s,q(s))),
$$
with initial condition $q_0$. Evaluating \eqref{eq64} at $t_0 - \epsilon, q(\epsilon)$ for $\epsilon > 0$, we obtain
$$
\langle U(t_0-\epsilon,q(\epsilon))- V,q(\epsilon)- \tilde{q}\rangle - \vartheta(t_0-\epsilon) \geq \langle U(t_0,q_0)- V,q_0- \tilde{q}\rangle - \vartheta(t_0).
$$
Rewriting this inequality yields
$$
\langle U(t_0-\epsilon,q(\epsilon)) -U(t_0,q_0),q(\epsilon)- \tilde{q}\rangle + \vartheta(t_0) - \vartheta(t_0-\epsilon) \geq \langle U(t_0,q_0)- V,q_0-q(\epsilon)\rangle.
$$
Since $U$ is a Lipschitz solution to \eqref{me:finite}, we deduce that
$$
\ba
\langle -\int_{t_0-\epsilon}^{t_0}&G(q(t_0-s),U(s,q(t_0-s)))ds,q(\epsilon)- \tilde{q}\rangle + \vartheta(t_0) - \vartheta(t_0-\epsilon)\\
& \geq \langle U(t_0,q_0)- V,q_0-q(\epsilon)\rangle.
\ea
$$
Dividing by $\epsilon$ and taking the limit $\epsilon \to 0$, we recover the relation that M-solutions must satisfy, thanks to the definition of $(q(s))_{s \in[0,t_0]}$, hence the claim follows.
\end{proof}

\begin{Rem}
If Lipschitz solutions are M-solutions, then one could wonder about the use of the weaker concept of M-solution. There are two main arguments in favor of the interest of M-solutions. The first one is that, using the stability explained in Propositions \ref{prop:stabfinite} and \ref{prop:stabc}, we can prove existence of M-solutions which are not Lipschitz continuous, as for instance in \citep{bertucci2023monotone}. Such solutions are particularly helpful when we consider MEs (say \eqref{me:finite}) with data (say $F$ and $G$) which are not Lipschitz continuous. The other type of case of interest for M-solutions is when the ME is so intricate that we are not able to write it properly, as is the case in MFG of optimal stopping, see \citep{bertucci2021monotone} for instance. In such situations, the underlying evolution of the forward quantity can be discontinuous, and thus the notion of Lipschitz solution quite meaningless, while the value function stays continuous and the notion of M-solution makes perfect sense.
\end{Rem}

\subsubsection{Existence of solutions on arbitrary long time intervals I: a priori estimates}
To justify that, in monotone regimes, we can prove existence of solutions on arbitrary long time intervals, a natural step is to try to prove a Lipschitz estimate on the solution to \eqref{me:finite}. If such an estimate can be justified, then from Proposition \ref{prop:lipn}, we know that the critical time of existence must be $+\infty$. In this section, we focus on proving formally an a priori estimate, and in Section \ref{sec:lastme} on how we can justify it.\\

In the following, I shall be quite formal, and simply assume that all the functions involved are sufficiently smooth to justify all the computations. What I want to establish is a quantitative estimate on the regularity of $U$, the solution to \eqref{me:finite}, based on the regularity of $F,G$ and $U_0$.\\

 The main idea of the following a priori estimate is that the propagation of monotonicity of Proposition \ref{prop:propagationN} already gave half of an estimate. Indeed, because $U$ is monotone, for any $q,p \in \bar\Omega, t \geq 0$,
 $$
 \langle U(t,p)-U(t,q),p-q\rangle \geq 0.
 $$
 Taking, if possible, $p = q + \eps \xi$ for $\xi \in \R^d$ and $\eps > 0$, dividing by $\eps$ and taking the limit $\eps \to 0$ yields the inequality
 $$
 \xi \cdot D_q U(t,q)\cdot \xi \geq 0.
 $$
 Hence, if we are able to show an estimate of the form
 $$
\forall T > 0, \exists A, \forall t \leq T,q \in \bar\Omega, \xi \in \R^d,\quad  \xi\cdot D_q U(t,q) \cdot\xi \leq A|\xi|^2,
 $$
 then we would have indeed proven an $L^\infty$ estimate on (the symmetric part) of $D_q U$. Because of the monotonicity, bounding from above quantities such as $\xi\cdot D_q U(t,q) \cdot\xi $ is harder than bounding them from below. Hence, I shall investigate in the following bounds from below on quantities of the form
 $$
  \xi\cdot D_q U(t,q) \cdot\xi - A|D_q U(t,q)\xi|^2,
 $$
 which imply uniform estimates on $D_q U$. Objects of the nature of $D_qU$ shall play a strong role in the proof. To lighten notation, I fix a convention that for any map $\phi: \bar\Omega \mapsto \R^d$, $x,y \in \R^d$
 $$
 x\cdot D_q \phi \cdot y := \sum_{i,j=1}^dx_i\partial_{q_i}\phi^jy_j,
 $$
 with the same convention if only $x$ or $y$ is used.
 \begin{Prop}\label{prop:estn}
 Assume that $U_0$, $F$ and $G$ are smooth Lipschitz functions, that $U$ is a Lipschitz solution associated to them, and that there exists $\alpha > 0$ such that
 \be\label{amon}
 \begin{pmatrix} &D_qG &D_q F\\&D_uG & D_u F \end{pmatrix} \geq \alpha Id, \quad D_q U_0 \geq \alpha Id.
 \ee
 Then $Z$ defined by $Z(t,q,\xi) = \xi\cdot D_q U(t,q)\cdot \xi - \beta |D_qU(t,q)\cdot \xi|^2$ is non-negative for some $\beta  >0$ which only depends on $\alpha$, $\|D_qU_0\|_\infty$, $\|D_qG\|_\infty$ and $\|D_qF\|_\infty$.
 \end{Prop}
 \begin{proof}
 Start with the case $\beta = 0$. Taking the scalar product of \eqref{me:finite} with $\xi$, differentiating this new relation with respect to $q$ and taking once again the scalar product with $\xi$ leads to
 $$
 \ba
  \partial_t Z& + \langle F(q,U), \nabla_q Z\rangle + (\xi\cdot D_q U)D_u F(q,U) (D_q U \cdot \xi)\\
  &+ \xi D_q F(q,U) (D_q U \cdot \xi) = \xi D_q G(q,U) \xi +  (\xi\cdot D_q U)D_u G(q,U) \xi.
 \ea
 $$
 Observe that, in this case $\beta = 0$,
\be\label{Zxi}
 \nabla_\xi Z(t,q,\xi) = \xi \cdot D_qU(t,q) + D_qU(t,q)\cdot \xi.
 \ee
Using this relation in the previous PDE leads to
$$
 \ba
  \partial_t Z& + \langle F(q,U), \nabla_q Z\rangle + (\nabla_\xi Z)D_u F(q,U) (D_q U \cdot \xi) - (\nabla_\xi Z)D_u G(q,U) \xi \\
   =&\, \xi D_q G(q,U) \xi -  ( D_q U\cdot \xi)D_u G(q,U) \xi - \xi D_q F(q,U) (D_q U \cdot \xi)\\
  & +(D_q U \cdot \xi)D_u F(q,U) (D_q U \cdot \xi) .
 \ea
$$
One can then recognize on the right hand side the matrix in the assumption \eqref{amon} as well as a transport equation on $Z$ on the left hand side. In the general case $\beta \ne 0$, a similar (and slightly more involved) computation, which is left as an exercise, leads to
 $$
\begin{aligned}
 \partial_t Z& + \langle F(q,U), \nabla_q Z\rangle + (\nabla_\xi Z) D_u F (q,U )( D_qU\cdot \xi)-(\nabla_\xi Z) D_u G (q,U ) \xi\\
 = &\,\xi D_qG(q, U) \xi - (D_qU\cdot\xi) D_uG(q, U) \xi - \xi D_q F(q,U)(D_q U \cdot \xi)\\
 & +  (D_q U \cdot\xi )D_u F(q, U) (D_q U \cdot\xi )\\
&- 2 \beta (D_q U \cdot \xi) D_q G(q, U)\xi + 2 \beta (D_q U \cdot \xi) D_q  F( q,U) (D_qU\cdot \xi). 
\end{aligned}
$$
 Using the $\alpha$-monotonicity of $(G,F)$, we obtain that the right-hand side is bounded from below by
 $$
 \ba
 \alpha&(|\xi|^2 + |D_q U\cdot \xi|^2) - 2 \beta (D_q U \cdot \xi) D_q G(q, U)\xi + 2 \beta (D_q U \cdot \xi) D_q  F( q,U) (D_qU\cdot \xi)\\
 &\geq \alpha(|\xi|^2 + |D_q U\xi|^2) - \beta K(|\xi|^2 + |D_q U\cdot\xi|^2),
 \ea
 $$
 where $K = 4 \max(\|D_qG\|_\infty,\|D_qF\|_\infty)$. Hence, if $\beta < \alpha /K$, $Z$ is a solution to
 $$
  \partial_t Z + \langle F(q,U), \nabla_q Z\rangle + (\nabla_\xi Z) D_u F (q,U )( D_qU\cdot \xi)-(\nabla_\xi Z) D_u G (q,U ) \xi \geq 0.
 $$
 Furthermore, if $\beta < \frac{\alpha}{\|D_qU_0\|^{2}}$, then $Z|_{t = 0} \geq 0$. Hence, arguing with a comparison principle, we indeed obtain that $Z \geq 0$.
 \end{proof}
 \begin{Rem}
The comparison principle used to conclude that $Z \geq 0$ is not detailed above. Since $\bar\Omega$ is compact, $Z$ smooth and quadratic in $\xi$, I leave it as an exercise to the interested reader. Remark also that the relation \eqref{Zxi} is the analogue of \eqref{relationwn} used earlier to show the propagation of monotonicity.
 \end{Rem}
 Note that in the previous result, I used more than simple monotonicity and Lipschitz regularity, a stronger form of monotonicity was quite necessary to get this a priori estimate. There are variants of this result where the requirements are still stronger than monotonicity, but weaker than the one used here, see \citep{lions20112012,bertucci2019some} for instance. There are counter-examples when we only assume simple monotonicity.

\subsubsection{Existence of solutions on arbitrary long time intervals II: justification of the a priori estimate}\label{sec:lastme}
In order to finally show existence of a Lipschitz solution to \eqref{me:finite} on arbitrary long time intervals, it remains to justify the previous a priori estimate, or in other words to show that a Lipschitz solution will actually verify the quantitative estimate $\|D_qU\|_\infty \leq \beta^{-1}$, even though it is not a smooth solution to \eqref{me:finite} as in the statement of Proposition \ref{prop:estn}. In order to do so, we are going to focus on a less differentiable quantity than the $Z$ of Proposition \ref{prop:estn} which is the quantity $W$ defined by
\be\label{defW4}
W(t,q,p) = \langle U(t,q)-U(t,p),q-p\rangle - \beta |U(t,q) - U(t,p)|^2.
\ee
In fact we are now going to establish the exact analogue of Proposition \ref{prop:estn} on this non-differentiated quantity.
\begin{Prop}
Consider $U$ a Lipschitz solution to \eqref{me:finite} on $[0,T]$. Assume that $U_0$ is Lipschitz, that $G$ and $F$ are Lipschitz continuous and that \eqref{amon} holds for some $\alpha > 0$, then $W$ defined by \eqref{defW4} is non-negative for some $\beta$ depending only on $\alpha$, $\|D_qU_0\|_\infty$, $\|D_qG\|_\infty$ and $\|D_qF\|_\infty$.
\end{Prop}
\begin{proof}
Assume that the result does not hold. Consider a point $(t_0,q_0,p_0)$ of minimum of $W$. Choosing $\beta$ sufficiently small, we know that $t_0 > 0$ and $U(t_0,q_0)\ne U(t_0,p_0)$. In particular, for any $s < t_0$
\be\label{eq1676}
\ba
&\langle U(t_0,q_0)-U(t_0,p_0),q_0-p_0\rangle - \beta |U(t_0,q_0) - U(t_0,p_0)|^2 \leq\\
 &\langle U(t_0-s,q(s))-U(t_0-s,p(s)),q(s)-p(s)\rangle - \beta |U(t_0-s,q(s)) - U(t_0-s,p(s))|^2,
\ea
\ee
where $q'(s) = -F(q(s),U(t_0-s,q(s))), q(0) =q_0$ and $p'(s) = -F(p(s),U(t_0-s,p(s))), p(0) =p_0$. Using the definition of Lipschitz solutions, we obtain that
$$
\ba
U(t_0,q_0)-U(t_0,p_0)& =U(t_0-s,q(s))-U(t_0-s,p(s))\\
& +\int_0^sG(q(s'),U(t_0-s',q(s'))) - G(p(s'),U(t_0-s',p(s')))ds'.
\ea
$$
Hence, plugging this information in \eqref{eq1676}, dividing by $s$ and taking the limit $s \to 0$ yields
$$
\ba
& \langle G(q_0,U(t_0,q_0)) - G(p_0,U(t_0,p_0)),q_0 - p_0\rangle \\
&+ \langle U(t_0,q_0)- U(t_0,p_0),F(q_0,U(t_0,q_0)) - F(p_0,U(t_0,p_0))\rangle\\
&\leq 2\beta \langle U(t_0,q_0) - U(t_0,p_0),G(q_0,U(t_0,q_0)) - G(p_0,U(t_0,p_0))\rangle
\ea
$$
using the monotonicity of $(G,F)$ we obtain that
$$
\alpha(|q_0-p_0|^2 + |U(t_0,q_0)-U(t_0,p_0)|^2) \leq 2\beta \langle U(t_0,q_0) - U(t_0,p_0),G(q_0,U(t_0,q_0)) - G(p_0,U(t_0,p_0))\rangle.
$$
Hence, the result follows, since the previous yields a contradiction if $\beta$ is below some threshold which does not depend on $U$, but simply on $G$ and $\alpha$.
\end{proof}
The previous result implies a quantitative Lipschitz bound on $U$ since it yields
$$
\sup_{t \leq T}\|D_qU(t,\cdot)\|_\infty \leq \beta^{-1}.
$$
Thus the maximal time of existence of Lipschitz solutions has to be $+\infty$. In particular, we have proven the following.
\begin{Cor}
Assume that $U_0$ is Lipschitz, that $G$ and $F$ are Lipschitz continuous and that \eqref{amon} holds for some $\alpha > 0$. Then, there exists a unique Lipschitz solution $U$ to \eqref{me:finite} on $[0,\infty)$.
\end{Cor}
Similar developments can also be made in the continuous state space case, following the same techniques, and I refer to the next section for references.

\subsection{Bibliographical comments}
The question of existence of solutions to the MFG MEs is of course deeply linked to the question of the notion of solution we are considering. Early results of existence of a solution through the methods of characteristics (which is not presented here and relies mostly on the formula \eqref{eq:626}) were presented by Lions in \citep{lions20112012}. This technique was then also used by Chassagneux, Crisan and Delarue \citep{chassagneux} and later by Cardaliaguet, Delarue, Lasry and Lions \citep{cdll}. The latter provided a detailed analysis of smooth solutions to the MEs by means of the associated forward-backward system. Specifically, the authors showed that under strong regularity requirements, we can differentiate the forward-backward system sufficiently many times to prove that the function defined through \eqref{eq:626} is indeed smooth, and thus the unique solution to the ME. Results were given in different monotone settings by Gangbo, Meszaros, Mou and Zhang \citep{gangbo2}, and MEs were notably studied in so-called anti-monotone regimes by Lions and Seeger \citep{lionsseeger} and Mou and Zhang \citep{mou2022mean}. See also Mou and Zhang \citep{mou2019wellposedness} which also passes through the characteristics. A method involving splitting was also introduced by Cardaliaguet, Cirant and Porretta in \citep{ccp}.

The notions of monotone solutions and Lipschitz solutions to the ME are weaker notions of solutions, which require weaker assumptions, and can be established by more abstract methods. The existence of monotone solutions through compactness arguments was presented in my works \citep{bertucci2021monotone,bertucci2023monotone}. With Lasry and Lions, we introduced the notion of Lipschitz solutions in \citep{bertucci2024lipschitz}, and showed existence of such solutions, also making the link with more regular notions of solutions.

The Lipschitz a priori estimate on the solution to the ME, which is crucial to establish the existence of solutions, was first established by Lions in \citep{lions20112012}, and later refined by Lasry, Lions and me in \citep{bertucci2019some} in finite dimensional cases, and by myself in \citep{bertucci2021master} and by Cardaliaguet et al in \citep{cdll} in the continuous state space case. A complete scheme of proof of the existence of Lipschitz solutions on large time intervals was written by Meynard and me in more involved cases in the series of papers \citep{meynard1,meynard2,meynard3}.

\newpage

\part{Mean field Hamilton-Jacobi-Bellman equations on the space of probability measures}
In this part, I focus on mean field Hamilton-Jacobi-Bellman (HJB) equations on spaces of probability measures and their link with mean field optimal control problems. I will start by showing how we can derive such equations from optimal control problems, omitting other areas where such PDEs appear, such as in large deviations for instance. I will then concentrate on the notion of viscosity solution to such equations, which is the correct one in this context. The notion of viscosity solution allows for an extremely robust proof of uniqueness of solutions, which involves, as usual in viscosity solution theory, the crucial comparison principle. I will then establish existence of solutions for cases arising from optimal control (I will also comment on more general cases). Finally, I will apply these methods to a problem of optimal transport toward a stochastic target.

\section{Derivation of an HJB equation}\label{sec:derivHJB}
In this section, I want to present what mean field optimal control problems are and how they give rise to HJB equations on the spaces of probability measures. Throughout this section, $\mo$ shall denote either $\T^d$ or $\R^d$. Furthermore, let me mention that the problem of interest in this section is clearly an optimal transport problem. However, the following presentation is quite different from usual ones in this theory.

\subsection{The optimal control problems}\label{sec:optcontr}
Starting from a measure $µ \in \mpo$ at time $t = 0$, we are interested in transporting it optimally in a certain time interval $[0,T]$. The optimality is with respect to a cost which is to be given just below.  As we saw in Section \ref{sec:variations}, there are several ways to consider variations on $\mpo$. Here, we focus on the case in which the controller chooses speeds, or velocities, which are going to transport the initial measure $µ$, as in Section \ref{sec:lagrange}. In particular, general Lagrangian evolutions, or randomization are going to be allowed.\\

To introduce the problem, assume that the controller chooses a (smooth) velocity field $\alpha : [0,T]\times \mo \mapsto \R^d, (t,x) \mapsto \alpha(t,x)$ such that $\|D_x \alpha\|_\infty < \infty$, then the evolution $(m_t)_{t \in [0,T]}$ of the initial measure $µ$ will be given through the continuity equation
\be\label{conteq}
\partial_t m + \text{div}(\alpha m) = 0 \text{ in } (0,T) \times \mo,
\ee
with initial condition $m_0 = µ$. In this case, we are going to consider costs of the form
$$
\tilde{\mathcal J}(\alpha,m) := \int_0^T\int_\mo L(x,\alpha(t,x),m_t)m_t(dx)dt + \mathcal G (m_T),
$$
where $L: \mo \times \R^d \times \mpo \mapsto \R_+$ and $\mathcal G: \mpo \mapsto \R_+ \cup \{+ \infty\}$ are given functions. The interpretation of the cost functional $\tilde{\mathcal J}$ is quite natural:
\begin{itemize}
\item $\mathcal G$ represents a standard terminal cost, which can be infinite to model constraints on the final measure $m_T$.
\item At time $t$, and position $x$, if the speed chosen is $\alpha(t,x)$, the cost is augmented by $L(x,\alpha(t,x),m_t)m_t(dx)dt$. The fact that this quantity is linear in $m_t(dx)$ plays a \emph{fundamental} role in the following, namely because it yields the shape of the HJB equation below. It is here that we see the mean field part of the problem. Indeed, if we forget about the dependence of $L$ in $m$, the cost is linear in $m_t(dx)$, the number of particles at $x$ at time $t$.
\end{itemize}

As we saw in Section \ref{sec:variations}, these kinds of dynamics are not entirely satisfactory as they do not allow us to go from any measure to any other measure under reasonable assumptions. The classical solution in optimal transport theory to this problem has been to consider the (very natural) possibility for the controller to choose a velocity field $\alpha$ which is discontinuous in $x$, namely since $\tilde{\mathcal J}$ does not require any regularity on $\alpha$. In this context, \eqref{conteq} may have multiple solutions and the control is then not only $\alpha$ but a pair $(\alpha,m)$. In the following I consider the setting of mixed flows of Section \ref{sec:lagrange}, which is more natural in these types of problems.\\

In what follows, I shall assume quite naturally that the controller can randomize its controls. Formally, this means that the controller has access to a way to draw randomly the control it chooses. This setting is natural in games, where it is called mixed strategies, or in optimal control, where it is sometimes called relaxed controls. Mathematically, the controller can choose a randomization space $K$ and a mixed vector field $b: [0,T] \times \mo \times K \mapsto \R^d$, and an associated weak solution\footnote{In general there might be multiple solutions to the continuity equation in this setting.} $\tilde m\in \mathcal C([0,T],\mathcal P(\mo\times K))$ to the continuity equation
\be\label{contpsi}
\partial_t \tilde m + \text{div}_x\left(b(t,x,k) \tilde m\right) = 0 \text{ in } (0,T)\times \mo\times K,
\ee
such that $(\pi_1)_\#\tilde m_0 = µ$. As presented in Section \ref{sec:lagrange}, from the point of view of dynamics, it is equivalent to choosing a standard probability space $(\Omega, \mathcal A, \mathbb P)$ and a process $(X_t,Z_t)_{t \in [0,T]}$ upon it, which has to be such that $\mathcal L(X_0) = µ$ and, almost surely, for all $t \in [0,T]$,
\be\label{sde3}
X_t = X_0 + \int_0^t Z_s ds.
\ee
In this Lagrangian setting, the cost is naturally given by
$$
\mathcal J_1(X,Z):= \mathbb E_{\mathbb P} \left[\int_0^T L(X_t,Z_t,\mathcal L(X_t))dt\right] + \mathcal G (\mathcal L(X_T)).
$$
In the previous, the last term is deliberately taken out of the expectation to highlight the fact that there is no extra randomness in the problem here, i.e. $\mathcal L(X_t)$ is deterministic. In the formulation of \eqref{contpsi}, the natural cost is written as
$$
\mathcal J_2(b,\tilde m) := \int_0^T\int_{\mo\times K} L(x,b(t,x,k),(\pi_1)_\#\tilde m_t)\tilde m_t(dx,dk)dt + \mathcal G ((\pi_1)_\#\tilde m_T).
$$
To pay this kind of cost for a mixed vector field is quite natural. The reader can refer to the enormous literature on mixed equilibria in game theory or on the not so small literature on relaxed control as well.

The optimal control problem is then given by
\be\label{problem}
\inf_{(b,\tilde m)} \mathcal J_2(b,\tilde m),
\ee
where the infimum is taken over all pairs $(b,\tilde m)$ such that $\tilde m \in \mathcal C([0,T],\mathcal P(\mo\times K))$, $(\pi_1)_\# \tilde m_0 = µ$, $b: [0,T]\times \mo\times K\mapsto \R^d$ is measurable and $(b,\tilde m)$ solves \eqref{contpsi} in the weak sense, that is, for all $t \in [0,T]$, and $\varphi \in \mathcal C^0_c([0,T]\times \mo\times K)$, $\mathcal C^1$ in $(t,x)$, 
$$
\langle \tilde m_t,\varphi(t,\cdot,\cdot) \rangle = \langle \tilde m_0,\varphi(0,\cdot,\cdot)\rangle + \int_0^t \int_{\mo \times K}(\partial_t \varphi(s,x,k) + \nabla_x\varphi(s,x,k)\cdot b(s,x,k) )\tilde m_s(dx,dk)ds.
$$
Several remarks are now in order. 
\begin{itemize}
\item Under reasonable growth conditions, see Remark \ref{rem:equiv} below, the problem \eqref{problem} is equivalent to 
\be\label{problem2}
\inf_{(X,Z)} \mathcal J_1(X,Z),
\ee
where the infimum is taken over pairs $(X_t,Z_t)_{t \in [0,T]}$ such that $\mathcal L(X_0) = µ$ and \eqref{sde3} is satisfied almost surely for all $t\in [0,T]$. This is an immediate consequence of Proposition \ref{prop:ambrosiov2} above.

\item The randomization set $K$ in \eqref{problem} or the probability space in \eqref{problem2} are in fact also a part of the decision of the controller. But since they can be fixed (provided that they are sufficiently rich), I shall omit this choice in the following.
\item I am here dealing with specific controlled trajectories, namely only the speed is controlled and there are no external forces. Generalizations are of course possible, and several other control trajectories are natural, for instance as in mean field optimal stopping or singular control.
\item When the cost $L$ is convex in its second argument, we immediately see that the problems \eqref{problem} and \eqref{problem2} are also equivalent to the infimum over all pair $(\alpha, m)$ which are measurable, with $m \in \mathcal C([0,T],\mpo)$ such that \eqref{conteq} is satisfied in a weak sense, of the cost $\tilde {\mathcal J}(\alpha,m)$. Indeed, one problem has a lower value than the other because it is a relaxation (hence more controls are available). Jensen's inequality and developments following Remark \ref{rem:lc} give the other inequality.
\end{itemize}
 Of the two problems \eqref{problem} and \eqref{problem2}, it is not clear that one has a strictly better form than the other. This is meant in terms of notation of course, since I just stated that they are equivalent. The formulation of \eqref{problem} is more insightful about the type of relaxation phenomena that appear for the controller, but it is quite a bit heavier in terms of notation. On the other hand, \eqref{problem2} uses the powerful language of the probability theory to convey the same idea in a more efficient way, but can be misleading about the nature of the problem, especially when there is an actual randomness in the problem.
 
 \subsection{Derivation of mean field optimal control from control of $N$-particle systems}
 I derive here the problems described above as the mean field limit of optimal control problems of systems of $N$ particles. Since the convergence of such problems is not a topic of this book, I am going to be quite brief, and simply state some ideas.
 
 Let $N \geq 1$ and consider $(x_0^i)_{1 \leq i \leq N}$ an initial configuration in $\mo^N$. The controller chooses an evolution of the state (in $\mo^N$) by means of a control $\beta: \{1,\dots,N\}\times [0,T]\mapsto \R^d$, whose action on the state is described by the ODE
 \be\label{edoN}
 \frac{dx^i_t}{dt} = \beta^i_t \text{ in } (0,T), \text{ for } 1 \leq i \leq N.
 \ee
 The associated cost is then given by 
 $$
\int_0^T L^N(x_t,\beta_t)dt + \mathcal G^N(x_T).
 $$
 When the costs $L^N$ and $\mathcal G^N$ satisfy specific structural conditions, we can expect to derive a problem of the form of \eqref{problem} from the previous. One of the first requirements is that $\mathcal G^N$ is symmetric. Another one is that, in $L^N$, the effects of $\beta^i$ are summed, as well as some symmetry condition on the $x^i$. Without trying to be as generic as possible, I will simply state a typical situation in which we expect to have a mean field limit. For some costs $L$ and $\mathcal G$ as in the previous section, assume that
 \be\label{structure}
 \ba
 L^N(x,\beta) = \frac{1}{N}\sum_{i =1}^NL(x^i,\beta^i, m^N_x),\\
 \mathcal G^N(x) = \mathcal G(m^N_x),
 \ea
 \ee
where $m^N_x = \frac 1N \sum_{i=1}^N\delta_{x^i}$. Following the formulation of the mean field optimal control problem through \eqref{problem2}, we can easily understand the type of limit we are looking for. Indeed, here, we want to exploit the fact that, formally, as $N \to \infty$, $((x^i_t,\beta^i_t)_{1 \leq i \leq N})_{t \in [0,T]}$ is going to converge in law toward some process $(X_t,Z_t)_{t \in [0,T]}$. Then, equation \eqref{sde3} becomes the limit of \eqref{edoN}, while the limit of the cost is simply $\mathcal J_1(X,Z)$. The previous can be intuited since, because of the condition \eqref{structure}, we already have that $L^N$ is a sort of expectation of the couple $(x^i,\beta^i)_{1 \leq i \leq N}$. Observe that from this point of view, problems \eqref{problem} and \eqref{problem2} are more natural than the minimization of $\tilde {\mathcal J}(\alpha,m)$, because there is no reason to impose the condition $\beta^i_t = \beta^j_t$ if $x^i_t=x^j_t$.

\subsection{Bellman's dynamic programming and HJB equation}
To study \eqref{problem} (or the equivalent \eqref{problem2}), I want to adopt Bellman's point of view, or in other words, instead of looking directly at \eqref{problem}, to consider the following family of problems indexed by $(t,µ) \in [0,T)\times \mpo$
\be\label{defU}
U(t,µ) := \inf_{(b,\tilde m)}\left\{\int_t^T\int_{\mo\times K} L(x,b(s,x,k),(\pi_1)_\#\tilde m_s)\tilde m_s(dx,dk)ds + \mathcal G ((\pi_1)_\#\tilde m_T)\right\},
\ee
where the infimum is taken over all pairs $(b,\tilde m)$ such that $\tilde m \in \mathcal C([t,T],\mathcal P(\mo\times K))$, $(\pi_1)_\#\tilde m_t = µ$, $b: [t,T]\times \mo\times K\mapsto \R^d$ is measurable and $(b,\tilde m)$ solves in a weak sense
$$
\partial_s \tilde m + \text{div}_x\left(b(s,x,k) \tilde m_s\right) = 0 \text{ in } (t,T)\times \mo\times K.
$$
Of course, $U$ could have been defined similarly with the probabilistic notation of \eqref{problem2}.

The function $U$ is called the value function of the problem, as it gives the value of the sub-problem starting in the state $(t,µ)\in [0,T)\times \mpo$. Observe that we expect\footnote{One has to be careful with such boundary conditions, because as usual in optimal control, they are not always satisfied.} $U$ to be given at time $T$ by $U(T,µ) = G(µ)$ for any $µ \in \mpo$. Bellman's main ideas are: i) that we can compute $U$ backwardly in time, and ii) that given $U$ we can find optimal controls in \eqref{problem}. The first idea relies mostly on the so-called dynamic programming principle, that is stated just below. I come back to the optimal controls later on.

\begin{Prop}\label{prop:dpp}
Consider $U$ defined by \eqref{defU}. For any $t < s < T$ and $µ\in \mpo$, 
\be\label{dpp}
U(t,µ) = \inf_{(b,\tilde m)} \left\{\int_t^s\int_{\mo\times K} L(x,b(u,x,k),(\pi_1)_\#\tilde m_u)\tilde m_u(dx,dk)du + U(s,(\pi_1)_\#\tilde m_s)\right\},
\ee
where the infimum is taken as in \eqref{defU}. In the previous, all quantities can be equal to $+\infty$.
\end{Prop} 
\begin{proof}
The fact that $U$ is greater than the $\inf$ simply follows from the definition and the splitting of the integral. To obtain the other inequality, let $t < s < T$, $µ \in \mpo$ and consider an $\eps$-optimal pair $(b^\eps,\tilde m^\eps)$ for the right hand side of \eqref{dpp} that we denote $V(t,µ)$. Consider now, an $\eps$-optimal pair $(a^\eps,\tilde n^\eps)$ for $U(s,(\pi_1)_\#\tilde m^\eps_s)$. If we are able to concatenate the pairs $(b^\eps,\tilde m^\eps)$ and $(a^\eps,\tilde n^\eps)$ into an admissible control for $U(t,µ)$, then the result holds since this would imply
$$
U(t,µ) \leq V(t,µ) + 2 \eps,
$$
with $\eps > 0$ arbitrary. Such a concatenation is always possible, even with the presence of the space $K$. Indeed, using the probabilistic formulation of the problem for instance, up to enlarging the probability space, we can always consider a concatenation of two such processes. Because we assumed that $K$ is always sufficiently large, we now deduce from standard probabilistic arguments (Theorem \ref{thm:existrv}), that we can find an admissible control $(b,\tilde m)$ for $U(t,µ)$ which has the same law as this concatenation.
\end{proof}

Using \eqref{dpp}, we can derive formally the HJB equation satisfied by $U$. Indeed, take $t < T, µ \in \mpo$ and $\eps > 0$ such that $t+ \eps < T$, from \eqref{dpp} we obtain
$$
0 = \inf_{(b,\tilde m)} \left\{\frac1\eps\int_t^{t+ \eps}\int_{\mo\times K} L(x,b(s,x,k),(\pi_1)_\#\tilde m_s)\tilde m_s(dx,dk)ds + \frac{U(t+ \eps,(\pi_1)_\#\tilde m_{t+\eps})-U(t,µ)}{\eps}\right\}.
$$
Arguing formally that $U$ is H-differentiable everywhere and that $b$ is such that we can indeed take the variation of $U$ along the continuity equation of the form \eqref{contpsi}, we obtain
$$
\ba
0 = \inf_{(b,\tilde m)} \bigg\{&\frac1\eps\int_t^{t+ \eps}\int_{\mo\times K} L(x,b(s,x,k),(\pi_1)_\#\tilde m_s)\tilde m_s(dx,dk)ds + \eps^{-1}o(\eps) + \partial_t U(t,µ)\\
 &+ \frac1\eps\int_t^{t+\eps}\int_{\mo\times K} D_µ U(s,µ)(x)\cdot b(s,x,k)\tilde m_s(dx,dk)ds\bigg\}
 \ea
$$
Pursuing formal computations, I proceed as if the $o(\eps)$ does not depend on $(b,\tilde m)$, exchanging the integral and the infimum and taking the limit $\eps \to 0$, we obtain that 
$$
0 = \partial_t U + \int_\mo \inf_{b\in \R^d}\left\{\int_{\R^d}D_µU(t,µ)(x)\cdot b + L(x,b,µ)\right\}µ(dx).
$$
If we set for $x\in \mo,p\in \R^d,µ\in \mpo$, $H(x,p,µ) := -\inf_{\alpha \in \R^d}\{L(x,\alpha,µ)+ \alpha \cdot p\}$, we then arrive at the HJB equation
\be\label{hjb1}
-\partial_t U(t,µ) + \int_\mo H(x,D_µU(t,µ)(x),µ)µ(dx) = 0 \text{ in } (0,T)\times \mpo.
\ee
These computations are the exact reason why it is the notion of H-differentiability of Part I that is going to be of interest in this third part. The nature of the underlying variations, the ones that are given in our mean field optimal control problem, dictates the notion of derivative that appears in the HJB equation.
\begin{Rem}
In this derivation of the HJB equation, the randomization of controls through the set $K$ seems entirely arbitrary. This is largely due to the fact that the function $U$ is assumed to be differentiable. As it is not the case in general, to make such a derivation more rigorous, which is done in Section \ref{sec:existvisc}, we shall need to use sub/super-differentials. Doing so, the use of this randomization will become clearer.
\end{Rem}

\subsection{Explicit solutions, dynamic formulation of Wasserstein distances and lack of smoothness}\label{sec:explicit}
The Wasserstein distances introduced in Part I can, in certain cases, be expressed through a dynamic formulation, quite in the spirit of \eqref{problem} and \eqref{problem2}. This will allow me to make certain remarks on what we can expect from the value function $U$ of such problems. I start with the following.

\begin{Prop}
For any $p \in (1,\infty)$, $T > 0$, $µ,\nu \in \mppo$
\be\label{static-dynamic}
\mathcal W_p^p(µ,\nu) = T^{p-1}\inf_{(X,Z)} \mathbb E\left[\int_{0}^T |Z_t|^pdt\right],
\ee
where the infimum is taken over $(X_t,Z_t)_{t \in [0,T]}$ such that $\mathcal L(X_0) = µ, \mathcal L(X_T) = \nu$ and, almost surely, \eqref{sde3} holds for all $t \in [0,T]$.

\end{Prop}
\begin{proof}
Consider the case $T =1$, $\mo =\R^d$. Using the probabilistic formulation of Wasserstein distances, we know that
$$
\mathcal W_p^p(µ,\nu) = \inf_{X_0 \sim µ, X_1 \sim \nu} \mathbb E [|X_0-X_1|^p].
$$
It remains to observe that the right hand side is smaller than $\inf_{(X,Z)}\mathbb E\left[\int_{0}^1 |Z_t|^pdt\right]$ as in the statement of the result, and that equality holds for $Z_t = X_1 - X_0$, because of the convexity of $x \mapsto |x|^p$. The case $T \ne 1$ follows from a change of variables and the case $\mo = \T^d$ is a direct extension.
\end{proof}

As slight generalizations of the previous, we can express more value functions as Wasserstein distances as follows.
\begin{Prop}
Take $p> 1$ and assume that for any $x\in \mo,\alpha \in \R^d, µ \in \mpo$, $L(x,\alpha,µ) = |\alpha|^p$ and that there exist $\eps > 0$ and $  \nu \in \mathcal P_p(\mo)$ such that for any $µ \in \mathcal P_p(\mo)$, $\mathcal G(µ) = \eps^{-1}W^p_p(µ,\nu)$.

Then, $U$ defined by \eqref{defU} is given for any $t \in (0,T), µ \in \mathcal P_p(\mo)$ by
$$
U(t,µ) = \frac{\mathcal{W}_p^p(µ,\nu)}{(T-t+\eps^{\frac{1}{p-1}})^{p-1}}.
$$
\end{Prop}
\begin{proof}
Performing the change of variable $s \leftarrow \lambda u + t$ and 
$b(s,x,k) = \lambda^{-1}b(u,x,k)$, for $\lambda = T-t$, we obtain that
$$
U(t,µ) = \inf_{(b,\tilde m)}\left\{\frac{1}{(T-t)^{p-1}}\int_0^1\int_{\mo\times K}|b(u,x,k)|^p\tilde m_u(dx,dk)du + \frac1\eps \mathcal{W}_p^p((\pi_1)_\#\tilde m_1,\nu)\right\},
$$
where the infimum is taken over the appropriate set (i.e. solutions of the continuity equation on $(0,1)$ starting from $µ$ at time $0$). Using the previous result on the second term, we end up with
$$
\ba
U(t,µ) =\frac{1}{(T-t)^{p-1}} \inf_{(b,\tilde m)}\bigg\{&\int_0^1\int_{\mo\times K}|b(u,x,k)|^p\tilde m_u(dx,dk)du+\\
& + \int_1^{1 + \kappa}\int_{\mo\times K}|b(u,x,k)|^p\tilde m_u(dx,dk)du\bigg\},
 \ea
$$
where $\kappa = \frac{\eps^{\frac{1}{p-1}}}{T-t}$ and the infimum is taken over appropriate $b$ and $\tilde m$, such that $(\pi_1)_\#\tilde m_0 = µ, (\pi_1)_\# \tilde m_{1 +\kappa} = \nu$. (An argument similar to the one used in Proposition \ref{prop:dpp} is implicitly used here.) Hence, we obtain
$$
U(t,µ) = \frac{1}{(T-t+\eps^{\frac{1}{p-1}})^{p-1}} \mathcal{W}_p^p(µ,\nu).
$$
\end{proof}
 
The main interest of the previous computation is to lead us to one of the main issues of problems such as \eqref{problem}, namely that Wasserstein distances (and more generally value functions of such problems) are not H-differentiable, see Remark \ref{rem:nondiff}. Hence, we shall need a weak notion of solution to \eqref{hjb1} in order to characterize the value function as the unique solution, in a sense to be given in the next section, of the associated HJB equation.

\subsection{Formal optimality conditions in deterministic settings}\label{sec:optcond}
I present here another way to address problems such as \eqref{problem}, which is particularly efficient when the problem at hand is deterministic. Since this other approach uses different tools from the ones introduced in this book, I will stay formal and only present it for the sake of completeness.\\

I start from almost the same problem as \eqref{problem}, except that I assume that the cost does not mix the dependencies in $\alpha$ and $m$, and work at the level of the non-relaxed problem:
$$
\inf_{\alpha,m} \int_0^T \left(\int_\mo L(x,\alpha_t) m_t(dx) + \mathcal F(m_t)\right)dt + \mathcal G(m_T),
$$
where the infimum is taken over pairs $(\alpha,m)$ which solve 
\be\label{ce12}
\ba
\partial_t m + \text{div}_x(\alpha m) = 0 \text{ in } (0,T)\times \mo,\\
m|_{t = 0} = m_0.
\ea
\ee
The idea is here to see \eqref{ce12} as a constraint and to introduce the associated Lagrange multiplier. Hence, we minimize over all $(\alpha,m)$ the following
$$
\ba
\inf_{\alpha,m} &\int_0^T \left(\int_\mo L(x,\alpha_t)m_t(dx) + \mathcal F(m_t)\right)dt  + \sup_{u}\int_0^T\int_\mo(\partial_t u + \alpha\cdot \nabla_x u)m_t(dx)dt\\
&+ \mathcal G(m_T) - \int_\mo u(T,x)m_T(dx) + \int_\mo u(0,x)m_0(dx),
\ea
$$
 where the supremum is taken over all $u \in \mathcal C^{1}$, and where I did not simply multiply the constraint \eqref{ce12} by the Lagrange multiplier but also perform an integration by parts.\\
 
Standard questions of optimization and calculus of variations can then be raised on this $\inf \sup$ problem. Arguing again formally, the optimality conditions of the previous problem can be written in terms of the following coupled PDEs  
$$
\ba
-\partial_t u + H(x,\nabla_x u) = \nabla_µ\mathcal F(m_t) &\text{ in } (0,T)\times \mo,\\
\alpha_t(x) = -D_pH(x,\nabla_x u(t,x)) &\text{ in } (0,T)\times \mo,\\
\partial_t m + \text{ div}_x(\alpha m) = 0 &\text{ in } (0,T)\times \mo,\\
m|_{t = 0} = m_0, \quad u|_{t= T} = \nabla_µ &\mathcal G(m_T).
\ea
$$
This kind of system is of the same structure as \eqref{mfg:sys} in MFG, and its study can allow, under certain convexity assumptions on the problem, one to obtain optimal controls.

\subsection*{Bibliographical comments}
HJB equations are a central object in the mathematical analysis of optimal control problems. I do not recall this general literature in finite dimension, and simply refer to Bardi and Capuzzo-Dolcetta \cite{bardiitalo} for deterministic problems and to Fleming and Soner \citep{flemingsoner} for stochastic problems. I focus here on the case of equations posed on the set of probability measures.

There are several independent starting points for the study of HJB equations on spaces of probability measures. An important one was the course of Lions \citep{lions20082009} in which he studied an HJB equation on $\mpt$ as an illustration of the calculus he introduced on this space through the Hilbertian lift. An equation was also derived and studied by Cardaliaguet and Quincampoix in \cite{cardaquincampoix} to study a zero sum game. Quincampoix then continued to work on HJB equations with Jimenez and Marigonda \cite{marigonda2018mayer,jimenez2020optimal}, associated with optimal control problems, using a quite geometric point of view. Motivated by questions of large deviations in statistical physics, Feng produced an impressive list of works on the topic with several co-authors \cite{feng1,feng2,feng3,feng4,feng5}. This last set of works is mainly concerned with HJB equations having singular terms which are the horizontal differential of some functionals, except for Ambrosio and Feng \citep{feng2} which deals with so-called metric cases. Such problems have also been considered independently as optimal control problems, by a long list of authors: Cosso, Gozzi, Kharroubi, Pham and Rosestolato \cite{cosso2021master}, Daudin and Seeger \cite{daudinseeger}, Soner and Yan \cite{soner2022}, Bayraktar, Ekren and Zhang \cite{bayraktar2025comparison}; or because of their links with MFG, for instance by Cecchin and Delarue \cite{cecchin2025weak}. 

The study of optimal control problems of the form of \eqref{problem} has also been carried out through usual optimality conditions, without passing through the HJB equation, as in Section \ref{sec:optcond}. It is for instance the case of all results on optimal transport, especially on its dynamic formulation. I do not expand too much on this literature either but refer the reader to \cite{santambrogio,bensoussan2013mean}.

\newpage

\section{Viscosity solutions to HJB equations and uniqueness of solutions}\label{sec:visc}
 This section introduces the notion of viscosity solutions to mean field HJB equations on $\mpo$. Before doing so, I recall some standard facts about viscosity solutions to first order HJB equations in finite dimension. Then, the notion of viscosity solutions is presented on $\mptd$ and then on $\mpprd$. Both times, a comparison result is established
 
\subsection{Viscosity solutions in finite dimension}

The goal of viscosity solution theory is, given PDEs for which a formal comparison principle holds, to define a weak notion of solution to these PDEs for which a comparison principle still holds. In other words, for equations for which, formally, we can say that super-solutions are above sub-solutions, viscosity solution theory aims at making rigorous such statements for weak solutions. 

In this section, I consider
\be\label{hjbn}
u + H(x,\nabla_x u) = 0 \text{ in } \mo,
\ee
where $H: \mo \times \R^d \mapsto \R$ is the Hamiltonian, $u : \mo \mapsto \R$ is the solution and $\mo$ is still either $\T^d$ or $\R^d$. I will sometimes allow myself to write $|x-y|$ instead of $d(x,y)$, where $d(\cdot,\cdot)$ is the distance on $\T^d$ and $|\cdot|$ the Euclidean norm. This equation is said to satisfy a comparison principle, in the following sense.

\begin{Prop}
Let $u$ and $v$ be two smooth functions such that
$$
u + H(x,\nabla_x u) \leq 0 \text{ in } \mo,
$$
$$
v + H(x,\nabla_x v) \geq 0 \text{ in } \mo.
$$
Assume also that $\limsup_{|x |\to \infty}u(x) - v(x) \leq 0$. Then $u \leq v$.
\end{Prop}
\begin{proof}
Consider a maximizing sequence $(x_n)_{n \geq 0}$ of $\sup_{x \in \mo} u(x) - v(x)$. If $|x_n|\to \infty$, then the result is proved. If not, up to a subsequence, $(x_n)_{n \geq 0}$ has a limit point $x_*$ which is a point of maximum of $u-v$. In particular, $\nabla_x u(x_*) = \nabla_x v(x_*)$. Evaluating the inequalities satisfied by $u$ and $v$ at $x_*$, we obtain
$$
0 \geq u(x_*) + H(x_*,\nabla_x u(x_*)) = u(x_*) + H(x_*,\nabla_x v(x_*)) \geq u(x_*)-v(x_*).
$$ 
Hence the result follows.
\end{proof}
The comparison principle implies the uniqueness of smooth solutions to \eqref{hjbn} in the compact case, and uniqueness of smooth solutions that share the same behavior at infinity in the non-compact case. Indeed, considering two such solutions $u_1$ and $u_2$, they each can play the role of either $u$ or $v$ in the previous proposition. Hence we deduce that $u_1\leq u_2\leq u_1$.
A possible way to understand the theory of viscosity solutions is to interpret it as the weakest definition of a notion of solution to \eqref{hjbn}, such that the comparison principle is preserved. Observe that, a priori, the roles played by $u$ and $v$ in the previous proposition are not symmetric. We must therefore define notions of super-solution and sub-solution, as we do now.
\begin{Def}\label{def:viscN}
A usc (resp. lsc) function $u : \mo \mapsto \R$ is a viscosity sub-solution (resp. super-solution) to \eqref{hjbn} if, for every smooth function $\phi : \mo \mapsto \R$, and every $x_0 \in \mo$ which is a point of maximum (resp. minimum) of $u - \phi$, we have
$$
u(x_0) + H(x_0,\nabla_x \phi(x_0)) \leq 0 \text{ (resp. } \geq 0). 
$$
A locally bounded function $u$ is a viscosity solution to \eqref{hjbn} if its lower and upper semi-continuous envelopes (resp. $u_*$ and $u^*$) are resp. viscosity super-solution and sub-solution.
\end{Def}
From this definition, we can easily compare a viscosity super-solution with a smooth sub-solution, or a viscosity sub-solution with a smooth super-solution. The main question is in fact whether one can compare a viscosity super-solution with a viscosity sub-solution. The standard technique to answer this is the so-called doubling of variables method, which is the main argument in the next proof.
\begin{Prop}
Assume that there exists a modulus of continuity $\omega$ such that for all $q \in \R^d$, $x,y\in \mo$
\be\label{hypH1}
|H(x,q) - H(y,q)|\leq \omega((1+ |q|)|x-y|).
\ee
Furthermore, if $\mo = \R^d$, assume that for any $R > 0$, there exists $C_R$ such that for all $x \in \R^d, q_1,q_2 \in B(0,R)$
\be\label{hypH2}
|H(x,q_1) - H(x,q_2)| \leq C_R|q_1-q_2|.
\ee
Let $u$ be a viscosity subsolution and $v$ a viscosity supersolution to \eqref{hjbn}, both bounded. Then $u \leq v$.
\end{Prop}
\begin{proof}
We start with the compact case. Consider a point of maximum $(x_\eps,y_\eps)$ of 
$$
w : (x,y)\mapsto u(x) - v(y) - \frac{1}{2\epsilon}d(x, y)^2.
$$
Such a point of maximum always exists because $w$ is a usc function on a compact set. From Lemma \ref{lemma:visc1} below, we know that $d(x_\eps,y_\eps) \to 0$ as $\eps \to 0$. Thus, taking $\eps$ sufficiently small, $d(x_\eps,y_\eps) = |x_\eps-y_\eps|$. Using the fact that $u$ and $v$ are respectively viscosity sub and super-solution, we find
$$
u(x_\eps) + H(x_\eps,\epsilon^{-1}(x_\eps - y_\eps)) \leq 0,
$$
$$
v(y_\eps) + H(y_\eps,\epsilon^{-1}(x_\eps - y_\eps)) \geq 0.
$$
Taking the difference of the inequalities, we obtain 
$$
\max (u-v) \leq u(x_\eps) - v(y_\eps) \leq H(y_\eps,\eps^{-1}(x_\eps-y_\eps)) - H(x_\eps,\eps^{-1}(x_\eps-y_\eps))
$$
Using \eqref{hypH1} we find $\max (u-v) \leq \omega(|x_\eps-y_\eps| + \eps^{-1}|x_\eps-y_\eps|^2)$. Thanks to Lemma \ref{lemma:visc1} below, we obtain the required result by taking the limit $\eps \to 0$.\\

Consider now the non-compact case. For $\eps > 0, \delta >0$, consider $(x_\eps,y_\eps)$ a point of maximum of
$$
w: (x,y) \mapsto u(x) - v(y) - \frac{1}{2\eps}|x-y|^2 - \delta\sqrt{1 + |x|^2}.
$$
Remark that such a point always exists because, since $u$ and $v$ are bounded and respectively usc and lsc, $w$ is a usc function such that $w(z) \to -\infty$ as $|z|\to \infty$. Using that $u$ and $v$ are viscosity sub and super-solution, we obtain
$$
u(x_\eps) + H\left(x_\eps,\epsilon^{-1}(x_\eps - y_\eps) + \delta \frac{x_\eps}{\sqrt{1+|x_\eps|^2}}\right) \leq 0,
$$
$$
v(y_\eps) + H(y_\eps,\epsilon^{-1}(x_\eps - y_\eps)) \geq 0.
$$
Hence, in this case we end up with
$$
\ba
\sup u-v &\leq H(y_\eps,\epsilon^{-1}(x_\eps - y_\eps)) - H\left(x_\eps,\epsilon^{-1}(x_\eps - y_\eps) + \delta \frac{x_\eps}{\sqrt{1+|x_\eps|^2}}\right)\\
& \leq \omega\left((1+\eps^{-1}|x_\eps-y_\eps|)|x_\eps - y_\eps|\right) + C_{1 + \eps^{-1}|x_\eps-y_\eps|}\delta,
\ea
$$
thanks to \eqref{hypH2}. Hence the result still holds, this time thanks to Lemma \ref{lemma:visc2}, by passing first to the limit $\delta \to 0$, and then $\eps \to 0$.
\end{proof}
\begin{Rem}
In this proof, the distinction between $d(\cdot,\cdot)$ and $|\cdot-\cdot|$ is not needed because in the limit $\eps \to 0$ $x$ and $y$ are close to each other. When working on $\mptd$, we shall not be able to make this simplification, and some care will be needed to handle $\T^d$ rigorously.
\end{Rem}
\begin{Rem}
The non-compact case relies on some additional estimate on the Hamiltonian $H$ and on the assumption that $u$ and $v$ are bounded. More general statements could be made, namely by considering different sets of assumptions for $H$, which allows for unbounded solutions.
\end{Rem}
We used the following lemma in the compact case.
\begin{Lemma}\label{lemma:visc1}
Let $E$ be a metric space. Let $u: E \mapsto \R$ be a usc function and $\Psi : E \mapsto \R_+$ a lsc function. For $\eps > 0$, let 
$$
M_\eps = \sup_{ x \in E} \left\{u(x)- \frac{1}{\eps} \Psi(x)\right\}.
$$
Assume $M_{\eps_0} < \infty$ for some $\eps_0> 0$ and $\inf_{\eps > 0} \{M_\eps\} > -\infty$. Consider $(x_\eps)_{\eps > 0}$ such that $\lim_{\eps\to 0} (M_\eps - u(x_\eps) + \frac1\eps \Psi(x_\eps)) = 0$. Then $\lim_{\eps \to 0}\frac{\Psi(x_\eps)}{\eps} = 0$ and if $x_*$ is a limit point of $(x_\eps)_{\eps > 0}$, $\lim_{\eps \to 0} M_\eps = u(x_*) = \sup_{\Psi(x) = 0}u(x)$.
\end{Lemma}
\begin{proof}
Denote $\alpha_\eps = M_\eps - u(x_\eps) + \frac1\eps \Psi(x_\eps)$. Since $\Psi \geq 0$, $M_\eps$ decreases as $\eps \to 0$, and thus $\lim_{\eps \to 0} M_\eps$ exists. Hence
$$
M_{2 \eps} \geq u(x_\eps) - \frac{1}{2\eps}\Psi(x_\eps) \geq u(x_\eps) - \frac{1}{\eps}\Psi(x_\eps) + \frac{1}{2\eps}\Psi(x_\eps) \geq M_\eps-\alpha_\eps + \frac{1}{2\eps}\Psi(x_\eps).
$$
Thus $\frac{\Psi(x_\eps)}{\eps} \leq 2 (M_{2 \eps}-M_{\eps} + \alpha_\eps) \to_{\eps \to 0} 0$ as $\lim_{\eps \to 0} M_\eps$ exists. Then, if $x_*$ is the limit of some $(x_{\eps_n})_{n \geq 0}$, 
$$
u(x_*) \geq \lim_{n \to \infty}u(x_{\eps_n}) \geq \lim_{n \to \infty} M_{\eps_n} -\alpha_{\eps_n}\geq \sup_{\Psi(x)= 0}u(x),
$$
which is the result.
\end{proof}
In the non-compact case, I used a slightly different version of the previous result.
\begin{Lemma}\label{lemma:visc2}
Let $E$ be a metric space. Let $u: E \mapsto \R$ be a usc function and $\Psi,\Phi : E \mapsto \R_+$ lsc functions. For $\eps,\delta > 0$, let
$$
M_\eps = \sup_{ x \in E} \left\{u(x)- \frac{1}{\eps} \Psi(x)\right\}.
$$
$$
M_{\eps,\delta} = \sup_{ x \in E} u(x)- \frac{1}{\eps} \Psi(x) - \delta \Phi(x).
$$
Assume $M_{\eps_0} < \infty$ for some $\eps_0> 0$ and $\inf_{\eps > 0} \{M_\eps\} > -\infty$. Consider $(x_{\eps,\delta})_{\eps,\delta > 0}$ such that $M_{\eps,\delta}=  u(x_{\eps,\delta}) - \frac1\eps \Psi(x_{\eps,\delta}) -\delta\Phi(x_{\eps,\delta})$. Then $(\delta \Phi(x_{\eps,\delta}))_{\delta > 0}$ is bounded for any $\eps>0$ and, $\lim_{\eps \to 0}\limsup_{\delta \to 0}\frac{\Psi(x_{\eps,\delta})}{\eps} = 0$.
\end{Lemma}
\begin{proof}
The first property is obvious from the non-negativity of $\Psi,\Phi$. For the second one, it follows from the previous proof. Indeed, a similar computation yields that
$$
\frac{\Psi(x_{\eps,\delta})}{\eps} \leq 2(M_{2\eps,\delta}-M_{\eps,\delta}).
$$
However, observe that for any $\eps > 0$, $\lim_{\delta \to 0}M_{\eps,\delta} = M_\eps$. Hence, the result follows since $(M_\eps)_{\eps >0}$ has a limit when $\eps \to 0$.
\end{proof}

To conclude these reminders on viscosity solutions in finite dimension, note that Definition \ref{def:viscN} can be replaced by the following:
\begin{Def}
A usc (resp. lsc) function $u : \mo \mapsto \R$ is a viscosity sub-solution (resp. super-solution) to \eqref{hjbn} if, for all $q \in \partial^+u(x)$ (resp. $q \in \partial^-u(x)$), we have
\be
u(x) + H(x,q) \leq 0  \text{ (resp. } \geq 0). 
\ee
\end{Def}

I do not go into the details of the equivalence between this notion and the previous one. Let me just mention that it relies on the fact that if $u - \phi$ has a local minimum at $x_0$, then $\nabla_x \phi(x_0) \in \partial^-u(x_0)$. And also on a form of converse: if $p \in \partial^-u(x_0)$, then there exists a $\mathcal C^1$ function $\phi$ such that $\phi(x_0) = p$ and $x \mapsto u(x) - \phi(x)$ has a local minimum at $x_0$. Such type of equivalence will not be at our disposal in the case of the space of probability measures.

\subsection{Viscosity solutions on the space of probability measures: the compact case}\label{sec:visccomp}
 In this section, I focus on the HJB equation
 \be\label{hjbt}
 U(µ) + \int_{\T^d}H(x,D_µU(µ)(x),µ)µ(dx) = 0 \text{ in } \mptd,
 \ee
 where $H: \T^d\times \R^d\times \mptd \mapsto \R$ is a function on which regularity assumptions shall be made later on, and $D_µ$ stands for the horizontal derivative on $\T^d$ (Section \ref{sec:horiz}). Motivated by Section \ref{sec:explicit}, there is a need to characterize non-differentiable solutions to \eqref{hjbt}, as we do not expect the typical solution to be differentiable.

 As mentioned earlier, in finite dimension, viscosity solutions can be defined equivalently with either super-differentials or smooth test functions. In this setting, it is not the case. The reason for this is quite apparent from Part I. Indeed, smooth functions have super-differentials which are deterministic, compared to the richness that one can obtain in general. Hence, because the functions we plan to characterize as viscosity solutions are not smooth, restricting ourselves to smooth test functions is equivalent to forgetting a potentially huge part of their super/sub-differentials. Note that the price to pay for such a choice is that the existence of such solutions shall then be harder to prove, compared to what should happen if I restricted the attention to smooth test functions.
 
For the integral term to be well defined in \eqref{hjbt}, a growth assumption on $H$ is necessary. We defined in Section \ref{sec:torus} elements in super/sub-differentials with an integrability condition. Thus, the following condition shall be sufficient 
 \be\label{star12}
\exists C > 0, \forall x \in \T^d, q \in \R^d, µ \in \mptd, \quad |H(x,q,µ)| \leq C(1+ |q|^2).
 \ee
  \begin{Rem}
 Note that even if HJB equations which are obtained through dynamic programming, as in the previous section, are such that $H$ is convex in its second argument, I consider here a more general situation and do not impose any condition of the sort. Such non-convex Hamiltonians are for instance involved in dynamic zero-sum games.
 \end{Rem}
 All topological statements on $\mptd$ are made with respect to the weak topology for the rest of this section. The previous discussion leads to the following definition.
 \begin{Def}
 A usc (resp. lsc) function $U: \mptd \mapsto \R$ is said to be a viscosity sub-solution (resp. super-solution) to \eqref{hjbt} if for any $µ \in \mptd, \psi \in \partial^+U(µ)$, (resp. $\psi \in \partial^-U(µ)$), it holds that
 \be\label{visc1}
 U(µ) + \int_{\T^d\times \R^d}H(x,z,µ)\psi_x(dz)µ(dx) \leq 0\quad  (\emph{resp. } \geq 0).
 \ee
 A locally bounded function $U: \mptd \mapsto \R$ is a viscosity solution to \eqref{hjbt} if $U^*$ is a sub-solution and $U_*$ a super-solution.
 \end{Def}
 \begin{Rem}
 Let me insist upon the fact that, a priori, a non-trivial choice has been made here. Indeed, evaluating the equation by replacing $D_µU(µ)$ by $\psi$ is not automatic, namely since those two objects do not live in the same space. I could have considered instead the term $\int_{\T^d}H(x,\int_{\R^d}z\psi_x(dz),µ)µ(dx)$. Of course, this does not look like a very good idea, namely because by doing so, we would lose all the richness of the extended super-differentials.
 \end{Rem}
 
 We have the following comparison result.
 \begin{Theorem}\label{thm:compstat}
 Assume that $H$ satisfies \eqref{star12}, and for a modulus of continuity $\omega$,
 \be\label{estH:comp}
 \ba
\forall x,y \in \T^d, &q \in \R^d, µ,\nu \in \mptd,\\
& |H(x,q,µ) - H(y,q,\nu)| \leq \omega((1+ |q|)(d(x,y) + \mathcal{W}_2(µ,\nu))).
\ea
 \ee
 Let $U$ and $V$ be respectively viscosity sub and super solution to \eqref{hjbt}, then $U \leq V$.
 \end{Theorem}
 \begin{proof}
 Let $\eps > 0$ and consider $(µ_\eps,\nu_\eps)$ a point of maximum of 
 $$
 (µ,\nu) \mapsto U(µ) - V(\nu) -\frac{1}{2\eps}\mathcal{W}_2^2(µ,\nu).
 $$
 Let $\gamma_\eps$ be an optimal coupling between $µ_\eps$ and $\nu_\eps$. Recalling Proposition \ref{prop:superdifftorus}, for any map $\zeta$ as in the statement of Proposition \ref{prop:superdifftorus}, $(\pi_1,\eps^{-1}\zeta)_\#\gamma_\eps \in \partial^+\mathcal{W}_2^2(\cdot,\nu_\eps)(µ_\eps)\subset\partial^+U(µ_\eps)$ and $(\pi_2,-\eps^{-1}\zeta)_\#\gamma_\eps \in \partial^+\mathcal{W}_2^2(µ_\eps,\cdot)(\nu_\eps)$, so that $(\pi_2,\eps^{-1}\zeta)_\#\gamma_\eps \in\partial^-V(\nu_\eps)$. Using that $U$ and $V$ are respectively sub and super-solution to \eqref{hjbt}, it follows that
 $$
 U(µ_\eps) + \int_{\T^{2d}}H\left(x,\frac{\zeta(x,y)}{\eps},µ_\eps \right)\gamma_\eps(dx,dy) \leq 0
 $$
 and
  $$
 V(\nu_\eps) + \int_{\T^{2d}}H\left(y,\frac{\zeta(x,y)}{\eps},\nu_\eps \right)\gamma_\eps(dx,dy) \geq 0.
 $$
 Taking the difference of the two previous inequalities, we can estimate
 $$
  \max (U-V) \leq \int_{\T^{2d}}H\left(y,\frac{\zeta(x,y)}{\eps},\nu_\eps \right)-H\left(x,\frac{\zeta(x,y)}{\eps},µ_\eps \right)\gamma_\eps(dx,dy).
 $$ Using \eqref{estH:comp}, it follows that
 $$
 \max (U-V) \leq \int_{\T^{2d}} \omega( (1 + \eps^{-1}|\zeta(x,y)|)(d(x,y) + \mathcal{W}_2(µ_\eps,\nu_\eps)))\gamma_\eps(dx,dy).
 $$
 Finally, using that $|\zeta(x,y)| = d(x,y)$ and Jensen's inequality (since $\omega$ is concave), it follows that
 $$
  \max (U-V) \leq \omega\left(2\mathcal{W}_2(µ_\eps,\nu_\eps) + \frac{\mathcal{W}_2^2(µ_\eps,\nu_\eps)}{\eps}  \right).
 $$
 Recalling Lemma \ref{lemma:visc1}, the right hand side goes to $0$ as $\eps \to 0$ and the result follows. 
 \end{proof}
  \begin{Rem}
Note that the only elements in the super-differential which are of interest in the previous proof have a bounded support. So that the condition 
$$
\forall R> 0, \exists C >0, \forall x \in \T^d, q \in B(0,R), µ \in \mptd, \quad |H(x,q,µ)| \leq C
$$
is sufficient to obtain a comparison principle.
 \end{Rem}
 A corollary of this comparison result is the following.
 \begin{Cor}
 If $H$ satisfies \eqref{star12} and \eqref{estH:comp}, there exists at most one viscosity solution to \eqref{hjbt}.
 \end{Cor}
 \begin{proof}
 The comparison principle states that for any two viscosity solutions $U$ and $V$, one has $U^*\leq V_*$, which yields the result by symmetry since $U_* \leq U^*$.
 \end{proof}

One of the most striking properties of viscosity solutions is their stability, which  allows one to  prove quite general existence result, notably by using approximation scheme. Even if such a program for existence of viscosity solutions is not exactly my objective here, I still state a stability result, both to give an idea of what one can expect in general, but also to highlight the particularity of the proof, on which I come back afterwards.

\begin{Theorem}\label{thm:stabvisc}
For $n \geq 1$, let $H_n: \T^d\times \R^d\times \mptd \mapsto \R$ and assume that $H$ satisfies \eqref{star12}, \eqref{estH:comp} as well as that there exists $C>0$ and a sequence $(\eps_n)_{n \geq 1}$ converging toward $0$ such that
\be\label{growthH}
\forall R > 0, \sup_{x,µ,|q| \leq R} |D_qH(x,q,µ)| \leq C(1 + R).
\ee
\be\label{convH}
\forall x \in \T^d, q \in \R^d, µ \in \mptd, \quad |H(x,q,µ)-H_n(x,q,µ)| \leq C(1+|q|^2)\eps_n.
\ee
Let $(U_n)_{n \geq 0}$ be a sequence of continuous viscosity solutions to \eqref{hjbt} when $H$ is replaced by $H_n$. Assume $(U_n)_{n \geq 0}$ converges uniformly to some $U: \mptd \mapsto \R$. Then $U$ is a viscosity solution to \eqref{hjbt}.
\end{Theorem}
\begin{proof}
I only prove that it is the case for the super-solution property, since the other property can be obtained in the same way.
 Let $µ \in \mptd$ and $\psi^* \in \partial^-U(µ)$ and write $\gamma(dx,dz) = µ(dx)\psi^*_x(dz)$. 

From the definition of sub-differentials, there exists a modulus of continuity $\omega$ such that for any $\eta \in \mathcal P(\T^d\times \R^d)$, such that $(\pi_1)_\#\eta = µ$, $\Gamma \in \mathcal P(\T^d\times \R^d\times \R^d)$ such that $(\pi_1,\pi_2)_\#\Gamma = \gamma$, $(\pi_1,\pi_3)_\#\Gamma = \eta$,
\be\label{eq:1290}
U(\exp_\#\eta) - U(µ) - \int_{\T^d\times \R^d\times \R^d}z\cdot v\, \Gamma(dx,dz,dv)  \geq - C_{\T^d}(\eta)\omega( C_{\T^d}(\eta)).
\ee
Without loss of generality, we assume that $\omega(r) \to \infty$ as $r \to \infty$. Let $W_n$ be the map defined by
\[
\begin{aligned}
W_n : & \{ \Gamma \in \mathcal P(\T^d\times \R^d\times \R^d) | (\pi_1,\pi_2)_\#\Gamma = \gamma\} \mapsto \R,\\
&\Gamma \mapsto U_n(\exp(\pi_1,\pi_3)_\#\Gamma) - \int_{\T^d\times\R^{d}\times \R^d}z\cdot v \,\Gamma(dx,dz,dv)\\
& \quad \quad \quad \quad + 2C_{\T^d}((\pi_1,\pi_3)_\#\Gamma)\omega( C_{\T^d}((\pi_1,\pi_3)_\#\Gamma)).
\end{aligned}
\]
Observe that $W_n$ is coercive in the sense that $W_n(\Gamma) \to \infty$ superlinearly as $M_2((\pi_3)_\#\Gamma) \to \infty$. Hence using Proposition \ref{cor2s}, we obtain that there exist $\beta_n \in \mathcal P^\gamma_2((\T^d\times \R^d)\times \R^d)$, $M_2((\pi_3)_\#\beta_n) \leq \frac{1}{n^2}$ such that $W_n(\Gamma) + \mathcal I^\gamma_2(\Gamma,\beta_n)$ has a point of minimum $\Gamma_n$. It follows that, noting $µ_n = \exp(\pi_1,\pi_3)_\#\Gamma_n$ and $\eta_n = (\pi_1,\pi_3)_\#\Gamma_n$,
\[
\ba
U_n(µ_n) &- \int_{\T^d\times \R^d\times \R^d}z\cdot v \,\Gamma_n(dx,dz,dv) + 2 C_{\T^d}(\eta_n)\omega( C_{\T^d}(\eta_n))\\
& \leq U_n(µ) - U(µ) + U(µ)  - \mathcal I^\gamma_2(\Gamma_n,\beta_n)\\
& \leq U_n(µ) -U(µ) + U(µ_n) - \int_{\T^d\times \R^d\times \R^d}z\cdot v \, \Gamma_n(dx,dz,dv) - \mathcal I^\gamma_2(\Gamma_n,\beta_n)\\
&\quad\quad +C_{\T^d}(\eta_n)\omega( C_{\T^d}(\eta_n))
\ea
\]
Thus, passing to the limit $n \to \infty$, we deduce that $ C_{\T^d}(\eta_n) \to 0$ because $U_n(µ) \to U(µ)$, $\mathcal I_2((\pi_3)_\#\Gamma_n,\nu_n) \to 0$ and $\sup_{µ'}|U(µ')- U_n(µ')|  \to_{n \to \infty}0$. The aim of the rest of the proof is to establish the following: i) at the point $\Gamma_n$, we can exhibit an element in the sub-differential of $U_n$, ii) using the previous convergence, we can show that it converges toward $µ(dx)\psi^*_x(dz)$, iii) this convergence yields the required viscosity solution property.

To lighten notation, I use the formalism of random variables, on a standard probability space $(\Omega,\mathcal A,\mathbb P)$. Let $(X_n,Z_n,V_n)$ be of law $\Gamma_n$, $\zeta_n$ be such that $((X_n,Z_n,V_n),(X_n,Z_n,\zeta_n))$ is optimal for $\mathcal I^\gamma_2(\Gamma_n,\beta_n)$ and observe that for any $V$ random variable in $\R^d$
$$
\ba
&U_n(\mathcal L(\exp(X_n, V_n))) - \mathbb E[Z_n\cdot V_n] + 2 \sqrt{\mathbb E[|V_n|^2]}\omega\left(\sqrt{\mathbb E[|V_n|^2]}\right) + \mathbb E[\zeta_n\cdot V]\\
& \leq U_n(\mathcal L(\exp(X_n, V_n+V)))- \mathbb E[Z_n\cdot (V_n + V)] + 2 \sqrt{\mathbb E[|V_n + V|^2]}\omega\left(\sqrt{\mathbb E[|V_n + V|^2]}\right)
\ea
$$
Since, on $\T^d$, $\exp(x,v + v') = \exp(\exp(x,v),v')$, we obtain that 
$$
\ba
&U_n(\mathcal L(\exp(X_n, V_n))) + \mathbb E\left[\left(Z_n + \zeta_n - 2\frac{V_n}{\|V_n\|_2}(\omega(\|V_n\|_2) + \omega'(\|V_n\|_2)\|V_n\|_2)\right)\cdot V\right]  + o(\|V\|_2)\\
& \leq U_n(\mathcal L(\exp(\exp(X_n, V_n),V))),
\ea
$$
where $\omega'$ is the derivative of $\omega$. In particular, we obtain that
$$
\mathcal L \left( \exp(X_n, V_n),Z_n + \zeta_n - 2\frac{V_n}{\|V_n\|_2}(\omega(\|V_n\|_2) + \omega'(\|V_n\|_2)\|V_n\|_2)\right) \in \partial^-U_n(µ_n).
$$
Denoting the random variable whose law is taken on the left-hand side by $(Y_n,Z_n + \xi_n)$, we know that $\|\xi_n\|_2 \to 0$ and $\mathcal L(Y_n,Z_n)(dx,dz) \to µ(dx)\psi^*_x(dz)$ as $n \to \infty$. On the other hand, since $U_n$ is a viscosity super-solution, we know that
\be\label{eq:4032}
U_n(µ_n) + \mathbb E[H_n(\exp(X_n,V_n),Z_n + \xi_n, µ_n)] \geq 0.
\ee
Using now the two estimates we have on the Hamiltonian, we deduce 
$$
\ba
|H(y,z+ q,\nu) - H(x,z,µ)| \leq &\sup_{|z'| \leq |z| + |q|}|D_qH(y,z',\nu)| |q|\\
& + \omega((1+ |z|)(d(x,y) + \mathcal{W}_2(µ,\nu))).
\ea
$$
We can thus pass to the limit in \eqref{eq:4032} thanks to \eqref{convH} to obtain the required
$$
U(µ) + \int_{\T^d\times \R^d} H(x,z,µ)µ(dx)\psi^*_x(dz) \geq 0.
$$
\end{proof}
\begin{Rem}
This proof of the stability of viscosity solutions is in fact a proof of the stability of sub/super differentials.
\end{Rem}
Note that a requirement on the continuity of $H$ in its second argument is needed in this kind of stability argument. It is classical and I do not know general proofs which work without it. Usually, stability proofs in viscosity solution theory are simpler, namely because of the equivalent definition of viscosity solution involving smooth test functions. In this case, because we need to work with sub/super-differentials, the proof is more technical (even though it still follows classical arguments and does not require additional assumptions).
 
\begin{Rem}
In this section, I made the choice of the $2$-Wasserstein distance. Since we are on a compact set, all Wasserstein distances give the same topology, however, their choice plays a role in some quantitative results or in the definition of super-differential for instance. All the arguments presented here can easily be extended to the case $p \in (1,\infty)$, by using $\mathcal W_p^p$ instead of $\mathcal W_2^2$.
\end{Rem}
 
\subsection{Viscosity solutions on Wasserstein spaces: the non-compact case}\label{sec:noncomp}

 As is slightly apparent in the previous section, integrability conditions of elements of the super-differential, as well as growth conditions on the Hamiltonian, play an important role in the theory of viscosity solutions on spaces of probability measures. The same goes in finite dimension for the growth of the Hamiltonian. In the non-compact case that I will now treat, such aspects are even more important. Here, I will focus on a particular setting: an equation set on $\mpprd$ for $p \in (1,\infty)$, with a time dependent problem and a Hamiltonian with an appropriate growth. The time dependent case is presented, as opposed to the stationary case of the previous section, simply for the sake of completeness and does not play a strong role in the following. However, in order to avoid quite challenging technicalities due to the presence of the boundary condition in time (the initial condition), the following results are restricted to continuous viscosity solutions.\\
 
 The main equation of interest is here
 \be\label{hjbtime}
 \partial_t U + \int_{\R^d}H(x,D_µU(t,µ)(x),µ)µ(dx) = 0 \text{ in } (0,T)\times \mpprd,
 \ee
 together with initial condition 
 $$
 U|_{t = 0} = U_0.
 $$
 The Hamiltonian $H$ satisfies for some $C>0$
 \be\label{ff7}
 \forall x,q\in \R^d, µ \in \mprd, \quad |H(x,q,µ)| \leq C(1+ |q|^{p'}),
 \ee
so that the equation is well defined on elements of the sub/super-differentials, together with the usual continuity requirement 
 \be\label{estH:comp2}
 \ba
\forall x,y \in \R^d, &\,q \in \R^d, µ,\nu \in \mpprd,\\
& |H(x,q,µ) - H(y,q,\nu)| \leq \omega((1+ |q|)(|x-y| + \mathcal{W}_p(µ,\nu))),
\ea
 \ee
 and an additional estimate on the growth of $D_qH$ to handle perturbation terms
 \be\label{estH:dq}
 \ba
 \forall x \in \R^d, &\,q \in \R^d, µ \in \mpprd,\\
& |D_qH(x,q,µ)| \leq C(1+ |q|^{p'-1}).
 \ea
 \ee
 Note that convexity in $q$ of $H$ would imply this last requirement, but I do not assume convexity here.
 The length $T > 0$ of the time interval is arbitrary.
 \begin{Def}\label{def:visctcb}
  A usc (resp. lsc) function $U: [0,T]\times \mpprd \mapsto \R$ is said to be a viscosity sub-solution (resp. super-solution) to \eqref{hjbtime} if for any $t \in (0,T], µ \in \mpprd, (\theta,\psi) \in \partial^+U(t,µ)$, (resp. $(\theta,\psi) \in \partial^-U(t,µ)$) it holds that
 \be\label{visc2}
 \theta + \int_{\R^{2d}}H(x,z,µ)\psi_x(dz)µ(dx) \leq 0\quad  (\emph{resp. } \geq 0).
 \ee
 A locally bounded function $U: [0,T]\times \mpprd \mapsto \R$ is a viscosity solution to \eqref{hjbtime} if $U^*$ is a sub-solution and $U_*$ a super-solution.
 \end{Def}
 In particular a continuous function is a viscosity solution if it is both a sub and a super-solution. Furthermore, super-differentials of functions of $[0,T]\times \mpprd$ are defined as follows: for $U : [0,T]\times \mpprd \mapsto \R$, we write $\R\times (\R^d \mapsto \mathcal P(\R^d))\ni(\theta,\psi) \in \partial^+U(t,µ)$ if there exists a modulus of continuity $\omega$ such that for all $(t',µ') \in [0,T]\times \mpprd$, for all $\Gamma \in \Pi(µ(dx)\psi_x(dz),µ')$,
 $$
 \ba
 U(t',µ') \leq& U(t,µ) + \theta(t'-t) + \int_{\R^{3d}}z\cdot(y-x)\Gamma(dx,dz,dy)\\
 & + (|t'-t| + \tilde C_p(\Gamma))\omega(|t'-t| + \tilde C_p(\Gamma)),
 \ea
 $$
 together with 
 $$
 \int_{\R^{2d}}|z|^{p'}µ(dx)\psi_x(dz) < \infty.
 $$
The following comparison principle holds.
 \begin{Theorem}\label{thm:compt}
Let $p \in (1,\infty)$. Assume that \eqref{ff7}, \eqref{estH:comp2} and \eqref{estH:dq} hold. Let $U,V \in \emph{BUC}([0,T]\times \mathcal P_p(\R^d))$ be respectively viscosity sub and super-solution to \eqref{hjbtime} such that $(U-V)|_{t = 0} \leq 0$. Then $U \leq V$.
 \end{Theorem}
 Recall that $BUC$ is the space of bounded uniformly continuous functions. The proof of this result is similar to the one of Theorem \ref{thm:compstat}, but has to deal with three additional difficulties: the time variable, the non-compactness of $\mpprd$ and the integrability condition of the elements of the super-differential which has to be in agreement with the growth of the Hamiltonian. Furthermore, an estimate of the form of \eqref{estH:dq} is required to deal with the various perturbations introduced in the proof. It is also important to note that this comparison result holds for sub and super-solutions which are in BUC. Such a regularity may seem arbitrary at the moment but I comment on it after the proof. 
 \begin{proof}[Proof of Theorem \ref{thm:compt}]
The following proof can seem quite heavy in terms of notation. I advise the unexperienced reader to make all the necessary simplifications resulting of the case $p=2$. For $\eps,\delta > 0$, define $$\Phi_{\eps,\delta}(t,s,µ,\nu):= U(t,µ) - V(s,\nu) - \frac{1}{2\eps}((t-s)^2 + \frac2p\mathcal{W}_p^p(µ,\nu)) - \delta \mathcal{W}_p^p(µ,\delta_0).$$ Assume that the result does not hold. Then, there exists $\lambda, \kappa,\bar \delta > 0$ such that, for any $\eps> 0$, $\delta \in (0,\bar \delta)$
 \be\label{eq:2000}
 \sup_{µ,\nu \in \mpprd\\t,s \in [0,T]} \Phi_{\eps,\delta} - \lambda t  \geq \kappa.
 \ee
 From Corollary \ref{cor1s}\footnote{I omit here the Stegall like perturbation in the $t$ variable as it is standard and does not raise any difficulty. Moreover, even if it is not done here, for viscosity solutions of equations of this type, we can always omit the doubling of variable in the time argument, as the term in $\partial_t$ is linear and constant.}, for any $\eps,\delta,\eta > 0$, there exists $\xi_1,\xi_2 \in \mathcal P_{p'}(\R^d), M_{p'}(\xi_1),M_{p'}(\xi_2) < \eta$, such that 
 $$
 (t,s,µ,\nu) \mapsto \Phi_{\eps,\delta}(t,s,µ,\nu) - \lambda t - \mathcal I_p(µ,\xi_1) - \mathcal I_p(\nu,\xi_2)
 $$
 has a point of strict maximum at some point $(t_{\eps,\delta,\eta},s_{\eps,\delta,\eta},µ_{\eps,\delta,\eta},\nu_{\eps,\delta,\eta})$.\\
 
 Assume first that
 \be\label{tspos}
 t_{\eps,\delta,\eta}\,,\,s_{\eps,\delta,\eta} > 0.
 \ee
 
 Consider the couplings $\gamma_\eps, \gamma^1_\delta(dx,dz) := µ_{\eps,\delta,\eta}(dx)\psi_x(dz)$ and $\gamma^2_\delta(dy,dz):= \nu_{\eps,\delta,\eta}(dy)\phi_y(dz)$, optimal in respectively $\mathcal W_p^p(µ_{\eps,\delta,\eta},\nu_{\eps,\delta,\eta}), \mathcal I_p(µ_{\eps,\delta,\eta},\xi_1)$ and $\mathcal I_p(\nu_{\eps,\delta,\eta},\xi_2)$.  From Propositions \ref{prop:superdiffWp}, \ref{prop:superI} and \ref{prop:propsuper}, we deduce the existence of elements in the super-differential of $U$ at $(t_{\eps,\delta,\eta},µ_{\eps,\delta,\eta})$ and in the sub-differential of $V$ at $(s_{\eps,\delta,\eta},\nu_{\eps,\delta,\eta})$ which are given by 
 
 $$
(\lambda + \eps^{-1}(t-s), (\pi_1,\eps^{-1}(\pi_1-\pi_2)|\pi_1-\pi_2|^{p-2} + p \delta \pi_1|\pi_1|^{p-2} - \pi_3)_\#\gamma_\eps(dx,dy)\psi_x(dz)) \in \partial_µ^+U(t_{\eps,\delta,\eta},µ_{\eps,\delta,\eta}),
 $$
  $$
 (\eps^{-1}(t-s),(\pi_2,\eps^{-1}(\pi_1-\pi_2) |\pi_1-\pi_2|^{p-2} + \pi_3)_\#\gamma_\eps(dx,dy)\phi_y(dz)) \in \partial_µ^-V(s_{\eps,\delta,\eta},\nu_{\eps,\delta,\eta}).
 $$
  Using the fact that $U$ and $V$ are respectively viscosity sub and super-solution, we deduce that
 $$
 \lambda + \frac{t-s}{\eps} + \int_{\R^{3d}}H\left(x,\frac{x-y}{\eps|x-y|^{2-p}} + p\delta x |x|^{p-2} - z,µ_{\eps,\delta,\eta}\right)\gamma_\eps(dx,dy)\psi_x(dz) \leq 0,
 $$
 as well as,
 $$
 \frac{t-s}{\eps} + \int_{ \R^{3d}}H\left(y,\frac{x-y}{\eps|x-y|^{2-p}} + z',\nu_{\eps,\delta,\eta}\right)\gamma_\eps(dx,dy)\phi_y(dz') \geq 0.
 $$
 Taking the difference of the two previous inequalities yields
 \be\label{star23}
 \ba
 \lambda \leq &\int_{\R^{4d}}H\left(y,\frac{x-y}{\eps|x-y|^{2-p}} + z',\nu_{\eps,\delta,\eta}\right)\\
 & -H\left(x,\frac{x-y}{\eps|x-y|^{2-p}} + p\delta x|x|^{p-2} - z,µ_{\eps,\delta,\eta}\right)\gamma_\eps(dx,dy)\psi_x(dz)\phi_y(dz').
 \ea
 \ee
 To estimate the integral of the difference of the two Hamiltonians, we use the assumptions \eqref{estH:comp2} and \eqref{estH:dq} to obtain
 $$
 \ba
 H\left(y,\frac{x-y}{\eps|x-y|^{2-p}} + z',\nu\right)& -H\left(x,\frac{x-y}{\eps|x-y|^{2-p}} + p\delta x|x|^{p-2} - z,µ\right)\\
 &= H\left(y,\frac{x-y}{\eps|x-y|^{2-p}} + z',\nu\right) -H\left(y,\frac{x-y}{\eps|x-y|^{2-p}},\nu\right)\\
 &\quad+H\left(y,\frac{x-y}{\eps|x-y|^{2-p}},\nu\right)-H\left(x,\frac{x-y}{\eps|x-y|^{2-p}},µ\right)\\
&\quad + H\left(x,\frac{x-y}{\eps|x-y|^{2-p}},µ\right) -H\left(x,\frac{x-y}{\eps|x-y|^{2-p}} + p\delta x|x|^{p-2} - z,µ\right)\\
& \leq C\left(1+ \eps^{-1}| x-y|^{p-1} + |z'|\right)^{p'-1}|z'|\\
&\quad + \omega\left( \left(1 + \eps^{-1}|x-y|^{p-1}\right)(|x-y|+ \mathcal{W}_p(µ,\nu))\right)\\
& \quad +C\left(1+ \eps^{-1}\left|x-y\right|^{p-1} + |z| + p\delta|x|^{p-1} \right)^{p'-1}(|z| + p \delta |x|^{p-1}).
 \ea
 $$
Injecting this inequality in \eqref{star23} and defining $\Gamma(dx,dy,dz,dz'):=\gamma_\eps(dx,dy)\psi_x(dz)\phi_y(dz')$ yields
$$
\ba
\lambda \leq& C \int_{\R^{4d}} \left(1+ \frac{|x-y|^{p-1}}{\eps}  + |z| + p\delta|x|^{p-1} + |z'|\right)^{p'-1}(|z| + p \delta |x|^{p-1} + |z'|)\Gamma(dx,dy,dz,dz')\\
&+\int_{\R^{4d}} \omega\left( \left(1 + \frac{|x-y|^{p-1}}{\eps}\right)(|x-y|+ \mathcal{W}_p(µ_{\eps,\delta,\eta},\nu_{\eps,\delta,\eta}))\right) \Gamma(dx,dy,dz,dz')
\ea
$$
Since $\omega$ is concave, by Jensen's inequality, we can take it outside of the integral and use the same estimate as in the compact case to bound the second integral by
$$
\omega\left ( \mathcal{W}_p(µ_{\eps,\delta,\eta},\nu_{\eps,\delta,\eta}) + \frac{\mathcal{W}_p^p(µ_{\eps,\delta,\eta},\nu_{\eps,\delta,\eta})}{\eps}\right) = o(1).
$$
Using H\"older's inequality, $(p'-1)p = p'$ and $(p'-1)(p-1)p = p$ on the first integral finally yields, for $K_p$ a constant depending on $p$ (and on $H$)
 $$
 \ba
&\lambda \leq \omega(\mathcal{W}_p(µ_{\eps,\delta,\eta},\nu_{\eps,\delta,\eta}) + \eps^{-1}\mathcal{W}_p^p(µ_{\eps,\delta,\eta},\nu_{\eps,\delta,\eta}))\\
& +  K_p\left( \eta + \delta^{p'}M_p(µ_{\eps,\delta,\eta})\right) \\
&+ K_p\left(1 +\frac{\mathcal{W}_p(µ_{\eps,\delta,\eta},\nu_{\eps,\delta,\eta})}{\eps^{p'-1}}\right)\left(\left(\int_{ \R^{4d}}(|z|^{p'}+|z'|^{p'})\Gamma(dx,dy,dz,dz')\right)^\frac{1}{p'} +  \delta \left(M_p(µ_{\eps,\delta,\eta})\right)^{\frac{1}{p'}}\right).
 \ea
 $$
 Observe that $\delta M_p(µ_{\eps,\delta,\eta})$ is bounded, and thus, thanks to Lemma \ref{lemma:visc2}, we can pass to the limit (in this order) $\eta \to 0, \delta\to 0, \eps\to 0$, to obtain $\lambda \leq 0$, which is not possible.\\

So, we deduce from the previous contradiction that \eqref{tspos} does not hold for sufficiently many $\eps,\delta,\eta$, in particular there exists a sequence $(\eps_n,\delta_n,\eta_n)_{n \geq 0}$, which goes toward $0$, such that either $t_n= 0$ or $s_n =0$, where $(t_n,s_n,µ_n,\nu_n)$ is the associated point of maximum. Recalling \eqref{eq:2000}, we then contradict $(U-V)|_{t = 0} \leq 0$ by passing to the limit $n \to \infty$. Indeed, consider for instance the case $s_n = 0$. It then follows that 
$$
U(t_n,µ_n) - V(0,\nu_n) \geq \frac{\kappa}{2},
$$
as well as $\mathcal{W}_p^p(µ_n,\nu_n) + (t_n-s_n)^2 \leq \eps_n$. Hence, using the uniform continuity of $U$, we deduce that for a certain modulus of continuity $\tilde \omega$, 
$$
U(0,\nu_n)-V(0,\nu_n) \geq \frac \kappa 2 - \tilde\omega(\sqrt{\eps_n}),
$$
which also leads to a contradiction. Hence the result.
 \end{proof}

 \begin{Rem}
 The assumption that $U$ and $V$ are in BUC is extremely helpful to deal with two problems at hand: i) the equation is posed on a non-compact space, ii) the problem is time dependent. It is in fact a very insightful exercise for the interested reader to verify that when one of these two problems is not present, one can simply assume that $U$ is usc and $V$ lsc. Nonetheless, here a difficulty truly exists. In practical applications, it is treated on a case by case basis, either by establishing that natural candidates for being viscosity solutions are in BUC, or by proving more information on the behavior of these candidates near infinity in space, to compactify the problem. Also, in certain situations, we can hope to prove the uniform continuity of the solution only close to the initial time (using a Cauchy-Lipschitz type result on small time).
 \end{Rem}

\subsection{Bibliographical comments}
The theory of viscosity solutions has been introduced by Crandall and Lions in \cite{crandall1983viscosity} to address first order Hamilton-Jacobi equations in finite dimension. This theory turned out to be the right one for such equations. I do not recall all the long literature on this topic but simply point out what I think are key elements. A major difficulty in the theory of viscosity solutions is the treatment of second order equations. A breakthrough was made by Jensen in \citep{jensen} and a complete presentation of the basics on viscosity solutions can be found in the survey of Crandall, Ishii and Lions \cite{crandall1992user}. Perhaps the two main studies which are absent from the previous (must-read) survey are the questions of growth of solutions in unbounded domains, for which I refer to Ishii \cite{ishii} or to Alvarez \cite{alvarez} in convex cases, and the half-relaxed limit of Barles and Perthame \cite{barlesperthame} their use in Barles and Souganidis \cite{barles1991convergence} which are key illustrations of the powerful stability of viscosity solutions.

Crandall and Lions published a series of papers \cite{crandall1,crandall2,crandall3,crandall4,crandall5,crandall6,crandall7} to provide a thorough study of HJB equations in Banach spaces, see also the book by Fabbri, Gozzi and Swiech \cite{fabbri} for more on this subject.

For HJB equations on sets of probability measures, the notion of viscosity solutions presented here (which I believe is the right one) stems from my work \citep{bertucci2024stochastic}. In the non-compact case, the comparison principle was done by Lions and me in \citep{bertucci2026stegall}, see also our work with Ceccherini Silberstein for more details on such a notion \citep{bertucci2025doubling}. Prior to that, several other notions were proposed. Lions proposed in \citep{lions20082009} to define a notion of viscosity solution through its lift, which works well when we are on $\mpt$, because the latter is then lifted into the Hilbert space of square-integrable random variables. Another notion, somehow closer to the one presented here, focused only on deterministic elements of the super-differential. In such a way, comparison principles can only be obtained for particular Hamiltonians, such as it was the case in Cardaliaguet and Quincampoix \citep{cardaquincampoix}, Marigonda and Quincampoix \citep{marigonda2018mayer} or in Jimenez, Marigonda and Quincampoix \citep{jimenez2020optimal}. The link between this last notion and the one of Lions (using the lift) was studied by Gangbo and Tudorascu in \citep{gangbotudorascu}.

Other notions of viscosity solutions were also considered in cases in which the HJB equations have additional singular terms but I do not comment on it here.

When the Hamiltonian $H$ is convex, using techniques of Alvarez \citep{alvarez}, we provided with Lions a more general version of a comparison result in the non-compact case in \citep{bertucci2026stegall}.

\newpage
 
\section{Regularity and viscosity properties of the value function, existence of viscosity solutions}\label{sec:existvisc}
In this section, I show how to establish various properties on the value function of mean field optimal control problems. In particular, I focus on regularity estimates and on the fact that the value function is a viscosity solution to the associated HJB equation. I shall not treat general questions of existence of solutions, i.e. for equations which are not given by an optimal control problem. For such cases, I give a simple existence result for equations that can be lifted on the space of random variables.

I will focus on the case of $\mptd$ for most of this section, mainly for two reasons: i) it allows for simpler notation than in the non-compact case, ii) it forces us to adopt a formalism on the space of measures as it prohibits simply lifting the equation on a Banach space of random variables with appropriate integrability conditions.

 Consider $U$ defined for $t \in [0,T), µ \in \mptd$ by
\be\label{defU3}
U(t,µ) := \inf_{(b,\tilde m)}\left\{\int_t^T \int_{\T^d\times K}L(x,b(s,x,k),(\pi_1)_\#\tilde m_s)\tilde m_s(dx,dk)ds + \mathcal G ((\pi_1)_\#\tilde m_T)\right\},
\ee
where the infimum is taken as in \eqref{defU}. Recalling Section \ref{sec:derivHJB}, an equivalent definition of $U$ through minimization over stochastic processes is also available. I shall make assumptions on $L$ and $\mathcal G$ incrementally, to highlight what is used in each proof. For the moment, I simply assume that $L, \mathcal G \geq 0$, and that $\mathcal G$ is bounded (and that they are measurable...). So that $U$ is well defined for every $t,µ$.\\

First I will establish some regularity estimates on $U$: boundedness, continuity, Lipschitz continuity. Then, I will pass to the verification of the fact that $U$ is indeed a viscosity solution to 
\be\label{hjb3}
-\partial_t U + \int_{\T^d}H(x,D_µU,µ)µ(dx) = 0 \text{ in } (0,T)\times \mptd,
\ee
with terminal condition $U(T,\cdot) = \mathcal G$. Stationary equations could also have been treated using arguments similar to the ones which are to follow.

 \subsection{Regularity estimates on the value function}
 \subsubsection{Boundedness}
 For this first case, it is quite immediate to verify the boundedness property. When dealing with problems which are stochastic and have terminal constraints, as in Section \ref{sec:stoot}, such bounds will turn out to be more difficult to obtain.
  \begin{Prop}\label{prop:boundeasy}
 For any $t \in [0,T],µ \in \mptd$,
 $$
 \inf \mathcal G \leq U(t,µ) \leq \mathcal G(µ) + (T-t)\int_{\T^d}L(x,0,µ)µ(dx).
 $$
 \end{Prop}
 \begin{proof}
 The first property follows from the non-negativity of $L$, the second from the sub-optimality of the control which consists in not moving.
 \end{proof}

 \subsubsection{Continuity through controllability}\label{sec:contcontr}
 
 Generally speaking, there are two main ways to prove that a value function is continuous. The first one relies on the controllability of the problem, while the second involves the use of the "same" control starting from two different measures. By controllability, it is understood that the controller can drive the state of the system, in any non-zero time interval, wherever it wants. This section presents the first approach and the second one is then presented in Section \ref{sec:parcont} just below.\\
 
 A first technical result consists in remarking that, in these kinds of problems, it is always almost optimal to wait a little bit before starting to make decisions. It is a time regularity result on the value function.
  \begin{Lemma}\label{lemma:density}
Assume that 
\be\label{hyp:dpl}
 \exists C >0, \forall x \in \T^d,z\in \R^d, µ \in \mptd,\quad D_zL(x,z,µ)\cdot z \leq C (1 + L(x,z,µ)).
\ee
Then, for $t < T, µ \in \mptd$ 
$$
\lim_{ \eps \to 0^+}U(t+\eps,µ) = U(t,µ).
$$
Moreover, for all $\eps > 0$ there exists an $\eps$-optimal control $(b,\tilde m)$ for $U(t,µ)$ such that $b \equiv 0$ on $[t, t + C \eps]$ with $C >0$ depending only on $U(t,µ)$ and $T-t$ (and on the data of the problem).
\end{Lemma}
The following proof is quite heavy in terms of notation. It might be interesting for the reader to first substitute for $L$ a function of the form $|z|^k$ for some $k >1$.
\begin{proof}
Consider the admissible control of not moving for a time $\eps$, we deduce that $U(t,µ) \leq \eps \int_{\T^d}L(x,0,µ)µ(dx) + U(t+\eps,µ)$, hence $U(t,µ) \leq \liminf_{ \eps \to 0^+}U(t+\eps,µ)$. 

 Consider $n \geq 1, \eps < T-t$ and a $n^{-1}$-optimal control $(b,\tilde m)$ for $U(t,µ)$. For $s \geq t + \eps$, define $(b',\tilde m')$ by 
 \[
 \ba
b'(s,\cdot) &= \frac{T-t}{T-t-\eps}b(\phi(s),\cdot)\\
\tilde m'_s &= \tilde m_{\phi(s)},
 \ea
 \]
where $\phi(s) = \frac{T-t}{T-t-\eps}(s-t - \eps) + t$. For $s \in [t,t+ \eps]$, set $b'(s,\cdot) = 0, \tilde m'_s = µ$. Evaluating the cost of such a pair yields

\[
\ba
 \int_t^T \int_{\T^d\times K}L&(x,b'(s,x,k),(\pi_1)_\#\tilde m'_s)\tilde m'_s(dx,dk)ds\\
 = (\phi')^{-1} &\int_{t}^T \int_{\T^d\times K}L\left(x,\phi'b(s,x,k),(\pi_1)_\#\tilde m_s\right)\tilde m_{s}(dx,dk)ds + \int_t^{t + \eps}\int_{\T^d}L(x,0,µ)µ(dx)ds\\
 = (\phi')^{-1}& \int_{t}^T \int_{\T^d\times K}L\left(x,b(s,x,k),(\pi_1)_\#\tilde m_s\right)\tilde m_{s}(dx,dk)ds+ \eps \int_{\T^d}L(x,0,µ)µ(dx)\\
 + (\phi')^{-1}&\int_{t}^T \int_{\T^d\times K}(L\left(x,\phi'b(s,x,k),(\pi_1)_\#\tilde m_s\right) - L(x,b(s,x,k),(\pi_1)_\#\tilde m_s))\tilde m_{s}(dx,dk)ds
 \ea
 \]
 It now remains to bound the difference of the two terms in $L$. From the bound \eqref{hyp:dpl}, it is an elementary analysis exercise, which requires at some point Gronwall's Lemma, to obtain
 $$
\exists M > 0, \forall x,z,µ,\, L(x, \phi' z ,µ) - L(x,z,µ) \leq M (1+ L(x,z,µ))(\phi'-1).
 $$
 Hence we deduce that the previous equality can be estimated by
 \[
 \ba
  \leq (\phi')^{-1}&(U(t,µ) +n^{-1} ) + \eps \int_{\T^d}L(x,0,µ)µ(dx)\\
  &+ \frac{\phi' -1}{\phi'}M \int_{t}^T \int_{\T^d\times K}( 1 +L(x,b(s,x,k),(\pi_1)_\#\tilde m_s))\tilde m_{s}(dx,dk)ds\\
   \leq (\phi')^{-1}&(U(t,µ) +n^{-1} ) + \eps \int_{\T^d}L(x,0,µ)µ(dx)+ \frac{\phi' -1}{\phi'}M((T-t) + U(t,µ) + n^{-1}) 
\ea
\]
where I have used the $n^{-1}$-optimality of $(b,\tilde m)$ in the last inequality. Since $\phi' \to 1$ as $\eps \to 0$, the first part of the claim follows. The second part can be observed simply by remarking that $$|\phi'- 1| \leq C \eps$$ as $\eps \to 0$, for $C$ depending only on $T-t$.
\end{proof}
\begin{Rem}
Here and in several proofs below, I use almost optimal controls, whereas it can in fact be proven that optimal controls exist under mild assumptions. The current presentation has the virtue to be easily extendable to stochastic cases, where the existence of optimal controls is more involved, or to zero-sum games. In any case, the important question of the existence of optimal controls is not addressed here.
\end{Rem}
The growth condition \eqref{hyp:dpl} implies in fact more than the inequality used in the previous proof. For instance, it implies that there exist constants $C > 0$, $k > 1$ such that
\be\label{eq:hyppoly}
\forall x \in \T^d, z \in \R^d, µ \in \mptd, \quad L(x,z,µ) \leq C( 1 + |z|^k),
\ee
where $k$ is the constant $C$ in \eqref{hyp:dpl} and $C$ in the previous estimate is another constant. 

The previous result of time regularity allows me to prove the following result of regularity in space.
\begin{Lemma}\label{lemma:contmu}
Assume that \eqref{hyp:dpl} holds. Then, for any $t < T$, there exists $C$ depending only on $\|U\|_\infty$ and $T-t$ such that
\[
\forall µ,µ' \in \mptd, \quad|U(t,µ) - U(t,µ')|\leq C \mathcal{W}_k(µ,µ'),
\]
where $k$ is the constant in \eqref{eq:hyppoly}.
\end{Lemma}
\begin{proof}
Take $t < T$ and $µ,µ' \in \mptd$. Thanks to Lemma \ref{lemma:density}, consider an $\eps$-optimal control $(b,\tilde m)$ for $U(t,µ)$ which is admissible for $U(t + C\eps,µ)$, for $C > 0$ depending only on $T-t$, $\|U\|_\infty$. Consider now the pair $(b',\tilde m')$ defined by 
\[
(b',\tilde m') = (b,\tilde m) \text{ for } s \in [t+C\eps,T),
\]
and $(b',\tilde m')_{s \in [t,t+ C\eps]}$ corresponds to an optimal trajectory for the deterministic optimal transport of $µ'$ toward $µ$ in time $C\eps$ for the $\mathcal{W}_k$ distance. It then follows from the $\eps$-optimality of $(b,\tilde m)$ that
\[
\ba
U(t,µ') &\leq U(t,µ) + \eps + \int_t^{t + C\eps} \int_{\T^d\times K}L(x,b'(s,x,k),(\pi_1)_\#\tilde m'_s)\tilde m'_s(dx,dk)ds\\
& \leq   U(t,µ) + \eps + K(C\eps +  \mathcal{W}_k^k(µ,µ')(C\eps)^{1-k}),
\ea
\]
where I have used \eqref{eq:hyppoly} in the last inequality. Taking $\eps = \mathcal{W}_k(µ,µ')$ yields the required result, if we restrict our attention to $µ',µ$ at $\mathcal{W}_k$ distance at most $(T-t)/C$ (because we cannot choose $\eps$ too large). Thus we have a local Lipschitz estimate, with a constant which depends only on time. The triangle inequality then implies that this Lipschitz estimate is in fact global.

\end{proof}
\begin{Rem}
Since the terminal cost $\mathcal G$ is not assumed to be continuous, it is natural that the estimate above deteriorates as $t \to T$. In fact, the previous result can be understood as a regularizing effect, even though such terminology is in general reserved for results proved at the level of the PDE, whereas here it is done at the level of the optimal control problem. It is a general feature of controllable problems, to have a Lipschitz continuous value function, and it is also well known that HJB equations have regularizing effects.
\end{Rem}

The previous leads easily to global Lipschitz continuity of the value function.

\begin{Prop}\label{prop:globalcont}
Assume that \eqref{hyp:dpl} holds. Then, for any $\eps > 0$, there exists $C> 0$ such that for all $0 \leq t,t' \leq T-\eps, µ,µ' \in \mptd$,
\be\label{eq:globalcont}
|U(t,µ) - U(t',µ')| \leq C (|t-t'| + \mathcal{W}_k(µ,µ')).
\ee
\end{Prop}
\begin{proof}
Consider $(t,µ), (t',µ') \in [0,T) \times \mptd$ and assume that $t' > t$. Considering any control for $U(t,µ)$ which does nothing in $[t,t']$, it follows that $U(t,µ) \leq U(t',µ) + (t'-t)\int_{\T^d}L(x,0,µ)dµ$. Thanks to Lemma \ref{lemma:density}, there exists a $C(t'-t)$-optimal control for $U(t,µ)$ which does nothing in the time interval $[t,t']$, for $C$ depending only on $\eps$ and $U(t,µ,w)$, i.e.
 \[
 U(t,µ) \geq U(t',µ) - C(t'-t).
 \]
Hence, Lipschitz regularity in $t$ is proven, and the result follows.
\end{proof}
 
 \subsubsection{Lipschitz regularity through parallel transport}\label{sec:parcont}
 I present here another way to obtain the regularity of the value function. This second method is, in my opinion, more elegant than the previous one, as it does not require a lot of playing around with the control problem. It will also make more apparent the use of mixed flows. The main result is the following.
 
\begin{Prop}
Assume that there exists $C > 0, k > 1$ such that $L$ and $\mathcal G$ satisfy
$$
\ba
&\quad \forall µ,\nu, \quad |\mathcal G(µ)-G(\nu)| \leq C\mathcal{W}_k(µ,\nu),\\
&\forall x,z,µ, \quad C^{-1}|z|^k - C \leq L(x,z,µ) ,\\
&\forall x,y,z,µ,\nu, \quad |L(x,z,µ) - L(y,z,\nu)| \leq C(1 + |z|^{k-1})(d(x,y) + \mathcal{W}_k(µ,\nu))
\ea
$$
Then, there exists $C_0$ depending only on $C$ and $T$ such that for all $t\in [0,T],µ,\nu \in \mptd$,
$$
|U(t,µ) - U(t,\nu)| \leq C_0\mathcal{W}_k(µ,\nu).
$$
\end{Prop}
\begin{Rem}\label{rem:equiv}
Because of the coercivity of the Lagrangian, all admissible controls are sufficiently integrable to apply Proposition \ref{prop:ambrosiov2}, and thus justify the equivalent reformulation of the value function with a minimization over processes.
\end{Rem}
\begin{proof}
Let $t \in [0,T], µ,\nu \in \mptd$ and $\gamma$ be an optimal coupling for $\mathcal{W}_k(µ,\nu)$. The proof consists in taking an almost optimal curve for $U(t,µ)$ and translating it along $\gamma$, in the spirit of Section \ref{sec:C1}. Using this translation for $U(t,\nu)$ shall then yield the result. To make this idea of translation more apparent, I will work at the level of random variables.\\

 Consider $E$ the space of $\T^d$ valued random variables (over a fixed standard probability space). Thanks to Proposition \ref{prop:ambrosiov2}, we can write down $U$ as
 \be\label{defUprocess}
 U(t,µ) := \inf_{(X_s,Z_s)_{s \in [t,T]}}\mathbb E\left[ \int_t^TL(X_s,Z_s,\mathcal L(X_s))ds + \mathcal G(\mathcal L(X_T))  \right],
 \ee
 where the infimum is taken over processes $(X_s)$ and $(Z_s)$, valued in respectively $E$ and $\R^d$, such that, almost surely, for almost every $s$
 $$
 X_s = \int_t^s Z_u du + X_t,
 $$
 and $\mathcal L(X_t) = µ$. Consider now $(X,Y)$ an optimal coupling between $µ$ and $\nu$ for the $\mathcal{W}_k$ distance and a process $(X_s,Z_s)_{s \in [t,T]}$ which is $\eps$-optimal in \eqref{defUprocess}. Up to the use of Theorem \ref{thm:proba} that I do not detail, we can assume that $X_t = X$. Then, one can consider the process 
 $$
 Y_s = \int_t^sZ_udu + Y.
 $$
 Using the admissibility of the processes $(Y_s,Z_s)_{s \in [t,T]}$ in $U(t,\nu)$ leads to the estimate
\be\label{eq:4428}
\ba
U&(t,\nu) - U(t,µ) \leq  \mathbb E[\mathcal G(\mathcal L(Y_T)) - \mathcal G(\mathcal L(X_T)) ]+ \eps +\\
&+\mathbb E\left [\int_t^TL(Y_s,Z_s,\mathcal L(Y_s)) - L(X_s,Z_s,\mathcal L(X_s)) ds\right]\\
\leq& \mathbb E\left[\int_t^TC(1 +|Z_s|^{k-1})(d(X_s,Y_s) + \mathcal{W}_k( \mathcal L(X_s),\mathcal L(Y_s)))ds\right]\\
 &+ C\mathcal{W}_k( \mathcal L(X_T),\mathcal L(Y_T)) + \eps,
\ea
\ee
where I used the regularity of both $L$ and $\mathcal G$ given by the statement of the result. Remark now that because of the growth condition on $L$, both $$I:=\mathbb E\left[\int_t^T |Z_s|^kds\right] <\infty$$ and $I \leq C( C+U(t,µ))$. Thus, using H\"older's inequality, we obtain that 
\be\label{eq:4437}
\ba
\mathbb E&\left[\int_t^T(1 +|Z_s|^{k-1})(d(X_s,Y_s) + \mathcal{W}_k( \mathcal L(X_s),\mathcal L(Y_s)))ds\right]\\
& \leq C((T-t)^\frac{k-1}{k} + I^\frac{k-1}{k})(\mathbb E[d(X_s,Y_s)^k]^\frac{1}{k} +  \mathcal{W}_k( \mathcal L(X_s),\mathcal L(Y_s))).
\ea
\ee
Recalling that for any $s \in [t,T], \mathcal{W}_k^k( \mathcal L(X_s),\mathcal L(Y_s)) \leq \mathbb E [d(X_s,Y_s)^k] = \mathbb E [d(X_t,Y_t)^k] = \mathcal{W}_k^k(µ,\nu)$, we obtain the required result by combining \eqref{eq:4428} and \eqref{eq:4437}.
\end{proof}
\begin{Rem}
The result is stated at the level of Lipschitz regularity. It extends to H\"older classes of regularity: if $\mathcal G$ and $L$ are simply H\"older continuous in $µ$ and $(x,µ)$, then H\"older regularity of $U$ can also be proven. In particular, the present proof does not show any regularizing effect, contrary to the previous one.
\end{Rem}
\begin{Rem}
The lower bound assumption on $L$ by a polynomial in $z$ can seem to be an arbitrary assumption. This assumption is here to compensate the fact that we allow the dependency of $L$ in $x$ and $µ$ to depend on $|z|$. I leave as an exercise to the interested reader to check that if $L(x,z,µ) = \tilde L(z) + f(x,µ)$, then no lower bound is needed and only Lipschitz continuity of $f$ is sufficient to obtain the result.
\end{Rem}

\subsection{The value function is a viscosity solution}\label{sec:propviscsol}
In this section, I show that the value function defined in \eqref{defU3} is indeed a viscosity solution to 
\be\label{hjbvisc}
-\partial_t U(t,µ) + \int_{\T^d} H(x,D_µU(t,µ)(x),µ)µ(dx) = 0 \text{ in } (0,T)\times \mptd,
\ee
with terminal condition 
\be
U(T,µ) = \mathcal G(µ) \text{ in } \mptd.
\ee
I shall use Definition \ref{def:visctcb}, with the obvious adaptation to this case in which time has been reversed and sub/super differentials with respect to the measure argument are taken as in Section \ref{sec:visccomp}. I focus on a case with quadratic growth to lighten notation and leave to the interested reader the extension to the case $p \ne 2$. Several requirements are in order on the costs. I now state them and they shall be in force for the rest of this section.
\begin{hyp}\label{hyp:final}
\begin{enumerate}
\item The growth condition \eqref{hyp:dpl} holds with $C =2$, which guarantees quadratic growth of the Lagrangian. 
\item The Hamiltonian $H: (x,q,µ)\mapsto H(x,q,µ)$ is such that $D_qH$ is a Lipschitz function in $x$, uniformly in $(q,µ)$ and such that, uniformly in $q$, it is uniformly continuous in $(x,µ)$.
\item There exists $C > 0$ such that
$$
|L(x,q,µ) - L(x,q,µ') | \leq C (1+ |q|)\mathcal W_2(µ,µ'),
$$
\item There exists $C >0$ such that
\be\label{eq:coerc}
C^{-1}|q|^2 - C \leq L(x,q,µ).
\ee
This implies that, for possibly another $C >0$,
$$
|H(x,q,µ)|\leq C(1 + |q|^2).
$$
Furthermore, from the convexity of $H$ in $q$, we obtain that $D_qH$ has at most a linear growth in $q$.
\end{enumerate}
\end{hyp}
In particular, under the previous assumptions, the value function is Lipschitz continuous away from $t= T$ thanks to Proposition \ref{prop:globalcont}. Furthermore, Remark \ref{rem:equiv} is also valid here.
\subsubsection{Sub-solution property}
 \begin{Prop}\label{prop:sub-solution}
Assume Hypothesis \ref{hyp:final} holds. Then, the value function $U$ is a viscosity sub-solution to \eqref{hjbvisc}.
 \end{Prop}
 \begin{proof}
Take $(t,µ) \in [0,T)\times \mptd$ and $(\theta,\psi) \in \partial^+U(t,µ)$. Consider a standard probability space $(\Omega,\mathcal{A},\mathbb{P})$ and a couple $(X_t,Z)$ of random variables on $\Omega$ such that $\mathcal{L}(X_t,Z)(dx,dz) = µ(dx)\psi(x,dz)$. Consider the ordinary differential equation 
 \[
 \frac{dX_s}{ds} = -D_qH(X_s,Z,µ)=: \xi_s \quad \text{ for } s \in (t,T),
 \]
 with initial condition $X_t$ at time $s = t$. Thanks to Hypothesis \ref{hyp:final}, it is well defined, as a standard ODE. Furthermore, 
 $$
 \ba
 \mathcal W&^2_2(\mathcal L(X_s),µ) \leq \mathbb E[d(X_s,X_t)^2] \leq \mathbb E\left[\left|\int_t^s\xi_udu\right|^2\right]\\
 &\leq C (s-t)\mathbb E\left[\int_t^sd(X_u,X_t)^2du\right] + (s-t)^2\mathbb E\left[|\xi_u|^2\right]\\
 & \leq C(s-t)^2,
 \ea
 $$
 where this last constant depends on $(X_t,Z)$, but not on $s$. In particular, note that $\|\int_t^s\xi_udu\|_2 \leq C(s-t)$. Take $\kappa \in (0,T-t)$ and for $s \in [t,t+\kappa]$, we want to use the control $(X_s,\xi_s)$. Evaluating such a cost leads to
 \[
 U(t,µ) \leq \mathbb E\left[\int_t^{t+ \kappa}L(X_s,\xi_s,\mathcal L(X_s))ds\right] + U(t+ \kappa,\mathcal L(X_{t + \kappa})).
 \]
 Using $(\theta,\psi) \in \partial^+U(t,µ)$ yields
 \[
 \ba
 0 \leq &\mathbb E\left[\int_t^{t+ \kappa}L(X_s,\xi_s,\mathcal L(X_s))ds\right] + \mathbb{E}\left[Z\cdot\int_t^{t+\kappa}\xi_sds\right]\\
 &+ \theta\kappa + (\kappa + C \kappa)\omega(\kappa + C\kappa),
 \ea
 \]
 where $\omega$ is the modulus of continuity appearing in the property of $(\theta,\psi) \in \partial^+U(t,µ)$. In particular, it does not depend on $\kappa$. I now use the regularity of the cost in its measure argument and the fact that $\|X_{t+\kappa}-X_t\|_{2} \leq C \kappa$ to obtain
 \[
 \ba
 0 \leq &\mathbb E\left[\int_t^{t+ \kappa}L(X_s,\xi_s,µ)ds\right] + \mathbb{E}\left[Z\cdot\int_t^{t+\kappa}\xi_sds\right]\\
 &+ \theta\kappa + (C+1)\kappa \omega((C+1)\kappa)\\
 & + C\int_t^{t+\kappa}\mathbb E[( 1+ |\xi_s|)\mathcal W_2(µ,\mathcal L(X_s))]ds.
 \ea
 \]
  Using the link between $L$ and $H$ and the estimate on $\mathcal W_2(\mathcal L(X_s),µ)$, we deduce
  \[
 \ba
 0 \leq &\mathbb E\left[\int_t^{t+ \kappa}-H(X_s,Z,µ)ds + Z\cdot D_qH(X_s,Z,µ)ds\right] + \mathbb{E}\left[Z\cdot\int_t^{t+\kappa}\xi_sds\right]\\
 &+ \theta\kappa  + (C+1)\kappa \omega((C+1)\kappa) + C(1 + \|Z\|_2)\int_t^{t+\kappa}(s-t)ds .
 \ea
 \]
 Simplifying, we obtain
 \[
  0 \leq \int_t^{t+ \kappa}\mathbb{E}_{\mathbb{P}}[-H(X_s,Z,µ)]ds+ \theta \kappa + (C+1)\kappa\omega((C+1)\kappa) +  C(1 + \|Z\|_2)\kappa^2.
 \]
 Dividing by $\kappa$ and taking the limit $\kappa \to 0$ yields the required result since $(X_s)_{s > t}$ converges toward $X_t$ as $s \to t$ in $\mathbb L^2$ and $H$ is continuous.
 \end{proof}
 
 \subsubsection{Super-solution property}
 Proving the viscosity super-solution property is slightly more technical and this can be understood quite simply. It requires similar developments as in the previous proof, but this time we cannot simply use the sub-optimality of a single control, but rather have to keep many of them, sufficiently many so that the infimum is not perturbed in fact. Proceeding in this way requires the inversion of the infimum and the limit $\kappa \to 0$, after having divided by $\kappa$ in the proof of Proposition \ref{prop:sub-solution}. Such an inversion has no reason to be true in general. There are several ways to bypass this difficulty. One consists in remarking that near optimal controls are in fact such that this inversion holds. Another consists in regularizing the problem, namely by restricting the set of controls so that the inversion is correct, and then pass to the limit in this regularization. I present the two approaches, even if the first one is more elegant, as the second one requires a Lemma which can be of independent interest, namely for questions of stability of viscosity solutions. In fact it involves the same techniques of stability that were presented in Theorem \ref{thm:stabvisc}. In both cases, the crucial assumption is the coercivity of the cost \eqref{eq:coerc}.\\
 
 \textbf{By regularity of almost optimal trajectories}\\
 From the coercivity of the cost, we can extract regularity of near optimal controls as the next result shows.
 
 \begin{Lemma}\label{lemma:giacomo}
 Assume Hypothesis \ref{hyp:final} holds. Let $\kappa,\eps > 0$, $(t,µ) \in [0,T)\times \mptd$ and consider an $\eps$-optimal control $(b,\tilde m)$ in 
 $$
 U(t,µ) = \inf_{(b,\tilde m)}\left\{\int_t^{t+\kappa}\int_{\T^d\times K}L(x,b(s,x,k),(\pi_1)_\#\tilde m_s)\tilde m_s(dx,dk)ds + U(t+\kappa,(\pi_1)_\#\tilde m_{t + \kappa})  \right\}.
 $$
Then, there exists a constant $C>0$ depending only on the constant appearing in \eqref{eq:coerc} and the Lipschitz continuity of $U$ on $[t,t+ \kappa]$ such that
 $$
\int_t^{t+\kappa}\int_{\T^d\times K}|b(s,x,k)|^2\tilde m_s(dx,dk)ds \leq C\kappa + \eps,
 $$
 which implies in particular that
 $$
 \mathcal{W}_2^2(µ,(\pi_1)_\#\tilde m_{t + \kappa}) \leq C\kappa^2 + \eps \kappa.
 $$
 \end{Lemma}
 \begin{proof}
The lower bound on the Lagrangian and the $\eps$-optimality of the control imply
 $$
 C^{-1}\int_t^{t+\kappa}\int_{\T^d\times K}|b(s,x,k)|^2\tilde m_s(dx,dk)ds \leq C\kappa + U(t,µ) - U(t+\kappa,(\pi_1)_\#\tilde m_{t + \kappa}) + \eps.
 $$
 Using the Lipschitz regularity of $U$, the previous leads to
 \be\label{eq:4000}
C^{-1} \int_t^{t+\kappa}\int_{\T^d\times K}|b(s,x,k)|^2\tilde m_s(dx,dk)ds \leq C (\kappa + \mathcal{W}_2(µ,(\pi_1)_\#\tilde m_{t + \kappa})) + \eps.
 \ee
 Writing, for any $\tilde C > 0$,
 $$
 \ba
 \mathcal{W}_2(µ,(\pi_1)_\#\tilde m_{t + \kappa}) &\leq \sqrt{\kappa} \left(\int_t^{t+\kappa}\int_{\T^d\times K}|b(s,x,k)|^2\tilde m_s(dx,dk)ds\right)^\frac12\\
 &\leq \frac{\tilde C^{-1}}{2} \int_t^{t+\kappa}\int_{\T^d\times K}|b(s,x,k)|^2\tilde m_s(dx,dk)ds + \tilde C\kappa,
 \ea
 $$
 and injecting it into \eqref{eq:4000} yields the required result.
 \end{proof}
 
The previous allows me to prove the super-solution property as follows.
\begin{Prop}
Assume Hypothesis \ref{hyp:final} holds. Then, the value function $U$ is a viscosity super-solution to \eqref{hjbvisc}.
\end{Prop}
\begin{proof}
Let $(t,µ) \in (0,T)\times \mptd$ and consider $(\theta,\psi) \in \partial^-U(t,µ)$. Let $\kappa \in (0,T-t)$ and take, thanks to Lemma \ref{lemma:giacomo} a $\kappa^2$-optimal control $(b,\tilde m)$ which satisfies the conclusion of the Lemma. Denote by $(X_s,Z_s)_{s \in [t,t+\kappa]}$ the associated process. Up to a use of Theorem \ref{thm:proba}, or a change of the probability space, we can always assume that there exists a random variable $Z^*$ such that $\mathcal L(X_t,Z^*) = µ(dx)\psi_x(dz)$. Using the same types of estimates as in the proof of Proposition \ref{prop:sub-solution}, combined with the estimate of the previous Lemma, we arrive at $\mathcal{W}_2(µ,(\pi_1)_\#\tilde m_{t + \kappa}) \leq \|d(X_{t+\kappa},X_t)\|_2\leq \|\int_t^{t+\kappa}Z_sds\|_2\leq C \kappa$. Moreover, evaluating the cost of this control leads to
$$
U(t,µ) \geq \mathbb E\left[\int_t^{t+\kappa}L(X_s,Z_s,\mathcal L(X_s))ds\right] + U(t+\kappa,\mathcal L(X_{t + \kappa}))  - \kappa^2
$$
Using the properties of sub-differentials, 
$$
\ba
0 \geq &\mathbb E\left[\int_t^{t+\kappa}L(X_s,Z_s,\mathcal L(X_s))ds\right] + \mathbb E\left[Z^*\cdot\int_t^{t+\kappa}Z_s ds\right] + \theta \kappa\\
& -(C+1) \kappa\omega((C +1) \kappa)- \kappa^2,
\ea
$$
where $\omega$ is the modulus of continuity appearing in $(\theta,\psi) \in \partial^-U(t,µ)$. Using that, for any $x,µ,q,\alpha$, $H(x,q,µ) \geq -L(x,\alpha,µ) -\alpha\cdot q$, together with the estimate on the trajectory, we find
$$
0 \geq \mathbb E\left[\int_t^{t+\kappa}-H(X_s,Z^*,\mathcal L(X_s))ds\right] + \theta \kappa-(C+1)\kappa\omega((C+1)\kappa)- \kappa^2.
$$
The continuity of $H$ in its first and third arguments, as well as $\|d(X_{t + \kappa} ,X_t)\|_2 \leq C\kappa$ imply
$$
0 \geq \mathbb E\left[\int_t^{t+\kappa}-H(X_t,Z^*,\mathcal L(X_t)) - \omega_H(\kappa) ds\right] + \theta \kappa-\kappa\omega( \kappa)- \kappa^2.
$$ 
Dividing by $\kappa$ and taking the limit $\kappa \to 0$ then yield the required result.
\end{proof}
 
\textbf{By truncation of the controls}\\

Let $M > 0$, $(t,µ) \in [0,T)\times \mptd$ and define 
$$
U_M(t,µ) = \inf_{(b,\tilde m)} \left\{\int_t^{T}\int_{ \T^d\times K}L(x,b(s,x,k), (\pi_1)_\#\tilde m_s)\tilde m_s(dx,dk)ds + \mathcal G((\pi_1)_\#\tilde m_{T})\right\},
$$
where the infimum is taken over controls as usual, except for the additional constraint that $b$ is valued in the ball $B(0,M)$. The proof that $U$ is a super-solution to \eqref{hjbvisc} then follows in 3 steps. First, show that $U_M$ is a viscosity super-solution to \eqref{hjbvisc} for any $M >0$. Second, prove that $\lim_{M \to \infty} U_M$ is also a super-solution to \eqref{hjbvisc}. Third that $\lim_{M \to \infty}U_M$ is equal to $U$. For the third step, some continuity on $L$ and $G$ is required, in order to show the convergence. Also, under these requirements, we obtain the continuity of $U_M$, which despite not being essential, is needed stricto sensu to talk about elements in the sub-differential (at least the lsc is needed). 
Such continuity requirements can take the form
\be\label{contG}
\text{The terminal cost $\mathcal G$ is continuous.}
\ee
\be\label{contL}
\ba
&L \text{ is continuous in $x,µ$, uniformly in $z$.}\\
\exists C > 0,& \forall x \in \T^d, z \in \R^d, µ \in \mptd,\quad |D_z L(x,z,µ)| \leq C(1+|z|).
\ea
\ee
They can be relaxed but I do not enter such questions here. I start with the following.
\begin{Lemma}
Assume Hypothesis \ref{hyp:final} and \eqref{contG}-\eqref{contL} hold. Then, for any $M > 0$, $U_M$ is a viscosity super-solution to the equation \eqref{hjbvisc}.
\end{Lemma}
\begin{proof}
Let $(\theta,\psi^*) \in \partial^-U_M(t,µ)$ for some $(t,µ) \in [0,T)\times \mptd$ and $\kappa \in (0,T-t)$. Consider a couple of random variables $(X,Z)$ of law $µ(dx)\psi^*_x(dz)$ on a standard probability space $(\Omega,\mathcal A, \mathbb P)$. Given the formulation of the control problem and Proposition \ref{prop:ambrosiov2}, we can express $U_M$ as follows.
$$
U_M(t,µ) = \inf_{(X_s,Z_s)_{s \in [t,T]}} \left\{\mathbb E\left[\int_t^{T}L(X_s,Z_s,\mathcal L(X_s))ds\right] + \mathcal G(\mathcal L(X_T))\right\},
$$
where the infimum is taken over process $(X_s,Z_s)$ such that almost surely: $X_s = X + \int_t^sZ_u du$ and almost surely $|Z_s|\leq M$. Since $X$ has already been chosen, note that we do not lose in generality with the initial condition on $(X_s)_{s \in [t,T]}$ thanks to Theorem \ref{thm:proba}. It then follows that for all $\kappa \in (0,T-t)$
\[
U_M(t,µ) = \inf_{(X_s,Z_s)} \left\{\mathbb E\left[\int_t^{t + \kappa}L(X_s,Z_s,\mathcal L(X_s))ds\right] + U_M(t+\kappa,\mathcal L(X_{t + \kappa}))\right\},
\]
From the bound on $|Z_s|$, it follows that $\|X_{t + \kappa} - X_t\|_2 \leq M \kappa$. This yields, since $(\theta,\psi^*) \in \partial^-U_M(t,µ)$,
\[
\ba
0 \geq \inf_{(X_s,Z_s)} \bigg\{&\mathbb E_{\mathbb P} \left[\int_t^{t + \kappa}L(X_s,Z_s,\mathcal L(X_s))ds\right] + \mathbb{E}_{\mathbb{P}}\left[Z\cdot\int_t^{t+\kappa}Z_s ds\right]\\
& + \theta \kappa - ((M+1)\kappa)\omega((M+1)\kappa )\bigg\},
\ea
 \]

Using the definition of $H$, we then obtain
\[
\ba
0 \geq \inf_{(X_s,Z_s)} \bigg\{&\mathbb{E}_{\mathbb{P}}\left[\int_t^{t + \kappa}-H(X_s,Z,\mathcal L(X_s))ds\right]  + \theta \kappa - ((M+1)\kappa) \omega((M+1)\kappa) \bigg\}.
\ea
 \]
 Hence, dividing by $\kappa$ and taking the limit $\kappa \to 0$ yields the required result, because $\omega$ does not depend on $(X_s,Z_s)$.

\end{proof}
I now produce an argument which is classical in the theory of viscosity solutions in finite dimension: the infimum of super-solutions to the HJB equation is itself a super-solution to the HJB equation.
\begin{Prop}\label{prop:infvisc}
Assume that Hypothesis \ref{hyp:final} and \eqref{contG}-\eqref{contL} hold. Then the function $V := \inf_{M >0} U_M$ is a super-solution to \eqref{hjbvisc}.
\end{Prop}
\begin{proof}
 Let $t \in (0,T), µ \in \mptd$ and $(\theta,\psi^*) \in \partial^-V(t,µ)$. I prove the result as if the time variable did not raise any issue, and leave this simple extension to the interested reader. Furthermore, the proof is very similar to the one of Theorem \ref{thm:stabvisc}.

From the definition of sub-differentials, there exists a modulus of continuity $\omega$ such that for any $\eta \in \mathcal P_2(\T^d\times \R^d)$, $\Gamma \in \mathcal P(\T^d\times \R^d\times \R^d)$ such that $(\pi_1,\pi_2)_\#\Gamma = µ(dx)\psi^*_x(dz)$ and $(\pi_1,\pi_3)_\#\Gamma = \eta$,
$$
V(\exp_\#\eta) - V(µ) - \int_{\T^d\times \R^d\times \R^d}z\cdot v\,\Gamma(dx,dz,dv)  \geq - C_{\T^d}(\eta)\omega( C_{\T^d}(\eta)).
$$
Let $V_M$ be the map defined by
\[
\Gamma \mapsto U_M(\exp_\#\eta) - \int_{\T^d\times\R^{d}\times \R^d}z\cdot v\,\Gamma(dx,dz,dv) + 2 C_{\T^d}(\eta)\omega( C_{\T^d}(\eta))
\]
over $\{\Gamma \in \mathcal{P}(\T^d\times\R^{d}\times \R^d) | (\pi_1,\pi_2)_\#\Gamma = \mu(dx)\psi^*_x(dz)\}$, where I used the notation $\eta = (\pi_1,\pi_3)_\#\Gamma$. Now, remark that arguing as we did using parallel transport, $U_M$ is continuous. Hence, so is $V_M$ and, up to the use of Proposition \ref{cor2s} that I do not detail once again but which is needed, we can consider a minimum $\Gamma_M$ of $V_M$. We then have, defining $\eta_M := (\pi_1,\pi_3)_\#\Gamma_M$,
\[
\ba
U_M(\exp_\#\eta_M) &- \int_{\T^d\times \R^d\times \R^d}z\cdot v\,\Gamma_M(dx,dz,dv) + 2 C_{\T^d}(\eta_M)\omega( C_{\T^d}(\eta_M))\\
& \leq U_M(µ) - V(µ) + V(µ)\\
& \leq U_M(µ) -V(µ) + V(\exp_\#\eta_M) - \int_{\T^d\times \R^d\times \R^d}z\cdot v\,\Gamma_M(dx,dz,dv)\\
&\quad\quad +  C_{\T^d}(\eta_M)\omega( C_{\T^d}(\eta_M))
\ea
\]
Thus, passing to the limit $M \to \infty$, we deduce that $C_{\T^d}(\eta_M) \to 0$ because $U_M(µ) \to V(µ)$ and $V(\exp_\#\eta_M) \leq U_M(\exp_\#\eta_M)$. The rest of the proof follows as in Theorem \ref{thm:stabvisc} and I do not detail it here.

\end{proof}

It then remains to show the following.
\begin{Lemma}\label{lemma:convtrunc}
Assume that Hypothesis \ref{hyp:final} and that \eqref{contG}-\eqref{contL} hold. Then, for any $(t,µ) \in [0,T)\times \mptd$, $\lim_{M \to \infty} U_M(t,µ) = U(t,µ)$.
\end{Lemma}
\begin{proof}
Let $(t,µ) \in [0,T)\times \mptd$ and consider an $\eps$-optimal control $(b,\tilde m)$ for $U(t,µ)$, and the associated process $(X_s,Z_s)_{s \in [t,T]}$. For $M > 0$, define the process $(Z^M_s)_{s \in [t,T]}$ by $Z_s^M = \mathbb 1_{ |Z_s| \leq M}Z_s$. By $\eps$-optimality of $(X_s,Z_s)_{s \in [t,T]}$, $\int_t^T\mathbb E[|Z_s|^2] ds < \infty$ follows from \eqref{eq:coerc}. In particular, the $L^2$ norm of $Z^M-Z$ goes to $0$ by dominated convergence. 

Define $(X^M_s)$ by $X^M_s = X_t + \int_t^sZ^M_udu$. Observe that $(X^M,Z^M)$ defines an admissible control for $U_M(t,µ)$. It then remains to estimate the difference between the cost of $(X,Z)$ in $U(t,µ)$ and the cost of $(X^M,Z^M)$ in $U_M(t,µ)$. 

Remark that 
\be\label{est:46}
\sup_{s \in [t,T]}\|d(X^M_s ,X_s)\|_2 \leq \int_t^T \|Z^M_s - Z_s\|_2ds \to_{M \to \infty} 0.
\ee
Define $m^M_s := \mathcal L(X^M_s)$ and 
$$
\ba
W_s :=& L(X_s,Z_s,m_s) - L(X^M_s,Z^M_s,m^M_s)\\
=& L(X_s,Z_s,m_s) - L(X_s,Z^M_s,m_s) + L(X_s,Z^M_s,m_s) - L(X^M_s,Z^M_s,m_s)\\
&+L(X^M_s,Z^M_s,m_s)- L(X^M_s,Z^M_s,m^M_s).
\ea
$$
Using assumption \eqref{contL}, we obtain
$$
\ba
\mathbb{E}[|L(X_s,Z_s,m_s) - L(X_s,Z^M_s,m_s)| ]&\leq \mathbb{E}[|Z^M_s-Z_s|\sup_{\theta \in [0,1]}|D_zL(X_s,\theta Z^M_s + (1- \theta)Z_s,m_s)|]\\
& \leq C \|Z^M_s - Z_s\|_2 (1 + \|Z_s\|_2).
\ea
$$
From the regularity of the cost, we deduce $| L(X_s,Z^M_s,m_s) - L(X^M_s,Z^M_s,m_s)|+|L(X^M_s,Z^M_s,m_s)- L(X^M_s,Z^M_s,m^M_s)|  \leq \omega( d(X_s ,X^M_s)+\mathcal{W}_2(m^M_s,m_s))$ for some modulus of continuity $\omega$. It then follows that $\int_t^T\mathbb{E}[|W_s|] ds\to 0$ as $M \to \infty$. The result then follows by taking $M$ sufficiently large in 
$$
U_M(t,µ) \geq U(t,µ) \geq U_M(t,µ) - \int_t^T\mathbb{E}[|W_s|] ds - |\mathcal G(\mathcal L(X_T)) - \mathcal G(\mathcal L(X_T^M))| - \eps,
$$
since $\eps$ is arbitrary and $\mathcal G$ continuous.
\end{proof}

\subsubsection{Behaviour at the final time and characterization of the value function}
I now conclude by showing that $U$ can be in fact characterized as the unique viscosity solution to \eqref{hjbvisc} with terminal condition $\mathcal G$. The behaviour of the value function at the terminal time is of importance here and I shall come back to this later on.

From the perspective of viscosity solutions, the most natural assumption to treat the condition at terminal time is to establish that the value function is continuous there, which the next result states.
\begin{Prop}
Assume that Hypothesis \ref{hyp:final} holds and that $\mathcal G$ is continuous. Then, the value function $U$ is continuous on $[0,T]\times \mptd$ and $U(T,µ) = \mathcal G(µ)$ for any $µ \in \mptd$.
\end{Prop}
\begin{proof}
Let $(t_n,µ_n)_{n \geq 0}$ be a $[0,T]\times \mptd$ valued sequence converging to some $(T,µ_*)$ for $µ_*\in \mptd$. From Proposition \ref{prop:boundeasy}, we know that for all $n \geq 0$, 
$$
U(t_n,µ_n) \leq (T-t_n)\int_{\T^d}L(x,0,µ_n)µ_n(dx) + \mathcal G(µ_n) \longrightarrow_{n \to \infty} \mathcal G(µ_*),
$$
because $\mathcal G$ is continuous. On the other hand, from the coercivity of the cost, considering a $(n+1)^{-1}$-optimal pair $(b_n,\tilde m_n)$ for $U(t_n,µ_n)$, we know that
$$
\ba
U(t_n,µ_n) &\geq \int_{t_n}^T \int_{\T^d\times K}L(x,b_n(s,x,k),(\pi_1)_\#(\tilde m_n)_s)(\tilde m_n)_s(dx,dk)ds\\
& \quad + \mathcal G((\pi_1)_\#(\tilde m_n)_T) - \frac{1}{n+1}\\
& \geq -C(T-t_n) + C^{-1}\int_{t_n}^T \int_{\T^d\times K}|b_n(s,x,k)|^2(\tilde m_n)_s(dx,dk)ds\\
& \quad + \mathcal G((\pi_1)_\#(\tilde m_n)_T) - \frac{1}{n+1}\\
& \geq -C(T-t_n) + C^{-1}\frac{\mathcal W_2^2((\pi_1)_\#(\tilde m_n)_T,µ_n)}{T-t_n}+ \mathcal G((\pi_1)_\#(\tilde m_n)_T) - \frac{1}{n+1}.
\ea
$$
Since $U$ is bounded, we deduce first that $\lim_{n \to \infty}\mathcal W_2^2((\pi_1)_\#(\tilde m_n)_T,µ_n) = 0$ and thus that $ \mathcal G((\pi_1)_\#(\tilde m_n)_T)$ converges toward $\mathcal G(µ_*)$ as $n \to \infty$. Hence the result is proved, since the previous implies $\liminf_{n \to \infty} U(t_n,µ_n) \geq \mathcal G(µ_*)$.
\end{proof}
We now have all the tools to prove the main result of this section.

\begin{Theorem}
Assume that Hypothesis \ref{hyp:final} holds and that $\mathcal G$ is continuous. Then, there exists a unique continuous viscosity solution to \eqref{hjbvisc} with terminal condition $\mathcal G$, it is given by the value function of the problem. 
\end{Theorem}
\begin{proof}
We already saw that the value function satisfies all these requirements. If another function $V$ satisfies them, then in particular $(V-U)|_{t = T} = 0$ and, from an immediate adaptation of Theorem \ref{thm:compt} to this equation on $\mptd$, it follows that $V \leq U \leq V$. Hence the result.
\end{proof}
\begin{Rem}
The case of discontinuous functions is much more involved, and I do not want to treat it here. Remark for instance that if $\mathcal G$ is discontinuous, then for two viscosity solutions $U$ and $V$ to \eqref{hjbvisc} satisfying the terminal condition $\mathcal G$, one does not have that $V^* \leq U_*$, in fact the reverse inequality holds. Hence, applying the comparison principle in such a case is more delicate.
\end{Rem}

\subsection{Existence of viscosity solutions for equations that can be lifted}\label{sec:existnonconvex}
The previous sections can be seen as an existence result for \eqref{hjbvisc}, when the Hamiltonian stems from an optimal control problem. There are various other ways to obtain the existence of viscosity solutions, in finite dimension at least: approximation by more regular equations, Perron's method, and representation through values of zero-sum games. Here I take some time to explain how we can build on an existing result for HJB equations with non-convex Hamiltonians on certain Banach spaces, which is based on the latter technique just mentioned, and which is due to Crandall and Lions \citep{crandall1}.

Let $p \in (1,\infty)$ and consider the Hamiltonian $H: \R^d\times \R\times\R^d \times \mpprd \mapsto \R$ and the problem of existence of a viscosity solution to
\be\label{hjb:exist}
U(µ) + \int_{\R^d}H(x, U(µ), D_µU(µ)(x),µ) µ(dx) = 0 \text{ in } \mpprd,
\ee
where $H$ is assumed to satisfy:
\begin{itemize}
\item There exist constants $C > 0$ and $l \leq \min(2,p')$ such that for all $x,y,q_1,q_2 \in \R^d$, for all $r,v \in \R$, for all $µ,\nu \in \mpprd$, 
$$
\ba
|H(x,r,q_1,µ) - H(y,v,q_2,\nu)| \leq & C(|q_1| + |q_2|)|x-y| + \\
 +C(|r-v|& + |q_1-q_2|+\min(1;\mathcal{W}_p(µ,\nu)))( 1 + |q_1|^{l-1} + |q_2|^{l-1}).
\ea
$$
\item For all $x,q,µ$, the map $r \mapsto H(x,r,q,µ)$ is non-decreasing.
\end{itemize}
In particular, $H$ is not supposed to be convex in $q$, hence it does not necessarily come from an optimal control problem as was the case in the first part of this Section \ref{sec:existvisc}. The current subsection can thus be thought of as presenting more general results of existence, in the case $\mo = \R^d$. The focus is put here on stationary equations, but similar results could also be proven on time dependent problems.\\

The existence of solutions to \eqref{hjb:exist} can be established through the existence of viscosity solutions to 
\be\label{hjblp}
\mathcal U(X) + \mathbb E [ H(X,\mathcal U(X), D \mathcal U(X),\mathcal L(X))] = 0 \text{ in } \mathbb L^p,
\ee
where $D \mathcal U(X)$ stands for the (Fréchet) differential of $\mathcal U: \mathbb L^p\mapsto \R$ at $X$. For a function $\mathcal U: \mathbb L^p \mapsto \R$, I denote by $\partial^+\mathcal U(X)$ and $\partial^-\mathcal U(X)$ the usual Fréchet super and sub-differentials of $\mathcal U$ at $X\in \mathbb L^p$, which are thus convex closed subsets of $\mathbb L^{p'}$. A viscosity solution to \eqref{hjb:exist} is defined similarly to what has been done above and a viscosity solution to \eqref{hjblp} is defined as follows.
\begin{Def}
Let $\mathcal U : \mathbb L^p \mapsto \R$ be a continuous function. It is a viscosity sub-solution (resp. super-solution) to \eqref{hjblp} if for any $Z \in \partial^+\mathcal U(X)$ (resp. $Z \in \partial^-\mathcal U(X)$), 
$$
\mathcal U(X) + \mathbb E [ H(X,\mathcal U(X), Z,\mathcal L(X))] \leq 0 \quad (\text{ resp. } \geq 0).
$$
It is a viscosity solution to \eqref{hjblp} if it is both a viscosity sub and super-solution.
\end{Def} 
Under the current requirements on $H$, Theorem 1 in \citep{crandall1} proves that there exists at most one viscosity solution to \eqref{hjblp}. Furthermore, \eqref{hjblp} is the lift of \eqref{hjb:exist} and the following result links the two equations.

\begin{Prop}
If $U : \mpprd\mapsto \R$ is a viscosity solution to \eqref{hjb:exist}, then its lift $\mathcal U$ is a viscosity solution to \eqref{hjblp}. If $\mathcal U: \mathbb L^p \mapsto \R$ is a continuous viscosity solution to \eqref{hjblp}, then for $X \in \mathbb L^p$, $\mathcal U(X)$ depends only on the law of $X$, i.e. it is equal to some $U(\mathcal L(X))$, and $U$ is a viscosity solution to \eqref{hjb:exist}.
\end{Prop}
\begin{proof}
The first assertion is a direct consequence of Proposition \ref{prop:superdifflift}. Indeed, for $Z \in \partial^+ \mathcal U(X)$, from Proposition \ref{prop:superdifflift} $\gamma := \mathcal L(X,Z) \in \partial^+ U(\mathcal L(X))$. Hence, since $U$ is a viscosity sub-solution, denoting $µ = \mathcal L(X)$ 
$$
U(µ) + \int_{\R^d\times \R^d}H(x, U(µ), z,µ) \gamma(dx,dz) \leq 0.
$$
But the previous is exactly 
$$
\mathcal U(X) + \mathbb E [ H(X,\mathcal U(X), Z,\mathcal L(X))] \leq 0.
$$
For the second one, consider a one-to-one bi-measurable, measure preserving mapping $\tau: \Omega\mapsto \Omega$ and define $\mathcal V(X) = \mathcal U(X\circ \tau)$. Consider $Z \in \partial^+ \mathcal V(X)$. Then, for all $Y$,
$$
\mathcal U(Y\circ \tau) \leq \mathcal U(X \circ \tau) + \mathbb E[Z\cdot(Y-X)] + o(\|Y-X\|),
$$
which implies that $Z\circ \tau \in \partial^+ \mathcal U(X\circ \tau)$. In particular, since $\mathcal U$ is a viscosity sub-solution to \eqref{hjblp},  
$$
\mathcal V(X) + \mathbb E [ H(X\circ \tau,\mathcal V(X), Z\circ \tau,\mathcal L(X))] \leq 0
$$
Since $$\mathbb E [ H(X\circ \tau,\mathcal V(X), Z\circ \tau,\mathcal L(X))]= \mathbb E [ H(X,\mathcal V(X), Z,\mathcal L(X))],$$ we obtain that $\mathcal V$ is also a viscosity sub-solution to \eqref{hjblp}. Arguing similarly, we deduce that $\mathcal V$ is also a viscosity super-solution, hence a viscosity solution. Since $H$ satisfies the requirements of Theorem 1 in \citep{crandall1}, the uniqueness of viscosity solutions to \eqref{hjblp} yields  $\mathcal V = \mathcal U$. Hence, since $\tau$ is any one-to-one measure preserving mapping and $\mathcal U$ is continuous, it follows that for any $X, Y$ which have same law, $\mathcal U(X) = \mathcal U(Y)$. In particular, $\mathcal U$ is the lift of some function $U : \mpprd \mapsto \R$, which is necessarily a viscosity solution to \eqref{hjb:exist}, thanks to Proposition \ref{prop:superdifflift}.
\end{proof}
\begin{Rem}
Note that the previous result does not imply that a viscosity sub-solution to \eqref{hjblp} is the lift of a viscosity sub-solution to \eqref{hjb:exist}, as there may exist sub-solutions which are not lifts. However, if we assume, in addition to the sub-solution property, such an invariance property, then the pull-down is indeed a viscosity sub-solution to \eqref{hjb:exist}.
\end{Rem}

The previous result allows one to extend instantly results of existence of \eqref{hjblp} to \eqref{hjb:exist}. There are numerous results of existence of first order HJB equations on Banach spaces, and I will not review all of them, nor prove the one I am going to state here, since the framework needed for the study of such equations is quite far from the tools developed and studied here. I simply mention that, in my opinion there are two general approaches to such problems. The first one, maybe the simpler, consists in approximating \eqref{hjblp} with finite dimensional elements. We can then use existence result for finite dimensional equations and pass to the limit in this approximation. The second one consists in studying particular Hamiltonians $H$ for which, \eqref{hjblp} arises as the so-called Hamilton-Jacobi-Isaacs equation of a differential zero-sum game. Existence is then obtained by verifying that the value function of the game is in fact a viscosity solution to the associated equation. By various approximation results, one can then extend the class of Hamiltonians for which one can prove existence. The first type of method is for instance presented in Section 3.7 in \citep{fabbri} or used in \citep{bertucci2025singular}; while the second one is for instance the object of \citep{crandall2} from which we extract the next result as a simpler version of their Theorem 1.1 (ii). Note that in both methods, stability properties of viscosity solutions are fundamental to pass to the limit.

\begin{Theorem}
Under the standing assumptions on $H$, there exists a (unique) viscosity solution to \eqref{hjb:exist}.
\end{Theorem}
I do not provide a self-contained proof of this result, and invite the readers who are not familiar with \citep{crandall2} to simply omit it.
\begin{proof}
I indicate how to verify that the assumptions (H0)-(H4) of Theorem 1.1 (ii) in \citep{crandall2} are satisfied. First, two auxiliary functions need to exist. We can take $d(X,Y) = \|X-Y\|_p$ and $\eta(X) = 2(1 + \mathbb E[|X|^p])^{\frac 1p}$. Next, their assumption (H1) is directly assumed here. Assumptions (H0) and (H2)-(H4) follow from the following computation, for $X,Y \in \mathbb L^p, P,Q \in \mathbb L^{p'}$, $r,s \in \R$, $µ= \mathcal L(X)$ and $\nu= \mathcal L(Y)$. 
$$
\ba
\mathbb E&[|H(X,r,P,µ) - H(Y,s,Q,\nu)|] \\
&\leq \mathbb E[C(|r-s| + |P-Q| + \mathcal{W}_p(µ,\nu))( 1+ |P|^{l-1} + |Q|^{l-1})  ] + C\mathbb E [|X-Y|(1 + |P| + |Q|)]\\
& \leq C \|X-Y\|_p(1 + \|P\|_{p'} + \|Q\|_{p'}) +C (\mathcal{W}_p(µ,\nu) + |r-s|)( 1 + \|P\|^{l-1}_{(l-1)} + \|Q\|^{l-1}_{(l-1)})\\
& \quad + C\|P-Q \|_{p'}(1 + \|P\|^{l-1}_{p'} + \|Q\|^{l-1}_{p'}).
\ea
$$
Since $\mathcal{W}_p(µ,\nu) \leq \|X-Y\|_p$, (H0) and (H2) in \citep{crandall2} follow. For (H4), we need to evaluate the previous for $P = Q = \lambda (X-Y)|X-Y|^{p-2}\|X-Y\|_p^{1-p}$ and $r =s$. Using $\|P\|_{l-1} \leq \|P\|_{p'} = \lambda$, we deduce (since $l-1 \leq 1$)
$$
\mathbb E[H(X,r,P,µ) - H(Y,r,Q,\nu)] \leq C ( 2 \lambda)d(X,Y) + C \min(1;d(X,Y)).
$$
The assumption (H4) is then verified by realizing that, setting $G(r) = C$ the requirement in Theorem 1.1 (ii) in \citep{crandall2} follows.

\end{proof}

\subsection{Bibliographical comments}
In the analysis of an optimal control problem, obtaining various regularity properties of the value function is now a classical topic. In finite dimension, I refer to the book of Bardi and Capuzzo-Dolcetta \citep{bardiitalo} for deterministic problems, and to the one of Fleming and Soner \citep{flemingsoner} for stochastic problems. The results I presented here are mostly taken from my work \citep{bertucci2024stochastic}. Usually a minimum of regularity needs to be established before trying to verify that the value function is indeed a viscosity solution to the associated HJB equation. 

The sub-solution property is quite often easier to prove (when looking at minimization problems), as is apparent here, namely because we need to have information on a single trajectory, that we can choose. The super-solution is usually more involved, as is the case here. A typical technique, which I choose to avoid here because it is not the most general one, consists in proving first that an optimal control exists. We are then in a similar position as for the sub-solution property, except for the fact that the minimizer is quite often simply given by some abstract result and we know little information about it. The two techniques I presented are the most general ones. The regularity of controls of finite cost is a powerful tool and I learned Lemma \ref{lemma:giacomo} from Giacomo Ceccherini Silberstein. The truncation of the controls has the interest of passing through Proposition \ref{prop:infvisc}, which follows the lines of Proposition 4.3 in Crandall, Ishii and Lions \citep{crandall1992user}, and which is a fundamental tool when one wants to establish some stability properties of viscosity solutions. I chose not to present many stability results as they seem to me easy to establish given the results I already presented here. I hope that this will not be held against this presentation.

I insist that, even though additional care is needed, the techniques presented here should be of interest for proving similar results in differential games. For instance, they could be of interest to prove that the value function of some zero-sum games is indeed the unique viscosity solution of the associated Hamilton-Jacobi-Isaacs equation.

The treatment of the terminal condition in time is always a difficulty that has to be taken with care.

Finally, using the formula given by the value function is a quite natural way of obtaining existence of a viscosity solution, even if restricted to HJB equations with convex Hamiltonians. Note that the result of existence of Crandall and Lions \citep{crandall2} that I used in Section \ref{sec:existnonconvex} relies on a similar idea, but with a representation through values of two-player zero-sum games instead.

\newpage

\section{A mean field optimal control problem with noise: stochastic optimal transport}\label{sec:stoot}
In this section, I introduce a variant of the problem considered in Section \ref{sec:optcontr}. There are two main differences here. The first one is that this new problem incorporates a stochastic component, and the second one is that it has a state constraint at the terminal time, which thus corresponds to a terminal cost $\mathcal G$ which is valued in $\R \cup \{+\infty\}$. \\

Several arguments already presented can be adapted to these cases, and I will not detail all these adaptations, focusing mostly on the challenging parts of this problem. I start by presenting the model, I then explain why the associated value function is indeed well defined by providing bounds. We will then see why it is indeed a viscosity solution to the associated PDE and I will conclude by explaining how it can be characterized as the unique viscosity solution with appropriate boundary requirements. All this section is restricted to the case $\mo = \T^d$, and I will see $\mptd$ equipped with the topology of weak convergence.

\subsection{The model}
\subsubsection{The optimization problem and the value function}
The problem of interest here is fairly similar to the one considered in Section \ref{sec:optcontr}, except that at time $t = T$, the state of the current measure $m_T$ needs to be equal to the value of a $\mptd$ valued, adapted, stochastic process $(\nu_t)_{t \in [0,T]}$ defined on a standard filtered probability space $(\Omega,\mathcal A, (\mathcal F_t)_{t \geq 0},\mathbb{P})$.

There are several such processes $(\nu_t)_{t \in [0,T]}$ that could be considered and I present here three of them.

\textbf{The reveal case}\\
The first, and simpler one, simply consists in a process which is constant equal to $\nu_0 \in \mptd$ in $[0,T/2)$ and which is constant in time afterwards, of law $\Theta \in \mathcal P(\mathcal P(\T^d))$, in $[T/2,T]$. In other words, at time $T/2$, a terminal target is randomly drawn according to the law $\Theta$. For instance, if $\Theta = \frac12(\delta_{\nu_1} + \delta_{\nu_2})$, then the target at terminal time is $\nu_1$ with probability $\frac 12$ or $\nu_2$ with probability $\frac12$.\\

\textbf{The jump case}\\
In this second case, there is a sequence of Poisson times $(\tau_i)_{i \geq 1}$ given a Poisson process of intensity $\lambda : [0,T]\mapsto \R_+$. The intensity $\lambda$ is assumed to be smooth. Hence, almost surely, the sequence of random times is finite. Then, at each of those times, the state of the target process $(\nu_t)_{t \geq 0}$ jumps from $\nu_{t_i^-}$ to $\mathcal T(\nu_{t_i^-})= \nu_{t_i^+}$, for a certain function $\mathcal T: \mptd \mapsto \mptd$. In this framework, we could also have considered a case in which at each of those times, instead of this jump function $\mathcal T$, the new value of the target process is drawn from a known distribution.\\

\textbf{The diffusion case}\\
In this third and last setting, the target process $(\nu_t)$ is given by 
$$
\nu_t := (\tau_{W_t})_\#\bar\nu
$$
for a given $\bar\nu \in \mptd$ and $(W_t)_{t \geq 0}$ is the strong solution to 
\be\label{SDE4}
dW_t = \sigma(t)dB_t,
\ee
with initial condition $W_0 = 0$, for $(B_t)_{t \geq 0}$ a standard $d$ dimensional Brownian motion and a continuous function $\sigma : [0,T]\mapsto \R$. More general SDEs could also be considered under standard regularity assumptions on the drift and volatility.\\

\textbf{The control problem and the value function}\\
 As usual in stochastic optimal control, the controls are allowed to be random, but they have to be adapted to the filtration $(\mathcal F_t)_{t \geq 0}$, or more generally to the filtration generated by the underlying randomness. Furthermore, I assume a sort of risk neutrality so that it is the expected value of the (now stochastic) cost which is minimized. In this case, this leads to the minimization problem
\be\label{sot}
\inf_{(b,\tilde m)} \mathbb E\left[ \int_0^T \int_{\T^d\times K}L(x,b(s,x,k))\tilde m_s(dx,dk)ds \right],
\ee
where $L: \T^d\times \R^d \mapsto \R$ is a given cost function on which assumptions shall be made later on, and where the infimum is taken over pairs $(b,\tilde m)$ as follows:
\begin{itemize}
\item $b : \Omega\times [0,T]\times \T^d \times K \mapsto \R^d$ is measurable and adapted to $(\mathcal F_t)_{t \geq 0}$.
\item $\tilde m: \Omega \times [0,T]\mapsto \mathcal P(\T^d\times K)$ is measurable, adapted to $(\mathcal F_t)_{t \geq 0}$ and, almost surely\footnote{In all this section, almost surely means $\mathbb P$ almost surely in  $\omega \in \Omega$.} is a continuous function of time with first marginal which equals at initial time $µ_0$, which is fixed.
\item Almost surely, $(b,\tilde m)$ satisfies the continuity equation \eqref{contpsi}.
\item Almost surely, the terminal target is reached, i.e.
$$
(\pi_1)_\#\tilde m_T = \nu_T.
$$
\end{itemize}
\begin{Rem}
To simplify notation, I do not consider cases in which the cost $L$ depends on the value of the controlled measure. Such cases can also be treated following the techniques presented above.
\end{Rem}

We are interested in the value function of the previous problem, which naturally depends on the value of the stochastic process $(\nu_t)_{t \in [0,T]}$ at the time in question. This value function is the map $U: [0,T)\times \mptd \times \mptd\mapsto\R$ defined by
\be\label{defUsto}
U(t,µ,\nu):=\inf_{(b,\tilde m)} \mathbb E\left[ \int_t^T \int_{\T^{d}\times K}L(x,b(s,x,k))\tilde m_s(dx,dk)ds \bigg| \nu_t = \nu\right],
\ee
where the infimum is taken over the same controls as above except for the fact that the time interval is now $[t,T]$, and that the initial condition of the continuity equation is $µ$.

In the diffusion case, we can reduce the dependence in $\nu$ to a dependence in a finite dimensional parameter, which is the value $w$ of the solution to \eqref{SDE4}. Hence, in this last case, we consider the value function defined by
$$
U(t,µ,w):=\inf_{(b,\tilde m)} \mathbb E\left[ \int_t^T \int_{\T^{d}\times K}L(x,b(s,x,k))\tilde m_s(dx,dk)ds \bigg| W_t = w\right],
$$
where the infimum is taken over the analogous space.

\subsubsection{The Hamilton-Jacobi-Bellman equations}
As in the deterministic case, we expect the value function to be given as the solution to an HJB equation. This equation still follows from an infinitesimal version of Bellman's dynamic programming principle. In this stochastic setting, this principle takes the form
$$
\ba
U(t,µ,\nu) = \inf_{(b,\tilde m)} \mathbb E\bigg[ &\int_t^{t + \eps} \int_{\T^d\times K}L(x,b(s,x,k))\tilde m_s(dx,dk)ds\\
& + U(t+\eps,(\pi_1)_\#\tilde m_{t + \eps},\nu_{t + \eps}) \bigg| \nu_t = \nu\bigg].
\ea
$$
Note that for the right-hand side to be well defined in general, some measurability is needed on $U$ for the expectation to have a meaning. In our setting, such requirement will follow directly from the continuity of the value function.\\

I do not reproduce the straightforward derivations of the HJB equations. Now, they rely on Ito's Lemma, and directly give the three following HJB equations.\\

\textbf{The reveal case}\\
$$
\ba
-\partial_t U + \int_{\T^d} H(x,D_µU(t,µ,\nu)(x))µ(dx) = 0 \text{ in } \left(0,\frac T2\right)\cup\left(\frac T2,T\right)\times \mptd^2,\\
U|_{t = \frac T2^-}(µ,\nu) = \int_{\mptd}U\left(\frac T2^+,µ,\nu'\right)\Theta(d\nu'),\\
U|_{t = T}(µ,\nu) = \begin{cases} 0 \text{ if } µ = \nu,\\ + \infty \text{ else.}\end{cases}
\ea
$$

\textbf{The jump case}\\
\be\label{hjbl}
\ba
&-\partial_t U + \int_{\T^d}H(x,D_µU(t,µ,\nu)(x))µ(dx)+ \\
&\quad + \lambda(t)(U(t,µ,\nu) - U(t,µ,\mathcal T (\nu))) = 0 \text{ in } (0,T)\times \mptd^2,\\
&U|_{t = T}(µ,\nu) = \begin{cases} 0 \text{ if } µ = \nu,\\ + \infty \text{ else.}\end{cases}
\ea
\ee

\textbf{The diffusion case}\\
$$
\ba
-\partial_t U + \int_{\T^d}H(x,D_µU(t,µ,w)(x))µ(dx) - \frac{\sigma^2(t)}{2}\Delta_w U = 0 \text{ in } (0,T)\times \mptd \times \T^d,\\
U(T,µ,w) = \begin{cases} 0 \text{ if } µ = (\tau_{w})_\#\nu,\\ + \infty \text{ else.} \end{cases}.
\ea
$$

For the rest of this section, I will focus on the jump case.

\subsection{Bounds on the value function}
The first step in the analysis of \eqref{sot} is to show that there exist admissible controls of finite cost. This will guarantee that the value function is not equal to $+\infty$ everywhere. In addition it will allow us to work nicely away from $t = T$. Indeed, if we prove that $U(T-\eps,µ,\nu)$ is bounded for $\eps > 0$, uniformly\footnote{In a non-compact case such as $\mpt$, it is the local boundedness which is of course sufficient.} in $µ$ and $\nu$, then we will be able to treat the problem on $[0,T-\eps]$ by using $U(T-\eps,\cdot,\cdot)$ as a terminal cost. Furthermore, in addition to obtaining bounds on the value function, I shall show a precise behaviour of the value function near $t = T$.\\

Because of the terminal constraint, to obtain bounds on the value function, a more precise control on the cost is needed. Here, I assume that the cost $L$ satisfies
\be\label{contrL}
\exists C >0, p \in (1,\infty), \forall x \in \T^d,z \in \R^d,\quad 0 \leq L(x,z) \leq C(1 + |z|^p).
\ee

Furthermore, let $U_{det}$ be the value of the associated deterministic control problem, that is $U_{det}$ is defined in the same way as $U$, but in the case $\lambda = 0$. This deterministic value $U_{det}$ shall be of importance in the study of $U$, so I start by showing that it is well defined under \eqref{contrL}.
 
\begin{Lemma}
Assume $L$ satisfies \eqref{contrL}. Then for all $(t,µ,\nu) \in [0,T)\times \mptd^2$, $U_{det}(t,µ,\nu)< \infty$.
\end{Lemma}
\begin{proof}
Let $\gamma$ be an optimal coupling for $\mathcal{W}_p(µ,\nu)$. Then, consider a couple $(X,Y)$ of law $\gamma$ and define $(X_s,Z_s) = (T-t)^{-1}((T-s)X +(s-t)Y,Y-X)$. The control associated with this process is indeed deterministic (it is $\mathcal F_0$ measurable), and admissible for $U_{det}(t,µ,\nu)$. Moreover, from the bound on $L$, we immediately obtain that
$$
U_{det}(t,µ,\nu) \leq C\left((T-t) +\frac{\mathcal{W}_p^p(µ,\nu)}{(T-t)^{p-1}}\right).
$$
\end{proof}

Hence, when $\lambda = 0$, the target problem is well defined. The next result gives some condition under which the same can be said in the case $\lambda \ne 0$. In order to state this result, define $$M(t):= \sup_{s \leq t,µ,\nu}\{U_{det}(s,µ,\nu)\}.$$

 \begin{Prop}\label{prop:behavjumps}
Assume \eqref{contrL} holds and that there exist $C >0$ and $\delta > -1$ such that for $T-t \leq C^{-1}$
\be\label{hypgrowth}
 \lambda(t)M(t) \leq C (T-t)^{\delta},
\ee
 and
\be\label{hyplambdatec}
 \lambda(t) \leq C (T-t)^{-1}\int_t^T\lambda.
\ee
Then there exists a continuous function $\beta : [0,T] \mapsto \mathbb{R}_+$, such that
 \[
 \sup_{µ,\nu\in \mptd}|U(t,µ,\nu) - U_{det}(t,µ,\nu)| \leq \beta(t) \underset{t \to T}{\longrightarrow} 0.
 \]
 \end{Prop}
 \begin{proof}
 I argue first as if the infimum in the deterministic problem is always reached. Notice that if $\lambda = 0$ in $L^1((T-\kappa,T),\R_+)$ for some $\kappa > 0$, then the result holds true trivially. Hence, I assume 
 \[
 \forall t \in [0,T), \int_t^T \lambda(s)ds > 0.
 \]
Consider the problem starting in $µ \in \mptd$ at time $t$ and where the target process is equal to $\nu$ at $t$. Let $n$ be the (random) number of jumps in $[t,T]$ and consider the sequence $t_0 = t < t_1 < t_2 < ... < t_n \leq T$ of random times at which the target process jumps. Note that this sequence is finite almost surely, possibly empty and that the event $\{t_n = T\}$ shall be ignored since it happens with probability $0$. Consider the random pair $(b,\tilde m)$ given by
 \begin{itemize}
 \item For $1\leq i \leq n$, $s \in (t_{i-1}, t_i)$, $(b(s,\cdot),\tilde m_s)$ is equal (up to a change of time) to an optimal control in the problem $U_{det}(t_{i-1},(\pi_1)_\#\tilde m_{t_{i-1}},\mathcal{T}^{i-1}\nu)$. In other words, we do as if there were no further jumps.
 \item For $s \in (t_n,T)$, $(b(s,\cdot),\tilde m_s)$ is given through an optimal control in $U_{det}(t_n,(\pi_1)_\#\tilde m_{t_n},\mathcal{T}^n\nu)$.
 \end{itemize}
 Such a pair is clearly admissible, since it is non-anticipative, thus adapted, and satisfies $(\pi_1)_\#\tilde m_T = \nu_T$ almost surely. We now estimate its cost.
 \[
 \begin{aligned}
& \mathbb{E}\left[ \int_t^T \int_{\mathbb{T}^d\times K}L(x,b(s,x,k))\tilde m_s(dx,dk)ds\right]\\
 &= \mathbb{P}(n > 0)\mathbb{E}\bigg[   \int_{t_n}^T\int_{\T^d\times K}L(x,b(s,x,k))\tilde m_s(dx,dk)ds \\
 &\quad \quad \quad \quad \quad + \sum_{i = 0}^{n-1} \int_{t_{i}}^{t_{i+1}} \int_{\mathbb{T}^d\times K}L(x,b(s,x,k))\tilde m_s(dx,dk)ds \bigg| n > 0\bigg] +  \mathbb{P}(n= 0)U_{det}(t,µ,\nu)\\
 & \leq \mathbb{P}(n > 0)\mathbb{E}\left[U_{det}(t_n, m_{t_n}, \mathcal{T}^n \nu)+  \sum_{i=0}^{n-1} U_{det}(t_i,m_{t_i},\mathcal{T}^i \nu) \bigg| n > 0\right] + \mathbb{P}(n= 0)U_{det}(t,µ,\nu)\\
 & \leq \mathbb{P}(n > 0)\mathbb{E}\left[ M(t_n) + \sum_{i=0}^{n-1} M(t_i)\bigg| n > 0\right] +  \mathbb{P}(n= 0)U_{det}(t,µ,\nu)\\
 & \leq \mathbb{P}(n> 0)\mathbb{E}\left[ (1 + n)M(t_n) | n > 0\right] + \mathbb{P}(n= 0)U_{det}(t,µ,\nu).
 \end{aligned}
 \]
 We can now compute
 \[
 \begin{aligned}
 \mathbb{E}\left[ (1 + n)M(t_n)| n > 0\right] = \sum_{i = 1}^{\infty}\mathbb{E}[(1+n)M(t_n)|n = i]\mathbb{P}( n = i| n > 0).
 \end{aligned}
 \]
 Since the $(t_n)_{n \geq 0}$ are given by a Poisson process, 
 \[
 \mathbb{P}(n > 0) \mathbb{P}(n = i | n > 0) = \mathbb{P}(n = i) = \frac{\left(\int_t^T\lambda(s)ds\right)^i}{i!}e^{-\int_t^T\lambda(s)ds},
 \]
 and also there exists $C >0$ such that for any $1\leq i < n$, the law of $t_i$ conditioned on $t_{i-1}$ has a density which is bounded by 
 \[
s \mapsto \mathbb{1}_{s \geq t_{i-1}} C\frac{\lambda(s)}{\int_{t_{i-1}}^T\lambda}.
 \]
Hence, we can estimate
\[
\begin{aligned}
\mathbb{E}[(1+n)M(t_n)|n = i] &\leq C^{i} (i+1)\int_t^T\int_{t_1}^T\dots\int_{t_{i-1}}^TM(t_i)\lambda(t_i)\frac{dt_i}{\int_{t_{i-1}}^T\lambda}\dots \frac{\lambda(t_2)dt_2}{\int_{t_1}^T\lambda}\frac{\lambda(t_1)dt_1}{\int_t^T\lambda}\\
& \leq C^i(i+1)\int_t^T\int_{t_1}^T\dots\int_{t_{i-2}}^T(T-t_{i-1})^{\delta +1}\frac{\lambda(t_{i-1})dt_{i-1}}{\int_{t_{i-1}}^T\lambda}\dots \frac{\lambda(t_2)dt_2}{\int_{t_1}^T\lambda}\frac{\lambda(t_1)dt_1}{\int_t^T\lambda}\\
& \leq \frac{C^{i}(i+1)(T-t)^{\delta +1}}{\int_t^T\lambda(s)ds} 
\end{aligned}
\]
From the previous estimate, it follows that
\[
\begin{aligned}
\mathbb{P}(n> 0)\mathbb{E}\left[ (1 + n)M(t_n) | n > 0\right] &\leq C\sum_{i = 1}^{\infty}\frac{C^i(i+1)}{i!} e^{-\int_t^T\lambda}\left(\int_t^T\lambda\right)^{i-1}(T-t)^{\delta +1}\\
& \leq C (T-t)^{ \delta+1}.
\end{aligned}
\]
The previous yields
\[
(1 - \mathbb{P}(n = 0))U_{\det}(t,µ,\nu) \leq C \int_t^TM(s)\lambda(s)ds \leq C (T-t)^{ \delta+1}.
\]
Thus, setting $\beta(t) = C (T-t)^{\delta+1}$ finally produces
\[
U(t,µ,\nu) \leq \mathbb{E}\left[ \int_t^T \int_{\mathbb{T}^d\times K}L(x,b(s,x,k))\tilde m_s(dx,dk)ds\right] \leq U_{det}(t,µ,\nu) + \beta(t).
\]

 Obtaining the lower bound is easier. Indeed, for $\epsilon > 0$, consider an $\epsilon$-optimal pair $(b,\tilde m)$ (which exists since the value is bounded from below). It then follows that
 \[
 \begin{aligned}
 U(t,µ,\nu) &\geq \mathbb{E}\left[ \int_t^T \int_{\mathbb{T}^d\times K}L(x,b(s,x,k))\tilde m_s(dx,dk)ds\right] - \epsilon\\
 &= \mathbb{P}(n = 0) \mathbb{E}\left[ \int_t^T \int_{\mathbb{T}^d\times K}L(x,b(s,x,k))\tilde m_s(dx,dk)ds\bigg| n = 0\right] - \epsilon \\
 & \quad + \mathbb{P}(n > 0)\mathbb{E}\left[ \int_t^T \int_{\mathbb{T}^d\times K}L(x,b(s,x,k))\tilde m_s(dx,dk)ds \bigg| n > 0\right]\\
 & \geq \mathbb{P}(n = 0) \mathbb{E}\left[\int_t^T \int_{\mathbb{T}^d\times K}L(x,b(s,x,k))\tilde m_s(dx,dk)ds \bigg| n = 0\right] - \epsilon\\
 & \geq  \mathbb{P}(n = 0)U_{det}(t,µ,\nu) - \epsilon.
 \end{aligned}
 \]
 Since $\epsilon$ is arbitrary, we deduce that the inequality also holds for $\epsilon = 0$. Hence, we deduce the lower bound
 \[
 U(t,µ,\nu) \geq U_{det}(t,µ,\nu) - \beta(t)
 \]
 following the same computation as in the part concerning the upper bound.\\
 
 We end the proof by remarking that if we are not able to consider optimal controls for the deterministic problem, then considering appropriate $\epsilon'$-optimal controls, and passing to the limit $\eps' \to 0$ yields the required estimates.
  \end{proof}

 \begin{Rem}
A few comments on the hypotheses of the result.
\begin{itemize}
\item The assumption \eqref{hyplambdatec} is purely technical (and thus quite regrettable). It is verified by any function such that $\lambda(t) \sim C (T-t)^{\alpha}$ for any $C> 0, \alpha > -1$. However it is not automatically verified, as for instance $\lambda(t) = -\frac{d}{dt}(e^{-(T-t)^{-2}})$ does not verify it. I do not know whether this can be removed.
\item On the contrary, \eqref{hypgrowth} is quite important in our proof! This assumption yields an integrability condition on the product $\lambda M$. Such an integrability condition is crucial. Note for example that if $\lambda$ is constant, then we require (among other things), that $M \in L^1_{loc}$, which is not the case for a quadratic cost. See the example below.\end{itemize}
 \end{Rem}
 The following example provides insight on what happens when the assumptions of Proposition \ref{prop:behavjumps} are not in force.
  \begin{Ex}
Let $L(x,\alpha) = |\alpha|^2$ be the quadratic cost, and let the jump rate $\lambda >0$ be constant. Furthermore, assume that the jumps are not trivial, i.e. that $\mathcal{T}\nu \ne \nu$ for some $\nu \in \mptd$. The deterministic problem is then simply the quadratic optimal transport and thus $M(t) = C(T-t)^{-1}$. Let $t > 0$ and assume that at this time, the target process is equal to $\nu$ such that $\mathcal T \nu \ne \nu$.
By conditioning on the number of remaining jumps, we can bound from below the cost of any admissible control $(b,\tilde m)$ by
\[
\begin{aligned}
& \mathbb{P}(n = 1)\mathbb{E}\left[\int_t^T \int_{\mathbb{T}^d\times K} |b(s,x,k)|^2\tilde m_s(dx,dk)ds \bigg| n =1\right]\\
& +\mathbb{P}(n = 0) \mathbb{E}\left[\int_t^T \int_{\mathbb{T}^d\times K} |b(s,x,k)|^2\tilde m_s(dx,dk)ds \bigg| n =0\right] .
\end{aligned}
\]
Let us denote by $\rho$ the density of the law of the jump $\tau_1$, conditioned on $\{n = 1\}$. By a slight abuse, I argue as if all admissible controls are Markovian, which allows us to consider $(µ_s)_{s \in [t,T]}$ the trajectory in the event $\{n = 0\}$. (In general, considering the mean of the trajectories on the event $\{n = 0\}$ could also yield something similar). By definition of the $2$-Wasserstein distance, we obtain that for any $t' \in [t,T)$
\be\label{eq1395}
\ba
 \mathbb{E}\left[\int_t^T \int_{\mathbb{T}^d\times K} |b(s,x,k)|^2\tilde m_s(dx,dk)ds \bigg| n =0\right]  &\geq  \mathbb{E}\left[\int_{t'}^T \int_{\mathbb{T}^d\times K} |b(s,x,k)|^2\tilde m_s(dx,dk)ds \bigg| n =0\right] \\
 &\geq\frac{W_2^2(µ_{t'},\nu)}{T-t'} 
 \ea
\ee
I then compute the expectation conditioned on $\{n =1\}$, namely by distinguishing on the time of the jump.
\[
 \mathbb{E}\left[\int_t^T \int_{\mathbb{T}^d\times K} |b(s,x,k)|^2\tilde m_s(dx,dk)ds \bigg| n =1\right] \geq \int_t^T\frac{W_2^2(µ_s, \mathcal{T}\nu)}{T-s}\rho(s)ds
\]
Integrating \eqref{eq1395} with respect to $\rho$, we deduce that
\[
\begin{aligned}
U(t,µ,\nu) \geq& \mathbb{P}(n = 0)\int_t^T\frac{W_2^2(µ_{s},\nu)}{T-s}\rho(s)ds\\
& +\mathbb{P}(n = 1)\int_t^T\frac{W_2^2(µ_{s},\mathcal{T}\nu)}{T-s} \rho(s)ds\\
\geq &\mathbb{P}(n = 1)\int_t^T\frac{W_2^2(µ_s,\nu) + W_2^2(µ_{s},\mathcal{T}\nu)}{T-s} \rho(s)ds,
\end{aligned}
\]
if $T-t$ is sufficiently small so that $\mathbb{P}(n= 1) \leq \mathbb{P}(n=0)$. The right hand side of the previous inequality is equal to $+ \infty$. Indeed, $\nu \ne \mathcal{T}\nu$ implies that the numerator in the last integral is bounded from below by a positive constant. Hence for any $µ \in \mptd$, $U(t,µ,\nu) = + \infty$, which thus extends to previous times by dynamic programming.
 \end{Ex}

The previous bound  allows one to  prove regularity of the value function away from $\{t = T\}$, namely by using the controllability of the problem, as in Section \ref{sec:contcontr}.
\begin{Prop}
Assume that Hypothesis \ref{hyp:final}, and \eqref{hypgrowth}, \eqref{hyplambdatec} hold. Then, for any $\eps > 0$, for any $\nu \in \mptd$ the function $(t,µ) \mapsto U(t,µ,\nu)$ is Lipschitz continuous in $[0,T-\eps]\times \mptd$.
\end{Prop}
\begin{proof}
This follows simply from the proof of Proposition \ref{prop:globalcont}. Indeed adapting exactly the same argument yields the required result, because we know that $U$ is uniformly bounded on $[0,T-\eps]\times \mptd$. Moreover, the dependence in $\nu$ does not affect the proof. 
\end{proof}

To obtain the regularity in $\nu$ (which is not particularly needed in the rest of the analysis), I pass by the dynamic programming principle.
\begin{Prop}
Under the assumptions of the previous Proposition, for any $t < t + \eps < T$, $µ,\nu \in \mptd$,
$$
\ba
U(t,µ,\nu) = \inf_{(b,\tilde m)}\mathbb E \bigg[\int_t^{t + \eps} &\int_{\T^d\times K}L(x,b(s,x,k))\tilde m_s(dx,dk)ds\\
& + U(t+\eps,(\pi_1)_\#\tilde m_{t+\eps},\nu_{t+ \eps})| \nu_t = \nu\bigg].
\ea
$$
\end{Prop}
I do not provide a proof of this result, which follows from standard techniques. I rather refer the reader to Claisse, Talay and Tan \citep{talay} for similar considerations, in more difficult situations in finite dimension, and to Djete, Possamai and Tan \citep{djete} in infinite dimension, namely because we have already proven that $U$ is continuous in $(t,µ)$, while, almost surely, the variable $\nu$ will only take a finite number of values, hence the measurability of $U$ in $\nu$ does not raise any issue.

I can now state the following continuity result in $\nu$.

\begin{Prop}
Assume that Hypothesis \ref{hyp:final}, and \eqref{hypgrowth}, \eqref{hyplambdatec} hold. Then, for any $\eps > 0$, the function $(t,µ,\nu) \mapsto U(t,µ,\nu)$ is continuous on $[0,T-\eps]\times \mptd^2$.
\end{Prop}
\begin{proof}
First note that we already established the continuity in $µ$. This implies that when $\lambda \equiv 0$, the Lipschitz regularity in $\nu$ is also established (it suffices to reverse time in this deterministic setting). Thus $U_{det}$ satisfies the conclusion of the Proposition. I now make use of the dynamic programming principle. Let $\theta \in (0,T)$ (which is going to be close to $T$) and $\beta$ be the function given by Proposition \ref{prop:behavjumps}, which is H\"older continuous. Let $t \in [0,T-\eps], µ \in \mptd, \nu, \nu' \in \mptd$. From the dynamic programming principle, we know that
$$
|U(t,µ,\nu) - U(t,µ,\nu')| \leq \sup_{m \in \mptd}|U(\theta,m,\nu) - U(\theta,m,\nu')|.
$$
Using the Lipschitz continuity of $U_{det}$ and Proposition \ref{prop:behavjumps} yields
$$
|U(t,µ,\nu) - U(t,µ,\nu')| \leq C\frac{\mathcal{W}_p(\nu,\nu')}{(T-\theta)^{p-1}} + 2\beta(T-\theta).
$$
Remark that, minimizing over $\theta$ alone yields some H\"older modulus of continuity for $U$ in $\nu$. Hence the result is proved because we already knew the Lipschitz regularity in the other variables.
\end{proof}

\subsection{The value function is a viscosity solution}
In this section, I explain how to adapt the proofs of Section \ref{sec:propviscsol} to this stochastic setting. In all this section, the value function is the function $U$ defined in \eqref{defUsto} in the jump case. Here, Hypothesis \ref{hyp:final}, and \eqref{hypgrowth}, \eqref{hyplambdatec} are in force.
Before showing that the value function is indeed a viscosity solution, the definition of viscosity solution needs to be (slightly) adapted to this new setting. Namely, since $U$ depends on two measure arguments, the notions of sub/super differentials need to be adapted. This is a bit arbitrary since the equation only involves the derivatives of $U$ with respect to its first measure argument. Hence, I do not justify the definition in this case, and hope that the reader can deduce it by herself or himself.
\begin{Def}\label{def:viscl}
A usc (resp. lsc) function $U: [0,T) \times \mptd^2\mapsto \R$ is said to be a viscosity sub-solution (resp. super-solution) of \eqref{hjbl} if for any $t\in[0,T), µ,\nu \in \mptd$ and $(\theta,\psi,\phi) \in \partial^+U(t,µ,\nu)$ (resp. $\in \partial^-U(t,µ,\nu)$)
\[
 -\theta +  \int_{\T^d\times \R^d}H(x,z)\psi_x(dz)µ(dx)+ \lambda(t)(U(t,µ,\nu) - U(t,µ,\mathcal{T}\nu)) \leq 0 \text{ (resp. } \geq 0 ). 
 \]
A viscosity solution to \eqref{hjbl} is a locally bounded function $U$ such that $U_*$ is a viscosity super-solution and $U^*$ is a viscosity sub-solution.
\end{Def}
Using the dynamic programming principle and the bounds obtained on the value function just above, we can argue as if we were dealing with a continuous bounded terminal cost. 
First, the property of viscosity sub-solution can be proven following exactly the same argument, as in the deterministic case.
 \begin{Prop}\label{prop:sjf}
The value function $U$ is a viscosity sub-solution to \eqref{hjbl}.
 \end{Prop}
 \begin{proof}
Take $(t,µ,\nu) \in [0,T)\times \mptd^2$ and $(\theta,\psi,\phi) \in \partial^+U(t,µ,\nu)$. Consider a standard probability space $(\Omega',\mathcal{A}',\mathbb{P}')$ (not the one on which the process $(\nu_t)$ is defined, but another one, on which the controller can randomize its controls upon), and a couple $(X_t,Z)$ of random variables on $\Omega'$ such that $\mathcal{L}(X_t,Z)(dx,dz) = µ(dx)\psi(x,dz)$. In particular, this couple is deterministic as a random variable over $\Omega$. As in the deterministic case, define the ODE
 \[
 dX_s = -D_qH(X_s,Z)ds \quad \text{ for } s \in (t,T),
 \]
 with initial condition $X_t$ at time $s = t$. Once again, this process is deterministic from the point of view of $\Omega$. Take $\kappa \in (0,T-t)$ and arguing as in the deterministic setup, we end up with 
 \[
 U(t,µ,\nu) \leq \mathbb E_{\mathbb P}\left[\mathbb E_{\mathbb P'}\left[\int_t^{t+ \kappa}L(X_s,-D_qH(X_s,Z))ds + U(t+ \kappa,\mathcal L(X_{t + \kappa}),\nu_{t + \kappa}) \right]\bigg| \nu_t = \nu\right].
 \]
 It suffices now to use the tower property on the first expectation to write it as the expectation of the expectation conditioned on the number of jumps $n$ of the target process in $[t,t+ \kappa]$. Since we have the equivalences, as $\kappa \to 0$, $\mathbb P(n = 0) \sim 1$, $\mathbb P(n = 1) \sim \lambda (t)\kappa$, $\mathbb P(n \geq 2) \leq C \kappa^2$, we can indeed conclude as in the proof of Proposition \ref{prop:sub-solution}.

 \end{proof}

 Proving the viscosity super-solution property can be done exactly as in the deterministic case, namely by considering first the value function $U^\delta_M$, for $M > 0$, $\delta \in (0,T)$, defined for $(t,µ,\nu) \in [0,T-\delta)\times \mptd^2$ by
$$
\ba
U^\delta_M(t,µ,\nu) = \inf_{(b,\tilde m)} \mathbb E\bigg[& \int_t^{T-\delta}\int_{ \T^{d}\times K}L(x,b(s,x,k))\tilde m_s(dx,dk)ds \\
&+ U(T-\delta,(\pi_1)_\#\tilde m_{T-\delta},\nu_{T-\delta}) | \nu_t = \nu\bigg],
\ea
$$
where the infimum is taken over the same controls as for $U$ in \eqref{defUsto}, except for the fact that, almost surely, in the time interval $[t,T-\delta]$, the control $b$ is bounded by $M$ everywhere. The fact that the restriction is only asked up to $T-\delta$ is necessary for the value function to be bounded. Indeed, otherwise, there will be a controllability issue which would imply that $U^\delta_M= +\infty$.

 In what follows, I am going to fix $\delta$ and pass to the limit $M \to \infty$ to show that $U$ is a viscosity solution on the time interval $(0,T-\delta)$. Then, since $\delta$ is arbitrary, this will lead to the required property. Showing that $U_M^\delta$ is a viscosity super-solution and that the $\lim_{M \to \infty} U^\delta_M$ is also a viscosity super-solution is the exact analogue of what has been done in the deterministic case, once we use the same conditioning as in the proof of Proposition \ref{prop:sjf} just above. Hence, to complete the proof of the fact that $U$ is a viscosity super-solution, it only remains to show that $\liminf_{M \to \infty} U^\delta_M = U$.
 
\begin{Lemma}
For any $(t,µ,\nu) \in [0,T-\delta]\times \mptd^2$, $\lim_{M \to \infty} U^\delta_M(t,µ,\nu) = U(t,µ,\nu)$.
\end{Lemma}
\begin{proof}
I proceed as in the deterministic case. Let $(t,µ,\nu) \in [0,T)\times \mptd^2$ and consider an $\eps$-optimal control $(b,\tilde m)$ for $U(t,µ,\nu)$, in the time interval $[t,T-\delta]$. Recall that this control is thus an adapted process to the filtration of the probability space $(\Omega,\mathcal A, (\mathcal F_t)_{t \geq 0},\mathbb P)$. For $M > 0$, I want to define an associated admissible control $(b^M, \tilde m^M)$ for $U_M(t,µ,\nu)$. This can be done exactly as previously, for each $\omega$. Indeed, for $\mathbb P$ almost every $\omega \in \Omega$, the control $(b,\tilde m)$ can be truncated (only on $[t,T-\delta]$), exactly as is being done in Lemma \ref{lemma:convtrunc}. Hence, for almost every $\omega$, we have a process $(X^{M,\omega}_s(\omega'), Z^{M,\omega}_s(\omega'))_{s \in [t,T- \delta]}$ defined on another probability space $(\Omega',\mathcal A',\mathbb P')$. A priori, this second probability space should depend on $\omega$ itself, but from standard probabilistic arguments, we can always assume that it is constant in $\omega$.\\

If we are able to prove that
\be\label{crucial}
\lim_{M \to \infty}\mathbb E_{\mathbb P} \left[ \int_t^{T-\delta}\mathbb E_{\mathbb P'}[|Z^{M,\omega}_s - Z^{\omega}_s|^p] ds\right] = 0,
\ee
then, it is a simple exercise to verify that the computation of Lemma \ref{lemma:convtrunc} yields the result. But once again, \eqref{crucial} follows from the dominated convergence theorem, and the double dependence in $\omega \in \Omega$ and $\omega' \in \Omega'$ does not perturb the argument. Hence, the result is proved. 
\end{proof}
I insist that in the previous proof, I looked at the optimal control problem at hand as a problem with a continuous terminal cost at time $t= T-\delta$, hence the truncation of the controls does not raise an issue with the constraint at terminal time.

Putting all the pieces of this section together, we have proven the following.
\begin{Theorem}
The value function $U$ is a viscosity solution to \eqref{hjbl} in $(0,T)\times \mptd^2$.
\end{Theorem}

\subsection{Characterization of the value function as a viscosity solution}
To conclude the present study of the value function $U$, I will show that it can be characterized as being the unique viscosity solution to \eqref{hjbl} satisfying the behaviour obtained in Proposition \ref{prop:behavjumps}. The standing assumptions are here the same as in the previous section, that is, Hypothesis \ref{hyp:final}, and \eqref{hypgrowth}, \eqref{hyplambdatec} hold. In addition, I focus on the case $p=2$ to lighten notation.

\begin{Theorem}
Assume that $\lambda$ is continuous and that $\mathcal T$ is Lipschitz continuous for $\mathcal{W}_2$. Then, there exists at most one viscosity solution $V$ to \eqref{hjbl} on $(0,T)\times \mptd^2$ satisfying 
\be\label{eq:unifter}
\lim_{t \to T} \|V(t,\cdot,\cdot)-U_{det}(t,\cdot,\cdot)\|_\infty = 0.
\ee
\end{Theorem}
\begin{proof}
The proof works in the same way as the comparison principles I presented in Section \ref{sec:visc} except for two features: i) the presence of the term in $\lambda$ and ii) the non-usual behaviour at the terminal time. I only highlight how those differences need to be handled. The structure of the proof is the same. We consider $V_1,V_2$, respectively sub and super-solution to \eqref{hjbl} satisfying both \eqref{eq:unifter}, and we focus on points of maxima of 
$$
\ba
(t,s,µ,µ',\nu,\nu') \mapsto V_1(t,µ,\nu) - V_2(s,µ',\nu') - &\frac{1}{2\eps}((t-s)^2 + \mathcal{W}_2^2(µ,µ') + \mathcal{W}_2^2(\nu,\nu'))\\
&-\rho (2T -t-s)\color{black}.
\ea
$$
We want to show that $V_1 \leq V_2$ and we assume it is not the case, hence, that there exists $ \kappa > 0, (t_0,µ_0,\nu_0) \in (0,T)\times \mptd^2$ such that
\be\label{contrlast}
V_1(t_0,µ_0,\nu_0) > V_2(t_0,µ_0,\nu_0) + \kappa.
\ee
 Consider a point of maximum $(t_{\epsilon},s_{\epsilon},µ_{\epsilon},µ_{\epsilon}',\nu_{\epsilon},\nu_{\epsilon}')$ of the previous function on $[0,T-\delta]^2\times \mptd^2$. Using \eqref{eq:unifter}, we know that $\|V_1(T-\delta,\cdot,\cdot)-V_2(T-\delta,\cdot,\cdot)\|_\infty \leq 2\beta(\delta)$, for some modulus of continuity $\beta$. Hence, taking $\delta, \eps$ sufficiently small, we contradict \eqref{contrlast} if $t_\eps$ or $s_\eps$ is equal to $T-\delta$. This argument allows us to handle the terminal time.
 
For the rest of the proof, I assume that I always have $t_\eps,s_\eps < T-\delta$, note that it is only possible if $\rho >0$ is sufficiently small relative to $\kappa$ and $t_0$. Using the properties of viscosity solutions of $V_1$ and $V_2$, and arguing as in the proof of Theorem \ref{thm:compt}, we arrive at
\[
\begin{aligned}
2 \rho + \lambda(t_{\epsilon})(V_1(t_{\epsilon},µ_{\epsilon},\nu_{\epsilon})-V_1(t_{\epsilon},µ_{\epsilon},\mathcal{T}\nu_{\epsilon})) - \lambda(s_{\epsilon})(V_2(s_{\epsilon},µ'_{\epsilon},\nu'_{\epsilon}) - V_2(s_{\epsilon},µ'_{\epsilon},\mathcal{T}\nu'_{\epsilon})) \leq o(1),
\end{aligned}
\]
where the right side term vanishes as $\epsilon \to 0$. We can then compute
\[
\begin{aligned}
\lambda&(t_{\epsilon})(V_1(t_{\epsilon},µ_{\epsilon},\nu_{\epsilon})-V_1(t_{\epsilon},µ_{\epsilon},\mathcal{T}\nu_{\epsilon})) - \lambda(s_{\epsilon})(V_2(s_{\epsilon},µ'_{\epsilon},\nu'_{\epsilon}) - V_2(s_{\epsilon},µ'_{\epsilon},\mathcal{T}\nu'_{\epsilon}))\\
&\geq -C |\lambda(t_{\epsilon}) - \lambda(s_{\epsilon})| + \lambda(t_{\epsilon})(V_1(t_{\epsilon},µ_{\epsilon},\nu_{\epsilon})-V_1(t_{\epsilon},µ_{\epsilon},\mathcal{T}\nu_{\epsilon}) -(V_2(s_{\epsilon},µ'_{\epsilon},\nu'_{\epsilon}) - V_2(s_{\epsilon},µ'_{\epsilon},\mathcal{T}\nu'_{\epsilon})))\\
&\geq -C |\lambda(t_{\epsilon}) - \lambda(s_{\epsilon})| + \lambda(t_{\epsilon})\frac{1}{2\epsilon}(\mathcal{W}_2^2(\nu_{\epsilon},\nu'_{\epsilon}) - \mathcal{W}_2^2(\mathcal{T}\nu_{\epsilon},\mathcal{T}\nu'_{\epsilon})\color{black})\\
&\geq - C |\lambda(t_{\epsilon}) - \lambda(s_{\epsilon})| - \lambda(t_{\epsilon})\frac{\tilde C^2\color{black}}{2\epsilon}\mathcal{W}_2^2(\nu_{\epsilon},\nu'_{\epsilon}),
\end{aligned}
\]
 where $\tilde C$ is the Lipschitz constant of $\mathcal T$. Remark that in the previous computation, $C$ only depends on the bounds on $V_1$ and $V_2$ on $\{t \leq T-\delta\}$. The limit of the last lower bound in the previous chain of inequalities is $0$ as $\epsilon \to 0$. Indeed, from Lemma \ref{lemma:visc1} $\lim_{\epsilon \to 0} s_{\epsilon}-t_{\epsilon} =\lim_{\epsilon \to 0} \epsilon^{-1}\mathcal{W}_2^2(\nu_{\epsilon},\nu'_{\epsilon})= 0$. Hence, using the continuity of $\lambda$ we obtain $\rho \leq 0$ which is a contradiction. Hence, $V_1 \leq V_2$.
 
 In particular, this comparison implies that the value function $U$ is the only viscosity solution to \eqref{hjbl} satisfying \eqref{eq:unifter}.

\end{proof}
\subsection*{Bibliographical comments}
I introduced and studied such models in \citep{bertucci2024stochastic} as a stochastic extension of the standard optimal transport problem.

\newpage
\section*{Acknowledgments}
Funded by the European Union (ERC, PaDiESeM, 101222038). Views and opinions expressed are however those of the author only and do not necessarily reflect those of the European Union or the European Research Council. Neither the European Union nor the granting authority can be held responsible for them.
I acknowledge a partial support from the Chair FDD (Institut Louis Bachelier).

\bibliographystyle{plain}
\bibliography{bib1}



\bigskip\bigskip\bigskip\bigskip
\uppercase{CEREMADE, CNRS, UMR7534, Université Paris Dauphine-PSL} 

\textit{E-mail address: } bertucci@ceremade.dauphine.fr
\end{document}